\documentclass[a4paper,11pt]{amsart}
\usepackage{amsmath,amssymb,amsthm,mathrsfs,comment,cases}
\usepackage[all]{xy}
\usepackage{graphicx}

\allowdisplaybreaks[3]
\usepackage{geometry}
\newtheorem{thm}{Theorem}[section]
\newtheorem{lemma}[thm]{Lemma}
\newtheorem{theorem}[thm]{Theorem}
\newtheorem{proposition}[thm]{Proposition}
\newtheorem{corollary}[thm]{Corollary}
\newtheorem{problem}[thm]{Problem}

\newtheorem{theoremintro}{Theorem}

\theoremstyle{definition}
\newtheorem{definition}[thm]{Definition}
\newtheorem{example}[thm]{Example}
\newtheorem{remark}[thm]{Remark}

\newtheorem{setting}[thm]{Setting}
\newtheorem{notation}[thm]{Notation}
\newtheorem{supply}[thm]{Supplementation}

\newtheorem{step}{Step}
\newtheorem{setting8}{Setting}

\newcommand{\QQ}{\mathbb Q}
\newcommand{\RR}{\mathbb R}
\newcommand{\RRR}{\mathbb {R}_{\geq 0}}
\newcommand{\R}{\mathcal R}
\newcommand{\NN}{\mathbb N}
\newcommand{\N}{\mathcal N}
\newcommand{\HH}{\mathbb H}
\renewcommand{\H}{\mathcal H}
\newcommand{\Z}{\mathbb Z}
\newcommand{\F}{\mathcal F}

\newcommand{\B}{\mathcal B}

\newcommand{\gD}{\Delta}

\newcommand{\gL}{\Lambda}
\newcommand{\gG}{\Gamma}
\newcommand{\gS}{\Sigma}
\newcommand{\id}{\mathrm{id}}
\newcommand{\SCyl}{\mathrm{SCyl}}
\newcommand{\rk}{\overline{\mathrm{rk}}}

\newcommand{\tn}[1]{\textnormal{#1}}
\newcommand{\ti}[1]{\textit{#1}}

\renewcommand{\:}{\colon}
\renewcommand{\tilde}{\widetilde}
\renewcommand{\hat}{\widehat}
\newcommand{\SC}{\mathrm{SC}}
\newcommand{\Cay}{\mathrm{Cay}}
\newcommand{\Stab}{\mathrm{Stab}}

\newcommand{\GC}{\mathrm{GC}}

\newcommand{\ol}[1]{\overline{#1}}
\newcommand{\Isom}{\mathrm{Isom}}
\newcommand{\An}{\mathrm{An}}
\newcommand{\Area}{\mathrm{Area}}
\newcommand{\CH}{\hat{CH}}
\newcommand{\Comp}{\mathrm{Comp}}

\begin{document}
\title{Subset currents on surfaces}

\author[D. Sasaki]{Dounnu Sasaki}
\address{Faculty of Science and Engineering, Waseda University, Okubo 3-4-1, Shinjuku, Tokyo 169-8555, Japan}
\email{dounnu-daigaku@moegi.waseda.jp}



\begin{abstract}
Subset currents on hyperbolic groups were introduced by Kapovich and Nagnibeda as a generalization of geodesic currents on hyperbolic groups, which were introduced by Bonahon and have been successfully studied in the case of the fundamental group $\pi_1 (\Sigma)$ of a compact hyperbolic surface $\Sigma$.
Kapovich and Nagnibeda particularly studied subset currents on free groups.
In this article, we develop the theory of subset currents on $\pi_1(\Sigma )$, which we call subset currents on $\Sigma$.
We prove that the space $\mathrm{SC}(\Sigma)$ of subset currents on $\Sigma$ is a measure-theoretic completion of the set of conjugacy classes of non-trivial finitely generated subgroups of $\pi_1 (\Sigma )$, each of which geometrically corresponds to a convex core of a covering space of $\Sigma$.
This result was proved by Kapovich-Nagnibeda in the case of free groups, and is also a generalization of Bonahon's result on geodesic currents on hyperbolic groups.
We will also generalize several other results of them. 
Especially, we extend the (geometric) intersection number of two closed geodesics on $\Sigma$ to the intersection number of two convex cores on $\Sigma $ and, in addition, to a continuous $\mathbb{R}_{\geq 0}$-bilinear functional on $\mathrm{SC}(\Sigma)$.
\end{abstract}

\maketitle

\tableofcontents

\section{Introduction}
Subset currents on hyperbolic groups were introduced by Kapovich and Nagnibeda \cite{KN13} as a generalization of geodesic currents on hyperbolic groups, which were introduced by Bonahon \cite{Bon86, Bon88b}. Geodesic currents have been successfully studied in the case of the fundamental group $\pi_1 (\Sigma)$ of a compact hyperbolic surface $\Sigma$ and also has been useful for several areas of mathematics, including geometric topology, geometric group theory, the study of Kleinian groups, ergodic theory and dynamics, etc.

Kapovich and Nagnibeda particularly studied the case of subset currents on free groups and proved that the space $\SC (F)$ of subset currents on a free group $F$ of finite rank can be thought of as a measure-theoretic completion of the set of conjugacy classes of (non-trivial) finitely generated subgroups of $F$, which is a generalization of Bonahon's result on geodesic currents on hyperbolic groups.
This result played a fundamental role in the study of subset currents on free groups.

In this article, we develop the theory of subset currents on $\pi_1(\Sigma )$ for a compact hyperbolic surface $\gS$, which we call subset currents on $\Sigma$.
Note that $\pi_1(\Sigma )$ is a free group of finite rank if $\gS$ has boundary.
We generalize the above result for free groups to the case of a closed hyperbolic surface $\gS$, in other words, to the case of a surface group $\pi_1(\gS)$.
We will also generalize several other results of Bonahon on geodesic currents on $\Sigma$ and that of Kapovich-Nagnibeda on subset currents on free groups.
Especially, we extend the (geometric) intersection number $i$ of two closed geodesics on $\gS$ to a continuous $\RRR$-bilinear functional $i_{\SC}$ on $\mathrm{SC}(\Sigma)$, which is also an extension of Bonahon's intersection number $i_{\GC}$ on $\GC(\gS)$.

In the case that $\gS$ has no boundary, the intersection number $i_{\GC}$ on $\GC(\gS)$ and the ``Liouville current'' $L(m) \in \GC(\gS)$ with respect to a hyperbolic metric $m$ or a ``non-positive curved metric'' $m$ on $\GC(\gS)$ have the interesting property that for any essential closed curve $c$ on $\gS$, $i_{\GC}(L(m), c)$ equals the $m$-length of $c$, which is the infimum of the length of a closed curve $c'$ on $(\gS ,m)$ free-homotopic to $c$.
This property plays the important role in the study of the length spectral rigidity problem for the class of the metric $m$, which asks whether the information of the all of the length of (simple) closed ``geodesics'' determines the metric or not (see \cite{Ota90,HP97,DLR10,BL18}).
More generally, it was also applied to the case of Kleinian surface groups and character varieties of higher dimensions in \cite{BC17,BCL17,BCLS18, BIPP17, BIPP19}. We also refer a new survey \cite{EU18}.

We generalize the above property of $i_{\GC}$ to the case of the intersection number $i_{\SC}$ and see that $i_{\SC}(L(m), S)$ equals the half of the sum of the $m$-length of the boundary components of $S$ for a ``simple compact surface'' $S$ on $\gS$ such as the convex core of a non-trivial finitely generated subgroup of $\pi_1(\gS)$.
By such results, subset currents can be considered as a natural good generalization of geodesic currents and expected to be useful similarly to geodesic currents.

We will give a detailed introduction in the following subsections.

\subsection{Background}
In general, the notion of geodesic currents can be defined on an infinite hyperbolic group $G$, which was introduced by Bonahon \cite{Bon88b}. A geodesic current on $G$ is a locally finite (i.e. finite on any compact subset) $G$-invariant Borel measure on the space $\partial_2G$ of 2-element subsets of the (Gromov) boundary $\partial G$.
The space $\GC(G)$ of geodesic currents on $G$, which is equipped with weak-$\ast$ topology, can be thought of as a completion of the space of conjugacy classes of infinite cyclic subgroups of $G$ with positive real weight in the following meaning. For an infinite-order element $g\in G$ we can define a counting geodesic current $\eta_g$ corresponding to the subgroup $\langle g\rangle$ by
\[ \eta_g:=\sum_{u\langle g\rangle \in G /\langle g \rangle }\delta_{u\gL (\langle g \rangle )},\]
where $\delta_{\gL (\langle g \rangle )}$ is the Dirac measure at the limit set $\gL(\langle g\rangle )$ of $\langle g\rangle $ on $\partial_2 G$. For $h\in G$ we can see that $\eta_{hgh^{-1}}=\eta_g$.
Bonahon \cite{Bon88b} proved that the set of all positive real weighted counting geodesic currents on $G$:
\[ \{ c\eta_g \mid c\in \RRR ,\ g\in G\setminus \{ \id \} \}, \]
where $c\eta_g$ is called a rational geodesic current on $G$, is a dense subset of $\GC(G)$. We call this property the denseness property of rational geodesic currents.

The study of geodesic currents was started from the case that a hyperbolic group $G$ is the fundamental group $\pi_1(\gS)$ of a orientable compact connected hyperbolic surface $\gS$ possibly with (geodesic) boundary. We always assume that a surface is orientable and connected, and can have boundary.
In this case geodesic currents can be understood in a more geometric way and has a lot of application.
If $\gS$ has the boundary, then $\pi_1(\gS)$ is a free group of finite rank, and if $\gS$ has no boundary, then $\pi_1(\gS)$ is called a surface group.
We write $\GC(\pi_1(\gS))$ simply as $\GC(\gS)$ and call $\GC(\gS)$ the space of geodesic currents on $\gS$ when we identify $\partial \pi_1(\gS)$ with the ideal boundary of the universal cover $\tilde{\gS}$ of $\gS$ by using the action of $\pi_1(\gS )$ on $\tilde{\gS}$.
In this setting, an element of $\partial_2 G$ corresponds to a geodesic line in $\tilde{\gS}$, and each conjugacy class of an infinite cyclic subgroup of $\pi_1(\gS)$ corresponds to a homotopy class of an unoriented closed curve on $\gS$ and also corresponds to an unoriented closed geodesic on $\gS$. 

For two closed curves $c_1,c_2$ on $\gS$, which are continuous maps from $S^1$ to $\gS$, the (geometric) intersection number $i$ of $c_1,c_2$ is the number of contractible components of the fiber product of $S^1$ and $S^1$ corresponding to $c_1,c_2$. If $c_1,c_2$ are simple and transversal, then $i(c_1,c_2)$ coincides with the cardinality of $c_1(S^1)\cap c_2(S^1)$. The intersection number $i$ of two homotopy classes of (unoriented) closed curves $[c_1],[c_2]$ is the minimum of $i(c_1',c_2')$ taken over all $c_1'\in [c_1],c_2'\in [c_2]$. 
For two non-trivial elements $g_1,g_2\in G$ we can define $i(g_1,g_2)$ to be the intersection number of homotopy classes of unoriented closed curves on $\gS$ corresponding to $g_1,g_2$. Note that if $c_1,c_2$ are closed geodesics on $\gS$, then $i(c_1,c_2)=i([c_1],[c_2])$. Such $c_1,c_2$ are said to be in minimal position. 
Bonahon \cite{Bon86} proved that there exists a unique continuous $\RRR$-bilinear functional $i_{\GC}$ from $\GC(\gS)\times \GC(\gS)$ to $\RRR$ such that for any non-trivial elements $g_1,g_2\in \pi_1(\gS)$ we have
\[ i_{\GC}(\eta_{g_1},\eta_{g_2})=i(g_1,g_2).\]
The uniqueness of $i_{\GC}$ is the result of the denseness property of rational geodesic currents. In this meaning, we say that $i_{\GC}$ is an extension of $i$.

In the case that $\Sigma$ has no boundary, Bonahon \cite{Bon88} proved that there exists an embedding $L$ from the Teichm\"uller space $\mathcal{T}(\gS)$ of $\gS$ to $\GC(\gS)$, and for $m\in \mathcal{T}(\gS)$ and a non-trivial $g\in \pi_1(\gS)$ the intersection number $i_{\GC}(L(m), \eta_g)$ coincides with the length of the $m$-geodesic corresponding to $g$, which we call the $m$-length of $g$. The geodesic currents $L(m)$ is called the Liouville current associated with $m$.
From the above, we obtain the $m$-length functional $\ell_m=i_{\GC}(L(m),\cdot )$ on $\GC (\gS )$ such that $\ell_m(\eta_g)$ equals the $m$-length of $g$ for every non-trivial $g\in \pi_1(\gS)$.

Even in the case that $\gS$ has boundary, we can consider the doubled surface $D\Sigma$ of $\Sigma$, which we can get by gluing two copies of $\Sigma$ together at corresponding boundary components and its hyperbolic structure $m_d$ is inherited from the hyperbolic structure $m$ of $\Sigma$. Then for the hyperbolic structure $m_d$ of $D\Sigma$ we obtain $L(m_d)\in GC(D\Sigma)$, which we will denote by $L(m)$ simply.
Note that $\GC(\Sigma)$ is naturally embedded into $\GC (D\Sigma)$ (see the above of Theorem \ref{thmintro:include map} for detail) and the restriction of $i_{\GC}(L(m),\cdot )$ to $\GC(\Sigma)$ will be the $m$-length functional $\ell_m\: \GC(\Sigma)\rightarrow \RRR$.
From an application viewpoint, the intersection number $i$ and the embedding $L$ often play important roles.

The notion of subset currents is also defined on an infinite hyperbolic group $G$. A subset current on $G$ is a locally finite $G$-invariant Borel measure on the space $\H (\partial G)$ of closed subsets of $\partial G$ containing at least $2$ points, which is endowed with the Vietoris topology. The Vietoris topology on $\H (\partial G)$ coincides with the topology induced by the Hausdorff distance. A geodesic current on $G$ can be considered as a subset current on $G$ since $\partial_2G$ is a $G$-invariant (closed) subspace of $\H (\partial G)$. Kapovich and Nagnibeda \cite{KN13} introduced the notion of subset currents on hyperbolic groups and particularly studied the space $\SC(F)$ of subset currents on a free group $F$ of rank $\geq 2$. For a non-trivial finitely generated subgroup $H$ of $F$ they defined a counting subset current $\eta_H$ by
\[ \eta_H:=\sum_{gH\in F/H} \delta_{g\gL(H)},\]
where $\delta_{\gL(H)}$ is the Dirac measure at the limit set $\gL(H)$ of $H$ on $\H (\partial F)$. We can see that $\eta_{gHg^{-1}}=\eta_H$ for $g\in F$.
They proved that the set $\SC_r(F)$ of all positive real weighted counting subset currents on $F$, which are called rational subset currents on $F$, is a dense subset of $\SC(F)$. In this meaning the space $\SC(F)$ can be thought of as a measure-theoretic completion of the set of conjugacy classes of non-trivial finitely generated subgroups of $F$.

Let $\gD$ be a finite connected graph whose fundamental group is isomorphic to $F$ and whose vertices have degree larger than or equal to $2$.
For a non-trivial finitely generated subgroup $H$ of $F$ we define a $\gD$-core graph $\gD_H$ to be the smallest subgraph of the covering space corresponding to $H$ such that the inclusion map is a homotopy equivalence map.
Some properties of counting subset currents tell us that the $\gD$-core graph $\gD_H$ is closely related to $\eta_H$. 
Let $H'$ be a $k$-index subgroup of $H$. Then we can see that $\eta_{H'}=k\eta_H$ by the definition. This property corresponds to the fact that there exists a $k$-fold covering map from $\gD_{H'}$ to $\gD_H$. Note that we have $\chi(\gD_{H'})=k\chi(\gD_H)$, where $\chi (\gD_H)$ is the Euler characteristic of $\gD_H$.

We define the reduced rank of a non-contractible connected graph to be the negative of the Euler characteristic and define the reduced rank of a contractible graph to be $0$. Moreover, we define the reduced rank $\rk$ of a free group $F_N$ of rank $N\in \NN \cup \{ 0\}$ to be $\max \{N-1, 0\}$. By the definition, the reduced rank of a connected graph whose fundamental group is isomorphic to $F_N$ equals the reduced rank of $F_N$.

A finitely generated subgroup of $F$ is also a free group of finite rank, and we can consider $\rk$ as a map from the set of finitely generated subgroups of $F$ to $\mathbb{Z}_{\geq 0}$.
Kapovich and Nagnibeda \cite{KN13} extended the reduced rank $\rk$ to a continuous $\RRR$-linear functional $\rk$ on $\SC(F)$.
In fact, they constructed $\RRR$-linear functionals $V^\#,E^\#$ from $\SC(F)$ to $\RRR$ satisfying the condition that for every non-trivial finitely generated subgroup $H$ of $F$, $V^\#(\eta_H)$ equals the number of vertices of $\gD_H$ and $E^\#(\eta_H)$ equals the number of edges of $\gD_H$. Then we can obtain the reduced rank functional $\rk$ as $E^\#-V^\# $.

For two finitely generated subgroups $H,K$ of $F$ we define the product $\N$ of $H$ and $K$ by
\[ \N (H,K):= \sum_{HgK \in H\backslash F/K}\rk (H\cap gKg^{-1}),\]
where $H\backslash F/K$ is the set of all double cosets of $H$ and $K$. By using this product $\N$ the Strengthened Hanna Neumann Conjecture can be written as follows: the inequality
\[ \N (H,K)\leq \rk (H)\rk (K)\]
holds for any two finitely generated subgroups $H$ and $K$ of $F$.
This conjecture was individually proved by Friedman \cite{Fri15} and Mineyev \cite{Min12}. Geometrically, the product $\N(H,K)$ equals the sum of the reduced rank of all connected components of the fiber product graph $\gD_H\times_\gD \gD_K$ when $H$ and $K$ are non-trivial.
In \cite{Sas15} the product $\N$ was extended to a continuous $\RRR$-bilinear functional $\N$ on $\SC (F)\times \SC(F)$.
As a corollary, we can obtain the following inequality:
\[ \N (\mu , \nu )\leq \rk (\mu )\rk (\nu )\]
for any two subset currents $\mu, \nu \in \SC (F)$. 

\subsection{Main results}
First, we develop the theory of subset currents on hyperbolic groups.
We prove that the space of subset currents on an infinite hyperbolic group $G$ is a locally compact, separable and completely metrizable space in Section \ref{sec:subset currents on hyperbolic groups}.
For a subgroup $H$ of $G$ we define a $G$-invariant measure $\eta_H$ on $\H (\partial G)$ by
\[ \eta_H:=\sum_{gH\in G/H} \delta_{g\gL (H)}.\]
If $H$ is a finite group, then we define $\eta_H$ to be the zero measure.
We prove that $\eta_H$ is a locally finite measure if and only if $H$ is a quasi-convex subgroup of $G$. In this case we call $\eta_H$ a counting subset current on $G$ and call a positive real weighted counting subset current on $G$ a rational subset current on $G$.

More generally, for a point $S\in \H (\partial G)$ we can define a $G$-invariant measure $\eta_S$ by taking the $G$-orbit of $S$.
Explicitly,
\[ \eta_S := \sum_{g\mathrm{Stab}(S)\in G/\mathrm{Stab}(S)} \delta_{g\gL (\mathrm{Stab}(S))},\]
where $\mathrm{Stab}(S)$ is the stabilizer of $S$ with respect to the action of $G$. Then we can see that $\eta_S$ is locally finite if and only if $\mathrm{Stab}(S)$ is a quasi-convex subgroup of $G$ and $S=\gL (\mathrm{Stab}(S))$.

Therefore, the set $\SC_r(G)$ of all rational subset currents on $G$ is a natural subset of $\SC(G)$ consisting of ``discrete measures''.
Hence we are interested in whether $\SC_r(G)$ is a dense subset of $\SC(G)$.
Note that the $\RRR$-linear subspace $\mathrm{Span}(\SC_r(G))$ of $\SC(G )$ generated by $\SC_r(G)$ is also a natural subspace of $\SC(G)$ consisting of ``discrete measures'', and we are also interested in whether $\mathrm{Span}(\SC_r(G))$ is a dense subset of $\SC(G)$.
Both of these problems are still open in general in contrary to the result of Bonahon on the space of geodesic currents on a hyperbolic group. The difficulty comes from the nature that constructing quasi-convex subgroups is much harder than finding generators of infinite cyclic subgroups. We say that an infinite hyperbolic group $G$ has the denseness property of rational subset currents if $\SC_r(G)$ is a dense subset of $\SC (G)$.

In the case of a free group $F$, Kapovich and Nagnibeda \cite{KN13} first proved that $\SC_r(F)$ is a dense subset of $\mathrm{Span}(\SC_r(F))$, and then proved that $\mathrm{Span}(\SC_r(F))$ is a dense subset of $\SC(F)$. Bonahon \cite{Bon88b} also divided the proof of the denseness property of rational geodesic currents for a hyperbolic group into such two steps.
We are not aware whether $\SC_r(G)$ is a dense subset of $\mathrm{Span}(\SC_r(G))$ for a general infinite hyperbolic group $G$.

From the viewpoint of the application of subset currents, solving either one of the two problems mentioned in the above for a surface group is important.
Actually, the former of the two problems was presented by Kapovich and Nagnibeda in \cite{KN13}.
In this article, we solve the problem and obtain the following theorem:

\begin{theoremintro}
Let $\gS$ be a compact hyperbolic surface. Then $\SC(\pi_1(\gS))$ has the denseness property of rational subset currents, that is,
the set 
\[ \{ r\eta_{H} \mid r\geq 0,\ H\leq \pi_1(\gS) : \text{a finitely generated subgroup} \}\]
is a dense subset of $\SC(\pi_1(\gS))$.
\end{theoremintro}

Note that a subgroup $H$ of $\pi_1 (\gS)$ is a quasi-convex subgroup of $\pi_1(\gS)$ if and only if $H$ is a finitely generated subgroup of $\pi_1 (\gS)$, that if $\gS$ has boundary, then $\pi_1(\gS)$ is a free group of finite rank, and that $\eta_{\{\id\}}=0$.

Our method of proving the denseness property for a surface group is partially based on the method of proving the denseness property for a free group of finite rank in \cite{Kap17}.
The point of our method is that we take a sequence of finite-rank free subgroups $\{ H_n\}$ of the surface group $\pi_1(\gS)$ ``approximating'' $\pi_1(\gS)$, and construct $\nu \in \mathrm{Span}(\SC_r(H_n))$ based on a given subset current $\mu \in \SC(\pi_1(\gS))$ for a sufficiently large $n$. From $\nu' \in \SC_r(H_n)$ sufficiently close to $\nu$, we can obtain a rational subset current on $\pi_1(\gS)$ sufficiently close to $\mu$.

As with geodesic currents on $\gS$, we write $\SC(\pi_1(\gS))$ simply as $\SC(\gS)$ and call $\SC(\gS)$ the space of subset currents on $\gS$ when we identify $\partial \pi_1(\gS)$ with the ideal boundary of the universal cover $\tilde{\gS}$ of $\gS$ by using the action of $\pi_1(\gS )$ on $\tilde{\gS}$.

From now on, we will talk about several continuous extensions of invariants of finitely generated subgroups (or pairs of finitely generated subgroups) of $\pi_1(\gS)$ to continuous $\RRR$-linear (or $\RRR$-bilinear) functionals on $\SC (\gS)$.
The outline of the strategy to prove the extensions is the same as that by Bonahon and Kapovich-Nagnibeda.
First, we construct an $\RRR$-linear functional on $\SC(\gS)$ associating a counting subset current for a non-trivial finitely generated subgroup of $\pi_1(\gS)$ with a certain invariant. Then we prove the continuity of the functional, which is the main part of the proof. Finally, we see that such a functional is unique by the denseness property of rational subset currents.
In this way we can obtain a concrete expression of the functional.

Since $\SC(\gS)$ is a completely metrizable space and the set $\SC_r(\gS)$ of rational subset currents on $\gS$ is a dense subset of $\SC(\gS)$, we can extend a continuous functional on $\SC_r(\gS)$ uniquely to a continuous functional on $\SC(\gS)$. We will also use this method in Section \ref{sec:an intersection functional}.

Let $\gG$ be a non-trivial torsion-free convex-cocompact Kleinian group acting on the $n$-dimensional hyperbolic space $\HH^n$ for $n\geq 2$. Then $\gG$ is a hyperbolic group, and we identify the boundary $\partial \gG$ with the limit set $\gL (\gG )\subset \partial \HH^{n}$. From the assumption, $\gG$ acts on the convex hull $CH(\gL( \gG ))$ of $\gL(\gG)$ cocompactly, which implies that the volume of the convex core $C_\gG:=\gG \backslash CH (\gL (\gG ))$ is finite. Then every non-trivial quasi-convex subgroup $H$ of $\gG$ also acts on the convex hull $CH(\gL (H))$ cocompactly. We prove that there exists a continuous $\RRR$-linear functional $\mathrm{Vol}$ on $\SC (\gG)$ such that for every non-trivial quasi-convex subgroup $H$ of $\gG$, $\mathrm{Vol}(\eta_H)$ equals the volume of the convex core $C_H$ corresponding to $H$.

In the case that $n=2$, $\gG$ is a free group of finite rank or a surface group, and the area of $C_H$ equals $-2\pi \chi (C_H)$ by the Gauss-Bonnet theorem. We define the reduced rank $\rk$ of a surface group to be the negative of the Euler characteristic of a closed surface whose fundamental group is isomorphic to the surface group. Then we obtain the following theorem, which is a generalization of the reduced rank functional on $\SC (F)$ in \cite{KN13}.

\begin{theoremintro}
Let $\gS$ be a compact hyperbolic surface. There exists a unique continuous $\RRR$-linear functional $\rk$ on $\SC (\gS)$ such that for every finitely generated subgroup $H$ of $\pi_1(\gS)$ we have 
\[ \rk (\eta_H ) =\rk (H).\]
\end{theoremintro}

From the definition of the reduced rank for surface groups, we can generalize the product $\N$ to the product of two finitely generated subgroups $H$ and $K$ of $\pi_1(\gS)$ for a closed hyperbolic surface $\gS$, that is,
\[ \N(H,K):=\sum_{HgK\in H\backslash \pi_1 (\gS )/K}\rk (H\cap gKg^{-1}).\]
In the case that $H$ and $K$ are non-trivial, the product $\N(H,K)$ equals the sum of the reduced rank of all connected components of the fiber product $C_H\times_\gS C_K$. The reduced rank of a non-contractible component is the negative of the Euler characteristic and the reduced rank of contractible component is $0$.
Then, we prove the following theorem, which is a generalization of the intersection functional $\N$ on $\SC (F)$.

\begin{theoremintro}
Let $\gS$ be a compact hyperbolic surface. There exists a unique symmetric continuous $\RRR$-bilinear functional $\N$ on $\SC (\gS )$ such that for any two finitely generated subgroups $H$ and $K$ of $\pi_1(\gS)$ we have 
\[ \N (\eta_H,\eta_K ) =\N (H, K).\]
\end{theoremintro}

As far as the author knows, the surface group version of the Strengthened Hanna Neumann Conjecture is still an open problem. By using the continuity of $\N$ and $\rk$ if we can prove the inequality for a dense subset of $\SC(\gS )$, then the conjecture is true for any two subgroups of $\pi_1(\gS)$ for a closed hyperbolic surface $\gS$. This gives us a new approach to the conjecture.

The intersection functional $\N$ on $\SC(\gS)$ has the property that for every $\mu \in \SC (\gS)$ we have
\[ \N (\eta_{\pi_1 (\gS )}, \mu )=\rk (\mu ).\]
In this meaning $\N$ can be thought of as a generalization of the reduced rank functional $\rk$.
We note that for any $\mu, \nu \in \SC (\gS)$ if $\nu \in \GC (\gS)$, then $\N(\mu,\nu )=0$.

Our method of proving the above theorem is based on the method of constructing the intersection functional $\N$ on $\SC (F)$ in \cite{Sas15}.
We will use the denseness property of rational subset currents for $\pi_1(\gS)$ in order to prove the existence of the functional $\N$.
Since the reduced rank of a contractible component is not the Euler characteristic, we need to count the number of contractible components of the fiber product $C_H\times_\gS C_K$. For this purpose we can use the \textit{(geometric) intersection number} $i_{\SC}$ on $\SC (\gS )\times \SC (\gS)$, which is a natural extension of the intersection number of closed geodesics (or curves).

The intersection number $i(H,K)$ of $H$ and $K$ is defined to be the number of contractible components of the fiber product $C_H\times_\gS C_K$.
Note that if $H$ and $K$ are infinite cyclic groups generated by $g_1,g_2\in \pi_1(\gS)$ respectively, then $i(H,K)=i(g_1,g_2)$ since $C_H$ and $C_K$ are closed geodesics and in minimal position.
We prove the following theorem:

\begin{theoremintro}
Let $\gS$ be a compact hyperbolic surface. There exists a unique continuous symmetric $\RRR$-bilinear functional $i_{\SC }$ on $\SC (\gS )$ such that for any non-trivial finitely generated subgroups $H$ and $K$ of $\pi_1(\gS)$ we have 
\[ i_{\SC} (\eta_H,\eta_K ) =i(H, K).\]
\end{theoremintro}

Note that $i(H,K)$ depends on $\gS$ if $\pi_ 1(\gS)$ is a free group, since there exist other compact hyperbolic surfaces that are not homeomorphic to $\gS$ but whose fundamental groups are isomorphic to $\pi_1 (\gS)$.

We also introduce the intersection number of two ``simple compact surfaces'' on an (orientable) compact surface $\gS$, where $\gS$ is not necessarily a hyperbolic surface and the definition is purely topological.
We overview the definition and its property here (see the beginning of Subsection \ref{subsec:Intersection number of surfaces} and Definition \ref{def:intersection number of simple compact surfaces} for detail).
Let $S$ be a compact surface or a circle $S^1$. A pair of $S$ and a continuous map $s$ from $S$ to $\gS$ is called a simple compact surface on $\gS$ if $s$ is locally injective and the restriction of $s$ to each component of the boundary $\partial S$ is not nullhomotopic and does not have a sub-arc forming a nullhomotopic closed curve on $\gS$.
If $S=S^1$, then we regard the boundary $\partial S$ as $S$.

For two simple compact surfaces $(S_1,s_1), (S_2,s_2)$ on $\gS$ we define the intersection number of $(S_1,s_1), (S_2,s_2)$, denoted by $i(s_1,s_2)$, to be the number of contractible components of the fiber product $S_1\times_\gS S_2$ corresponding to $s_1,s_2$.
When we consider the intersection number, we always assume that $s_1$ and $s_2$ are ``transverse'', which means that the restriction of $s_1$ and $s_2$ to any components of their boundaries intersect transversely or virtually coincide if they intersect.
We say that two closed curves $c_1,c_2$ on $\gS$ virtually coincide if there exist a closed curve $c$ on $\gS$ and $m_1,m_2\in \NN$ such that $c_i$ equals $c^{m_i}$ up to reparametrization for $i=1,2$.
We define the intersection number of two homotopy classes $[s_1],[s_2]$ of simple compact surfaces to be the minimum of $i(s_1',s_2')$ taken over $s_1'\in [s_1]$ and $s_2'\in [s_2]$ that are transverse.
If $i(s_1,s_2)=i([s_1],[s_2])$, then we say that $s_1,s_2$ are in minimal position.

In the case that $\gS$ is a compact hyperbolic surface, we can see that for any simple compact surface $(S,s)$ (not a cylinder) on $\gS$ there exists a finitely generated subgroup $H$ of $\pi_1(\gS)$ such that the pair of the convex core $C_H$ and the natural projection from $C_H$ to $\gS$ induced by the universal covering belongs to the homotopy class $[s]$.
We also introduce the notion of an immersed bigon formed by $s_1,s_2$ and generalize the well-known bigon criterion for two closed curves on $\gS$ to two simple compact surfaces on $\gS$.

\begin{theoremintro}
Let $(S_1,s_1),(S_2,s_2)$ be transverse simple compact surfaces on a compact surface $\gS$.
If $s_1$ and $s_2$ do not form an immersed bigon, then $s_1,s_2$ are in minimal position.
If either $S_1$ or $S_2$ is $S^1$, then the converse is also true.
\end{theoremintro}

Let $\gS$ be a compact hyperbolic surface.
From the above theorem, we can see that for any two non-trivial finitely generated subgroups $H$ and $K$ of $\pi_1(\gS)$, $C_H$ and $C_K$ are in minimal position, that is,
\[ i(H,K)=i(C_H,C_K)=i([C_H], [C_K]).\]
Suppose that $C_H$ is a surface and $C_K$ is a closed geodesic on $\gS$, that is $H$ is non-cyclic and $K$ is infinite cyclic (see Example \ref{example:surface and closed geodesic intersection} and Figure \ref{intersection number formula} for detail).
Then we see that each contractible components of the fiber product $C_H \times_\gS C_K$ is an arc, whose two endpoints are on the boundary of $C_H$. Therefore we have the following formula:
\begin{equation}\label{intro:eq}
i(C_H,C_K)=\frac{1}{2}\sum_{c\in \partial C_H }i(c, C_K ) \tag{$\ast$}
\end{equation}
where $\partial C_H$ means the set of all boundary components of $C_H$.
For non-trivial $g\in \pi_1(\gS)$ we interpret $\partial C_{\langle g \rangle}$ as $\{ C_g , C_{g^{-1}} \}$, and then Equation (\ref{intro:eq}) also holds even if $H$ is also cyclic.

We generalize the above correspondence between $C_H$ and ``the half'' of $\partial C_H$ to the projection from $\SC(\gS)$ to $\GC (\gS)$.
For a non-cyclic finitely generated subgroup $H$ of $\pi_1(\gS)$, every component of the boundary of the convex core $C_H$ is a closed geodesic on $\gS$, and for each closed geodesic $c$ on $\gS$ we can obtain a counting geodesic current $\eta_c$ on $\GC (\gS )$, which equals $\eta_g$ for $g\in \pi _1(\gS )$ freely homotopic to $c$.
Then we can obtain a projection $\B$ from $\SC (\gS )$ onto $\GC (\gS )$:

\begin{theoremintro}
Let $\gS$ be a compact hyperbolic surface. There exists a unique continuous $\RRR$-linear map
\[ \B\: \SC (\gS )\rightarrow \GC (\gS )\]
such that for every non-cyclic finitely generated subgroup $H$ of $\pi_1(\gS)$ we have
\[ \B (\eta_H )=\frac{1}{2}\sum_{c\in \partial C_H} \eta_c \]
and the restriction of $\B$ to $\GC (\gS)$ is the identity map.
\end{theoremintro}

Note that if $\partial C_H$ is empty, then $B(\eta_H)$ is the zero measure.
For a non-trivial $g\in \pi_1(\gS)$, $B(\eta_g)=\frac{1}{2}(\eta_g +\eta_{g^{-1}})\ =\eta_g$.

Concerning the projection $\B$ we can obtain the following theorem, which generalizes Equation (\ref{intro:eq}) in the above:

\begin{theoremintro}
For any $\mu,\nu \in \SC (\gS )$ the following inequality holds:
\[ i_{\SC }(\mu,\nu )\leq i_{\GC}(\B (\mu ),\B (\nu )).\]
If either $\mu$ or $\nu$ belongs to $\GC (\gS )$, then the equality holds.
\end{theoremintro}

Recall that the hyperbolic structure $m$ of $\Sigma$ (or $m_d$ of the doubled surface $D\Sigma$ if $\Sigma$ has boundary) is associated to the geodesic current $L(m)$.
Note that in the case that $\Sigma $ has boundary, $\SC (\gS)$ is naturally embedded in $\SC (D\gS)$ (see the above of Theorem \ref{thmintro:include map} for detail).
From the above theorem, since $L(m)$ belongs to $\GC (\gS)$ for $m\in \mathcal{T}(\gS)$, we can generalize the $m$-length functional $\ell_m$ on $\GC (\gS)$ to the $m$-length functional $\ell_m$ on $\SC (\gS)$ by defining
\[ \ell_m(\mu ):=i_{\SC}( L(m), \mu )=i_{\GC}(\B (L(m)),\B (\mu ))=i_{\GC}(L(m), \B (\mu ))\]
for $\mu \in \SC (\gS)$. 
As a result, we obtain the following theorem:

\begin{theoremintro}
Let $\gS$ be a compact hyperbolic surface and let $m$ be its hyperbolic structure.
There exists a unique continuous $\RRR$-linear functional $\ell_m$ on $\SC (\Sigma)$ such that for any non-trivial finitely generated subgroup $H$ of $\pi_1(\gS)$ we have
\begin{equation}\label{intro:eq2}
\ell_m (\eta_H )=\frac{1}{2}\sum_{c\in \partial C_H} \ell_m (c) ,\tag{$\ast \ast$}
\end{equation}
where $\ell_m(c)$ is the $m$-length of $c$.
\end{theoremintro}

By using Equation (\ref{intro:eq}), we can give a more direct proof of Formula (\ref{intro:eq2}). Since $\GC(\gS)$ has a denseness property of rational geodesic currents, there is a sequence $c_n\eta_{g_n} \in \GC (\gS)\ (c_n>0, g_n\in \pi_1(\gS))$ converging the Liouville current $L(m)$.
Then we have
\begin{align*}
\ell_m(\eta_H)
&=i_{\SC}(L(m), \eta_H )=\lim_{n\rightarrow \infty }i_{\SC}(c_n\eta_{g_n}, \eta_H )\\
&=\lim_{n\rightarrow \infty} \frac{1}{2}\sum_{c\in \partial C_H}  i_{\GC}(c_n\eta_{g_n}, \eta_c )=\frac{1}{2}\sum_{c\in \partial C_H}i_{\GC}(L(m), \eta_c) 
=\frac{1}{2}\sum_{c\in \partial C_H}\ell_m (c).
\end{align*}

Assume that $\Sigma$ has no boundary. Bonahon's result with respect to the embedding of the Teichm\"uller space $\mathcal{T}(\gS)$ to $\GC (\gS)$ by sending a hyperbolic metric $m$ to the Liouville current corresponding to $m$ was extended to all negatively curved Riemannian metrics by Otal in \cite{Ota90}, to negatively curved cone metrics by Hersonsky and Paulin in \cite{HP97}, to singular flat metrics (or nonpositively
curved Euclidean cone metrics) coming from quadratic differentials by Duchin-Leininger-Rafi in \cite{DLR10}, and to all singular flat metrics by Bankovic-Leininger in \cite{BL18}.
For any such metric $m$ on $\gS$, we can obtain an associated geodesic current $L_m\in \GC (\gS)$, and for non-trivial $g\in \pi_1(\gS)$, the intersection number $i_{\GC}(L_m ,\eta_g)$ equals the $m$-length of $g$.  Note that the argument in the above direct proof of Formula (\ref{intro:eq2}) uses only the property that $i_{\GC}(L(m), \eta_c)$ equals the $m$-length of $c$ for any essential closed curve $c$ on $\gS$.
Therefore we have the following theorem:

\begin{theoremintro}
Let $\gS$ be a orientable closed surface of genus $\geq 2$ and let $m$ be a metric on $\gS$ mentioned in the above.
There exists a unique continuous $\RRR$-linear functional $\ell_m$ on $\SC (\Sigma)$ such that for any non-trivial finitely generated subgroup $H$ of $\pi_1(\gS)$ we have
\[
\ell_m (\eta_H )=\frac{1}{2}\sum_{c\in \partial C_H}\ell_m (c) ,
\]
where $\ell_m(c)$ is the $m$-length of $c$.
\end{theoremintro}

Consider two infinite quasi-convex subgroups $H$ and $J$ of an infinite hyperbolic group $G$. Assume that $J$ is a subgroup of $H$.
Then $\H (\partial J) $ is considered as a subspace of $\H (\partial H)$ and we have a continuous $\RRR$-linear map $\iota_J^H$ from $\SC (J)$ to $\SC (H)$ by defining
\[ \iota_J^H (\mu ):= \sum_{h J\in H/J} h_\ast (\mu )\]
for $\mu \in \SC (J)$, where $h_\ast (\mu )$ is the push-forward of $\mu$ by the self-homeomorphism $h$ on $\H (\partial H)$.
We write $\iota_H^{G}$ simply as $\iota_H$. For a quasi-convex subgroup $K$ of $H$ we denote by $\eta_K^H$ the counting subset current on $H$ corresponding to $K$. Then we can see that
\[ \iota_H (\eta_K^H )=\eta_K \in \SC (G ).\]
In general, $\iota_J^H$ is not injective but if $hJh^{-1}\cap J =\{ \id \}$ for any $hJ\ (\not=J) \in H/J$, then the map $\iota_J^H$ is injective.

For example, for a compact hyperbolic surface $\gS$ with boundary, its fundamental group $\pi_1(\gS)$ is considered as a subgroup of the fundamental group $\pi_1 (D\gS)$ of the doubled surface $D\gS$.
Then $\iota_{\pi_1(\gS)}^{\pi_1(D\gS)}\: \SC(\gS)\rightarrow \SC(D\gS)$ is injective (see Example \ref{exa:doubled surface inclusion} for detail). By using this map, $\SC(\gS)$ (or $\GC (\gS)$) is considered as a subspace of $\SC(D\gS)$ (or $\GC (D\gS)$ respectively).
In the context of geodesic currents, it is well-known that $\GC (\gS)$ can be considered as a subspace of $\GC (D\gS)$ as above, and the map $\iota_H$ introduced in the above is a natural extension.

We note that the map introduced in the above will play a fundamental role in proving the denseness property for a surface group.
Moreover, by using the map we can obtain the following theorem:

\begin{theoremintro}\label{thmintro:include map}
Let $H$ be a finite index subgroup of an infinite hyperbolic group $G$. If $H$ has the denseness property of rational subset currents, then $G$ also has the denseness property of rational subset currents.
\end{theoremintro} 

\subsection{Future study}
Consider the automorphism group $\mathrm{Aut}(G)$ of an infinite hyperbolic group $G$. The group $\mathrm{Aut}(G)$ acts on the boundary $\partial G$ continuously, which induces a continuous action on $\H (\partial G)$. Moreover, by considering the push-forward of subset currents by $\varphi \in \mathrm{Aut}(G)$ we have a continuous $\RRR$-linear action of $\mathrm{Aut}(G)$ on $\SC (G)$.
Since a subset current is $G$-invariant, the action of the inner automorphisms is trivial. Hence we have a continuous $\RRR$-linear action of the outer automorphism group $\mathrm{Out}(G)$ on $\SC (G)$.
For a quasi-convex subgroup $H$ of $G$ and $[\varphi] \in \mathrm{Out}(G)$ we have
\[ [\varphi ] (\eta_ H)=\eta_{\varphi (H)}\]
since $\varphi( \gL (H))=\gL (\varphi (H))$.
The action of $\mathrm{Out}(G)$ on $\SC (G)$ fixes the point $\eta_G$.
We also note that $\GC (G)$ is an $\mathrm{Out}(G)$-invariant subspace of $\SC (G)$.

Define an equivalence relation $\sim$ on $\SC(G)\setminus \{ 0 \}$ as follows: for $\mu_1,\mu_2\in \SC(G)\setminus \{ 0 \}$, $\mu_1 \sim \mu_2$ if $\mu_1=c\mu_2$ for some $c>0$.
Set $\mathbb{P} \SC (G):=(\SC(G)\setminus \{ 0 \})/\sim$ and endow $\mathbb{P} \SC (G)$ with the quotient topology. Then the action of $\mathrm{Out}(G)$ on $\SC (G)$ induces the continuous action of $\mathrm{Out}(G)$ on $\mathbb{P} \SC (G)$, which can be thought of as a generalization of the action of $\mathrm{Out}(G)$ on the set of all conjugacy classes of non-trivial quasi-convex subgroups of $G$.
Note that $\mathbb{P} \SC (G)$ is compact since $\SC (G)$ is locally compact.

From the Dehn-Nielsen-Baer theorem, the mapping class group $\mathrm{MCG}(\gS)$ of a closed surface $\gS$ is isomorphic to a $2$-index subgroup of $\mathrm{Out}(\pi_1(\gS))$.


For a closed hyperbolic surface $\gS$, we can see that our maps $\rk$, $i_{\SC}$ and $\N$ on $\SC (\gS )$ are $\mathrm{Out}(\pi_1 (\gS ))$-invariant, especially, $\mathrm{MCG}(\gS)$-invariant.
The projection $\B\: \SC (\gS)\rightarrow \GC (\gS)$ is $\mathrm{Out}(\pi_1 (\gS ))$-equivariant.
We are interested in investigating $\mathrm{MCG}(\gS)$ by using $\SC (\gS ), \mathbb{P}\SC (G)$ and functionals on $\SC (\gS)$.
We remark that if $\gS$ has boundary, then $\rk$ and $\N$ are $\mathrm{Out}(\pi_1 (\gS ))$-invariant but $i_{\SC}$ is not $\mathrm{Out}(\pi_1 (\gS ))$-invariant, and $\B$ is not $\mathrm{Out}(\pi_1 (\gS ))$-equivariant.
In the case of a free group $F$ of finite rank, Kapovich and Nagnibeda made an observation on the action of $\mathrm{Out}(F)$ on $\mathbb{P}\SC (F)$ and presented two problems in \cite[Subsection 10.6]{KN13}.


\subsection{Organization of article}
In Section \ref{sec:subset currents on hyperbolic groups}, we will introduce subset currents on a hyperbolic group $G$ and develop a general theory on the space $\SC(G)$. We also give a short introduction to the background of measure theory related to subset currents.

In Section \ref{sec:Volume functionals for Kleinian groups}, we will prove the existence of the volume functional $\mathrm{Vol}$ on $\SC (\gG )$ for a non-trivial torsion-free convex-cocompact Kleinian group $\gG$ on $\HH^n$ for $n\geq 2$ (see Theorem \ref{thm:volume functional}). As a corollary, we obtain the reduced rank functional $\rk$ on $\SC (\gS)$ for a compact hyperbolic surface $\gS$ (see Corollary \ref{cor:reduced rank functional}).

In Section \ref{sec:relation between subgroups}, we will introduce the natural continuous $\RRR$-linear map $\iota_H$ from $\SC (H)$ to $\SC (G)$ for a quasi-convex subgroup $H$ of a hyperbolic group $G$. By using the map $\iota_H$ we prove that if a hyperbolic group $G$ has the denseness property of rational subset currents, then the finite index extension of $G$ also has the denseness property of rational subset currents (see Theorem \ref{thm:denseness and finite index}). We present a method of extending a functional on $\SC(H)$ to a functional on $\SC(G)$ if $H$ is a finite index subgroup of $G$ in Subsection \ref{subsec:finite index extension}.

In Section \ref{sec:intersection number of subset currents}, first, we will review several facts on the intersection number of two closed curves on a compact surface $\gS$, and then introduce the intersection number of two simple compact surfaces on $\gS$. We prove the bigon criterion for two simple compact surfaces on $\gS$ as a generalization of the bigon criterion for two (simple) closed curves on $\gS$ (see Theorem \ref{thm:bigon criterion 3}). Finally, we prove the existence of the intersection number $i_{\SC }$ on $\SC (\gS)$ (see Theorem \ref{thm:intersection number of subset currents}). During the proof, we introduce some new techniques for proving the continuity of a functional on $\SC (\gS)$. In the beginning of Subsection \ref{conti ext of int number}, we review several facts on geodesic currents on hyperbolic groups.

In Section \ref{sec:an intersection functional}, we will introduce the product $\N$ of two finitely generated subgroups of $\pi_1(\gS)$ for a compact hyperbolic surface $\gS$. Our proof of the bigon criterion for two simple compact surfaces on $\gS$ gives a geometric characterization of $\N$ and also gives us an idea for extending $\N$ to an $\RRR$-bilinear functional on $\SC (\gS )$. Our proof of the continuity of $\N$ on $\SC (\gS )$ is partially based on the proof of the continuity of $i_{\SC }$.

In Section \ref{sec: projection B}, we will prove the existence of the continuous $\RRR$-linear projection $\B$ from $\SC (\gS )$ onto $\GC (\gS )$ for a compact hyperbolic surface $\gS$ (see Theorem \ref{thm:projection from sc to gc}).
By using the projection $\B$ and the denseness property of rational subset currents for $\pi_1(\gS)$, we obtain the inequality on the intersection number on $\SC (\gS )$ and $\GC (\gS )$ (see Theorem \ref{thm:projection and inequality}). As a corollary, we obtain the $m$-length functional $\ell_m$ on $\SC (\gS)$ for an element $m$ of the Teichm\"uller space of $\gS$.

In Section \ref{sec:denseness property of rational subset currents}, our goal is proving that surface groups have the denseness property of rational subset currents (see Theorem \ref{thm:surface group has denseness property}). In Subsection \ref{subsec: denseness property of free groups}, we will give a proof of the denseness property for a free group $F$ of finite rank based on the proof by Kapovich in \cite{Kap17}. In the proof we give some new ideas for understanding the denseness property such as SC-graphs on $F$. In Subsection \ref{subsec:approximation by a sequence of sugroups}, we will give a sequence of finitely generated subgroups $H_n$ of $F$ so that the union of the image of $\SC (H_n)$ by the natural map $\iota_{H_n}$ taken over all $n$ is a dense subset of $\SC (F)$. Finally, in Subsection \ref{subsec:denseness property of surface groups}, we will prove that for a closed hyperbolic surface $\gS$, the surface group $\pi_1(\gS)$ has the denseness property. Several methods for this proof have been introduced in Subsection \ref{subsec: denseness property of free groups} and \ref{subsec:approximation by a sequence of sugroups} in advance but also we generalize some of those methods. Especially, we use a sequence of finitely generated subgroups of $\pi_1(\gS )$, each of which is a free group of finite rank.
A lot of constants are involved in the proof, and we need to be careful of the relation between constants.
We note that we will use the denseness property for surface groups in several sections before Section \ref{sec:denseness property of rational subset currents}.

\subsection{Acknowledgements}
I would like to show my greatest appreciation to Prof. Katsuhiko Matsuzaki who provided helpful comments and suggestions.
The author also would like to thank the referee for her/his careful reading of the manuscript and valuable comments.
The author is partially supported by JSPS KAKENHI Grant Number JP16J02814 and JP19K14539.

\section{Subset currents on hyperbolic groups}\label{sec:subset currents on hyperbolic groups}
First we define the hyperspace of a topological space, which consists of compact subsets. Later, we consider only the case where the topological space is the (Gromov) boundary of a hyperbolic group, which is compact metrizable. The hyperspace is used for considering limit sets of subgroups of the hyperbolic group.

\begin{definition}[See {\cite[Subsection 4.F]{Kec95}}]
Let $X$ be a topological space. We will denote by $\hat{\H}(X)$ the set of all compact subsets of $X$ including the empty set $\emptyset$ with the \ti{Vietoris topology}, which is generated by the sets of the form
\[ \{ S\in \hat{\H}(X)\mid S\subset U\} \text{ and }\{ S\in \hat{\H}(X)\mid S\cap U\not= \emptyset \}\]
for an open subset $U \subset X$. We call $\hat{\H}(X)$ the \ti{hyperspace} of $X$ consisting of compact sets.
\end{definition}

\begin{theorem}[See {\cite[Theorem 4.26]{Kec95}}]\label{thm:hyperspace is compact metrizable}
If $X$ is a compact metrizable space, then so is $\hat{\H}(X)$. In particular, $\hat{\H}(X)$ is separable.
\end{theorem}

\subsection{Space of subset currents on hyperbolic group}\label{subsec:space of subset currents}

Let $G$ be an infinite hyperbolic group. We do not consider the case that $G$ is a finite group throughout this article. Fix a finite generating set of $G$ and denote by $\Cay (G)$ the Cayley graph of $G$ with respect to the generating set. When we want to emphasize a generating set $A$ of $G$, we will denote by $\Cay (G,A)$ the Cayley graph of $G$ with respect to $A$. Recall that the vertex set $V(\Cay (G,A))=G$, the edge set $E(\Cay (G,A))=G\times A$ and $(g,a)\in E(\Cay (G,A))$ is an edge connecting $g$ to $ga$. Then we have the natural action of $G$ on $X$ from the left. We consider a connected graph as a metric space by endowing the graph with the path metric such that every edge has length $1$.

Since the boundary $\partial G$ of $G$ is compact metrizable,
the space $\hat{\H}(\partial G)$ is compact metrizable by Theorem \ref{thm:hyperspace is compact metrizable}. Now, we consider an open subspace
\[ \H(\partial G):=\{ S \in \hat{\H}(\partial G)\mid \#S \geq 2\} \]
of $\hat{\H}(\partial G)$. Then $\H(\partial G)$ is a locally compact separable metrizable space.

Let $d_{\partial G}$ be a metric on $\partial G$ that is compatible with its topology. Then we can define the Hausdorff distance $d_{\mathrm{Haus}}$ on $\H(\partial G)$ as
\[ d_{\mathrm{Haus}}(S_1,S_2):=\max\, \{ \max_{s\in S_1}d_{\partial G}(s,S_2),\max_{s\in S_2}d_{\partial G}(S_1,s)\} \quad (S_1,S_2\in \H(\partial G) ).\]
We can see that the Hausdorff distance is compatible with the subspace topology on $\H(\partial G)$ given by the Vietoris topology.
When we consider the topology of $\H (\partial G)$, the Hausdorff distance $d_{\mathrm{Haus}}$ is convenient.
Note that $d_{\mathrm{Haus}}$ actually can be defined on $\hat{\H}(\partial G)\setminus \{ \emptyset \}$.

Since $G$ acts on $\partial G$ continuously, the action extends to the continuous action on $\H (\partial G)$.
Now, we introduce the main object, subset currents.

\begin{definition}[Subset currents on $G$]
A \ti{subset current} on $G$ is a $G$-invariant locally finite Borel measure on $\H (\partial G)$.
The space of subset currents on $G$ is denoted by $\SC (G)$. We give $\SC (G)$ the weak-$\ast$ topology. (See Subsection \ref{subsec:measure} for the definitions of measure-theoretic terminology.)
\end{definition}

\begin{remark}
For a finite hyperbolic group $G$, since the boundary $\partial G$ is empty, we define $\SC (G)$ to be the set consisting of the zero measure. In the case that $G$ is an infinite cyclic group, the boundary $\partial G$ consists of two points and $G$ acts on $\partial G$ trivially. Hence $\SC (G)$ is the $\RRR$-linear space generated by the Dirac measure $\delta_{\partial G}$ at $\partial G$.
\end{remark}

For $S\in \H(\partial G)$ the \ti{weak convex hull} $WC(S)\subset \Cay(G)$ of $S$ is the union of all geodesic lines connecting two points of $S$. A geodesic line in a metric space is an isometric embedding of $\RR$ into the metric space.
For each vertex $g \in V(\mathrm{Cay}(G))$ we consider a subset
\[ A_{g}:=\{ S\in \H(\partial G) \mid WC(S)\ni g \}.\]
Since for any $g\in G$ and $S\in \H(\partial G)$ we have $gWC(S)=WC(g(S))$, 
\[ G(A_\id )=\bigcup_{g\in G}gA_\id =\bigcup_{g\in G}A_{g}=\H(\partial G).\]

\begin{lemma}\label{lem: A g is compact subset}
The set $A_{g}$ is a compact subset of $\H (\partial G)$ for every $g\in G$.
\end{lemma}
\begin{proof}
Recall that the space $\hat{\H}(\partial G)$ is compact. Therefore, it suffices to show that the set $A_\id$ is closed in $\hat{\H}(\partial G)$.
Consider a sequence $\{S_n\}\subset A_\id$ converging to $S\in \hat{\H}(\partial G)$. It is clear that $S\not =\emptyset$.
For each $n\in \NN$ take $\xi_n, \zeta_n \in S_n $ such that there exists a geodesic line $\gamma_n$ containing $\id$ and joining $\xi_n$ to $\zeta_n$.
We can take convergent subsequences $\{\xi_{k_n}\}$ and $\{\zeta_{k_n}\}$ by the compactness of $\partial G$.
Since $S_n$ converges to $S$ in the Hausdorff distance $d_{\mathrm{Haus}}$, the sequences $\{\xi_{k_n}\}$ and $\{\zeta_{k_n}\}$ converge to $\xi ,\zeta \in S$, respectively. From the Ascoli-Arzel\`a theorem there exists a subsequence of $\gamma_{k_n}$ that converges uniformly on compact subsets to a geodesic line $\gamma$ joining $\xi$ to $\zeta$. Since each $\gamma_{k_n}$ contains the vertex $\id$, so is $\gamma$.
Therefore, $WC(S)$ contains the vertex $\id$, which proves our claim.
\end{proof}

From the above lemma, we can see that $G$ acts on $\H(\partial G)$ cocompactly.
By applying Theorem \ref{thm:M_G(X)} in Subsection \ref{subsec:measure} to $\SC(G)$, we have the following theorem.

\begin{theorem}\label{thm:SC(G)}
The space $\SC (G)$ is a locally compact, separable and completely metrizable space.
\end{theorem}

We assume some background knowledge on the properties of limit sets of subgroups of hyperbolic groups. 

For each subgroup $H$ of $G$ we have the limit set $\gL (H) \subset \partial G$, which is the set of accumulation points of $H$ in $G\cup \partial G$. We usually consider the case that $\gL(H)\not =\emptyset$, which implies $\gL(H)\in \H(\partial G)$. We define a measure $\eta_H$ on $\H(\partial G)$ as
\[ \eta_H:= \sum_{gH\in G/H} \delta_{g\gL(H)},\]
where $\delta_{g\gL(H)}$ is the Dirac measure at $g\gL(H)\in \H(\partial G)$. It is easy to check that $\eta_H$ is $G$-invariant.
In the case that $H$ is a finite group, especially the trivial group $\{ \id \}$, the limit set $\gL (H)$ is empty, and so we define $\eta_H$ to be the zero measure.

A subgroup $H$ of $G$ is called \ti{quasi-convex} if $H$ is a quasi-convex subset of $\mathrm{Cay }(G)$, that is, there exists $k>0$ such that any geodesic connecting two points of $H$ is included in the $k$-neighborhood of $H$.
It is known that a subgroup $H$ of $G$ is quasi-convex if and only if $H$ acts on the weak convex hull $WC(\gL(H))$ cocompactly (see \cite{Swe01}).
The following lemma is a generalization of \cite[Lemma 4.4]{KN13} in the case of hyperbolic groups.

\begin{lemma}\label{lem:continug subset currents}
Let $H$ be a subgroup of $G$. The measure $\eta_H$ is locally finite if and only if $H$ is quasi-convex.
\end{lemma}
\begin{proof}
We generalize the compact subset $A_\id \subset \H(\partial G)$.
For $r\geq 0$ we define $A(\id ,r )$ to be a subset of $\H (\partial G)$ consisting of $S\in \H (\partial G)$ such that $WC(S)$ intersects the closed ball $B(\id, r)$ centered at $\id $ with radius $r$.
Note that $A(\id, 0)=A_\id$, and $A_\id$ can be not open.
We can see that if $r$ is sufficiently large compared with the hyperbolic constant of $\Cay(G)$, then the interior $\mathrm{Int}(A(\id,r))$ includes $A_\id$, and so
\[ G(\mathrm{Int}(A(\id ,r)))=\H (\partial G).\]
Therefore, any compact subset of $\H(\partial (G))$ is covered by a finite union of $g(A(\id, r))$ $(g\in G)$.
Moreover, any compact subset of $\H(\partial (G))$ is covered by a finite union of $A_g \ (g\in G)$ since we have
\[ A(\id ,r) =\bigcup_{g\in G\cap B(\id,r )} A_g ,\]
which is a finite union. As a result, a $G$-invariant measure $\mu$ on $\H (\partial G)$ is locally finite if and only if $\mu (A_\id)$ is finite.

For the measure $\eta_H$ we have
\begin{align*}
\eta_H (A_\id )
&= \# \{ g H\in G/H\mid WC(g\gL (H))\ni \id \}\\
&= \# \{ g H\in G/H\mid g WC(\gL (H))\ni \id \}\\
&= \# \{ g H\in G/H\mid WC(\gL (H))\ni g^{-1} \}.
\end{align*}
For $g_1H, g_2H\in G/H$ with $g_1H\not =g_2H$, there is no $h\in H$ such that $hg_1^{-1}=g_2^{-1}$. Therefore, $\eta_H(A_\id)$ equals the number of vertices of the quotient graph $H\backslash WC(\gL(H))$ of $WC(\gL(H))$ by $H$.
Hence, $\eta_H(A_\id)$ is finite if and only if $H$ acts on $WC(\gL(H))$ cocompactly, which completes the proof.
\end{proof}

In general, for any $S\in \H(\partial G)$ we can obtain a $G$-invariant Borel measure (not necessarily locally finite)
\[ \eta_S:=\sum_{gH\in G/H}\delta_{gS} \]
on $\H(\partial G)$, where $H:=\mathrm{Stab}_G(S)=\{ g\in G\mid g(S)=S \}$, the stabilizer of $S$.
For any $G$-invariant Borel measure $\mu$ on $\H(\partial G)$, if $\mu$ has an atom $S$, that is, $\mu (\{ S\})>0$, then
$\mu(E)\geq \mu (\{ S\}) \eta_S(E)$
for every Borel subset $E\subset \H(\partial G)$. Therefore, if $\mu$ is locally finite, then so is $\eta_S$.

\begin{theorem}\label{thm:counting sc}
Let $S\in \H(\partial G)$. The measure $\eta_S$ is locally finite if and only if $H:=\mathrm{Stab}_G(S)$ is quasi-convex and $S=\gL(H)$.
In particular, if a subset current $\mu\in \SC(G)$ has an atom $S$, then $H$ is quasi-convex and $S=\gL(H)$.
\end{theorem}
\begin{proof}
The ``if part'' follows by Lemma \ref{lem:continug subset currents}. We prove the ``only if" part.
Assume that $\eta_S$ is locally finite. From the proof of Lemma \ref{lem:continug subset currents}, $\eta_S(A_\id)$ equals the number of vertices of the quotient graph $H\backslash WC(S)$, which implies that $H$ acts on $WC(S)$ cocompactly. Note that for every $\xi\in S$ there exists a sequence of $WC(S)$ converging to $\xi$ and we can take the sequence from $H(x)$ for some $x\in WC(S)$. Therefore $S=\gL(H)$ and $H$ is quasi-convex.
\end{proof}

\begin{definition}
We call $\eta_H$ the \ti{counting subset current} for a quasi-convex subgroup $H$ of $G$.
A subset current $\mu\in \SC(G)$ is called \ti{rational} if there exists a quasi-convex subgroup $H$ of $G$ and $c\in \RRR$ such that $\mu =c\eta_H$.
We denote by $\SC_r(G)$ the set of all rational subset currents on $G$.
\end{definition}

Counting subset currents have the following properties:
\begin{proposition}\label{prop:property}
For two quasi-convex subgroups $H_1,H_2$ of $G$,
\begin{enumerate}
\item if $H_1$ is a $k$-index subgroup of $H_2$, then $\eta_{H_1}=k\eta_{H_2}$;
\item if $H_1$ is conjugate to $H_2$, then $\eta_{H_1}=\eta_{H_2}$.
\end{enumerate}
\end{proposition}
\begin{proof}
(1). Assume that $H_1$ is a $k$-index subgroup of $H_2$.
Note that $\gL (H_1)=\gL (H_2)$.
Take a complete system of representatives $R$ of $G/H_2$.
Then a map sending $(g,hH_1)\in R\times H_2/H_1$ to $ghH_1\in G/H_1$ is a bijective map.
Hence
\begin{align*}
\eta_{H_1}=\sum_{gH_1\in G/H_1}\delta_{g\gL (H_1)}
&=\sum_{g\in R}\sum_{hH_1\in H_2/H_1}\delta_{gh\gL (H_1)}\\
&=\sum_{g\in R}k\delta_{g\gL (H_2)}=k\eta_{H_2}.
\end{align*}

(2). Assume that $H_1=g_0H_2g_0^{-1}$ for $g_0\in G$.
Note that $\gL(H_1)=g_0\gL (H_2)$.
Take a complete system of representatives $R$ of $G/H_2$.
Then $g_0Rg_0^{-1}$ is a complete system of representatives of $G/H_1$ since
\[ G=\bigsqcup_{g\in R} gH_2=\bigsqcup_{g\in R}g_0gH_2g_0^{-1}=\bigsqcup_{g\in R}(g_0gg_0^{-1})H_1.\]
Hence
\[ \eta_{H_1}=\sum_{g\in R}\delta_{g_0gg_0^{-1}\gL (H_1)}=\sum_{g\in R}\delta_{g_0g\gL (H_2)}=\sum_{gH_2\in G/H_2}\delta_{g\gL (H_2)}=\eta_{H_2},\]
which is the required equation.
\end{proof}

Kapovich and Nagnibeda \cite{KN13} proved the following theorem, which played a fundamental role in their study of the space of subset currents on a free group. Kapovich \cite{Kap17} gave another proof to the following theorem.

\begin{theorem}[See {\cite[Theorem 5.8]{KN13}} and \cite{Kap17}]\label{thm:F_N dense}
For a free group $F$ of finite rank, the set $\SC_r(F)$ of all rational subset currents on $F$ is a dense subset of $\SC(F)$.
\end{theorem}

Note that a subgroup $H$ of $F$ is quasi-convex if and only if $H$ is finitely generated.
By Proposition \ref{prop:property} (2) and Theorem \ref{thm:F_N dense},
we can thought of $\SC(F)$ as a measure-theoretic completion of the set of all conjugacy classes of finitely generated subgroups of $F$.

We say that an infinite hyperbolic group $G$ has \ti{the denseness property of rational subset currents} if the set of all rational subset currents on $G$ is a dense subset of $\SC(G)$.
Recall that the space $\SC(G)$ is separable for any hyperbolic group $G$.
If $G$ has the denseness property of rational subset currents, then we have a concrete countable dense subset 
\[ \{ q\eta_H\mid q\in \QQ_{\geq 0}\ \text{$H\leq G:$ a quasi-convex subgroup }\}\]
of $\SC(G)$.

In Subsection \ref{subsec:denseness property of surface groups}, we will prove that surface groups have the denseness property of rational subset currents (see Theorem \ref{thm:surface group has denseness property}).
In Subsection \ref{subsec: iota J H}, we will prove that for a hyperbolic group $G$ and a finite index subgroup $H$ of $G$, if $H$ has the denseness property of rational subset currents, then $G$ also has the denseness property of rational subset currents (see Theorem \ref{thm:denseness and finite index}).

From the above, it is natural to propose the following problem. 
\begin{problem}\label{problem: denseness property for hyperbolic group}
Does any infinite hyperbolic group $G$ have the denseness property of rational subset currents?
\end{problem}

Note that from the viewpoint of the application, it is sufficient to see that the $\RRR$-linear subspace $\mathrm{Span}(\SC_r(G))$ generated by $\SC_r(G)$ is a dense subset of $\SC (G)$. In the case that $G$ is a free group $F$ of finite rank, Kapovich-Nagnibeda \cite{KN13} first proved that $\SC_r(F)$ is a dense subset of $\mathrm{Span}(\SC_r(F))$, and then they proved that $\mathrm{Span}(\SC_r(F))$ is a dense subset of $\SC(F)$, which implies that $\SC_r(F)$ is a dense subset of $\SC(F)$.
However, for a general infinite hyperbolic group $G$, we do not know whether $\SC_r(G)$ is a dense subset of $\mathrm{Span}(\SC_r(G))$ or not.

Let $G$ be an infinite hyperbolic group with denseness property of rational subset currents.
The denseness property of rational subset currents has a lot of application.
For example, if we have an $\RRR$-linear functional on $\SC (G)$ that is a continuous extension of an invariant of a quasi-convex subgroup of $G$, then we can see that the functional is unique.
We will use this argument frequently in this article for the case that $G$ is the fundamental group of a compact hyperbolic surface.

In addition, if we have a continuous $\RRR$-linear functional $\phi$ on $\mathrm{Span}(\SC_r(G))$, then $\phi$ is uniquely extended to a continuous $\RRR$-linear functional on $\SC(G)$ by the Hahn-Banach theorem, which we will use in the proof of Proposition \ref{prop: angle bilinear functional}.

\subsection{Measure theory background}\label{subsec:measure}
In this subsection, we give an introduction to the space of measures. Most of the contents are well-known in the measure theory (see \cite[Section 8]{Bog07} for more detail). First, we consider the space of locally finite measures with weak-$\ast$ topology, and then we consider a group action additionally.

Let $(X,d)$ be a locally compact second countable metric space. We consider the space $M (X)$ of locally finite Borel measures on $X$ in this subsection. Our goal is to see that the space $M(X)$ with the weak-$\ast$ topology is second countable and completely metrizable.

First we recall some definitions from the measure theory.

\begin{definition}
A Borel measure $\mu$ on $X$ is called \ti{locally finite} if $\mu (K)$ is finite for any compact subset $K\subset X$.
A Borel measure $\mu$ on $X$ is called \ti{regular} if for any Borel subset $E\subset X$,
\[ \mu (E) =\mathrm{inf}\{ \mu (O)\mid O\subset X\: \tn{open and } E \subset O\} \]
and if for any Borel subset $E\subset X$ with $\mu (E)<\infty$,
\[ \mu (E) =\sup \{ \mu (K)\mid K\subset X\: \tn{compact and } E \supset K\} .\]
\end{definition}

Since $X$ is a locally compact second countable metric space, locally finite Borel measures are regular (see \cite[2.18 Theorem]{Rud86}).

\begin{definition}
Let $C_c (X)$ be the space of compactly supported continuous functions from $X$ to $\RR$ with
the topology of uniform convergence on compact sets. This means that $f_n$ converges to $f$ in $C_c(X)$ if there exists a compact subset $K\subset X$ such that $\mathrm{supp} f_n,\ \mathrm{supp} f \subset K$, and $f_n$ converges to $f$ uniformly.
With this topology, for any $\mu \in M(X)$ the functional 
\[ f\in C_c(X) \mapsto \int f d\mu \]
is continuous. We often denote $\int fd\mu$ briefly by $\mu (f)$.

A sequence $\{ \mu_n\} \subset M(X)$ converges to $\mu \in M(X)$ in the \ti{weak-$\ast$ topology} if and only if for any $f\in C_c(X)$ we have 
$\mu_n(f) \rightarrow \mu (f) \ (n\rightarrow \infty )$.
\end{definition}

\begin{proposition}\label{prop:c c(x) is separable}
The space $C_c(X)$ is separable.
\end{proposition}
\begin{proof}
If $X$ is compact, then we can see that $C_c(X)=C(X)$ is separable from the Stone-Weierstrass Theorem. In a general case, we take a sequence of compact subsets $K_n\subset X\ (n\in \NN )$ satisfying the condition that
\begin{equation}
X=\bigcup_{n\in \NN } K_n \ \textnormal{and}\ K_n\subset \mathrm{Int}(K_{n+1})\ \tn{for any }n\in \NN. \tag{$\ast $}\label{eq:cpt subsets}
\end{equation}
This implies $ X=\bigcup_{n\in \NN } \mathrm{Int}(K_n)$.
Then we have
\[ C_c(X)=\bigcup_{n\in \NN} \{ f\in C_c(X)\mid \mathrm{supp}f\subset K_n\}.\]
Since $\{ f\in C_c(X)\mid \mathrm{supp}f\subset K_n\} \subset C(K_n)$ is separable for every $n\in \NN$, so is $C_c(X)$.
\end{proof}

Now, we define a metric $d_M$ on $M(X)$ as follows. Fix a sequence of compact subsets $K_n\subset X\ (n\in \NN )$ satisfying the above Condition (\ref{eq:cpt subsets}).
Then take a countable dense subset $C=\{ \phi_n \mid n\in \NN\} \subset C_c(X)$ such that $C$ contains a compactly supported continuous function $\chi_n$ with $\chi_n \geq 0$ and $\chi_n(x)=1$ for any $x\in K_n$ for each $n\in \NN$, which implies that $\mu (\chi_n)\geq \mu (K_n)$ for any $\mu \in M(X)$.
Moreover, from the proof of Proposition \ref{prop:c c(x) is separable}, we can assume that for any $f\in C_c(X)$ with $\mathrm{supp} f\subset K_{n}$ for some $n\in \NN$ there exists a sequence $\{ f_j\}$ of $C$ such that $\{ f_j\}$ converges to $f$ and $\mathrm{supp} f_j\subset K_n$ for each $j$.
For $\mu,\nu\in M(X)$ we define
\[ d_M(\mu ,\nu ):=\sum_{n\in \NN } 2^{-n}\min \, \{ | \mu(\phi_n)-\nu(\phi_n) |, 1\}.\]
\begin{theorem}
The metric $d_M$ on $M(X)$ is compatible with the weak-$\ast$ topology.
\end{theorem}
\begin{proof}
For $\mu_n ,\mu\in M(X)\ (n\in \NN)$, it is easy to see that $d_M(\mu_n ,\mu )\rightarrow 0 \ (n\rightarrow \infty )$ if and only if $\mu_n(f)\rightarrow \mu(f)\ (n\rightarrow \infty)$ for any $f\in C$. Assume that $d_M(\mu_n ,\mu )\rightarrow 0 \ (n\rightarrow \infty )$.
It is sufficient to prove that $\mu_n(f)\rightarrow \mu(f)$ for any $f\in C_c(X)$. We can take $k \in \NN$ such that there exists a sequence $\{ f_j\} \subset C$ converging to $f$ uniformly and $\mathrm{supp} f_j, \mathrm{supp}f \subset K_k$.
Since $\mu_n (\chi_k)\rightarrow \mu(\chi_k)\ (n\rightarrow \infty)$, the sequence $\{ \mu_n(\chi_k )\}_{n \in \NN}$ is bounded and so is $\{ \mu_n(K_k)\}_{n\in \NN}$. Therefore, for any $\varepsilon>0$ and a sufficiently large $j\in \NN$ we have
\begin{align*}
|\mu_n(f)-\mu(f)|
\leq &|\mu_n(f)-\mu_n(f_j) |+|\mu_n (f_j)-\mu(f_j)| +|\mu (f_j)-\mu (f)|\\
\leq &\| f-f_j\| _{{\infty}}\mu_n (K_k)+|\mu_n (f_j)-\mu(f_j)|+\| f-f_j\|_{{\infty}}\mu (K_k)\\
\leq &2\varepsilon +|\mu_n (f_j)-\mu(f_j)|,
\end{align*}
where $\| \cdot \|_{\infty}$ is the sup norm. Hence if $n\in \NN$ is sufficiently large, then $ |\mu_n(f)-\mu(f)|\leq 3\varepsilon $.
This completes the proof.
\end{proof}

\begin{theorem}
The metric space $(M(X),d_M)$ is complete.
\end{theorem}

\begin{proof}Let $\{\mu_n\}$ be a Cauchy sequence in $(M(X),d_M)$. For any $f\in C$ we can see that $\{ \mu_n (f)\}$ is also a Cauchy sequence. Since $\RR$ is complete, we obtain a map
\[ \Phi \: C\rightarrow \RR ;\ f\mapsto \lim_{n\rightarrow \infty} \mu _n(f).\]
Then we extend the map $\Phi$ to a positive linear functional from $C_c(X)$ to $\RR$ by using the denseness of $C$ in $C_c(X)$.
Finally from the Riesz-Markov-Kakutani representation theorem, there exists a unique locally finite measure $\mu$ such that we have
\[ \Phi (f)=\int f d\mu \ \tn{for any }f\in C_c(X).\]
The measure $\mu$ is the limit of the Cauchy sequence $\{ \mu_n\}$.
\end{proof}

To see that $M(X)$ is separable, we decompose $X$ into ``small'' subsets by using the condition that $X$ is a locally compact second countable metric space, whose property is similar to that of the Euclidean space. Note that on a metric space being separable is equivalent to being second countable, and we use both words according to each situation.

For each $n\in \NN$ we take a family of Borel subsets $\{ E_{\lambda}^n\} _{\lambda \in \Lambda_n}$ satisfying the following conditions:
\begin{enumerate}
\item $X$ is a disjoint union of $\{ E_{\lambda}^n\} _{\lambda \in \Lambda_n}$;
\item for any compact subset $K\subset X$ only finitely many $E_{\lambda}^n$ intersect $K$, which in particular implies that $\Lambda _n$ is countable;
\item the diameter of $E_{\lambda}^n$ is smaller than $1/n$.
\end{enumerate}
For each $E_{\lambda}^n$ we fix $p_{\lambda}^n\in E_{\lambda}^n$. Since $\Lambda_n$ is countable for each $n\in \NN$, the set $P:=\{ p_{\lambda }^n\mid n\in \NN,\ \lambda \in \Lambda_n\}$ is also countable. For each $(p ,q)\in P\times \QQ_{\geq 0}$ we consider a measure $q\delta_p \in M(X)$, where $\delta_p$ is the Dirac measure at $p$, that is, for any Borel subset $E\subset X$, if $E\ni p$, then $\delta_p (E)=1$; if $E\not \ni p$, then $\delta_p(E)=0$.
Now, set
\[ D:= \bigcup_{k\in \NN} \left\{ \sum_{i=1}^k\ q_i\delta_{p_i} \Big| \ (p_i, q_i)\in P \times \QQ_{\geq 0}\right\},\]
which is countable.

\begin{theorem}
The set $D$ is a dense subset of $M(X)$. Hence $M(X)$ is separable.
\end{theorem}
\begin{proof}
Take an arbitrary $\mu \in M(X)$. For each $n\in \NN$ set
\[ \mu_n :=\sum_{\lambda \in \Lambda_n} \mu (E_\lambda ^n )\delta_{p_\lambda ^n}.\]
Then take $q_\lambda^n\in \QQ_{\geq 0}$ such that
\[ \sum_{\lambda \in \Lambda_n}| \mu (E_\lambda ^n )-q_{\lambda}^n| <\frac{1}{n}\]
and set
\[ \mu_n':=\sum_{\lambda \in \Lambda_n} q_\lambda ^n\delta_{p_\lambda ^n}.\]
Next, recall that the sequence of compact subsets $\{ K_n\}$ of $X$ satisfies the Condition (\ref{eq:cpt subsets}).
For each $n\in \NN$ the restriction of $\mu_n'$ to $K_n$, denoted by $\nu_n$, is contained in $D$ since only finitely many $E_\lambda^n\ (\lambda \in \Lambda_n)$ intersect $K_n$.

Now, we prove that the sequence $\{ \nu_n \}$ converges to $\mu$.
Take an arbitrary $f\in C_c(X)$. For a sufficiently large $n\in \NN$ the support of $f$ is included in $K_n$, and so
\[ \int f d\nu_n =\int f d\mu_n'=\sum _{\lambda \in \Lambda_n}q_\lambda^n f(p_\lambda ^n).\]
Hence
\begin{align*}
| \nu_n(f)-\mu_n(f)| 
\leq & \left| \sum_{\lambda \in \Lambda_n}\Bigl(q_\lambda^n f(p_\lambda^n)-\mu(E_\lambda^n)f(p_\lambda^n)\Bigr) \right|\\
\leq & \| f\|_{\infty} \sum_{\lambda \in \Lambda_n}|q_\lambda^n-\mu(E_\lambda^n)|\\
\leq & \| f\|_{\infty} \frac{1}{n}\ \rightarrow 0\quad (n\rightarrow \infty).
\end{align*}
From the above, it is sufficient to prove that $\mu_n(f)$ converges to $\mu(f)$.
Note that
\[ \mu_n(f)=\sum_{\lambda \in \Lambda_n}\mu(E_\lambda^n)f(p_\lambda^n) =\int \sum_{\lambda \in \Lambda_n}f(p_\lambda^n)\chi_{E_\lambda^n}d\mu, \]
where $\chi_{E_\lambda^n}$ is the characteristic function of $E_\lambda^n$. Since $f$ is continuous and the diameter of $E_\lambda^n$ tends to $0$, the function $\sum_{\lambda \in \Lambda_n}f(p_\lambda^n)\chi_{E_\lambda^n}$ converges pointwise to $f$.
Therefore $\mu_n(f)$ converges to $\mu(f)$ by the bounded convergence theorem.
\end{proof}

Let $G$ be a group acting on $X$ continuously and cocompactly, that is, there exists a compact subset $K\subset X$ such that 
$G(K):= \bigcup_{g\in G} g(K)=X$. We define an action of $G$ on $M(X)$ by pushing forward, namely, for $\mu \in M(X)$ and $g\in G$ we define $g_\ast (\mu) \in M(X)$ to be the push-forward of $\mu$ by $g$, explicitly,
\[ g_\ast (\mu) (E):= \mu (g^{-1}(E))\]
for any Borel subset $E\subset X$.
A measure $\mu \in M(X)$ is said to be $G$\ti{-invariant} if we have $g_\ast(\mu) =\mu$ for any $g\in G$.
Set 
\[M_G(X):= \{ \mu \in M(X)\mid \mu \: G\text{-invariant} \}.\]
We will prove that the space $M_G(X)$ is locally compact, separable and completely metrizable.
The topological property of $M_G(X)$ is similar to that of the space of probability measures on a compact metric space with weak-$\ast$ topology.
As a side note, a locally compact second countable Hausdorff space is completely metrizable in general.

\begin{lemma}
For $\mu \in M(X)$ the following are equivalent:
\begin{enumerate}
\item $\mu$ is $G$-invariant;
\item for any $f\in C_c(X)$ and $g\in G$
 \[ \int f d\mu = \int f \circ gd\mu ;\] 
\item for any $f\in C$ and $g\in G$
 \[ \int f d\mu = \int f \circ gd\mu .\] 
\end{enumerate}
\end{lemma}
\begin{proof}
(1)$\Rightarrow$(2):
For the characteristic function $\chi_{E}$ of a Borel subset $E\subset X$ and for $g\in G$, we have
\[ \int \chi_{E}d\mu =\mu (E)=\mu (g^{-1}(E))=\int \chi_ {g^{-1}(E)}d\mu =\int \chi_E\circ gd\mu.\]
Recall that any $f\in C_c(X)$ can be approximated by step functions, each of which is a finite sum of constant multiplication of characteristic functions.
Hence, (2) follows.

(2)$\Rightarrow$(1):
First, we check that $\mu(J)=g_\ast(\mu) (J)$ for any $g\in G$ and any compact subset $J\subset X$.
The characteristic function $\chi_J$ can be approximated 
by a sequence $\{ f_n\} \subset C_c(X)$, that is to say,
\[ \int |\chi_J-f_n|d\mu \rightarrow 0 \ (n\rightarrow \infty ). \]
Therefore, $\chi_{g^{-1}(J)}$ is approximated by the sequence $\{ f_n\circ g\}$, and so we have
\begin{align*}
\mu (g^{-1}(J))=\int \chi_{g^{-1}(J)}d\mu =&\lim_{n\rightarrow \infty }\int f_n\circ g d\mu \\
=&\lim_{n\rightarrow \infty }\int f_nd\mu =\int \chi_{J}d\mu=\mu (J).
\end{align*}
Since $\mu$ is regular, we have $\mu(E)=g_\ast (\mu)(E)$ for any $g\in G$ and any Borel subset $E\subset X$ with $\mu (E)<\infty$.
In general, we can decompose a Borel subset $E\subset X$ into a disjoint union of countably many Borel subsets $\{ E_n \}_{n\in \NN}$ with $E=\sqcup_n E_n$ and $\mu (E_n)<\infty$ for any $n\in \NN$ since $X$ is a union of countably many compact subsets. Then,
\[ \mu(E)=\sum_{n\in \NN} \mu(E_n)=\sum_{n\in \NN} g_\ast(\mu) (E_n)=g_\ast(\mu) (E)\]
for any $g\in G$.

(2)$\Rightarrow$(3): Obvious.

(3)$\Rightarrow$(2): This follows from the denseness of $C$ in $C_c(X)$.
\end{proof}

\begin{proposition}
The space $M_G(X)$ is a closed subspace of $M(X)$. Hence $M_G(X)$ is also a separable complete metric space.
\end{proposition}
\begin{proof}
Let $\mu_n \in M_G(X)\ (n\in \NN)$ and assume that $\mu_n$ converges to $\mu \in M(X)$.
For any $f\in C_c(X)$ and $g\in G$ we have
\[ \mu_n(f\circ g)\rightarrow \mu (f\circ g )\ \tn{and}\ \mu_n(f)\rightarrow \mu (f)\ (n\rightarrow \infty) .\]
Hence $\mu (f\circ g)=\mu (f)$, which implies $\mu \in M_G(X)$.
\end{proof}

Recall that $G$ acts on $X$ cocompactly. Since $X$ is locally compact, any compact subset of $X$ is included in a finite union of relatively compact open subsets of $X$. Therefore, we can take a compact subset $K\subset X$ such that $G(\mathrm{Int}(K))=X$.

\begin{proposition}\label{prop:cocompact action and measure of compact space}
Let $K$ be a compact subset of $X$ with $G(\mathrm{Int}(K))=X$. Let $\mu ,\nu \in M_G(X)$. If the restriction of
$\mu$ to $K$ equals that of $\nu$, then $\mu =\nu$.
\end{proposition}
\begin{proof}
Since $\mu$ and $\nu$ are regular, it is sufficient to show that $\mu (J)=\nu (J)$ for any compact subset $J\subset X$.
From the assumption we can take $g_1,\dots g_m\in G$ such that
\[ J\subset \bigcup_{i=1}^m g_i(K).\]
Then by using $g_1(K),\dots ,g_m(K)$, we divide $J$ into $J_1,\dots ,J_m$ such that $J_i \subset g_i(K)$ and $J$ is a disjoint union of $J_1,\dots ,J_m$.
Hence
\begin{align*}
\mu (J)=\sum_{i=1}^m\mu (J_i)&=\sum_{i=1}^m\mu (g_i^{-1}(J_i))\\
&=\sum_{i=1}^m\nu (g_i^{-1}(J_i))=\sum_{i=1}^m\nu (J_i)=\nu(J).
\end{align*}
This completes the proof.
\end{proof}

By Proposition \ref{prop:cocompact action and measure of compact space}, we see that the property of $M_G(X)$ is similar to that of $M(Y)$ for a compact metric space $Y$.
To prove the local compactness of $M_G(X)$ we use the following lemma:

\begin{lemma}[See {\cite[2.13 Theorem]{Rud86}}]\label{lem:Rudin}
Let $K$ be a compact subset of $X$. Suppose $V_1,\dots ,V_n$ are open subsets of $X$ and
\[ K\subset V_1\cup \cdots \cup V_n.\]
Then there exist continuous functions $h_1,\dots ,h_n\in C_c(X)$ such that $h_i\geq 0$, $\mathrm{supp}h_i\subset V_i$ and
$ h_1(x)+\cdots +h_n(x)=1$
for any $x\in K$.
\end{lemma}
The collection $\{ h_1,\dots ,h_n\}$ is called a \ti{partition of unity} on $K$, subordinate to the cover $\{ V_1,\dots ,V_n\}$.

\begin{theorem}\label{thm:M_G(X)}
The space $M_G(X)$ is a locally compact, separable and complete metric space.
\end{theorem}
\begin{proof}
We need to prove only that $M_G(X)$ is locally compact.
Take any $\mu \in M_G(X)$. Recall that we included the function $\chi_k\ (k\in \NN)$ with respect to the compact subset $K_k$ in the set $C=\{ \phi_n \mid n\in \NN\}$ when we defined the metric $d_M$ on $M(X)$.
Take a sufficiently large $k\in \NN$ such that $G(\mathrm{Int} (K_k))=X$ and then take $k_0 \in \NN$ such that $\phi_{k_0}=\chi_k$. 
Now, we take $\varepsilon >0$ with $\varepsilon <2^{-k_0}$ and prove that the closed ball
\[ B(\mu ,\varepsilon ):=\{ \nu \in M_G(X) \mid d_M(\mu ,\nu )\leq \varepsilon \}\]
is compact.
For any $\nu \in B(\mu ,\varepsilon )$ we have
\[ 2^{-k_0}\min \{| \mu (\chi _{k})-\nu (\chi _{k})| ,1\} \leq d_M(\mu ,\nu )\leq \varepsilon <2^{-k_0}\]
by the definition of $d_M$. Thus
\[ | \mu (\chi _{k})-\nu (\chi _{k})|<1,\]
which implies
\[ \nu (K_k)\leq \nu (\chi_{k})< 1+\mu (\chi_{k}).\]
Put $K:=K_k$ and $M:=1+\mu (\chi_{k})$.

Now, we take any sequence $\{ \mu_n\}\subset B(\mu ,\varepsilon)$ and prove that $\{ \mu_n\}$ contains a convergent subsequence.
Set 
\[ C_K:=\{ f\in C \mid \mathrm{supp} f\subset K\}.\]
Note that $C_K$ is countable. For each $f\in C_K$
\[ |\mu_n(f) |=\left| \int f d \mu_n \right| \leq \| f\|_{\infty} \mu_n (K) \leq \| f\|_{\infty} M, \]
which implies that the sequence $\{ \mu_n(f)\}$ is bounded and has a convergent subsequence.
From the diagonalization argument, we can take a subsequence $\{\mu_{\phi (n)}\}$ of $\{ \mu_n\}$ such that $\{ \mu_{\phi (n)}(f)\}$ is a convergent sequence for any $f\in C_K$. Then we obtain a map $\Phi \: C_K\rightarrow \RR$ as
\[ \Phi (f):=\lim_{n\rightarrow \infty } \mu_{\phi(n)}(f)\ (f\in C_K).\]
By the choice of $C$, for any $f\in C_c(X)$ with $\mathrm{supp}f\subset K$ there is a sequence in $C_K$ converging to $f$. Hence we can extend $\Phi$ to a positive linear functional on $\{ f\in C_c(X) \mid \mathrm{supp}f\subset K\}$ such that $\Phi (f)=\lim_{n\rightarrow \infty}\mu_{\phi (n)}(f)$. Finally, we extend $\Phi$ to a positive linear functional on $C_c(X)$ as follows.
For every $f\in C_c(X)$ take $g_1,\dots ,g_m\in G$ such that
\[ \mathrm{supp}f\subset g_1(\mathrm{Int}(K))\cup \cdots \cup g_m(\mathrm{Int}(K)).\]
By using Lemma \ref{lem:Rudin}, take a partition of unity $\{ h_1,\dots ,h_m\}$ on $\mathrm{supp}f$, subordinate to the cover $\{ g_1(\mathrm{Int}(K)), \dots , g_m(\mathrm{Int}(K))\}$. Then $f=f_1+\cdots +f_m$ for $f_i:=fh_i\ (i=1,\dots ,m)$. Note that
$\mathrm{supp}f_i\subset \mathrm{supp}h_i\subset g_i (\mathrm{Int}(K))$,
and so $\mathrm{supp} (f_i\circ g_i)= g_i^{-1}(\mathrm{supp}f_i)\subset \mathrm{Int}(K)$.
Now, we define $\Phi (f)$ as
\[ \Phi (f):=\sum_{i=1}^m \Phi (f_i\circ g_i).\]
To see that $\Phi (f)$ does not depend on the choice of $h_i$ and $g_i$, we check that the following equality holds:
\[ \sum_{i=1}^m \Phi (f_i\circ g_i)=\lim_{n\rightarrow \infty} \mu_{\phi (n)}(f).\]
Actually, we have
\[
\sum_{i=1}^m \Phi (f_i\circ g_i)=\sum_{i=1}^m \lim_{n\rightarrow \infty }\mu_{\phi (n)}(f_i\circ g_i)
=\sum_{i=1}^m \lim_{n\rightarrow \infty }\mu_{\phi (n)}(f_i)
=\lim_{n\rightarrow \infty} \mu_{\phi (n)}(f).
\]
From the above equality, we can see that $\mu_{\phi (n)}(f)$ converges to $\Phi (f)$ for any $f\in C_c(X)$.
From the Riesz-Markov-Kakutani representation theorem we obtain $\nu \in M(X)$ where $\mu_{\phi (n)}$ converges.
Since $M_G(X)$ is a closed subspace of $M(X)$, $\nu$ is contained in $B(\mu ,\varepsilon)$, which completes the proof.
\end{proof}

\begin{remark}
In general, a topological vector space $V$ over $\RR$ is locally compact if and only if $V$ is finite dimensional. However, the space $M_G(X)$ is not necessarily finitely generated. Especially, the space $\SC (G)=M_G(\H (\partial G))$ of subset currents on a hyperbolic group $G$ with $\# \partial G =\infty$ is not finitely generated. The reason is that $G$ includes a quasi-convex subgroup $H$ that is isomorphic to a free group of rank $\geq 2$, which includes a decreasing sequence $\{ H_n \}$ of finitely generated subgroups such that $\{ \gL (H_n)\}$ is a strictly decreasing sequence.
\end{remark}

\section{Volume functionals on Kleinian groups}\label{sec:Volume functionals for Kleinian groups}
First, we recall some fundamental notions on Kleinian groups.
Let $\HH^n$ be the $n$-dimensional hyperbolic space for $n\geq 2$ and $d_{\HH^n}$ the distance function on $\HH^n$.
We usually consider the Poincar\'e ball model of $\HH^n$.
We will denote by $\Isom (\HH^n)$ the group of orientation-preserving isometries of $\HH^n$.
The action of $\Isom (\HH^n)$ extends to the boundary $\partial \HH^n$, which is homeomorphic to $(n-1)$-dimensional sphere $S^{n-1}$. A \ti{Kleinian group} is a discrete subgroup of $\Isom (\HH^n)$. 
It is known that a subgroup $\gG$ of $\Isom (\HH^n)$ is discrete if and only if $\gG$ acts on $\HH^n$ properly discontinuously, that is, for any compact subset $K$ of $\HH^n$ there are at most finitely many $g\in \gG$ such that $gK\cap K\not= \emptyset$.
Here, we remark that our definition of Kleinian group includes Fuchsian groups, which is a discrete subgroup of $\Isom (\HH^2)$.
The \ti{limit set} of a Kleinian group $\gG$, denoted by $\gL (\gG)$, is the set of accumulation points of the orbits $\gG (x)$ in $\partial \HH^n$ for $x\in \HH^n$, which is independent of the choice of $x$.
More generally, the limit set of a subset $X$ of $\HH^n$, denoted by $X(\infty )$, is the set of accumulation points of $X$ in $\partial \HH^n$.
For a closed subset $S\subset \partial \HH^n$ containing at least two points, we define the \ti{convex hull} $CH(S)$ of $S$ to be the smallest convex closed subset of $\HH^n$ including all geodesic lines connecting two points of $S$.

For a Kleinian group $\gG$ we denote by $CH_\gG$ the convex hull of the limit set $\gL(\gG)$.
Note that $\gG$ acts on $CH_\gG$.
The quotient space $\gG\backslash CH_\gG$, denoted by $C_\gG$, is called the \ti{convex core} of $\gG$.
We say that a Kleinian group $\gG$ is \ti{convex-cocompact} if the convex core $C_\gG$ is compact.
A group is said to be \ti{torsion-free} if it does not have any non-trivial element with finite order.
It is known that a Kleinian group $\gG$ is \ti{torsion-free} if and only if $\gG$ acts on $\HH^n$ freely.
A Kleinian group $\gG$ is said to be \ti{non-elementary} if $\gL(\gG)$ contains infinitely many points.
Note that if a Kleinian group $\gG$ is finite, then $\gL (\gG), CH_\gG$ and $C_\gG$ are empty.

In this section, we consider only a non-elementary, torsion-free and convex-cocompact Kleinian group.
Let $\gG$ be a Kleinian group satisfying the above condition. Since $\gG$ acts on $CH_\gG$ properly discontinuously and cocompactly by isometry, $\gG$ is a (Gromov) hyperbolic group by the \v{S}varc-Milnor Lemma. We identify the limit set $\gL (\gG)$ with the boundary of $\partial \gG$.
We can see that a subgroup $H$ of $\gG$ is quasi-convex if and only if $H$ is convex-cocompact. 

Let $m_{\HH^n}$ be the measure on $\HH^n$ induced by the Riemannian metric on $\HH^n$, which implies that $m_{\HH^n}$ is invariant with respect to the action of $\Isom (\HH^n)$. Note that the set of measurable subsets for $m_{\HH^n}$ coincides with that for the restriction of Lebesgue measure to the Poincar\'e ball model of $\HH^n$.
A measurable subset $A\subset CH_\gG$ is called a (\ti{geometric}) \ti{fundamental domain} for the action of $\gG$ on $CH_\gG$ if the boundary $\partial A$ of $A$ in $CH_\gG$ has measure zero with respect to $m_{\HH^n}$, $\gG (A)=CH_\gG$ and $g(A)\cap A$ is included in $\partial A$ or empty for any non-trivial $g\in \gG$.
We define the volume of $C_\gG$ to be $m_{\HH^n}(A)$ for a fundamental domain $A$ for the action of $\gG$ on $CH_\gG$, which is independent of the choice of $A$.
Actually, the following lemma follows:
\begin{lemma}\label{lem:volume of fundamental domain}
Let $A$ be a fundamental domain for the action of $\gG$ on $CH_\gG$. Let $B$ be a measurable subset of $CH_\gG$ satisfying the condition that $\gG(B)=CH_\gG$ and $g(B)\cap B$ has measure zero for any non-trivial $g\in \gG$. Then we have
\[ m_{\HH^n} (A)=m_{\HH^n}(B).\]
\end{lemma}
\begin{proof}
From the assumption, for any measurable subset $X$ of $CH_\gG$ and any finite subset $\gG_0$ of $\gG$ we have
\[ m_{\HH^n}( X )\geq m_{\HH^n}\left( \left( \bigcup_{g\in \gG_0 }g(B) \right) \cap X \right) =\sum_{g\in \gG_0 } m_{\HH^n}(gB\cap X). \]
Hence by taking a limit on $\gG_0$ we have
\[ m_{\HH^n}( X )\geq \sum_{g\in \gG } m_{\HH^n}(gB\cap X). \]
Since the opposite inequality is obvious, we have
\[ m_{\HH^n}( X ) =\sum_{g\in \gG } m_{\HH^n}(gB\cap X). \]
Therefore
\begin{align*}
m_{\HH^n}( A ) 
&=\sum_{g\in \gG } m_{\HH^n}(gB\cap A)\\
&=\sum_{g\in \gG } m_{\HH^n}(B\cap g^{-1}A)\\
&=m_{\HH^n}( B ),
\end{align*}
which is our claim.
\end{proof}

A measurable subset $B$ of $CH_\gG$ satisfying the condition in the above lemma is called a \ti{measure-theoretic} fundamental domain for the action of $\gG$ on $CH_\gG$.

Since $\gG$ acts on $CH_\gG$ cocompactly, we can take a fundamental domain $\F_\gG$ for the action of $\gG$ on $CH_\gG$ such that $\F_\gG$ is convex and bounded. The Dirichlet domain centered at any point $x\in CH_\gG$,
\[ \{ z\in CH_\gG \mid d_{\HH^n}(x,z)\leq d_{\HH^n}(g(x),z)\text{ for any }g\in \gG \}, \]
is a compact convex geometric fundamental domain.

Recall that $\H (\partial \gG )$ is the hyperspace of $\partial \gG$ consisting of all closed subsets of $\partial \HH$ containing at least two points.
We define a function 
$f_\gG \:\H(\partial \gG)\rightarrow \RRR $
by
\[ f_\gG (S):=m_{\HH^n}(CH(S)\cap \F_\gG) \]
for $S\in \H (\partial \gG)$.

\begin{proposition}\label{prop:f is continuous}
The function $f_\gG$ is a continuous function with compact support.
\end{proposition}

For the moment, we assume that the above proposition follows. Then we can define the continuous $\RRR$-linear functional $f_\gG^\ast \: \SC (\gG)\rightarrow \RRR$ by
\[ f_\gG^\ast (\mu ):=\int f_\gG d\mu \]
for $\mu \in \SC (\gG)$. Now, we check that $f_\gG^\ast (\eta_H )$ equals the volume of the convex core $C_H$ for a non-trivial quasi-convex subgroup $H$ of $\gG$. Let $R \subset \gG $ be a complete system of representatives of $\gG /H$. Then we have
\begin{align*}
f_\gG^\ast (\eta_H )= \int f_\gG d\eta_H 
&=\sum_{gH\in \gG /H } f(g\Lambda (H)) \\
&=\sum_{g\in R } m_{\HH^n}( g CH_H\cap \F_\gG ) \\
&=\sum_{g\in R } m_{\HH^n}(CH_H \cap g^{-1}(\F_\gG )) \\
&=m_{\HH^n }\left( CH_H \cap \bigcup_{g\in R} g^{-1}(\F_\gG ) \right) .
\end{align*}
In the last of the above equation we used the property that for any non-trivial $g\in \gG$ the intersection $g(\F_\gG) \cap \F_\gG$ has measure zero.
Note that $R^{-1}=\{ g^{-1}\mid g\in R\}$ is a complete system of representatives of $H\backslash \gG$.
Then it is sufficient to prove that
\[ A:=CH_H\cap \bigcup_{g\in R} g^{-1}(\F_\gG ) \]
is a measure-theoretic fundamental domain for the action of $H$ on $CH_H$. 
First, we check that $H(A)=CH_H$. Take any $x\in CH_H$. Then there exist $g\in \gG,\ g_0\in R$ and $h\in H$ such that $g(x)\in \F_\gG$ and $g=g_0h^{-1}$. Therefore
\[ x\in g^{-1}(\F_\gG)=h g_0^{-1}(\F_\gG) \subset h(A).\]
This concludes that $H(A)=CH_H$.
For a non-trivial $h\in H$ we have
\[ h(A)\cap A=CH_H \cap \left( \bigcup_{g_1,g_2\in R} hg_1^{-1}(\F_\gG )\cap g_2^{-1}(\F_\gG ) \right). \]
If $g_2hg_1^{-1}=\id$ for $g_1,g_2\in \R$, then $g_2h=g_1$ and so $h=\id$, a contradiction.
Hence $g_2hg_1^{-1}$ is not the identity element for any $g_1,g_2\in R$.
Therefore we have
\[ hg_1^{-1}(\F_\gG )\cap g_2^{-1}(\F_\gG )\subset g_2^{-1}(\partial \F_\gG)\]
and so
\[ h(A)\cap A \subset CH_H\cap \gG (\partial \F_\gG) .\] 
This implies that $h(A)\cap A$ has measure zero.

Therefore, $f_\gG^\ast (\eta_H )$ equals the volume of the convex core $C_H$ for every non-trivial quasi-convex subgroup $H$ of $\gG$.

From the above argument, we obtain the following main theorem in this section, where $\mathrm{Vol}=f_\gG^\ast$.

\begin{theorem}\label{thm:volume functional}
Let $\gG$ be a non-elementary, torsion-free and convex-cocompact Kleinian group in $\mathrm{Isom}(\HH^n)$.
There exists a continuous $\RRR$-linear functional
\[ \mathrm{Vol} \: \SC (\gG )\rightarrow \RRR \]
such that for every non-trivial quasi-convex subgroup $H$ of $\gG$, $\mathrm{Vol} (\eta_H )$ equals the volume of the convex core $C_H$.
\end{theorem}

Now, we prepare some lemmas for proving Proposition \ref{prop:f is continuous}.

\begin{lemma}\label{lem:boundary has measure zero}
Let $X$ be a convex subset of $\HH^n$. Then the boundary $\partial X$ has measure zero with respect to $m_{\HH^n}$ and $X$ is particularly measurable. 
\end{lemma}
\begin{proof}
Recall that $\HH^n$ is the Poincar\'e ball model of the $n$-dimensional hyperbolic space, that is, $\HH^n$ is the unit open ball of $\RR^n$.
We can assume that $X$ contains $0$ without loss of generality since the action of $\Isom (\HH^n)$ on $\HH^n$ is transitive and $m_{\HH^n}$ is $\Isom (\HH^n)$-invariant.
Let $m_L$ be the Lebesgue measure on $\RR^n$. It is sufficient to see that $\partial X$ has measure zero with respect to $m_L$.

First, we consider the case that $X$ contains $0$ as an interior point. Since $X$ is convex, for any $x\in X$ there exists a unique geodesic joining $0$ to $x$, which is also a geodesic in $\RR^n$. Therefore $X$ is a star-like domain centered at $0$ in $\RR^n$.
For $t\geq 0$ set
\[ X_t:=\{ tx\in \RR^n \mid x\in X\} .\]
Since $0$ is an interior point of $X$, there exists a small open ball $U\subset X$ containing $0$. For $t_0\in [0,1)$ and $x\in X_{t_0}$ there exists $t>1$ such that $tx \in X$ and the convex hull of $U\cup \{ tx\}$ in $\HH^n$ contains $x$ as an interior point, so $x$ is an interior point of $X$, which implies that for any $t_0\in [0,1)$ the set $X_{t_0}$ is included in the interior $\mathrm{Int}(X)$ of $X$.
Then for any $t>1$
\[ \mathrm{Int}(X)_t:=\{ tx \mid x\in \mathrm{Int}(X)\} \]
includes $X$ since for any $x\in X$ we have $x/t\in \mathrm{Int}(X)$ and
\[ x=t\left( \frac{1}{t}x\right) \in \mathrm{Int}(X)_t.\]
Note that $\mathrm{Int}(X)_t$ is similar to $\mathrm{Int}(X)$ in $\RR^n$. Hence we have
\[ m_L(\mathrm{Int}(X)_t)=t^n m_L(\mathrm{Int}(X)).\]
Therefore for any $t>1$ we have
\[ \partial X\subset \mathrm{Int}(X)_t\setminus \mathrm{Int}(X). \]
As a result, we obtain
\begin{align*}
m_L(\partial X)
&\leq m_L(\mathrm{Int}(X)_t)-m_L(\mathrm{Int}(X))\\
&=(t^n-1)m_L(\mathrm{Int}(X) )\rightarrow 0\quad (t\rightarrow 1).
\end{align*}
This implies that the boundary $\partial X$ has measure zero with respect to the Lebesgue measure $m_L$.
The equation $X=(\partial X\cap X)\cup \mathrm{Int}(X)$ implies the measurability of $X$.

If $X$ does not contain any interior points and contains $0$, then $X$ is included in a hyperplane of $\RR^n$, which implies that both $X$ and $\partial X$ have measure zero.
\end{proof}

From the proof of the above lemma, we see that we can apply some techniques of convex geometry in Euclidean spaces to $\HH^n$ by using $0$ as a base point.
Let $d_{\RR^n}$ be the Euclidean metric on $\RR^n$.

A \ti{hyperplane} of $\HH^n$ is a totally geodesic codimension-1 submanifold. Here, ``totally geodesic'' means that for any two different points in the submanifold the geodesic line passing through the two points is included in the submanifold. Any hyperplane of $\HH^n$ is isometric to $\HH^{n-1}$ and its limit set is homeomorphic to $S^{n-2}$.

Any hyperplane divides $\HH^n$ into two connected components, and the union of the hyperplane and one of the connected components is called a \ti{half-space} of $\HH^n$. In this case the hyperplane is the boundary of the half-space in $\HH$.
The following property of a convex set is well-known in $\RR^n$ and also follows in $\HH^n$:
for a convex subset $X$ of $\HH^n$, $x\in \partial X$ and an exterior point $y$ of $X$ there exists a half-space $U$ of $\HH^n$ such that $U\supset X$, $x\in \partial U$ and $y\not\in U$.

From the above property we can see that for any closed convex subset $X$ of $\HH^n$ the intersection of all half-spaces including $X$ coincides with $X$.
Therefore for any $S\in \H (\partial \gG)$ the convex hull $CH(S)$ coincides with the intersection of all half-spaces whose limit sets include $S$.

Recall that $\H (\partial \gG)$ is a metric space with a Hausdorff distance $d_{\mathrm{Haus}}$, which is compatible with the Hausdorff distance induced by the restriction of $d_{\RR^n}$ to $\partial \gG$.
In this section we use the Hausdorff distance $D$ induced by $d_{\RR^n}$ instead of $d_{\mathrm{Haus}}$.
Note that we can consider the distance $D$ for any two non-empty subsets of $\HH\cup \partial \HH$.

Take $S,S'\in \H (\partial \gG)$ such that $\# S=\# S'=2$.
Then we can see that for any $\varepsilon>0$ there exists $\delta>0$ such that if $D(S,S')<\delta $, then $D(CH(S),CH(S'))<\varepsilon$.
This property follows since for $S\in \H (\partial \gG)$ with $\#S=2$, $CH(S)$ is the intersection of $\HH^n$ and the circle in $\RR^n$ intersecting $\partial \HH^n$ orthogonally at each point of $S$.
See the argument below Lemma \ref{lem:taking convex hull is continuous} for the reason of the uniform continuity.

For a hyperplane $H$ of $\HH^n$ the union of all geodesic lines connecting two points of the limit set $H(\infty )$ coincides with $H$ itself. Therefore for a hyperplane $H$ of $\HH^n$ and $\varepsilon >0$ there exists $\delta >0$ such that if a hyperplane $H'$ satisfies the condition that $D(H(\infty ), H'(\infty))<\delta$, then $D(H,H')<\varepsilon$.

Consider the hyperplane $H:=(\RR^{n-1}\times \{ 0\}) \cap \HH^n$.
For $a >0$ we call the set
\[ H_a :=(\RR^{n-1}\times [-a ,a ])\cap \HH^n \]
the $[a]$-neighborhood of $H$.
Then we can see that if a hyperplane $H'$ of $\HH^n$ is included in $H_a$, then we have
\[ D(H',H)\leq D(CH((\RR^{n-1}\times \{ a\} )\cap \partial \HH^n ), H) .\]
Moreover, for a hyperplane $H'$ of $\HH^n$ if either $D(H,H')$ or $D(H(\infty),H'(\infty))$ is sufficiently small, then $H'$ is included in $H_a$.
Note that for any hyperplane $H'$ of $\HH^n$ there exists $\phi \in \Isom (\HH^n)$ such that $\phi (H')=H$.
Then we define the $[a]$-neighborhood $H'_a$ of $H'$ to be $\phi^{-1} (H_a)$, and $H'_a$ has the same property as $H_a$.


\begin{lemma}\label{lem:exterior converging}
Let $S\in \H (\partial \gG)$ and $\{ S_k\}_{k\in \NN}$ a sequence in $\H( \partial \gG)$ converging to $S$. For any exterior point $x$ of $CH(S)$ there exists $N\in \NN$ such that if $k\geq N$, then $x$ is an exterior point of $CH(S_k)$.
\end{lemma}
\begin{proof}
We can assume that $\partial \gG =\partial \HH$ without loss of generality since $\H (\partial \gG)$ can be considered as a subspace of $\H (\partial \HH)$.
Take a half-space $U$ of $\HH^n$ such that $CH(S)\subset U$ and $x\not\in U$.
Note that $x$ is also an exterior point of $U$ and we can take the $[a]$-neighborhood $(\partial U)_a$ of $\partial U$ such that $x\not\in (\partial U)_a$. Then there exists a half-space $U'$ of $\HH^n$ such that $x$ is an exterior point of $U'$, and $U$ and $(\partial U)_a$ are included in $U'$.
Therefore if we take a sufficiently large $N\in \NN$ and $k\geq N$, then $S_k\subset U'(\infty)$, which implies that $x$ is an exterior point of $CH(S_k)$.
\end{proof}

\begin{lemma}\label{lem:interior converging}
Let $S\in \H (\partial \gG)$ and $\{ S_k\}_{k\in \NN}$ a sequence in $\H( \partial \gG)$ converging to $S$. For any interior point $x$ of $CH(S)$ there exists $N\in \NN$ such that if $k\geq N$, then $x$ is an interior point of $CH(S_k)$.
\end{lemma}
\begin{proof}
We can assume that $\partial \gG =\partial \HH$ without loss of generality.
We also assume that $x=0$ in $\HH^n$.
Note that an isometry of $\HH^n$ fixing $0$ is the restriction of an isometry of $\RR^n$ to $\HH^n$.
Take $r>0$ such that the open ball $B(0,r)$ centered at $0$ with radius $r$ with respect to $d_{\HH^n}$ is included in $CH(S)$.
Take a half-space $U$ such that $U\supset B(0,r)$ and the hyperbolic distance from $0$ to $\partial U$ equals $r$, that is, $\partial U$ is tangent to the boundary of $B(0,r)$.
Then by considering the $[a]$-neighborhood of $\partial U$, there exists $\varepsilon >0$ such that for any half-space $U'$ of $\HH^n$ if $D(U(\infty),U'(\infty))<\varepsilon$, then $U' \supset B(0,r/2)$.
Moreover, for this $\varepsilon$, we can see that for any two half-spaces $U_1,U_2$ of $\HH^n$, if $U_1\supset B(0,r)$ and $D(U_1(\infty),U_2(\infty))<\varepsilon$, then $U_2\supset B(0,r/2)$.

Assume that $k$ is sufficiently large and $D(S,S_k)<\varepsilon$.
Take any half-space $V$ of $\HH^n$ such that $V\supset CH(S_k)$.
By considering the $[a]$-neighborhood of $\partial V$, there exists a half-space $V'$ of $\HH^n$ such that $D(V(\infty),V'(\infty))<\varepsilon$ and $V'\supset CH(S)$. Since $V'\supset B(0,r)$, we have $V\supset B(0,r/2)$. This implies that $CH(S_k)\supset B(0,r/2)$.
\end{proof}

\begin{lemma}\label{lem:A(K) is relatively compact}
For a bounded subset $K$ of $CH_\gG $ the set
\[ A(K) :=\{ S\in \H (\partial \gG )\mid CH(S)\cap K\not= \emptyset \}\]
is a relatively compact subset of $\H (\partial \gG)$.
Moreover, for any compact subset $E$ of $\H (\partial \gG )$ there exists a bounded subset $K$ of $CH_\gG$ such that $E\subset A (K)$.
\end{lemma}
\begin{proof}
From Lemma \ref{lem: A g is compact subset}, for a Cayley graph $\Cay (\gG)$ with respect to a finitely generating set and $g\in \gG$ the set
\[ A_g=\{ S\in \H (\partial \gG )\mid WC(S)\ni g\} \]
is a compact subset of $\H (\partial \gG)$. Take $x_0\in \HH^n$. Then we have a quasi-isometry 
\[ \theta \: \Cay (\gG)\rightarrow CH_\gG;\ g\mapsto g(x_0).\]
Recall that $\theta $ induces a homeomorphism $\partial \theta $ from $\partial \gG$ to $\gL (\gG)$, which is independent of the choice of $x_0$. We identify $\partial \gG$ with $\gL (\gG)$ by this homeomorphism.
Take a quasi-inverse $\theta '$ to $\theta $. Since $K$ is bounded, $\theta '(K)$ is also bounded in $\Cay (\gG)$.
For $S\in A(K)$ we can see that the weak convex hull $WC(S)$ in $\Cay (\gG)$ intersects the $c$-neighborhood of $\theta '(K)$ for some $c>0$ by the property of quasi-isometry.
Hence the set $A(K)$ is included in a union of $A_{g_1},\dots ,A_{g_m}$ for some $g_1,\dots g_m\in \gG$.
Since $A_g$ is compact for any $g\in \gG$, the set $A(K)$ is relatively compact.

From the proof of Lemma \ref{lem:continug subset currents}, for any compact subset $E$ of $\H (\partial \gG )$ there exist finitely many elements $g_1,\dots , g_m \in \gG$ such that $E$ is included in the union of $A_{g_1}, \dots , A_{g_m}$.
Then by considering $\theta ( \{ g_1,\dots ,g_m\} )$ we can take a bounded subset $K$ of $\HH ^n$ such that for any $g_i$ and any $S\in A_{g_i}$ the convex hull $CH (S)$ in $\HH^n$ must intersects $K$. This implies that $E \subset A (K)$.
\end{proof}

\begin{proof}[Proof of Proposition \ref{prop:f is continuous}]
The support of $f$ is included in the closure of $A(\F_\gG)$, which is compact since $\F_\gG$ is bounded.
Now, we prove the continuity of $f$.
Let $S\in \H (\partial \gG)$. Let $\{ S_k\}$ be a sequence in $\H( \partial \gG)$ converging to $S$.
It is sufficient to see that $m_{\HH^n}(CH(S_k)\cap \F_\gG)$ converges to $m_{\HH^n}(CH(S)\cap \F_\gG)$.
By the bounded convergence theorem it is sufficient to see that the characteristic function of $CH(S_k)\cap \F_\gG$ converges pointwise to the characteristic function of $CH(S)\cap \F_\gG$ almost everywhere.
Actually, from Lemmas \ref{lem:boundary has measure zero}, \ref{lem:exterior converging} and \ref{lem:interior converging} this claim follows.
\end{proof}

If the dimension $n$ is $2$, then we can obtain a stronger result than Lemmas \ref{lem:exterior converging} and \ref{lem:interior converging}. We will write $\HH^2$ simply as $\HH$.

\begin{lemma}\label{lem:taking convex hull is continuous}
Let $S\in \H (\partial \HH)$. For any $\varepsilon >0$ there exists $\delta >0$ such that if $D(S,S')<\delta$ for $S'\in \H (\partial \HH)$, then $D(CH(S),CH(S'))<\varepsilon$.
\end{lemma}
\begin{proof}
In this proof we use only the Euclidean metric $d$ in $\RR^2$ and the Hausdorff distance $D$ induced by $d$.
However, we will use the term ``geodesic'' as a geodesic in $\HH$.

First of all, we consider the description of $CH(S)$ in the case that $S\not=\partial \HH$. Since $S$ is a closed subset of $\partial \HH= S^1$, the complement $S^c=\partial \HH \setminus S$ is a union of at most countably many open intervals $\{ I_\lambda \}_{\lambda \in \gL}$ of $\partial \HH$, that is,
\[ S=\partial \HH \setminus \bigsqcup_{\lambda \in \gL}I_\lambda. \]
For each $I_\lambda$ we consider the interior $\mathrm{Int}(CH(\ol{I_\lambda}))$ in $\HH$, which equals the union of all geodesic line connecting two points of $I_\lambda$. Then we can see that
\[ CH(S)=\HH \setminus \bigsqcup_{\lambda \in \gL}\mathrm{Int}(CH(\ol{I_\lambda})) .\]
Note that the boundary $\partial CH(S)$ coincides with the union of all geodesic line connecting the two points of $\partial I_\lambda$ taken over $\lambda \in \gL$.

Fix $\varepsilon>0$.
First, we consider the case that $\#S=2$.
Then we can take $\delta>0$ such that for any $S'\in \H (\partial \HH)$ with $\# S'=2$ and $D(S,S')<\delta $ we have $D(CH(S),CH(S'))<\varepsilon$.
Now, we do not assume that $\#S'=2$. Then $\partial \HH \setminus S'$ is a disjoint union of countably many intervals $\{ I_\lambda \}_{\lambda \in \gL }$.
Since $D(S,S')<\delta$, there exists two $\lambda _1,\lambda_2\in \gL$ such that
\[ D(S, \partial \HH \setminus (I_{\lambda_1}\cup I_{\lambda_2}))<\delta.\]
Then we can see that $D(CH(S),CH(S'))<\varepsilon$. 

Next, we consider the case that $S=\partial \HH$. Take $\delta>0$ such that if the diameter of an open interval $I\subset \partial \HH$ is smaller than $2\delta$, then the diameter of $CH(\ol {I})$ is smaller than $\varepsilon$. Then for $S'\in \H (\partial \HH)$ with $D(S,S')<\delta $, the complement $S^c$ never includes an open interval with diameter $>2\delta$ . Therefore $D(CH(S),CH(S'))<\varepsilon$.

Finally, we consider the case that $S\not =\partial \HH$ and $\#S\geq 3$.
Take open intervals $\{ I_\lambda \}_{\lambda \in \gL}$ of $\partial \HH$ such that $\partial \HH \setminus S$ is a disjoint union of $\{ I_\lambda \}$.
Take $\delta>0$ satisfying the following two conditions:
\begin{enumerate}
\item for any $S_1,S_2\in \H (\partial \HH)$ with $\# S_1=\#S_2=2$, if $D(S_1,S_2)<\delta $, then we have $D(CH(S_1),CH(S_2))<\varepsilon$;
\item if the diameter of an open interval $I\subset \partial \HH$ is smaller than $2\delta$, then the diameter of $CH(\ol {I})$ is smaller than $\varepsilon$.
\end{enumerate}
Take $S'\in \H(\partial \HH)$ with $D(S,S')<\delta$ and open intervals $\{ I'_\lambda \}_{\lambda \in \gL'}$ of $\partial \HH$ such that $\partial \HH \setminus S'$ is a disjoint union of $\{ I'_\lambda \}$. 

Take $x\in CH(S)$. First, we consider the case that $d(x,\partial CH(S))<\varepsilon $. Then there exists $\lambda \in \gL$ such that $d(x, CH (\partial I_\lambda ))<\varepsilon $. If the diameter of $I_\lambda$ is smaller than $2\delta$, then the diameter of $CH(\ol{I_\lambda})$ is smaller than $\varepsilon$ and there exists $\xi \in S'$ such that $CH(\partial I_\lambda)$ is included in the $(\delta +\varepsilon)$-neighborhood of $\xi$. This implies that $x$ belongs to the $(\delta +2\varepsilon)$-neighborhood of $CH(S')$.
If the diameter of $I_\lambda$ is larger than or equal to $2\delta$, then by the Condition (1) there exists $\lambda'\in \gL'$ such that $D(CH(\partial I_\lambda ),CH(\partial I'_{\lambda '}))<\varepsilon$, which implies that $x$ is contained in the $2\varepsilon$-neighborhood of $CH(S')$.

Next, we consider the case that $d(x,\partial CH(S))>\varepsilon $. Assume that $x\not\in CH(S')$, that is, there exists $\lambda'\in \gL'$ such that $x\in \mathrm{Int}(CH(\ol{I'_{\lambda'}}))$. If the diameter of $I'_{\lambda'}$ is smaller than $2\delta$, then $x$ is contained in the $\varepsilon $-neighborhood of $CH(S')$ by the Condition (2).
If the diameter of $I'_{\lambda'}$ is larger than or equal to $2\delta$, then there exists $\lambda \in \gL$ such that $D(CH(\partial I_\lambda ),CH(\partial I'_{\lambda '}))<\varepsilon$. Since $x\not\in \mathrm{Int}(CH(\ol{I_\lambda}))$, we have
\[ x\in \mathrm{Int}(CH(\ol{I'_{\lambda'}}))\setminus \mathrm{Int}(CH(\ol{I_\lambda})),\]
which implies that $d(x,CH(\partial I_\lambda))<\varepsilon$. This is a contradiction. Hence $x\in CH(S')$.

Therefore, in any cases $CH(S)$ is included in the $(\delta +2\varepsilon )$-neighborhood of $CH(S')$.
By the same way as the above we can see that $CH(S')$ is included in $(\delta +2\varepsilon )$-neighborhood of $CH(S)$.
This completes the proof.
\end{proof}

From Lemmas \ref{lem:A(K) is relatively compact} and \ref{lem:taking convex hull is continuous}, we see that if $Y$ is a bounded open subset of $\HH$, then $A(Y)$ is a relatively compact open subset of $\H (\partial \HH)$; if $Y$ is a compact subset of $\HH$, then $A(Y)$ is also a compact subset of $\H (\partial \HH)$. 

Recall that $\hat{\H }(\partial \HH)$ is the hyperspace of $\partial \HH$ consisting of all closed subsets of $\partial \HH$.
The closure $\ol{\H (\partial \HH )}$ in $\hat{\H }(\partial \HH)$ coincides with $\hat{\H }(\partial \HH)\setminus \{ \emptyset \}$, which is compact.
We define a map $\Phi$ from $\ol{\H (\partial \HH )}$ to $\hat{\H } (\HH \cup \partial \HH)$ by
\[
  \Phi (S):=
\begin{cases}
    CH(S)\cup S & (\# S \geq 2)\\
    S & (\# S=1)
\end{cases} ,
\]
for $S\in \ol{\H (\partial \HH )}$.
Note that we have $D(CH(S_1),CH(S_2))=D(\Phi (S_1), \Phi (S_2))$ for $S_1,S_2\in \H (\partial \HH )$.
Hence Lemma \ref{lem:taking convex hull is continuous} implies that $\Phi$ is continuous at every $S\in \H (\partial \HH )$.
It is easy to see that $\Phi$ is continuous at any point in $\ol{\H (\partial \HH )}$ from the proof of Lemma \ref{lem:taking convex hull is continuous}.
Moreover, $\Phi$ is uniformly continuous on $\ol{\H (\partial \HH )}$ since $\ol{\H (\partial \HH )}$ is compact.
As a result, we obtain the following proposition:

\begin{proposition}
For any $\varepsilon>0$ there exists $\delta >0$ such that for $S_1,S_2\in \H (\partial \HH)$ if $D(S_1,S_2)<\delta$, then $D(CH(S_1),CH(S_2))<\varepsilon$.
\end{proposition}

In the case that the dimension $n$ is $2$, the area of the convex core $C_\gG$ equals $-2\pi \chi (C_\gG)$ from the Gauss-Bonnet theorem for the Euler characteristic of $C_\gG$.
We define the Euler characteristic $\chi (G)$ of a group $G$ to be the Euler characteristic of a $K(G,1)$-space if we can take a $K(G,1)$-space as a finite-dimensional CW-complex. Here, we can see that $C_\gG$ is a finite-dimensional CW-complex and a $K(\gG,1)$-space since the universal cover $CH_\gG$ of $C_\gG$ is contractible.
Then we obtain the following corollary from Theorem \ref{thm:volume functional}, where $\chi=-\mathrm{Vol}/2\pi$.

\begin{corollary}
Let $\gG$ be a torsion-free convex-cocompact Fuchsian group. Then there exists a unique continuous $\RRR$-linear functional
\[ \chi \: \SC (\gG )\rightarrow \RR_{\leq 0}=\{ r\in \RR \mid r\leq 0\} \]
such that for every non-trivial quasi-convex subgroup $H$ of $\gG$ we have
\[ \chi (\eta_H )=\chi (H).\]
\end{corollary}

Note that a torsion-free convex-cocompact Fuchsian group is isomorphic to a surface group or a free group of finite rank since $C_\gG$ is a compact hyperbolic surface possibly with boundary or a closed geodesic, and that a subgroup $H$ of a surface group or a free group of finite rank is quasi-convex if and only if $H$ is finitely generated. The uniqueness of the functional $\chi$ is a result of the denseness property of rational subset currents for $\gG$ (see Theorem \ref{thm:subset currents on hyperbolic surface has denseness property}).
We also remark that in the above corollary our claim is independent of the action of $\gG$ on $\HH$.

For a non-trivial free group $F$ of finite rank the \ti{reduced rank} $\rk (F)$ of $F$ is defined to be $-\chi (F)$, which coincides with $\mathrm{rank}(F)-1$. We define the reduced rank of the trivial group to be $0$.
In the same way, for a surface group $\gG$ we define the reduced rank $\rk (\gG)$ of $\gG$ to be $-\chi (\gG)$.
Then we have the following corollary. Note that in the case that $\gG$ is a free group of finite rank the following corollary was proved in \cite[Theorem 8.1]{KN13}.

\begin{corollary}\label{cor:reduced rank functional}
Let $\gG$ be a surface group or a free group of finite rank. Then there exists a unique continuous $\RRR$-linear functional 
\[ \rk\: \SC (\gG )\rightarrow \RRR\]
such that for every finitely generated subgroup $H$ of $\gG$ we have
\[ \rk(\eta_H)=\rk (H).\]
\end{corollary}

We call $\rk$ the \ti{reduced rank functional} on $\SC (\gG)$.

Let $H$ be a quasi-convex subgroup of $\gG$ and $K$ a finite index subgroup of $H$. Then $\eta_K =[H:K]\eta_H$, where $[H:K]$ is the index of $K$ in $H$. Since $\rk$ is $\RRR$-linear, we have
\[ \rk (K)=\rk (\eta_K)=\rk ([H:K]\eta_H)=[H:K]\rk(\eta_H)=[H:K]\rk(H).\]
This property comes from the property that $C_K$ is a $[H:K]$-fold covering of $C_H$.

\section{Subgroups, inclusion maps and finite index extension}\label{sec:relation between subgroups}

Let $G$ be an infinite hyperbolic group. Since a quasi-convex subgroup $H$ of $G$ is also a hyperbolic group, we want to consider a relation between $\SC (G)$ and $\SC (H)$, especially, in the case that $H$ is a finite index subgroup of $G$.
We assume that $H$ is also an infinite group.
First, we identify the boundary $\partial H$ of $H$ with the limit set $\gL (H)$ in $\partial G$. Then the space $\H(\partial H)$ is a closed subspace of $\H (\partial G)$. Note that if $H$ is a finite index subgroup of $G$, then $\partial H=\partial G$.
Now we consider an infinite quasi-convex subgroup $J$ of $H$ and identify $\partial J$ with $\gL (J)$ in $\partial G$. 
For $\mu \in SC(J)$ we consider $\mu$ as a measure on $\H (\partial H)$, whose support is included in $\H (\partial J)$. Recall that the support of a measure $\mu$ is the smallest closed subset such that the restriction of $\mu$ to the exterior of the closed subset is the zero measure.

\subsection{Natural continuous $\RRR$-linear maps between subgroups}\label{subsec: iota J H}
We can define a natural continuous $\RRR$-linear map $\iota_J^H$ from $SC(J)$ to $SC(H)$ as follows.
Since $H$ acts on $\H (\partial H)$, we define the push-forward $h_\ast (\mu )$ of $\mu\in \SC(J)$ by $h\in H$ as
\[ h_\ast (\mu )(E):=\mu (h^{-1}(E))\]
for every Borel subset $E$ of $\H (\partial H)$. Note that the support of $h_\ast (\mu )$ is included in $h(\H( \partial J))=\H (h \gL (J))\subset \H(\partial H)$. Since $\mu$ is $J$-invariant, $h_\ast (\mu)=\mu$ for $h\in J$.
Now, we define a measure $\iota_J^H(\mu)$ on $\H (\partial H)$ by
\[ \iota_J^H(\mu ):=\sum_{hJ\in H/J} h_\ast (\mu ).\]

\begin{lemma}
Let $H,J$ be infinite quasi-convex subgroups of an infinite hyperbolic group $G$ with $J\subset H$. For any $\mu \in \SC (J)$ the measure $\iota_J^H (\mu)$ is an $H$-invariant locally finite Borel measure on $\H (\partial H)$, that is, $\iota_J^H(\mu)$ is a subset current on $H$. Moreover, the map
\[ \iota_J^H\: \SC(J)\rightarrow \SC (H)\]
is a continuous $\RRR$-linear map.
\end{lemma}
\begin{proof}
First we check that $\iota_J^H (\mu)$ is $H$-invariant. For $g\in H$ we have
\[
g_\ast (\iota_J^H (\mu )) 
= \sum_{hJ\in H/J} g_\ast (h_\ast (\mu ))
=\sum_{hJ\in H/J} (gh)_\ast (\mu )
=\sum_{hJ\in H/J} h_\ast (\mu ),
\]
which implies that $\iota_J^H (\mu)$ is $H$-invariant.

From Lemma \ref{lem: A g is compact subset}, by considering the Cayley graph $\Cay (H)$ of $H$ with respect to a finite generating set of $H$ and the vertex $\id \in \Cay (H)$, the set
\[ A_\id^H=\{ S\in \H (\partial H)\mid WC(S)\ni \id \}\]
is a compact subset of $\H (\partial H)$ and any compact subset $K$ of $\H (\partial H)$ is included in a finite union of $h_1A_\id^H,\dots ,h_mA_\id^H$ for some $h_1,\dots ,h_m\in H$.

Now, for the local finiteness of $\iota_J^H(\mu)$, it is sufficient to see that $\iota_J^H (\mu)(A_\id^H)$ is finite.
Since $J$ is a quasi-convex subgroup of $H$, the counting subset current
\[ \eta _J^H:=\sum_{hJ\in H/J} \delta_{h\gL (J)}\]
on $H$ is locally finite. Hence there are at most finitely many $h_1J,\dots ,h_mJ\in H/J$ such that 
\[ h_1\gL (J),\dots , h_m\gL(J) \in A_\id^H.\]
For $h\in H$ satisfying the condition that $h\gL(J)\not\in A_\id^H$, that is, $WC(h\gL(J))\not\ni \id$, we can see that $h\H (\partial J)\cap A_\id^H=\emptyset$ since for any $S\in h\H (\partial J)$ the weak convex hull $WC(S)$ is included in $WC(h\gL(J))$. Note that $A_\id^H\cap h_i(\H (\partial J))$ is a compact subset of $h_i(\H (\partial J))$ for $i=1,\dots ,m$.
Therefore we have
\[ \iota_J^H (\mu )(A_\id^H ) =\sum_{i=1}^m (h_i)_\ast (\mu )(A_\id^H )=\sum_{i=1}^m \mu ((h_i^{-1}A_\id^H) \cap \H (\partial J))<\infty.\]

Finally, we check that $\iota_J^H$ is continuous. Take $\mu_n ,\mu \in \SC(J)\ (n\in \NN)$ such that $\mu_n$ converges to $\mu$ by taking $n\rightarrow \infty$. Take any compactly supported continuous function $f\:\H(\partial H)\rightarrow \RR$. Since the intersection of a compact subset of $\H (\partial H)$ and $\H (\partial J)$ is compact, the restriction of $f$ to $h\H(\partial J)$ is a continuous function with compact support for any $h\in H$.
From the above argument, there are at most finitely many $h_1J,\dots h_mJ\in H/J$ such that the support of $f$ intersects each of $h_1\H(\partial J),\dots ,h_m\H (\partial J)$. Therefore
\begin{align*}
\int fd\iota_J^H (\mu_n )
=&\sum_{hJ\in H/J}\int f d h_\ast (\mu_n )\\
=&\sum_{i=1}^m\int f \circ h_i d(\mu_n )\\
\underset{n\rightarrow \infty }\rightarrow &\sum_{i=1}^m\int f \circ h_i d(\mu )
=\sum_{hJ\in H/J}\int f d h_\ast (\mu )=\int fd\iota_J^H (\mu ).
\end{align*}
This implies that $\iota_J^H (\mu_n )$ converges to $\iota_J^H (\mu )$.
\end{proof}

\begin{notation}
For a quasi-convex subgroup $K$ of $H$ we denote by $\eta_K^H$ the counting subset currents on $H$ corresponding to $K$, that is,
\[ \eta_K^H=\sum_{hK\in H/K} \delta_{h\gL(K)}.\]
In the case that $H=G$, we usually write $\eta_K$ in place of $\eta_K^G$.
\end{notation}

Since $\gL (J)$ is $J$-invariant, the Dirac measure $\delta_{\gL (J)}=\eta_J^J$ is a subset current on $J$. Then we can see that
\[ \iota_J^H ( \eta_J^J )=\sum_{hJ\in H/J}h_\ast \delta_{\gL (J)}=\sum_{hJ\in H/J}\delta_{h\gL (J)}=\eta_J^H.\]

For simplicity of notation, we write $\iota_H$ in place of $\iota_H^G$.
Then we can see that the composition $\iota_H\circ \iota_J^H$ equals $\iota_J$. Actually, for $\mu \in \SC(J)$,
\begin{align*}
\iota_H\circ \iota_J^H (\mu )
&=\sum_{gH\in G/H}g_\ast \left( \sum_{hJ\in H/J} h_\ast (\mu )\right) \\
&=\sum_{gH\in G/H, hJ\in H/J}g_\ast (h_\ast (\mu ))\\
&=\sum_{gH\in G/H, hJ\in H/J}(gh)_\ast (\mu ).
\end{align*}
Let $\{ g_i \}, \{ h_j\}$ be complete systems of representatives of $G/H$ and $H/J$ respectively. Then $\{ g_ih_j\}$ is a complete system of representatives of $G/J$. Hence
\[ \iota_H\circ \iota_J^H (\mu )=\sum_{i,j}(g_ih_j)_\ast (\mu )=\sum_{gJ\in G/J}g_\ast (\mu )=\iota_J (\mu ).\]
Then, we can see that
\[ \iota_J ( \SC (J ))= \iota_H\circ \iota_J^H (\SC (J ))\subset \iota_H (\SC (H )).\]
Moreover, we have
\[ \iota_H ( \eta_J^H )=\iota_H\circ \iota_J^H (\eta_J^J )=\eta_J .\]
Since $\iota_H$ is $\RRR$-linear, $\iota_H$ maps a rational subset current on $H$ to a rational subset current on $G$.
As a result of the above, we obtain the following theorem:

\begin{theorem}
Let $H$ be an infinite quasi-convex subgroup of an infinite hyperbolic group $G$.
Then $\iota_H$ is a continuous $\RRR$-linear map from $\SC(H)$ to $\SC (G)$ satisfying the condition that for every quasi-convex subgroup $J$ of $H$ we have
\[ \iota_H (\eta_J^H) =\eta_J.\]
If $H$ has the denseness property of rational subset currents, then such a map is unique. 
Moreover, if $gHg^{-1}\cap H=\{ \id \}$ for any $gH \ (\not=H )\in G/H$, then $\iota_H$ is injective.
\end{theorem}
\begin{proof}
We prove only the last claim.
Assume that $gHg^{-1}\cap H=\{ \id \}$ for any $gH \ (\not=H )\in G/H$.
Let $\mu \in \SC (H)$. 
Note that for $g\in G$ the support of $g_\ast (\mu )$ is included in $g(\H (\partial H ))=\H (g\gL (H))$ and
\[ g\gL (H ) \cap \gL (H ) =\gL (gHg^{-1}\cap H).\]
Hence from the assumption the support of $g_\ast (\mu )\ (gH \in G/H)$ is pairwise disjoint.
Hence $\iota_H$ is injective.
\end{proof}

\begin{example}\label{exa:doubled surface inclusion}
We give a example of a case when $\iota_H$ is injective.
Let $\gS$ be a compact hyperbolic surface and $\gS_0$ a compact connected subsurface of $\gS$. Assume that any boundary component of $\gS_0$ is a simple closed geodesic of $\gS$.
Then the fundamental group $H:=\pi_1(\gS_0)$ of $\gS_0$ is considered as a subgroup of $G:=\pi_1(\gS)$ and we obtain the map $\iota_H\: \SC(H)\rightarrow \SC(G)$.
In this setting, we can see that $gHg^{-1}\cap H=\{ \id \}$ for any $gH \ (\not=H )\in G/H$ and so $\iota_H$ is injective.
One of the ways to prove this claim is using the cubic commutative diagram in the above of Definition \ref{def:size of cpt comp} and considering the fiber product $C_H\times_{\gS}C_H$, which is homeomorphic to $C_H$.

The inclusion map $\iota_H$ is especially useful in the case that $\gS$ is the doubled surface $D\gS_0$ of $\gS_0$.
For a compact hyperbolic surface $S$ with boundary, the doubled surface $DS$ of $S$ is obtained by gluing two copies of $S$ together at corresponding boundary components. In the context of geodesic currents, this case is already known (see \cite[Subsection 2.6]{DLR10}) and the map $\iota_H$ introduced in this section is a natural extension.
\end{example}

Let $H$ be a finite index subgroup of $G$. We denote by $[G:H]$ the index of $H$ in $G$. A subset current on $G$ can be considered as a subset current on $H$ since $\H (\partial G)=\H (\partial H)$.
Therefore $\SC (G)$ can be considered as an $\RRR$-linear subspace of $\SC (H)$.
Moreover, for $\mu \in \SC (G)$ we have
\[ \iota_H (\mu )=\sum_{gH\in G/H}g_\ast (\mu )=\sum_{gH\in G/H}\mu =[G:H] \mu. \]
Then we have the following theorem.
\begin{theorem}\label{thm:denseness and finite index}
Let $H$ be a finite index subgroup of an infinite hyperbolic group $G$.
Then $\iota_H$ is surjective.
Moreover, if $H$ has the denseness property of rational subset currents, then $G$ also has the denseness property of rational subset currents.
\end{theorem}
\begin{proof}
Take any $\mu \in \SC (G)\ (\subset \SC(H))$.
Then we see that
\[ \iota_H \left(\frac{1}{[G:H]}\mu \right)=\frac{1}{[G:H]} [G:H]\mu =\mu ,\]
which implies that $\iota_H\: \SC(H)\rightarrow \SC(G)$ is surjective.

By considering $\mu$ as a subset current on $H$ we can take a sequence of rational subset currents $\{ \mu_n \}$ on $H$ such that $\{ \mu _n \}$ converges to $\mu$. Since $\iota_H$ is continuous, $\{ \iota_H(\mu_n)\}$ converges to $\iota_H( \mu )=[G:H] \mu$. Since $\{ \iota_H (\mu_n)\}$ is a sequence of rational subset currents on $G$, the sequence
\[ \left\{ \frac{1}{[G:H]}\iota_H(\mu_n) \right\}\]
is a sequence of rational subset currents on $G$ converging to $\mu$.
\end{proof}

\begin{remark}
Recall that $\mathrm{Span}(\SC_r(G))$ is the $\RRR$-linear subspace of $\SC(G)$ generated by the set $\SC_r(G)$ of all rational subset currents on $G$.
Even if we consider $\mathrm{Span}(\SC_r(H))$ and $\mathrm{Span}(\SC_r(G))$ instead of $\SC_r(H)$ and $\SC_r(G)$ in the above theorem, the same statement follows.
\end{remark}

\subsection{Finite index extension of functionals}\label{subsec:finite index extension}

Let $G$ be an infinite hyperbolic group.
From the previous subsection, for a finite index subgroup $\gG$ of $G$ we can consider $\SC (G)$ as an $\RRR$-linear subspace of $\SC(\gG)$. By using this fact, we provide a method for extending functionals on $\SC (\gG)$ to functionals on $\SC (G)$. Especially, we will consider the case that $\gG$ is a free group of finite rank or a surface group.
In this case we can see that the extension is unique from the denseness property of rational subset currents.

Assume that the hyperbolic group $G$ has a finite index subgroup $\gG$ that is isomorphic to a free group of finite rank or a surface group. For example, a finitely generated Fuchsian group satisfies this property.
From Theorem \ref{thm:subset currents on hyperbolic surface has denseness property} and Theorem \ref{thm:denseness and finite index}, the set of all rational subset currents on $G$ is dense in $\SC (G)$.

\begin{supply}\label{supply:finite index normal subgroup}
(1). Let $H$ be a group. Let $J,K$ be finite index subgroups of $H$.
Then the following formula is well-known:
\[ [J:J\cap K]=[JK:K],\]
where $JK$ may not be a subgroup of $H$ but $JK$ can be represented as a disjoint union of cosets of $K$.
From the above formula we can see that $J\cap K$ is also a finite index subgroup of $H$. Actually, we have
\[ [H:J\cap K]=[H:J][J:J\cap K]=[H:J][JK:K]\leq [H:J][H:K] .\]

(2). Consider the conjugacy class of $\gG$ in $G$,
\[ \mathrm{Conj}(\gG) :=\{ g\gG g^{-1}\mid g\in G \} .\]
Then we have a surjective map $\phi$ from $G/\gG$ to $\mathrm{Conj}(\gG )$, which is defined by $\phi (g\gG ):= g\gG g^{-1}$ for $g\gG \in G/\gG$.
Since $\gG$ is a finite index subgroup of $G$, the cardinality of $\mathrm{Conj}(\gG )$ is also finite.
Note that $g\gG g^{-1}$ for any $g\in G$ is also a finite index subgroup of $G$. Actually, if $G$ is a disjoint union of $g_1\gG,\dots ,g_m\gG$, then $G=gGg^{-1}$ is a disjoint union of
\[ (gg_1g^{-1})g\gG g^{-1},\dots ,(gg_mg^{-1})g \gG g^{-1}.\]
Let $\gG_0$ be the intersection of all $g\gG g^{-1} \in \mathrm{Conj}(\gG)$. Then $\gG_0$ is a finite index normal subgroup of $G$.
Since $\gG_0$ is also a finite index subgroup of $\gG$, the group $\gG_0$ is isomorphic to a free group of finite rank or a surface group.
Therefore, we can take $\gG$ as a finite index normal subgroup of $G$. Then we obtain the exact sequence:
\[ \{ \id \}\rightarrow \gG \rightarrow G \rightarrow G/\gG \rightarrow \{ \id \} ,\]
which implies that $G$ is a finite extension of $\gG$ by $G/\gG$.
\end{supply}

\begin{lemma}\label{lem:finitely generated and quasi-convex}
Let $G$ be a hyperbolic group with a finite index subgroup $\gG$ that is isomorphic to a free group of finite rank or a surface group.
A subgroup $H$ of $G$ is quasi-convex if and only if $H$ is finitely generated.
\end{lemma}
\begin{proof}
The ``only if'' part is known from the property of quasi-convexity.
Assume that $H$ is finitely generated.
The intersection $H\cap \gG$ is a finite index subgroup of $H$ since 
\[ [H:H\cap \gG ]=[H\gG :\gG]\leq [G:\gG ]<\infty .\]
Hence we see that $H\cap \gG$ is finitely generated. Note that a finite index subgroup of a finitely generated subgroup is also finitely generated.
Since $H\cap \gG$ is a finitely generated subgroup of $\gG$, $H\cap \gG$ is a quasi-convex subgroup of $\gG$ and also a quasi-convex subgroup of $G$.
Then $H$ is quasi-isometric to $H\cap \gG$ in $G$, which implies that $H$ is a quasi-convex subgroup of $G$.
\end{proof}

For a quasi-convex subgroup $H$ of $G$, the intersection $H\cap \gG$ is a finite index subgroup of $H$ from the proof of Lemma \ref{lem:finitely generated and quasi-convex}.
Recall that for a finite index subgroup $J$ of $\gG$ the reduced rank of $J$ equals $[\gG:J]\rk(\gG)$, that is,
\[ \rk(\gG )=\frac{1}{[\gG :J]}\rk (J).\]
Now, we define the reduced rank $\rk (H)$ of $H$ by
\[ \rk(H):= \frac{1}{[H:H\cap \gG]} \rk (H\cap \gG ).\]
We check that this definition is independent of the choice of $\gG$.
Take a finite index subgroup $\gG'$ of $G$ isomorphic to a free group of finite rank or a surface group. 
Then we have
\begin{align*}
[H:H\cap \gG][H\cap \gG:H\cap \gG\cap \gG' ]
&=[H:H\cap \gG\cap \gG'] \\
&=[H:H\cap \gG'][H\cap \gG':H\cap \gG\cap \gG' ]
\end{align*}
and so
\begin{align*}
&\frac{1}{[H:H\cap \gG]} \rk (H\cap \gG )\\
=&\frac{1}{[H:H\cap \gG]}\frac{1}{ [H\cap \gG:H\cap \gG\cap \gG' ]} \rk (H\cap \gG \cap \gG' )\\
=&\frac{1}{[H:H\cap \gG \cap \gG']} \rk (H\cap \gG \cap \gG' )\\
=&\frac{1}{[H:H\cap \gG']}\frac{1}{ [H\cap \gG':H\cap \gG\cap \gG' ]} \rk (H\cap \gG \cap \gG' )\\
=&\frac{1}{[H:H\cap \gG']} \rk (H\cap \gG' ).
\end{align*}

Recall that $\gG$ is a finite index subgroup of $G$ isomorphic to a free group or a surface group and we have the reduced rank functional $\rk_\gG$ on $\SC (\gG)$ from Corollary \ref{cor:reduced rank functional}.
We define the reduced rank functional $\rk_G$ on $\SC(G)$ by
\[ \rk_G (\mu ):= \frac{1}{[G:\gG]}\rk_\gG (\mu )\]
for $\mu \in \SC (G)$, that is,
\[\rk_G=\frac{1}{[G:\gG]} \rk_\gG|_{\SC (G)} .\]
Then we have the following theorem:

\begin{theorem}
Let $G$ be a hyperbolic group with a finite index subgroup $\gG$ that is isomorphic to a free group of finite rank or a surface group.
Then the following equality holds on $\SC(\gG)$:
\[ \rk_G\circ \iota_\gG =\rk_\gG.\]
Moreover, $\rk_G$ is a unique continuous $\RRR$-linear functional on $\SC (G)$ satisfying the condition that for every quasi-convex subgroup $H$ of $G$ we have
\[ \rk_G (\eta_H^G )=\rk (H).\]
\end{theorem}
\begin{proof}
First, we consider the case that $\gG$ is a normal subgroup of $G$.
Take a quasi-convex subgroup $H$ of $\gG$ and a complete system of representatives $\{ \gamma _i\}$ of $\gG/H$.
For $g\in G$ the set $\{ g\gamma_i g^{-1}\}$ is a complete system of representatives of $\gG /(gHg^{-1})$ since
\[ \gG =g\gG g^{-1}=g\left( \bigsqcup_{i}\gamma_i H\right) g^{-1}=\bigsqcup_{i} (g\gamma _ig^{-1})gHg^{-1}.\]
Then we have
\[
g_\ast (\eta_H^\gG )
=\sum_{i}g_\ast (\delta_{\gamma_i \gL (H)})
=\sum_{i}\delta_{g\gamma_i \gL (H)}
=\sum_{i}\delta_{g\gamma_i g^{-1} \gL (gHg^{-1})}
=\eta_{gHg^{-1}}^\gG.
\]
Note that $\rk (gHg^{-1})=\rk (H)$.
Therefore 
\begin{align*}
\rk_G\circ \iota_\gG (\eta_H^\gG)
&=\frac{1}{[G:\gG ]}\rk_\gG \left( \sum_{ g\gG \in G/\gG}g_\ast (\eta_H^\gG ) \right)\\
&=\frac{1}{[G:\gG ]}\rk_\gG \left( \sum_{ g\gG \in G/\gG}\eta_{gHg^{-1}}^\gG \right)\\
&=\frac{1}{[G:\gG ]}\sum_{ g\gG \in G/\gG}\rk (gHg^{-1})\\
&=\frac{1}{[G:\gG ]}[G:\gG ]\rk (H)=\rk_\gG (\eta_H^\gG ).
\end{align*}
From the denseness property of rational subset currents for $\gG$ we have $\rk_G\circ \iota_\gG =\rk_\gG$.

From now on, we do not assume that $\gG$ is a normal subgroup of $G$.
We can take a normal subgroup $\gG_0$ of $G$ such that $\gG_0$ is a finite index normal subgroup of $\gG$ from Supplementation \ref{supply:finite index normal subgroup}.
Note that we have $\rk_{\gG}\circ \iota_{\gG_0}^\gG=\rk_{\gG_0}$ from the above argument.
Hence
\[
 (\rk_G\circ \iota_\gG) \circ \iota_{\gG_0}^\gG = \rk_G\circ \iota_{\gG_0}= \rk_{\gG_0}= \rk_\gG \circ \iota_{\gG_0}^\gG.
\]
Since the map $\iota_{\gG_0}^\gG$ from $\SC (\gG_0)$ to $\SC (\gG)$ is surjective by Theorem \ref{thm:denseness and finite index}, we obtain the required equality
\[ \rk_G\circ \iota_\gG =\rk_\gG.\]

Take a quasi-convex subgroup $H$ of $G$. Then $\eta_H^G= \frac{1}{[H:H\cap \gG]}\eta_{H\cap \gG}^G$, and we have
\begin{align*}
\rk_G( \eta_H^G)
&=\frac{1}{[H:H\cap \gG ]}\rk_G (\eta_{H\cap \gG}^G)\\
&=\frac{1}{[H:H\cap \gG ]}\rk_G \circ \iota_{\gG} (\eta_{H\cap \gG}^\gG )\\
&=\frac{1}{[H:H\cap \gG ]}\rk_\gG (\eta_{H\cap \gG}^\gG )\\
&=\frac{1}{[H:H\cap \gG ]}\rk (H\cap \gG )=\rk (H).
\end{align*}
This completes the proof.
\end{proof}


\section{Intersection number}\label{sec:intersection number of subset currents}

\if{o}

compact surface と言ったら連結として closed diskは省いておくべき！
いくつかの定理や補題の証明をすべて後半に回しているため，その説明はしておくべきかもしれない．
\fi

Let $\gS$ be a non-contractible compact surface. A surface which we refer to may have boundary. We always assume that a surface is orientable and connected.
In this section, our goal is to generalize the notion of the (geometric) intersection number of two closed curves on $\gS$ to the intersection number of two ``simple compact surfaces'' on $\gS$ by using the fiber product.
Moreover, we extend the intersection number of two simple compact surfaces to a continuous $\RRR$-bilinear functional on $\SC (\gS)$ in the case that $\gS$ is a compact hyperbolic surface in Subsection \ref{conti ext of int number}.
As a side note, the condition of the compactness of $\gS$ is not essential in the study of the intersection number of simple compact surfaces in Subsection \ref{subsec:Intersection number of surfaces}.

\subsection{Intersection number of closed curves}\label{subsec:intersection number of closed curves}
In this subsection, we review the notion of the intersection number of closed curves on $\gS$.

A continuous map $c\:S^1\rightarrow \gS$ is called a closed curve on $\gS$. For two closed curves $c_1,c_2$ on $\gS$ we will denote by $c_1\times_\gS c_2$ the fiber product corresponding to $c_1$ and $c_2$. Explicitly,
\[ c_1\times_\gS c_2:=\{ (x,y)\in S^1\times S^1 \mid c_1(x)=c_2(y)\} . \]
\begin{supply}\label{supply:fiber product}
Let $X,Y$ and $Z$ be topological spaces. Let $f\:X\rightarrow Z,\ g\:Y\rightarrow Z$ be continuous maps. In the topological category, the \ti{fiber product} $X\times_ZY$ corresponding to $f,g$ is defined to be 
\[ X\times_Z Y:=\{ (x,y)\in X\times Y\mid f(x)=g(y)\} ,\]
equipped with the subspace topology of $X\times Y$. If $Z$ is Hausdorff, then $X\times_Z Y $ is closed since $X\times_Z Y$ is the preimage of the diagonal component of $Z\times Z$ with respect to the map
\[ f\times g\: X\times Y\rightarrow Z\times Z;\ (x,y)\mapsto (f(x),g(y)).\]
Therefore, if $Z$ is Hausdorff and $X,Y$ are compact, then $X\times_Z Y$ is compact.

If $f,g$ are injective, then the map
\[ \phi\: X\times_Z Y\rightarrow f(X)\cap g(Y);\ (x,y)\mapsto f(x)\]
is a bijective continuous map. 
Therefore, if $c_1,c_2$ are simple closed curves, then $c_1\times_\gS c_2$ is homeomorphic to $c_1(S^1)\cap c_2(S^1)$.
More generally, if $f,g$ are embedding maps, then $\phi$ is a homeomorphism. In fact, the maps $f^{-1}|_{f(X)\cap g(Y)}, g^{-1}|_{f(X)\cap g(Y)}$ are continuous maps from $f(X)\cap g(Y)$ to $X$ and $Y$, which induce a continuous map from $f(X)\cap g(Y)$ to $X\times_Z Y$, and this map is the inverse map of $\phi$.
\end{supply}
We will denote by $[c]$ the homotopy class of a closed curve $c$ on $\gS$. We say that a closed curve $c$ is nullhomotopic if $c$ is homotopic to a constant map.

\begin{definition}[Intersection number of two closed curves]\label{def:intersection number of closed curves}
Let $c_1,c_2$ be closed curves on $\gS$. The \ti{intersection number} $i(c_1,c_2)$ of $c_1,c_2$ is the number of contractible components of the fiber product $c_1\times_\gS c_2$. We define the intersection number $i([c_1],[c_2])$ of two homotopy classes $[c_1],[c_2]$ by
\[ i([c_1],[c_2]):=\min_{c_1'\in [c_1],c_2'\in [c_2]}i(c_1',c_2').\]
If $c_1'\in [c_1]$ and $c_2'\in [c_2]$ attain the minimum of the intersection number of two homotopy classes, we say that $c_1'$ and $c_2'$ are \ti{in minimal position}.
\end{definition}

Note that a closed curve on $\gS$ has an orientation induced by an orientation of $S^1$ but we usually do not care about the orientation since it does not influence the intersection number.

For a closed curve $c$ on $\gS$ and $m\in \NN$, we have the closed curve $c^m$ on $\gS$, which can be considered as an $m$-fold covering of $c$.
For another closed curve $c'$ on $\gS$ we have $i(c^m,c')=m \cdot i(c, c')$.
We say that two closed curves $c_1,c_2$ on $\gS$ \ti{virtually coincide} if there exist a closed curve $c$ on $\gS$ and $m_1,m_2\in \NN$ such that $c_i$ coincides with $c^{m_i}$ up to reparametrization for $i=1,2$.

We usually consider only the case that two closed curves on $\gS$ intersect transversely or virtually coincide if they intersect .
When we say that two closed curves on $\gS$ are transverse, we allow them to virtually coincide.

From the above definition of the intersection number, it is natural to ask when two closed curves are in minimal position. The bigon criterion is one of the well-known answer.

\begin{definition}[Bigon and immersed bigon]\label{def:bigon}
A \ti{bigon} is a closed disk $D$ with two subsets $e_1,e_2$ of $\partial D$, called \ti{edges}, satisfying the condition that each of $e_1$ and $e_2$ is homeomorphic to a closed interval, $\partial D=e_1\cup e_2$ and $e_1\cap e_2$ is two points, called \ti{vertices}.

Let $I_1,I_2$ be closed intervals of $\RR$. Let $f_i$ be a continuous map from $I_i$ to a 2-dimensional manifold $M$ possibly with boundary ($i=1,2$). We say that $f_1$ and $f_2$ form an \ti{immersed bigon} if there exists a locally injective continuous map $b$ from a bigon $D$ into $M$ such that there exists a homeomorphism $b_i$ from the edge $e_i$ of $D$ to $I_i$ and $f_i\circ b_i$ coincides with the restriction of $b$ to $e_i$ for $i=1,2$. In this case we say that $f_1$ and $f_2$ form an immersed bigon $b$.
If $b$ is an embedding map, then we say that $f_1$ and $f_2$ form a \ti{bigon} $b$.

\[
\xymatrix{
 & &I_i  \ar[d]^{f_i}\\
e_i\ar@<0.8ex>[urr]^{b_i}\ar@{^{(}->}[r]& D\ar[r]^{b} &M}
\]

A \ti{sub-arc} of a continuous map $f$ from a 1-dimensional manifold $I$ possibly with boundary to a topological space is the restriction of $f$ to a subset of $I$ that is homeomorphic to a closed interval.
We say that a sub-arc of a closed curve $c$ forms a closed curve if the image of the endpoints of the sub-arc is one point.

Let $c_1,c_2$ be closed curves on $\gS$. Let $p:\RR \rightarrow S^1$ be a universal covering of $S^1$.
We say that $c_1$ and $c_2$ form an immersed bigon if there exist sub-arcs $p_1,p_2$ of $p$ such that $c_1\circ p_1$ and $c_2\circ p_2$ form an immersed bigon.
We say that $c_1$ and $c_2$ form a bigon if there exist sub-arcs $f_1,f_2$ of $c_1,c_2$ such that $f_1$ and $f_2$ form a bigon.
\end{definition}

\begin{example}
See Figure \ref{immersed bigon but not bigon}.
Two closed curves on a closed surface of genus 2 form an immersed bigon but do not form a bigon. The points $p,q$ are the vertexes of the immersed bigon. The intersection number of those closed curves is $2$ but they are not in minimal position.
By ``enlarging'' the inner simple closed curve, the intersection number of those closed curves will be $0$.

\begin{figure}[h]
\begin{center}
\includegraphics[width=5.0cm]{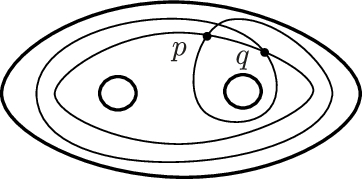}
\vspace{-0.3cm}
\caption{Two closed curves form an immersed bigon but do not form a bigon.}\label{immersed bigon but not bigon}
\end{center}
\end{figure}

\end{example}

Let $c_1,c_2$ be transverse closed curves on $\gS$ such that sub-arcs $f_1,f_2$ of $c_1,c_2$ form a bigon.
We can modify $f_1$ by a homotopy in the bigon such that $f_1$ and $f_2$ coincides, and then we can modify a neighborhood of $f_1$ by a homotopy such that $f_1$ and $f_2$ are disjoint.
Therefore, if two transverse closed curves form a bigon, then we can reduce the intersection number by a homotopy.
The following lemma says that the converse is also true in the case that closed curves are simple.

\begin{lemma}[The bigon criterion 1 (see {\cite[Proposition 1.7]{FM12})}]\label{lem:bigon criterion 1}
Let $c_1,c_2$ be transverse simple closed curves on $\gS$. Then two simple closed curves $c_1,c_2$ do not form a bigon if and only if $c_1,c_2$ are in minimal position.
\end{lemma}

In the case that a closed curve $c$ is not simple, $c$ can have a sub-arc which forms a nullhomotopic closed curve on $\gS$.
Such a nullhomotopic closed curve is easy to reduce (but difficult to deal with), so we usually assume that a non-simple closed curve do not have a sub-arc forming a nullhomotopic closed curve.

\begin{lemma}[The bigon criterion 2]\label{lem:bigon criterion 2}
Let $c_1,c_2$ be transverse closed curves on $\gS$. Assume that no sub-arc of $c_i$ forms a nullhomotopic closed curve on $\gS$ for $i=1,2$. Then $c_1,c_2$ do not form an immersed bigon if and only if $c_1,c_2$ are in minimal position.
\end{lemma}

We can obtain Lemma \ref{lem:bigon criterion 2} as a corollary to Theorem \ref{thm:bigon criterion 3}, which we will prove later.

Recall that any non-nullhomotopic closed curve on a surface with a Riemannian metric of constant curvature 0 or $-1$ is homotopic to a closed geodesic on the surface. Especially, in the case that the constant curvature of the surface is $-1$, which we call \ti{a hyperbolic surface}, such a closed geodesic is unique. We assume that every boundary component of a compact hyperbolic surface is totally geodesic.
When we consider a geodesic on $\gS$, we always assume that $\gS$ is equipped with a Riemannian metric with constant curvature 0 or $-1$, which implies that $\gS$ is not a sphere.
The following theorem is well-known, which is a direct corollary to Lemma \ref{lem:bigon criterion 2}.

\begin{theorem}\label{thm:intersection number of geodesics}
Two closed geodesics on $\gS$ are in minimal position.
\end{theorem}

We can see that our definition of the intersection number works effectively in the case that two closed geodesics coincide (cf. Example \ref{exa:non-simple closed geodesic and intersection number}). The half of $i([c],[c])$ is called the \ti{self-intersection number} of a closed curve $c$ on $\gS$, which coincides with the half of $i(c',c')$ for a closed geodesic $c'$ homotopic to $c$ if $c$ is not nullhomotopic.

\begin{example}\label{exa:non-simple closed geodesic and intersection number}
We see an example of a closed geodesic with self-intersection, which means that the self-intersection number of the closed geodesic is positive.
\begin{figure}[h]
\begin{center}
\includegraphics[width=9cm]{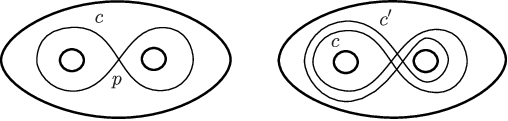}
\vspace{-0.3cm}
\caption{The closed geodesic $c$ has one self-intersection in the left figure. In the right figure, the intersection number of $c$ and $c'$, which is homotopic to $c$, equals $2$.}\label{geodesic and self-intersection}
\end{center}
\end{figure}
In the left figure of Figure \ref{geodesic and self-intersection}, the closed geodesic $c$ on a closed hyperbolic surface $\gS$ of genus $2$ is not simple. Take $a,b\in S^1$ with $a\not=b$ such that $c(a)=c(b)=p$. Then we can see that the fiber product $c \times_\gS c$ equals
\[ \{ (x,x)\in S^1\times S^1 \mid x\in S^1\} \sqcup \{ (a,b)\} \sqcup \{ (b,a)\} ,\]
which includes two contractible components $\{ (a,b)\}$ and $\{ (b,a)\}$.
Hence the intersection number of $c$ and $c$ equals $2$ and the self-intersection number of $c$ equals $1$.
In the right figure of Figure \ref{geodesic and self-intersection}, two closed curves $c$ and $c'$ are in minimal position since they do not form any immersed bigon. The closed curve $c'$ is homotopic to $c$ and the intersection number of $c$ and $c'$ equals $2$.
\end{example}

\subsection{Intersection number of surfaces}\label{subsec:Intersection number of surfaces}

Now, we generalize the notion of the intersection number of two closed curves on a compact orientable surface $\gS$ to the intersection number of ``two simple compact surfaces'' on $\gS$, and we prove ``the bigon criterion'' for this intersection number (see Theorem \ref{thm:bigon criterion 3}).

\begin{definition}[Simple compact surfaces]
Let $S$ be a compact surface possibly with boundary or $S=S^1$.
A \ti{simple compact surface} on $\gS$ is a pair of $S$ and a continuous map $s\: S\rightarrow \gS$ satisfying the following condition:
\begin{enumerate}
\item $s$ is locally injective;
\item the restriction of $s$ to each component of the boundary $\partial S$ of $S$ is not nullhomotopic and does not have a sub-arc forming a nullhomotopic closed curve on $\gS$.
\end{enumerate}

If $S=S^1$, then we regard the boundary $\partial S$ as $S$. Here, we remark that in the case that $S=S^1$, a simple compact surface $(S,s)$ on $\gS$ may not be a simple closed curve on $\gS$. We will often write $s$ in place of $(S,s)$ for simplicity.

A simple compact surface $(S_1,s_1)$ is said to be homotopic to a simple compact surface $(S_2,s_2)$ if there exist a homeomorphism $f\:S_1\rightarrow S_2$ and a continuous function $F\: S_1\times [0,1]\rightarrow \gS$ such that $F(\cdot, 0)=s_1$ and $F(f^{-1}(\cdot),1)=s_2$.
Being homotopic is an equivalence relation and the equivalence class of a simple compact surface $(S,s)$, called a homotopy class, is denoted by $[S,s]$ or $[s]$. Note that if $S=S^1$, then changing the orientation of $(S,s)$ does not change the homotopy class of $(S,s)$.

Let $(S_1,s_1),(S_2,s_2)$ be two simple compact surfaces on $\gS$.
We say that $(S_1,s_1) $ and $(S_2,s_2)$ are transverse if the restriction of $s_1$ and $s_2$ to any components of their boundaries intersect transversely or virtually coincide if they intersect. We consider only the case that two simple compact surfaces are transverse.
\end{definition}

\begin{definition}[(Geometric) Intersection number of simple compact surfaces]\label{def:intersection number of simple compact surfaces}
Let $(S_1,s_1),(S_2,s_2)$ be two simple compact surfaces on $\gS$.
The \ti{intersection number} $i(s_1,s_2)$ of $s_1$ and $s_2$ is the number of contractible components of the fiber product $S_1\times_\gS S_2$ corresponding to $s_1,s_2$. We define the intersection number $i([s_1],[s_2])$ of two homotopy classes $[s_1],[s_2]$
to be the minimum of $i(s_1',s_2')$ taken over $s_1'\in [s_1]$ and $s_2'\in [s_2]$ that are transverse.
If two transverse simple closed surfaces $s_1'\in [s_1]$ and $s_2'\in [s_2]$ attain the minimum of the intersection number of two homotopy classes, we say that $s_1'$ and $s_2'$ are \ti{in minimal position}.
\end{definition}

Note that for a simple compact surface $(S,s)$ on $\gS$, the surface $S$ can not be a closed disk by the definition.
We usually do not consider the case that $\gS$ is a sphere because we see that the intersection number of simple compact surfaces on a sphere always equals 0 in Lemma \ref{prop:simple compact surface on sphere}.

We will give some examples of the calculation of the intersection number of simples compact surfaces in Example \ref{example:surface and closed geodesic intersection} and \ref{example:intersection number 0}.

In the definition of a simple compact surface $(S,s)$ on $\gS$, the required property of the continuous map $s$ seems to be strict.
However, in the following example, we will see that if $s$ does not have this property, then the definition of the intersection number does not work well.

\begin{example}\label{exa:simple condition is necessary}
First, we consider a simple model of the fiber product of two 2-dimensional manifolds over a 2-dimensional manifold. Set $X:=[-1,1]\times \RR, Y:=\RR \times [-1,1]$ and $Z:=\RR^2$. The fiber product corresponding to the inclusion maps from $X,Y$ to $Z$ is homeomorphic to $X\cap Y=[-1,1]\times [-1,1]$, which implies that the number of contractible components of the fiber product is one. Let $i_Y$ be the inclusion map from $Y$ to $Z$. We define a continuous map $f\:X\rightarrow Z$ to be
\[
f(x) := \begin{cases}
 x & (\| x \| \geq 1) \\
 x + (0,2(1-\| x\| )) & (\| x \| \leq 1)
 \end{cases}
\]
for $x\in X$, where $\| \cdot \|$ is the Euclidean norm.
We can see that the following map $F\: X\times [0,1]\rightarrow Z$ is a homotopy from the inclusion map to $f$:
\[
F(x,t) := \begin{cases}
 x & (\| x \| \geq 1) \\
 x + (0,2t(1-\| x\| )) & (\| x \| \leq 1)
 \end{cases}
\]
for $(x,t)\in X\times [0,1]$.
We consider the fiber product $X\times_Z Y$ corresponding to $f,i_Y$ and want to say that $X\times_Z Y$ is connected and not contractible.
Then we can see that that we can reduce the number of contractible components of the fiber product of two dimensional spaces by a homotopy which deforms a ``local'' part of one of the spaces. 

Let $p_X$ be the natural projection from $X\times_Z Y$ onto $X$, that is, $p_X$ maps $(x,y)\in X\times_Z Y$ to $x\in X$.
We can see that $p_X(X\times_Z Y)=f^{-1}(Y)$, which includes the unit circle $S^1=\{ x\in X\mid \| x\| =1\}$ but does not contain $(0,0)$ since $f(0,0)=(0,2)\not\in Y$.

Now, we consider a closed curve $c:S^1\rightarrow X \times_Z Y$ defined by $c(x)=(x,x)$ for $x\in S^1$. Then we can see that $p_X\circ c$ is not nullhomotopic in $f^{-1}(Y)$ since $(0,0)\not\in f^{-1}(Y)$, which implies that $c$ is not nullhomotopic in $X \times_Z Y$.

Finally, we check that $X \times_Z Y$ is connected. Take any $(x,y)\in X \times_Z Y$. If $|| x || \geq 1$, then $x=y$ and a line segment joining $1/|| x|| x$ to $x$ induces a path joining $1/ ||x|| (x,x)$ to $(x,x)$.
We consider the case that $|| x|| <1$. Note that $x\not = (0,0)$ and $y=f(x)$. Set
\[ a(t):=\frac{1}{(1-t)||x|| +t} \]
for $t\in [0,1]$. Then we can see that the path $ (a(t)x, f(a(t)x))\in X\times_Z Y$ for $t\in [0,1]$ joining $1/||x|| (x,x)\in c(S^1)$ to $(x,y)$.
Therefore $X \times_Z Y$ is path-connected.

See Figure \ref{intersection number of simple compact surfaces} and its caption.
The pairs $(S,s_1)$ and $(S,s_2)$ are simple compact surfaces on $\gS$.
Then we can see that the fiber product of $(S,s_1)$ and $(S,s_2)$, which is homeomorphic to $s_1(S)\cap s_2(S)$, includes two contractible components, which implies that the intersection number of $(S,s_1)$ and $(S,s_2)$ equals $2$.
However, we can modify $s_1$ (or $s_2$) and obtain $s_1'$ by the same way as above, then the intersection number of $(S,s_1')$ and $(S,s_2)$ will be $0$. This contradicts Theorem \ref{thm:bigon criterion 3}, since Theorem \ref{thm:bigon criterion 3} implies that $(S,s_1)$ and $(S,s_2)$ are in minimal position.
The reason is that $s_1'$ is not locally injective.
From the above, we see that the condition of local injectivity for a simple compact surface is necessary.

\begin{figure}[h]
\begin{center}
\includegraphics[width=9cm]{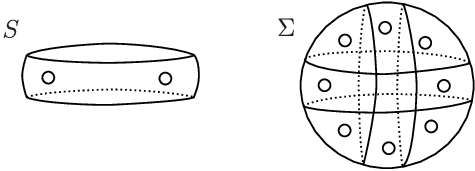}
\vspace{-0.3cm}
\caption{The left figure shows a compact surface $S$ of genus $2$ with $2$ boundary components. The right figure shows two inclusion maps $s_1$ and $s_2$ from $S$ to a closed surface $\gS$ of genus $8$. We see that $i(s_1,s_2)=2$.}\label{intersection number of simple compact surfaces}
\end{center}
\end{figure}

\end{example}

We define the notion of an immersed bigon formed by two simple compact surfaces on $\gS$ in order to characterize the condition that they are in minimal position.
\begin{definition}[Bigon formed by simple compact surfaces]\label{def:surface bigon}
Let $(S_1,s_1)$ and $(S_2,s_2)$ be simple compact surfaces on $\gS$.
We say that $s_1$ and $s_2$ form an immersed bigon if there exist components $B_1,B_2$ of $\partial S_1, \partial S_2$ such that $s_1|_{B_1}$ and $s_2|_{B_2}$ form an immersed bigon.
\end{definition}

Proving the following theorem is our goal in this subsection.

\begin{theorem}[The bigon criterion 3]\label{thm:bigon criterion 3}
Let $(S_1,s_1),(S_2,s_2)$ be transverse simple compact surfaces on a compact surface $\gS$.
If $s_1$ and $s_2$ do not form an immersed bigon, then $s_1,s_2$ are in minimal position.
If either $S_1$ or $S_2$ is $S^1$, then the converse is also true.
\end{theorem}

The following theorem, which is a corollary to the above theorem, is a generalization of Theorem \ref{thm:intersection number of geodesics}.

\begin{theorem}\label{thm:geodesic criterion}
Let $\gS$ be a torus or a compact hyperbolic surface.
Let $(S_1,s_1),(S_2,s_2)$ be simple compact surfaces on $\gS$. If the restriction of $s_i$ to each component of $\partial S_i$ is a closed geodesic on $\gS$ for $i=1,2$, then $s_1$ and $s_2$ are in minimal position.
\end{theorem}
\begin{example}\label{example:surface and closed geodesic intersection}
Consider the case that $\gS$ is a compact hyperbolic surface.
Recall that for a non-trivial finitely generated subgroup $H$ of the fundamental group $\pi_1(\gS)$ of $\gS$ we have the convex core $C_H$ and the canonical map $p_H\: C_H \rightarrow \gS$ induced by the universal covering map (see the beginning part of Section \ref{sec:Volume functionals for Kleinian groups}).
Then $(C_H,p_H)$ is a simple compact surface on $\gS$ satisfying the condition that the restriction of $p_H$ to each component of $\partial C_H$ is a closed geodesic on $\gS$.
Moreover, we will prove that any simple compact surface $(S,s)$ on $\gS$ that is not a cylinder is homotopic to a convex core $(C_H,p_H)$ for a finitely generated subgroup $H$ of $\pi_1(\gS)$ in Proposition \ref{prop:characterize simple cpt surf}.
We will also see that for non-trivial finitely generated subgroups $H$ and $K$ of $\pi_1(\gS)$ the ``shape'' of a contractible component $J$ of the fiber product $C_H\times_\gS C_K$ can be classified into the following three cases:

\begin{itemize}
\item $J$ is a point if and only if both of $H$ and $K$ are cyclic;
\item $J$ is a geodesic segment if and only if $H$ is cyclic and $K$ is not cyclic, or vice versa;
\item $J$ is a compact hyperbolic polygon (with geodesic boundary) if and only if both of $H$ and $K$ are non-cyclic.
\end{itemize}

Consider the case that finitely generated subgroups $H$ and $K$ of $\pi_1(\gS)$ are as shown in Figure \ref{intersection number formula}.
Then the fiber product $C_H\times_\gS C_K$ is identified with the intersection of $p_H(C_H)\cap p_K(C_K)$, which is a sub-arc $\alpha$ of $C_K$.
Hence the intersection number $i(C_H ,C_K)$ equals $1$.
From Figure \ref{intersection number formula} we can see that two endpoints of the sub-arc $\alpha$ are the intersection of the boundary component $c$ of $C_H$ and $C_K$. Note that the number of the endpoints of $\alpha$ is the double of the intersection number $i(C_H, C_K)$.

In general, consider the case that $C_H$ is a compact surface with boundary and $C_K$ is a closed geodesic. Then each contractible component of $C_H\times_\gS C_K$ is an arc, whose two endpoints lie on a boundary component of $C_H$ and can be identified with the intersection of a boundary component $C_H$ and $C_K$.
Hence we have the following equation:
\[ i(C_H, C_K)=\frac{1}{2}\sum_{c\in \partial C_H }i(c, C_K),\]
where $\partial C_H$ is the set of boundary components of $C_H$.
In Section \ref{sec: projection B}, we generalize the correspondence between $C_H$ and ``the half'' of $\partial C_H$ to the projection from $\SC(\gS)$ to $\GC (\gS)$, and also generalize the above equality to the equality of subset currents in Theorem \ref{thm:projection and inequality}.

If both $C_H$ and $C_K$ are surfaces, then each contractible component of $C_H\times_\gS C_K$ is a (hyperbolic) polygon and the number of the vertices is even and greater than or equal to $4$ (see Lemma \ref{lem:intersection is contractible}). Hence in general the following inequality holds (see Figure \ref{intersection number of simple compact surfaces} for the case that a contractible component is a square and see Theorem \ref{thm:projection and inequality} for detail):
\[ i(C_H, C_K)\leq \frac{1}{4}\sum_{c\in \partial C_H, \gamma \in \partial C_K}i(c, \gamma ).\]

\begin{figure}[h]
\begin{center}
\includegraphics[width=9cm]{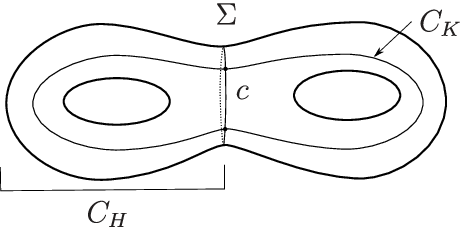}
\vspace{-0.3cm}
\caption{$C_H$ is a subsurface of $\gS$ with the boundary component $c$ and $C_K$ is a closed geodesic on $\gS$.
The maps $p_H$ and $p_K$ are inclusion maps. }\label{intersection number formula}
\end{center}
\end{figure}
\end{example}

\begin{example}\label{example:intersection number 0}
We give another example of the calculation of the intersection number of simple compact surfaces in the case that $\gS$ is a compact hyperbolic surface.
In Subsection \ref{conti ext of int number}, we will consider the intersection number of $C_H$ and $C_K$ for non-trivial finitely generated subgroups $H$ and $K$ of $\pi_1(\gS)$ and extend it to a continuous $\RRR$-bilinear functional $i_{\SC}$ on $\SC (\gS)$, that is,
\[ i_{\SC} (\eta_H ,\eta_K) =i(C_H ,C_K).\]
By the $\RRR$-bilinearlity of $i_{\SC}$, we see that if $H'$ is an $h$-index subgroup of $H$ and $K'$ is a $k$-index subgroup of $K$, then $\eta_{H'}=h\eta_H, \eta_{K'}=k\eta_K$ and
\[ i(C_{H'},C_{K'})=i_{\SC} (\eta_{H'},\eta_{K'})=hk i_{\SC}(\eta_H ,\eta_K )=hk i(C_H, C_K).\]
Especially, if $H$ is a finite-index subgroup of $\pi_1(\gS)$, then $i(C_H, C_K)=0$, since 
\[C_{\pi_1(\gS)}\times_\gS C_K=\gS \times_\gS C_K \cong C_K. \]
Note that in general, for any simple compact surface $(S,s)$ on $\gS$ the fiber product $S\times_\gS \gS$ is homeomorphic to $S$, and so the intersection number $i(s ,\gS)$ equals $0$.
\end{example}

The following lemma, which is intuitively obvious, plays a fundamental role in proving the bigon criterion 3 (Lemma \ref{thm:bigon criterion 3}).
\begin{lemma}\label{lem:bigon easy}
Let $M$ be a contractible 2-dimensional manifold possibly with boundary. Let $I_1,I_2$ be intervals of $\RR$. Let $f_i$ be an embedding map from $I_i$ to $M$ for $i=1,2$. Assume that $f_1,f_2$ are transverse, $f_1(I_1)$ divides $M$ into two connected components $M_1,M_2$ and there exist $a,b,c \in I_2$ with $a<b<c$ such that $f_2(a), f_2(c)\in M_1$ and $f_2(b)\in M_2$. Then there exist sub-arcs of $f_1,f_2$ that form a bigon.
\end{lemma}
\begin{proof}
By the assumption $f_2|_{[a,c]}$ intersects $f_1$ transversely. Then we can take a sub-interval $[a',c']$ of $[a,c]$ containing $b$ such $f_2((a',c'))\subset M_2$ and $f_2(a'),f_2(c')\in f_1(I_1)$, which implies that the union of a sub-arc of $f_1$ and $f_2|_{[a',c']}$ forms a simple closed curve $c$ in $M$. From the Jordan curve theorem, $c$ divides $M$ into two subsets such that one of the subsets is contractible. From the Riemann mapping theorem and Carath\'eodory's theorem, there exists an embedding map $b$ from a closed disk $D$ into $M_1\cup c(S^1)$ such that $b(\partial D)$ coincides with $c(S^1)$.
This completes the proof.
\end{proof}

The following lemma is useful to understand a simple compact surface on $\gS$.
\begin{lemma}\label{lem:locally injective and homeomorphic}
Let $S$ be a compact surface possibly with boundary and $s$ a continuous map from $S$ to $\gS$. If $s$ is locally injective, then the restriction of $s$ to $S\setminus \partial S$ is a local homeomorphism and $s(S\setminus \partial S) \cap \partial \gS =\emptyset$.
\end{lemma}
\begin{proof}
Take $x\in S\setminus \partial S$ and a compact neighborhood $U\subset S\setminus \partial S$ of $x$ such that $U$ is homeomorphic to a closed disk and $s|_U$ is injective. Since $U$ is compact, the map $s|_U \: U\rightarrow s(U)$ is homeomorphic. 
Since $U\setminus \{ x\}$ is non-contractible, so is $s(U)\setminus \{ s(x)\}$, which implies $s(x)\not\in \partial \gS$. Then we can assume that $s(U)$ does not intersect $\partial \gS$.
Since $\partial U$ is homeomorphic to $S^1$, so is $s(\partial U)$. By applying the Jordan curve theorem to $s(\partial U)$ we can see that $s(\partial U)$ divides $\gS$ into two regions $\gS_1,\gS_2$. Assume that $\gS_1$ contains $s(x)$. Then $s(\mathrm{Int}(U))$ coincides with $\gS_1$. Therefore $\mathrm{Int}(U)$ is homeomorphically mapped to $\gS_1$ by $s$, which is an open neighborhood of $s(x)$. Hence our claim follows.
\end{proof}

From the above lemma, we can obtain the following proposition immediately:
\begin{proposition}\label{prop:simple compact surface on sphere}
Let $\gS$ be a sphere and $(S,s)$ a simple compact surface on $\gS$. Then $S$ is also a sphere and $s$ is a homeomorphism from $S$ to $\gS$. Moreover, the intersection number of any two simple compact surfaces on $\gS$ equals zero.
\end{proposition}
\begin{proof}
By the definition of simple compact surfaces on $\gS$, the compact surface $S$ does not have a boundary. By Lemma \ref{lem:locally injective and homeomorphic}, $s$ is a local homeomorphism, which implies that $s$ is a covering map since $S$ is compact. 
Moreover, $s$ must be a homeomorphism from $S$ to $\gS$ since $\gS$ is simply-connected.
Therefore, the fiber product of any two simple compact surfaces on $\gS$ is homeomorphic to a sphere, which implies that the intersection number of these simple compact surfaces equals $0$.
\end{proof}

By the above proposition, any two simple compact surfaces on a sphere are always in minimal position.
From now on, we assume that $\gS$ is not a sphere.

The following lemma related to a bigon and an immersed bigon will be used in Lemma \ref{lem:cubic diagram and immersed bigon}.
\begin{lemma}\label{lem:immersed bigon}
Let $M$ be $\gS$ or the universal cover $\tilde{\gS}$ of $\gS$. Let $b$ be a locally injective continuous map from a closed disk $D$ to $M$. If the restriction of $b$ to the boundary $\partial D$ of $D$ is injective, then so is $b$. Hence $b$ is an embedding map.
\end{lemma}
\begin{proof}
We can assume that $M$ does not have boundary since $M$ can be embedded into a 2-dimensional orientable manifold without boundary whose universal cover is contractible.
It is sufficient to prove that the map $b\: D\rightarrow b(D)$ is a local homeomorphism.
In fact, if $b\: D\rightarrow b(D)$ is a local homeomorphism, then we can see that $b\: D\rightarrow b(D)$ is a covering map since $D$ is compact.
Note that $b|_{\partial D}$ is injective. Hence $b\: D\rightarrow b(D)$ is a homeomorphism.

First, we consider the case that $M$ is $\tilde{\gS}$.
Since $M$ does not have boundary, $M$ is homeomorphic to $\RR^2$.
From the Jordan curve theorem, $b(\partial D)$ divides $M$ into the interior region $M_1$ and the exterior region $M_2$ of $M$.
We prove that $b(\mathrm{Int}(D))=M_1$. Assume that $b(x)\in M_2$ for some $x\in \mathrm{Int}(D)$. Since $M_2$ is path-connected, we can take a path $\ell \:[0,1]\rightarrow M_2$ such that $\ell (0)=b(x)$ and $\ell (1)\not \in b(D)$. Let $t$ be the maximum of $\ell ^{-1}(b(D))$.
Then $t\in (0,1)$, $\ell(t)\in b(D)$ and there exists $y\in \mathrm{Int}(D)$ such that $b(y)=\ell(t)$. By Lemma \ref{lem:locally injective and homeomorphic}, $b(y)$ is an interior point of $b(D)$, which contradicts that $t$ is the maximum of $\ell^{-1}(b(D))$. Therefore $b(\mathrm{Int}(D))\subset M_1$.

To see that $M_1\subset b(\mathrm{Int}(D))$, assume that there exists $z\in M_1$ such that $z\not\in b(D)$.
Since $\partial D$ is a contractible closed curve in $D$, $b(\partial D)$ is also a contractible closed curve in $b(D)$, which contradicts that $z\not\in b(D)$. Hence $b(\mathrm{Int}(D))=M_1$.

Take any $x\in \partial D$. We prove that the map $b\: D\rightarrow b(D)$ is a local homeomorphism at $x$. Take an open neighborhood $V$ of $x$ in $D$ such the restriction of $b$ to $\ol{V}$ is a homeomorphism onto $b(\ol{V})$ and $\ol{V}$ is homeomorphic to a closed disk.
Then we can take a contractible open neighborhood $U$ of $b(x)$ such that $(U\setminus M_2)\subset b(V)$ since $b(D)=b(\partial D)\sqcup M_1$. Set $W:=b^{-1}(U)\cap V$. Then $W$ is an open neighborhood of $x$ and 
\[ b(W)=U\cap b(V)=U\setminus M_2=U\cap b(D)\]
is an open subset of $b(D)$. Hence $b$ is a local homeomorphism onto $b(D)$.

In the case that $M$ is $\gS$, we take a lift $\tilde{b}\: D\rightarrow \tilde{\gS}$ of $b$ with respect to the universal covering $\pi\: \tilde{\gS}\rightarrow \gS$.
Then $\tilde{b}|_{\partial D}$ is injective and $\tilde{b}$ is a local homeomorphism onto $b(D)$. This implies that $b=\pi \circ \tilde{b}$ is a local homeomorphism onto $b(D)$.
\end{proof}

The following lemma, which characterizes simple compact surfaces on $\gS$, will play a fundamental role in proving Theorem \ref{thm:bigon criterion 3}.

\begin{lemma}\label{lem:good compact surface}
Let $(S,s)$ be a simple compact surface on $\gS$. Then there exist a covering $s' \: S'\rightarrow \gS$ and an embedding map $f$ from $S$ into $S'$ such that $f$ is a homotopy equivalence and $s=s'\circ f$. Moreover, the embedding map $f$ lifts to an embedding map from the universal cover $\tilde{S}$ of $S$ into the universal cover $\tilde{\gS}$ of $\gS$.
\end{lemma}
\begin{proof}
Let $p\: \tilde{S}\rightarrow S$ be the universal covering of $S$ and $\pi \: \tilde{\gS}\rightarrow \gS$ the universal covering of $\gS$. Take a base point $x$ of $\gS$ such that $x\in s(S)$.
Take base points $\tilde{y}\in \tilde{S}$, $y\in S$ and $\tilde{x}\in \tilde{\gS}$ such that $p(\tilde{y})=y,\ s(y)=x$ and $\pi(\tilde x)=x$.
Then we have a lift $\tilde s \:(\tilde{S}, \tilde y)\rightarrow (\tilde{\gS},\tilde x)$ of the map $s\circ p\: (\tilde S, \tilde y)\rightarrow (\gS ,x)$ with respect to $\pi$. Then we obtain the following commutative diagram of based topological spaces:

\[
\xymatrix{
(\tilde S, \tilde y) \ar[d]_p\ar[r]^{\tilde s} & (\tilde \gS ,\tilde x) \ar[d]^{\pi}\\
(S, y) \ar[r]_s & (\gS ,x)\ar@{}[lu]|{\circlearrowright} }
\]

If $\partial S=\emptyset$, then $s\: S\rightarrow \gS$ is a covering map and our statement follows immediately.
If $S=S^1$, then we can see that $\tilde{S}$ is homeomorphic to $\RR$ and the lift $\tilde s$ is an embedding map since no sub-arc of the closed curve $s$ forms a nullhomotopic closed curve. In this case, there exists $g\in \pi_1( \gS , x)$ corresponding to $(S,s)$ such that $\langle g \rangle $ acts on $\tilde{s}(\tilde S)$. Let $S'$ be the quotient space $\langle g \rangle \backslash \tilde{\gS}$ and $s'$ the covering map from $S'$ to $\gS$ induced by the universal covering map $\pi$. Then $\tilde{s}$ induces an embedding map $f$ from $S$ to $S'$. Hence our claim follows.

From now on, we assume that $S\not= S^1$ and $\partial S \not=\emptyset$.

\underline{Step 1. Construct $S'$:} 
Let $B_1,\dots B_m$ be all connected components of $\partial S$. We can consider $c_j:=s|_{B_j}$ as a closed curve on $\gS$ since $B_j$ is homeomorphic to $S^1$ for $j=1,\dots , m$. For each $B_j$ we can take a component $\tilde{B_j}$ of $\partial \tilde{S}$ such that the restriction of $\tilde{s}$ to $\tilde{B_j}$ is a universal covering of $B_j$.
We will apply the same argument as that for $S=S^1$ to each $c_j$.
Set $\tilde{c_j}:=\tilde{s}|_{B_j}$, which is an embedding map from $\tilde{B_j}$ into $\tilde \gS$.
We endow $B_j$ with an orientation such that the left-hand side of $B_j$ is the interior of $S$, which induces the orientation of $\tilde{B_j}$ and $\tilde{c_j}(\tilde{B_j})$.
Let $U_j\subset \tilde{\gS}$ be the right-hand side of $\tilde{c_j}(\tilde{B_j})$ including $\tilde{c_j}(\tilde{B_j})$. Note that if $c_j(B_j)\subset \partial \gS$, then $U_j=\tilde{c_j}(\tilde{B_j})$. Since $\tilde{c_j}$ is a lift of $c_j$, there exists $g_j\in \pi_1(\gS, x)$ corresponding to $c_j$ such that $\langle g_j \rangle$ acts on $\tilde{c_j}(\tilde{B_j})$ and also acts on $U_j$.
Set $L_j:=\langle g_j \rangle \backslash U_j$. Now, we obtain $S'$ by gluing $S$ to $L_j$ along $B_j$ and $\langle g_j \rangle \backslash \tilde{c_j}(\tilde{B_j})$ for $j=1,\dots ,m$. Since $U_j$ is a subset of $\tilde{\gS}$, the universal covering map induce the map $\pi_j$ from $L_j$ to $\gS$. Then by gluing those maps $\pi_1, \dots , \pi_m$ and $s$, we obtain a continuous map $s'$ from $S'$ to $\gS$.

\underline{Step 2. Prove the map $s'\:S'\rightarrow \gS$ is a covering map:}
Take $z\in \gS$.
We prove that there exists a connected open neighborhood $W$ of $z$ such that the restriction of $s'$ to every connected component of $s'^{-1}(W)$ is a homeomorphism onto $W$.

First, we consider the case that $z\in s(S)$ and $s^{-1}(z)\cap \partial S=\emptyset$.
In this case $z\not\in \partial \gS$. Since $S$ is compact and $s$ is locally injective, $s^{-1}(z)$ is a finite set. In fact, if $s^{-1}(z)$ is an infinite set, then $s^{-1}(z)$ has an accumulation point $w$, which contradicts the assumption that $s$ is locally injective.
We can take a contractible open neighborhood $V$ of $z$ such that $s^{-1}(V)\cap \partial S =\emptyset$.
Then the restriction of $\pi$ to each connected component of $\pi^{-1}(V)$ is a homeomorphism onto $V$ and $\pi^{-1}(V)\cap \tilde{s}(\partial \tilde{S})=\emptyset$. Hence the restriction of $\pi_j$ to each connected component of $\pi_j^{-1}(V)\subset L_j$ is a homeomorphism onto $V$ and $\pi_j^{-1}(V)\cap \langle g_j \rangle \backslash \tilde{c_j}(\tilde{B_j})=\emptyset$ for every $j=1,\dots ,m$.

For each $u\in s^{-1}(z)$ we can take a connected open neighborhood $V_u$ of $u$ such that the restriction of $s$ to $V_u$ is homeomorphic to an open subset of $\gS$ not intersecting $\partial \gS$. Let $M$ denote the complement of the union of all $V_u$ for $u\in s^{-1}(z)$ in $S$. Since $S$ is compact, so is $M$.
If we take a connected open neighborhood $W$ of $z$ included in $V$, $\gS\setminus s(M)$ and $s(V_u)$ for every $u$, then $W$ satisfies the required condition. 

In the case that $z\not\in s(S)$, if the contractible open neighborhood $V$ as above is sufficiently small, then $V$ does not intersect $s(S)$ and satisfies the required condition.

Finally, we consider the case that $z\in s(S)$ and $s^{-1}(z)$ intersects $\partial S$. Note that if $z\in \partial \gS$, then $s^{-1}(z)\subset \partial S$. For each $u\in s^{-1}(z)\setminus \partial S$ we can take a connected open neighborhood $V_u$ of $u$ in $S$ such that the restriction of $s$ to $V_u$ is homeomorphic to an open subset of $\gS$.
For $v\in s^{-1}(z)\cap \partial S$ take a lift $\tilde{v}\in \tilde{B_j}$ of $v$ when $v\in B_j$.
Since $s$ is locally injective, so is $\tilde s$.
Hence there is an open neighborhood $W_{\tilde v}$ of $\tilde v$ in $\tilde S$ and an open neighborhood $W$ of $\tilde s(\tilde v)$ in $\tilde \gS$ such that $\tilde s$ maps $W_{\tilde v}$ homeomorphically to $\tilde s(\tilde S)\cap W$, and $W$ is homeomorphically projected onto an open subset $O_v$ of $\gS$ by $\pi$.
We also have an open subset $W_j $ of $L_j$ by projecting $W\cap U_j$ onto $L_j$.
Now we can see that $W_v:=p(W_{\tilde v})\cup W_j$ in $S'$ is an open neighborhood of $v$ in $S'$ and $s'$ maps $W_v$ homeomorphically to $O_v$. Let $M$ be the complement of the union of all $V_u$ for $u\in s^{-1}(z)\setminus \partial S$ and all $W_v$ for $v\in s^{-1}(z)\cap \partial S$ in $S$.
Then we can see that $M$ is a compact subset of $S$ and $s(M)$ is a closed subset of $\gS$.
Now, take a contractible open neighborhood $O$ of $z$ included in $s(V_u)$ for every $u\in s^{-1}(z)\setminus \partial S$, $O_v$ for every $v\in s^{-1}(z)\cap \partial S$ and $\gS \setminus s(M)$. Then $O$ satisfies the required condition, that is, the restriction of $s'$ to each connected component of $s'^{-1}(O)$ is a homeomorphism onto $O$.

\underline{Step 3. Prove that $f$ and $\tilde{s}$ have the stated properties:}
The inclusion map $f$ from $S$ to $S'$ is an embedding map since $S$ is compact. Each $L_j=\langle g_j \rangle \backslash U_j$ is homotopy equivalent to $S^1$ and so the inclusion map $f$ is a homotopy equivalence. We get a universal covering map $\pi'$ from $(\tilde \gS,\tilde x) $ to $(S', y)$, which is a lift of the covering map $\pi \: \tilde \gS \rightarrow \gS$ with respect to $s'$.
Now, we check that the map $\tilde s\:\tilde S\rightarrow \tilde \gS$ is a lift of $f$, that is, $f\circ p=\pi'\circ \tilde s$. Then we get the following commutative diagram of based topological spaces.

\[
\xymatrix{
(\tilde S, \tilde y) \ar[dd]_p\ar[rr]^{\tilde s} & & (\tilde \gS ,\tilde x) \ar[dd]^{\pi}\ar[ld]_{\pi'}\\
	& (S' ,y) \ar[rd]^{s'}& \\
(S, y) \ar[rr]_s \ar@{^{(}->}[ru]^f & & (\gS ,x) }
\]

Take $y_0\in \tilde S$ and a path $\ell$ from $\tilde y$ to $y_0$. Then $f\circ p(y_0)$ is the terminal point of the lift of $s\circ p\circ \ell$ to $(S',y)$, and $\pi' \circ \tilde s(y_0)$ is the terminal point of the lift of $\pi \circ \tilde s \circ \ell$ to $(S',y)$. Since $s\circ p=\pi \circ \tilde s$, we have $f\circ p(y_0)=\pi'\circ \tilde s(y_0)$. Therefore, $\tilde s$ is a lift of $f$.

Finally, we check that the map $\tilde s$ is an embedding map. First, we check the injectivity of $\tilde s$. Let $y_1,y_2\in \tilde S$ and assume that $\tilde s(y_1)=\tilde s(y_2)$. Let $\ell$ be a path from $y_1$ to $y_2$. Since $\tilde s(y_1)=\tilde s(y_2)$, we have a nullhomotopic closed curve $\pi' \circ \tilde s\circ \ell $ in $S'$, which equals $f\circ p \circ \ell$. Since $f$ is injective and a homotopy equivalence, $p\circ \ell$ is also a nullhomotopic closed curve in $S$, which implies that $y_1=y_2$.
To see that the inverse map $\tilde s^{-1}\: \tilde s (\tilde S)\rightarrow \tilde S$ is continuous, take $x_0\in \tilde s (\tilde S)$ and an open neighborhood $V$ of $\tilde s^{-1}(x_0)$. We can assume that the restriction of $p$ to $V$ is a homeomorphism onto an open subset of $S$. Take a small open neighborhood $W$ of $f\circ p(\tilde s^{-1}(x_0))=\pi '(x_0)$ such that $W\cap S\subset p(V)$ and there exists an open neighborhood $\tilde W$ of $x_0$ such that the restriction of $\pi'$ to $\tilde W$ is homeomorphic to $W$. Then we can see that $\tilde {s}^{-1}(\tilde W\cap \tilde s (\tilde S))= p^{-1}(W\cap S)\cap V \subset V$, which concludes that $\tilde{s}^{-1}$ is continuous. This completes the proof.
\end{proof}

\begin{remark}\label{rem:injective fundamental group}
Under the setting in the above lemma, we can also see that the map $\tilde{s}\:\tilde{S}\rightarrow \tilde{\gS}$ is a proper map, that is, for any compact subset $K$ of $\tilde{\gS}$ the preimage $\tilde{s}^{-1}(K)$ is a compact subset of $\tilde{S}$ since the image $\tilde{s}(\tilde{S})$ is a closed subset of $\tilde{\gS}$.

If either $S$ does not have a boundary or $S$ is a surface whose boundary is mapped to a boundary of $\gS$ by $s$, then the map $s$ itself is a covering map.

The map $s\: S\rightarrow \gS$ induces an injective group homomorphism $s_{\#}$ from the fundamental group $\pi_1(S)$ of $S$ to $\pi_1(\gS)$. By identifying $\pi_1(S)$ with $s_{\#}(\pi_1(S))$ we can see that the map $\tilde{s}\:\tilde{S}\rightarrow \tilde{\gS}$ is a $\pi_1(S)$-equivariant embedding and we can identify $S'$ with the quotient space $\pi_1(S)\backslash \tilde{\gS}$.
Moreover, we can classify a simple compact surface on $\gS$ that is not homeomorphic to a cylinder
by using non-trivial finitely generated subgroups of $\pi_1(\gS)$ (see Proposition \ref{prop:characterize simple cpt surf} for the case that $\gS$ is a compact hyperbolic surface).

Consider the case that $S$ is a cylinder and $S=S^1\times [0,1]$. Then $s|_{S^1\times \{ 0\} }$ is homotopic to $s|_{S^1\times \{ 1\} }$ and the property of $(S,s)$ is the same as that of the closed curve $(S^1\times \{ 0\} ,s|_{S^1\times \{ 0\} })$ on $\gS$ (see Proof of Theorem \ref{thm:bigon criterion 3} from p.\pageref{proof:proof of bigon criterion 3}).

If $\gS$ is a cylinder, then both $\pi_1(\gS)$ and $\pi_1(S)$ are isomorphic to $\mathbb{Z}$. Hence $S$ is homeomorphic to $S^1$ or a cylinder.

Consider the case that $\gS$ is a torus and $S$ is neither a cylinder nor $S^1$. Since $\pi_1(\gS)$ is isomorphic to $\mathbb{Z}^2$ and a non-trivial subgroup of $\mathbb{Z}^2$ is isomorphic to $\mathbb{Z}$ or $\mathbb{Z}^2$, $\pi_1(S)$ is isomorphic to $\mathbb{Z}^2$.
Then $S$ is also a torus and $s$ is a finite-fold covering map.

From the above, we can say that the case that $\gS$ is a compact hyperbolic surface and $\pi_1(S)$ is non-cyclic is essentially new when we consider the intersection number of simple compact surfaces on $\gS$.
\end{remark}

\begin{proposition}\label{prop:characterize simple cpt surf}
Let $\gS$ be a compact hyperbolic surface.
For any simple compact surface $(S,s)$ on $\gS$ that is not a cylinder, there exists a finitely generated subgroup $H$ of the fundamental group of $\gS$ such that the convex core $(C_H,p_H)$ is homotopic to $(S,s)$.
\end{proposition}
\begin{proof}
The notation in this proof is based on the proof of Lemma \ref{lem:good compact surface}.
We consider the universal cover $\tilde{\gS}$ of $\gS$ as a closed convex subspace of the hyperbolic plane $\HH$.
From Lemma \ref{lem:good compact surface}, there exists a covering $S'$ of $\gS$ and a homotopy equivalent embedding map $f$ from $S$ to $S'$.
Let $H$ be a subgroup of the fundamental group of $\gS$ corresponding to the covering $S'$ of $\gS$. Since $S$ is a compact surface or $S^1$, $H$ is finitely generated.

In the case that $S\not=S^1$, since $S$ and $C_H$ have the same genus and the same number of boundary components, there exists a homeomorphism $\phi$ from $S$ to $C_H$.
Even if $S=S^1$, we have a homeomorphism $\phi$ from $S$ to $C_H$.
Note that if $S$ is a cylinder, $S$ is not homeomorphic to $C_H$ since $C_H$ is homeomorphic to $S^1$.
The homeomorphism $\phi$ extends to an $H$-equivariant homeomorphism $\tilde{\phi}$ from $\tilde{S}$ to the convex hull of the limit set of $H$, which is the universal cover of $C_H$ and included in $\tilde{\gS}$. Note that we also have an $H$-equivariant embedding $\tilde{s}\: \tilde{S}\rightarrow \tilde{\gS}\subset \HH$.

Now, we define a homotopy $F\: \tilde{S}\times [0,1]\rightarrow \tilde{\gS}$ from $\tilde{s}$ to $\tilde{\phi}$ by the rule that for $(x,t)\in \tilde{S}\times [0,1]$, $F(x,t)$ is the point on the geodesic from $\tilde{s}(x)$ to $\tilde{\phi}(x)$ in $\HH$ that divides the length of the geodesic in $t:(1-t)$.
Note that $\tilde{\gS}$ is a convex subspace of $\HH$.
Since $H$ acts on $\tilde{\gS}$ by isometry, $F$ is $H$-equivariant, that is, for any $(x,t)\in \tilde{S}\times [0,1]$ and $h\in H$ we have $F(hx,t)=hF(x,t)$. Therefore $F$ induces a homotopy $F'\:S\times [0,1]\rightarrow \gS$ such that for $(x,t)\in S\times [0,1]$ and $\tilde{x}\in \tilde{S}$ with $p(\tilde{x})=x$, $F'(x,t)=\pi (F(\tilde{x},t))$.
We can see that $F'(\cdot, 0)=s$ and $F'(\cdot,1)=p_H\circ \phi$ since for $(x,1)\in S\times [0,1]$ and $\tilde{x}\in \tilde{S}$ with $p(\tilde{x})=x$ we have
\[ F'(x,1)=\pi \circ \tilde{\phi}(\tilde{x})=p_H\circ \phi(x).\]
Therefore $(S,s)$ is homotopic to $(C_H,p_H)$ by the homotopy $F'$ and the homeomorphism $\phi$.
\end{proof}

Let $(S,s)$ be a simple compact surface on $\gS$.
Let $(T,t)$ be a simple compact surface on $\gS$ homotopic to $(S,s)$.
We identify $S$ with $T$ for simplicity of notation.
Let $F\:S\times [0,1]\rightarrow \gS$ be a homotopy from $s$ to $t$.
Consider the universal covering $p\:\tilde{S}\rightarrow S$ of $S$ and a lift $\tilde{s}\:\tilde{S}\rightarrow \tilde{\gS}$ of $s$.
Then $F':=F (p(\cdot ),\cdot )$ is a homotopy from $s\circ p$ to $t\circ p$.
Since we have a lift $\tilde{s}$ of $s\circ p$ with respect to $\pi \:\tilde \gS\rightarrow \gS$, there exists a unique lift $\tilde{F}$ of $F'$ from the homotopy lifting property (see the following commutative diagram).
\[
\xymatrix{
\tilde S \ar[d]\ar[r]^{\tilde{s}} & \tilde \gS \ar[d]^{\pi}\\
\tilde S \times [0,1]\ar@{.>}[ru]^{\tilde{F}}\ar[r]_(.60){F'} & \gS }
\]
Here the map from $S$ to $S\times[0,1]$ maps $x\in S$ to $(x,0)\in S\times [0,1]$.
Since $\pi\circ \tilde{F}(x, 1)=F'(x,1)=F(p(x),1)=t\circ p(x)$ for $x\in \tilde{S}$, the map $\tilde{t}:=\tilde{F}(\cdot,1)\: \tilde{S}\rightarrow \tilde \gS$ is a lift of $t$.

Consider a Riemannian metric of constant curvature $0$ or $-1$ on $\gS$, which induces a Riemannian metric on $\tilde{\gS}$.
Then the fundamental group $\pi_1(\gS)$ of $\gS$ acts on $\tilde{\gS}$ isometrically and we have the following lemma:

\begin{lemma}\label{lem:finite distance}
The Hausdorff distance between $\tilde{s}(\tilde{S})$ and $\tilde{t}(\tilde{S})$ is finite.
\end{lemma}
\begin{proof}
Let $d$ be the distance function on $\tilde{\gS}$.
Let $H$ be the subgroup of $\pi_1(\gS)$ corresponding to $S$.
Then both $\tilde{s}$ and $\tilde{t}$ are $H$-equivariant maps by Remark \ref{rem:injective fundamental group}. This implies that $H$ acts on both $\tilde{s}(\tilde{S})$ and $\tilde{t}(\tilde{S})$ cocompactly.
Take $x\in \tilde{s}(\tilde{S})$ and $y \in \tilde{t}(\tilde{S})$. Then there exists $R>0$ such that 
\[ H(B(x,R))\supset \tilde{s}(\tilde{S}) \text{ and }H(B(y,R))\supset \tilde{t}(\tilde{S}),\]
where $B(x,R)$ is the closed ball centered at $x$ with radius $R$. 
Hence for any $z\in \tilde{s}(\tilde{S})$ there exists $h\in H$ such that $d(z,hx)\leq R$.
Then we have
\[ d(z, \tilde{t}(\tilde{S}))=d(h^{-1}z,\tilde{t}(\tilde{S}))\leq d(h^{-1}z,x)+d(x,y)\leq R+d(x,y).\]
We can also apply this argument to $w\in \tilde{t}(\tilde{S})$ and have $d(\tilde{s}(\tilde{S}),w)\leq R+d(x,y)$. Therefore the Hausdorff distance between $\tilde{s}(\tilde{S})$ and $\tilde{t}(\tilde{S})$ is finite.
\if{0}
By the compactness of $S_i\times [0,1]$ there exist a constant $C\geq 0$ such that for any $(x,r_0)\in S_i\times [0,1]$ there exists a path with length smaller than $C$ from $F(x,0)$ to $F(x,r_0)$ that is homotopic to $F(x, r)\ (r\in [0,r_0])$ fixing the starting point $F(x,0)$ and the terminal point $F(x,1)$. 
Explicitly, for $[x,r_0]\in S_i$ define $f(x,r_0)$ to be the length of a geodesic that is homotopic to $F(x, r)\ (r\in [0,r_0])$ fixing the starting point $F(x,0)$ and the terminal point $F(x,1)$. Then we can see that $f$ is a continuous function on $S_i\times [0,1]$ and take $C$ to be bigger than the maximum of $f$. 

For $x\in \tilde{S_i}$ the distance between $\tilde{s_i}(x)=\tilde{F_i}(x,0)$ and $\tilde{t_i}(x)=\tilde{F_i}(x,1)$ is smaller than $C$ by considering the lift of a geodesic connecting $s_i\circ p_i(x)$ to $t_i\circ p_i(x)$ that is homotopic to $F_i(x,r)\ (r\in [0,1])$ fixing the starting point and the terminal point.

Now, we check that $\tilde{t_i}$ is $H_i$-equivariant. Let $x\in \tilde{S_i}$ and $h\in H_i$. Take a path $\ell$ from $x$ to $h(x)$. Then the closed curve $t_i\circ p_i\circ \ell$ on $\gS$ is homotopic $s_i\circ p_i \circ \ell$ by using the restriction of $F_i$. Since the terminal point of the lift of $s_i\circ p_i\circ \ell$ to $\tilde{s_i}(x)$ equals $h\tilde{s_i}(x)$, we can see that the terminal point of the lift of $t_i\circ p_i\circ \ell$ to $\tilde{t_i}(x)$ equals $h\tilde{t_i}(x)$, which implies that $\tilde{t_i}(h(x))=h\tilde{t_i}(x)$. This completes the proof.
\fi
\end{proof}


We are going to construct the cubic commutative diagram in Proposition \ref{prop:cubic diagram}, which will be used for proving Theorem \ref{thm:bigon criterion 3}. Our construction of the cubic diagram is originally based on that in \cite{Min11}, which was used for studying the Strengthened Hanna Neumann Conjecture. The product $\N$, which will be studied in Section \ref{sec:an intersection functional}, is a certain term appearing in the inequality of the conjecture, and we will also use the cubic diagram for studying the product $\N$.

\begin{setting}[Setting for the proof of the Bigon criterion 3]
We will use the following setting in the rest of this subsection.
Let $\gS$ be a compact cylinder, a torus or a compact hyperbolic surface possibly with boundary.
Let $(S_1,s_1),(S_2,s_2)$ be transverse simple compact surfaces on $\gS$.
Let $G$ be the fundamental group of $\gS$. From Lemma \ref{lem:good compact surface} we can take a covering $s_i'\:S_i'\rightarrow \gS$ such that there is an embedding map $f_i$ from $S_i$ to $S_i'$ with $s_i=s_i'\circ f_i$ ($i=1,2$).
Let $H_i$ be a subgroup of $G$ corresponding to the covering space $S_i'$. We identify $S_i'$ with $H_i\backslash \tilde \gS$ and $\pi_1(S_i)$ with $H_i$.
Let $p_i\: \tilde{S_i}\rightarrow \gS$ be the universal covering and $\tilde{s_i}\: \tilde{S_i}\rightarrow \tilde{\gS}$ a lift of $s_i$ for $i=1,2$.
Then $\tilde{s_i}$ is an $H_i$-equivariant embedding map. Then we obtain the following commutative diagram.
\[
\xymatrix{
\tilde{S_i}\ar[r]^{\tilde{s_i}}\ar[d]_{p_i} &\tilde{\gS} \ar[d]^{\pi}\\
S_i\ar[r]_{s_i} &\gS \ar@{}[lu]|{\circlearrowright}}
\]

Let $\Lambda_i \subset G $ be a complete system of representatives of $G/H_i$.
Assume that the identity element $\id$ belongs to $\gL_i$.
We endow $\Lambda_i$ with the discrete topology and define 
$\hat{S_i}$ to be the direct product $\tilde{S_i}\times \Lambda_i$, which is equipped with the direct product topology. We define a continuous map $\hat{s_i}$ from $\hat{S_i}$ to $\tilde \gS$ by
\[ \hat{s_i}(x,g):=g\circ \tilde{s_i}(x)=g(\tilde{s_i}(x))\quad ((x,g)\in \hat{S_i}).\]
Note that $g\circ \tilde{s_i}$ is also a lift of $s_i$ for any $g\in \Lambda_i$. For $g,g'\in G$ if $gH_i=g'H_i$, then
\[ g\circ \tilde{s_i}(\tilde{S_i})=g'\circ \tilde{s_i}(\tilde{S_i})\]
since $\tilde{s_i}$ is $H_i$-equivariant. Therefore, $(\hat{S_i}, \hat{s_i})$ can be considered as the space consisting of ``all lifts'' of $s_i$.

We define a continuous action of $G$ on $\hat{S_i}$ such that $\hat{s_i}\: \hat{S_i}\rightarrow \tilde \gS$ is a $G$-equivariant map. Let $g\in G$ and $(x,g_0)\in \hat{S_i}$. Note that the natural action of $G$ on $G/H_i$ induces the action of $G$ on $\gL_i$.
We can choose $g_0'\in \Lambda_i $ such that $gg_0= g_0'h$ for some $h\in H_i$. Then we define $g(x,g_0)$ to be $(hx,g_0')$.
We can see that $\hat{s_i}$ is $G$-equivariant from the following equation:
\[ \hat{s_i}(g(x,g_0))=\hat{s_i}(hx,g_0')=g_0'h(\tilde{s_i}(x))=gg_0(\tilde{s_i}(x))=g(\hat{s_i}(x,g_0)).\]
Moreover, $(g_1g_2)(x,g_0)=g_1(g_2(x,g_0))$ for any $g_1,g_2\in G$ and $(x,g_0)\in \hat{S_i}$. Therefore, we get an action of $G$ on $\hat{S_i}$. 

The stabilizer of the connected component 
\[ (\tilde{S_i}, g_0):=\{ (x, g_0)\in \hat{S_i}\mid x\in \tilde{S_i}\} \subset \hat{S_i}\]
coincides with $g_0H_i g_0^{-1}$ for $g_0\in \gL_i$. Especially, the stabilizer of $(\tilde{S_i}, \id)$ is $H_i$.
For $g\in G$ and $g_0, g_0'\in \gL_i$ with $gg_0H_i =g_0'H_i$, we have $g(\tilde{S_i}, g_0 )=(\tilde{S_i} ,g_0' )$.
As a result, for any two connected components of $\hat{S_i}$ there exists $g\in G$ such that $g$ maps one component to the other component. Therefore, the quotient space $G\backslash \hat{S_i}$ can be identified with $H_i\backslash \tilde{S_i}=S_i$. By this identification we get the canonical projection $\hat{p_i}$ from $\hat{S_i}$ to $S_i$. From the construction of $\hat{s_i}$ and $\hat{p_i}$, we have the following commutative diagram.

\[
\xymatrix{
\hat{S_i}\ar[r]^{\hat{s_i}}\ar[d]_{\hat{p_i}} &\tilde{\gS} \ar[d]^{\pi}\\
S_i\ar[r]_{s_i} &\gS \ar@{}[lu]|{\circlearrowright}}
\]
\end{setting}

\begin{remark}\label{rem:expression of hat S i}
Set
\[ \mathfrak{S}_i:=\{ (gH_i, x)\in G/H_i\times \tilde{\gS} \mid x \in g\circ \tilde{s_i}(\tilde{S_i})\} \]
for $i=1,2$. Let $\sigma_i$ be the canonical projection from $\mathfrak{S}_i$ onto $\tilde{\gS}$, that is, $\sigma_i(gH_i,x)=x$.
We see that $(\hat{S_i},\hat{s_i})$ is ``isomorphic'' to $(\mathfrak{S}_i,\sigma_i)$.
Actually, we define a map $\tau_i\: \hat{S_i}\rightarrow \mathfrak{S}_i$ by
\[ \tau_i (x ,g_0):=(g_0H_i, g_0\circ \tilde{s}_i (x)) \]
for $(x,g_0)\in \hat{S_i}$. Then $\tau_i$ is a homeomorphism. Moreover, we can define a natural diagonal action of $G$ on $\mathfrak{S}_i$ by
\[ g(g'H_i ,x):=(gg'H_i,gx)\]
for $g\in G$ and $(g'H_i,x )\in \mathfrak{S}_i$.
Then we can check that $\tau_i$ is a $G$-equivariant homeomorphism.
Actually, for $g\in G$ and $(x, g_0)\in \hat{S_i}$ take $h\in H_i$ and $g_0'\in \gL_i$ such that $gg_0=g_0'h$. Then $g(x,g_0)=(hx, g_0')$.
Hence
\begin{align*}
\tau_i(g(x,g_0))&=\tau_i (hx,g_0')=(g_0'H_i ,g_0' \circ \tilde{s_i}(hx) )=(gg_0H_i, g_0'h \circ \tilde{s_i}(x))\\
&=(gg_0H_i, g(g_0'\circ \tilde{s_i}(x)))=g(g_0H_i,g_0'\circ \tilde{s_i}(x))=g\tau _i(x,g_0).
\end{align*}
We can identify $(\hat{S_i},\hat{s_i})$ with $(\mathfrak{S}_i,\sigma_i)$ through $\tau_i$, and $(\mathfrak{S}_i,\sigma_i)$ is convenient to consider the action of $G$. Moreover, $(\mathfrak{S}_i,\sigma_i)$ can be associated with the counting subset current $\eta_{H_i}$ on $G$ naturally. We will consider this association more concretely in the next subsection.
\end{remark}

We say that $g_1\circ \tilde{s_1}$ and $g_2\circ \tilde{s_2}$ form a bigon for $(g_1,g_2)\in \Lambda_1\times \Lambda_2$ if there exist boundary components $\tilde{B_1}, \tilde{B_2}$ of $\tilde{S_1},\tilde{S_2}$ such that sub-arcs of $(g_1\circ \tilde{s_1})|_{\tilde{B_1}}$ and $(g_2\circ \tilde{s_2})|_{\tilde{B_2}}$ form a bigon. From the following lemma, we can see that considering $\hat{S_1}, \hat{S_2}$ is useful for finding an immersed bigon formed by $s_1$ and $s_2$.

\begin{lemma}\label{lem:cubic diagram and immersed bigon}
Two simple compact surfaces $s_1$ and $s_2$ form an immersed bigon if and only if $g_1\circ \tilde{s_1}$ and $g_2\circ \tilde{s_2}$ form a bigon for some $(g_1,g_2)\in \Lambda_1\times \Lambda_2$.
\end{lemma}
\begin{proof}
\underline{If part:} Assume that $g_1\circ \tilde{s_1}$ and $g_2\circ \tilde{s_2}$ form a bigon $b\: D\rightarrow \tilde{\gS}$ for some $(g_1,g_2)\in \Lambda_1\times \Lambda_2$. 
Take components $\tilde{B_1}$ and $\tilde{B_2}$ of $\partial \tilde{S_1}$ and $\partial \tilde{S_2}$ such that sub-arcs of $g_1\circ \tilde{s_1}|_{\tilde{B_1}}$ and $g_2\circ \tilde{s_2}|_{\tilde{B_2}}$ form the bigon $b$. Then we can see that $\pi \circ b\: D\rightarrow \gS$ is an immersed bigon formed by $s_1|_{p_1(\tilde{B_1})}$ and $s_2|_{p_2(\tilde{B_2})}$.

\underline{Only if part:} Assume that the restriction of $s_1,s_2$ to boundary components $B_1,B_2$ of $S_1, S_2$ form an immersed bigon $b\: D\rightarrow \gS$. Take a boundary component $\tilde{B_i}$ of $\tilde{S_i}$ such that $p_i(\tilde{B_i})=B_i$ for $i=1,2$. Then $p_i|_{\tilde{B_i}}$ is a universal covering of $B_i$.
By the definition of an immersed bigon formed by closed curves, there exists a closed sub-interval $I_i$ of $\tilde{B_i}$ for $i=1,2$ such that $s_1\circ p_1|_{I_1}$ and $s_2\circ p_2|_{I_2}$ form the immersed bigon $b$. Let $b_i$ be a homeomorphism from the edge $e_i$ of $D$ to $I_i$ such that $s_i\circ p_i \circ b_i$ coincides with the restriction of $b$ to $e_i$ for $i=1,2$.
Take a lift $\tilde{b}\: D\rightarrow \tilde{\gS}$ of $b$ with respect to the universal covering $\pi \: \tilde{\gS}\rightarrow \gS$. 
Then $\tilde{b}|_{e_i}$ is a lift of $s_i\circ p_i \circ b_i$ and there exists $\gamma_i \in G$ such that $\tilde{b}|_{e_i}$ coincides with $\gamma_i\circ \tilde{s_i} \circ b_i$. 
Take $g_i \in \gL_i$ and $h_i \in H_i$ such that $\gamma_i =g_i h_i$.
Then $\gamma_i \circ \tilde{s_i} \circ b_i= g_i \circ \tilde{s_i} \circ h_i \circ b_i$ since $\tilde{s_i}$ is $H_i$-equivariant.
This implies that $g_1 \circ \tilde{s_1}|_{h_1 (I_1)}$ and $g_2 \circ \tilde{s_2}|_{h_2 (I_2)}$ form the immersed bigon $\tilde{b}$.
Note that $\tilde{s_i}$ is an embedding map. Hence the restriction of $\tilde{b}$ to $\partial D$ is injective, which implies that $\tilde{b}$ is an embedding map by Lemma \ref{lem:immersed bigon}.
Therefore, $g_1\circ \tilde{s_1}$ and $g_2\circ \tilde{s_2}$ form the bigon $\tilde{b}$.
\end{proof}

In order to consider the intersection of $g_1\circ \tilde{s_1}(\tilde{S_1})$ and $g_2\circ \tilde{s_2}(\tilde{S_2})$ for every $(g_1,g_2)\in \gL_1\times \gL_2$, we take the fiber product $\hat{S_1}\times _{\tilde \gS}\hat{S_2}$ corresponding to $\hat{s_1},\hat{s_2}$. Explicitly,
\[ \hat{S_1}\times _{\tilde \gS}\hat{S_2}:= \{ ((x_1,g_1),(x_2,g_2))\in \hat{S_1}\times \hat{S_2}\mid \hat{s_1}(x_1,g_1)=\hat{s_2}(x_2,g_2)\},\]
which can be identified with the formal union of the fiber product of connected components of $\hat{S_1}$ and $\hat{S_2}$. Therefore we have
\[ \hat{S_1}\times _{\tilde \gS}\hat{S_2}=\bigsqcup_{(g_1,g_2)\in \Lambda_1\times \Lambda_2}(\tilde{S_1},g_1)\times_{\tilde{\gS}}(\tilde{S_2}, g_2).\]
Since the restriction of $\hat{s_i}$ to each connected component of $\hat{S_i}$ is an embedding map, the fiber product $(\tilde{S_1},g_1)\times_{\tilde{\gS}}(\tilde{S_2}, g_2)$ is homeomorphic to the intersection $g_1\circ\tilde{s_1}(\tilde{S_1})\cap g_2\circ \tilde{s_2}(\tilde{S_2})$ for any $(g_1,g_2)\in \Lambda_1\times \Lambda_2$ (see Supplementation \ref{supply:fiber product}). Therefore, we have
\begin{align*}
 \hat{S_1}\times _{\tilde \gS}\hat{S_2}
&\cong \bigsqcup_{(g_1,g_2)\in \Lambda_1\times \Lambda_2}g_1\circ\tilde{s_1}(\tilde{S_1})\cap g_2\circ \tilde{s_2}(\tilde{S_2})\\
&= \bigsqcup_{(g_1 H_1,g_2 H_2)\in G/H_1\times G/H_2}g_1\circ\tilde{s_1}(\tilde{S_1})\cap g_2\circ \tilde{s_2}(\tilde{S_2}).\\
&\cong \{ (g_1 H_1,g_2 H_2,x)\in G/H_1\times G/H_2 \times \tilde{\gS}\ |\\
&\quad	\quad \quad x\in g_1\circ\tilde{s_1}(\tilde{S_1})\cap g_2\circ \tilde{s_2}(\tilde{S_2}) \}.
\end{align*}
Here, we remark that $g_1\circ\tilde{s_1}(\tilde{S_1})\cap g_2\circ \tilde{s_2}(\tilde{S_2})$ can be empty.

Let $\phi_i$ be the canonical projection from $\hat{S_1}\times _{\tilde \gS}\hat{S_2}$ to $\hat{S_i}$. The action of $G$ on $\hat{S_1}$ and $\hat{S_2}$ induces the action of $G$ on $\hat{S_1}\times _{\tilde \gS}\hat{S_2}$ such that $\phi_i$ is a $G$-equivariant map.
Explicitly, for $g\in G$ and $((x_1,g_1),(x_2,g_2))\in \hat{S_1}\times _{\tilde \gS}\hat{S_2}$, we define 
\[ g((x_1,g_1),(x_2,g_2)):=(g(x_1,g_1),g(x_2,g_2)).\]
Note that $(g(x_1,g_1),g(x_2,g_2))$ belongs to $\hat{S_1}\times _{\tilde \gS}\hat{S_2}$ since 
\[ \hat{s_1}(g(x_1,g_1))=g\hat{s_1}(x_1,g_1)=g\hat{s_2}(x_2,g_2)=\hat{s_2}(g(x_2,g_2)).\]
We will prove that the quotient space $G\backslash \hat{S_1}\times _{\tilde \gS}\hat{S_2}$ is homeomorphic to $S_1\times_{\gS}S_2$ in Proposition \ref{prop:cubic diagram}, which will plays an essential role in proving Theorem \ref{thm:bigon criterion 3}.

\begin{lemma}\label{lem:intersection is contractible}
If the intersection $g_1\circ\tilde{s_1}(\tilde{S_1})\cap g_2\circ \tilde{s_2}(\tilde{S_2})$ is not empty for $(g_1,g_2)\in \Lambda_1\times \Lambda_2$, then any connected component of $g_1\circ\tilde{s_1}(\tilde{S_1})\cap g_2\circ \tilde{s_2}(\tilde{S_2})$ is contractible. Moreover, for any compact connected component $M$ of $g_1\circ\tilde{s_1}(\tilde{S_1})\cap g_2\circ \tilde{s_2}(\tilde{S_2})$ with interior points, the number of boundary components of $g_1\circ\tilde{s_1}(\tilde{S_1})$ surrounding $M$ equals that of $g_2\circ \tilde{s_2}(\tilde{S_2})$. Therefore $M$ can be considered as a polygon with even sides.
\end{lemma}
\begin{proof}
Let $M$ be a connected component of $g_1\circ\tilde{s_1}(\tilde{S_1})\cap g_2\circ \tilde{s_2}(\tilde{S_2})$.
In the case that either $S_1$ or $S_2$ is $S^1$, our claim follows obviously.
If $S_1$ (or $S_2$) does not have boundary, then $M$ coincides with $g_2\circ\tilde{s_2}(\tilde{S_2})$ (or $g_1\circ \tilde{s_1}(\tilde{S_1})$ respectively) and $M$ is contractible.

We consider the case that neither $S_1$ nor $S_2$ is $S^1$ and both $S_1$ and $S_2$ have boundaries.
We can assume that $\tilde{\gS}$ does not have boundary by embedding $\tilde{\gS}$ into $\RR^2$ or $\HH^2$.
Then $M$ is a connected subspace of $\tilde{\gS}$ surrounded by the boundaries $g_1\circ \tilde{s_1}(\partial \tilde{S_1})$ and $g_2\circ \tilde{s_2}(\partial \tilde{S_2})$.
Each component $\tilde{B}$ of $g_i\circ \tilde{s_i}(\partial \tilde{S_i})$ is homeomorphic to $\RR$ and divides $\gS$ into two contractible components since there exists $u\in H_i$ with infinite-order such that $\langle u \rangle$ acts on $\tilde{B}$. Hence we can see that $M\setminus \partial M$ is a simply-connected region since the interior region of any simple closed curve on $M\setminus \partial M$ is included in $M\setminus \partial M$.
Note that $s_1$ and $s_2$ are transverse.
Then we can see that $M$ is a 2-dimensional manifold with boundary. By the Riemann mapping theorem, $M\setminus \partial M$ is contractible, which implies that $M$ is contractible.

Now, we assume that $M$ is compact and has some interior points. We can see that $M$ is surrounded by finite boundary components of $g_1\circ\tilde{s_1}(\tilde{S_1})$ and $g_2\circ \tilde{s_2}(\tilde{S_2})$ from Lemma \ref{lem:finite distance}. Since $\tilde{S_i}$ is a 2-dimensional manifold with boundary, any boundary component of $g_i\circ\tilde{s_i}(\tilde{S_i})$ does not intersect the other boundary components of $g_i\circ\tilde{s_i}(\tilde{S_i})$ for $i=1,2$, which implies that the number of boundary components of $g_1\circ\tilde{s_1}(\tilde{S_1})$ surrounding $M$ equals that of $g_2\circ\tilde{s_2}(\tilde{S_2})$. This completes the proof.
\end{proof}

The following proposition, which is corresponding to Remark \ref{rem:expression of hat S i}, is useful to understand the fiber product $\hat{S_1}\times _{\tilde \gS}\hat{S_2}$.

\begin{proposition}\label{prop:characterize fiber product}
Set 
\begin{align*}
Z:=\{ & (g_1 H_1,g_2 H_2,x)\in G/H_1\times G/H_2 \times \tilde{\gS} \ | \\ 
		&\qquad x\in g_1\circ\tilde{s_1}(\tilde{S_1})\cap g_2\circ \tilde{s_2}(\tilde{S_2}) \}.
\end{align*}
Define a map $\theta$ from $\hat{S_1}\times _{\tilde \gS}\hat{S_2}$ to $Z$ as
\[ \theta ((x_1,g_1),(x_2,g_2)):=(g_1 H_1,g_2 H_2, g_1\circ \tilde{s_1}(x_1)) \]
for $((x_1,g_1),(x_2,g_2))\in \hat{S_1}\times _{\tilde \gS}\hat{S_2}$. Then $\theta$ is a homeomorphism.

Define a natural action of $G$ on $Z$ as
\[ g (g_1 H_1,g_2 H_2,x):= (gg_1H_1,gg_2H_2,gx)\]
for $(g_1 H_1,g_2 H_2,x)\in G/H_1\times G/H_2 \times \tilde{\gS}$ and $g\in G$. Then $\theta$ is a $G$-equivariant map.
Moreover, the map $\hat{s_i}\circ \phi_i\circ \theta^{-1} $ is the projection for $i=1,2$, that is,
\[ \hat{s_i}\circ \phi_i \circ \theta^{-1} (g_1 H_1,g_2 H_2,x)=x.\]
This implies that the following diagram is commutative.

\[
\xymatrix{
\hat{S_1}\times_{\tilde{\gS}}\hat{S_2} \ar[r]^{\quad \theta}\ar[d]_{\phi_i} & Z \ar[d]^{{\rm projection}}\\
\hat{S_i}\ar[r]_{s_i} &\tilde{\gS} }
\]
\end{proposition}
\begin{proof}
For $(g_1,g_2)\in \Lambda_1\times \Lambda_2$, the restriction of $\theta$ to $(\tilde{S_1},g_1)\times_{\tilde{\gS}}(\tilde{S_2}, g_2)$ is a homeomorphism onto
\[ \{ (g_1 H_1,g_2 H_2,x)\mid x\in g_1\circ\tilde{s_1}(\tilde{S_1})\cap g_2\circ \tilde{s_2}(\tilde{S_2}) \} \]
since $(\tilde{S_1},g_1)\times_{\tilde{\gS}}(\tilde{S_2}, g_2)$ is mapped homeomorphically to $g_1\circ\tilde{s_1}(\tilde{S_1})\cap g_2\circ \tilde{s_2}(\tilde{S_2})$ by $\tilde{s_1}\circ \phi_1$, which maps $((x_1,g_1),(x_2,g_2))$ to $g_1\circ \tilde{s_1}(x)$. Recall that $\Lambda_i$ is a complete system of representatives of $G/H_i$. Therefore $\theta$ is a homeomorphism.

To see that $\theta$ is $G$-equivariant, take $((x_1,g_1),(x_2,g_2))\in \hat{S_1}\times _{\tilde \gS}\hat{S_2}$ and $g\in G$. 
Take $g_i'\in \Lambda_i $ such that $gg_i= g_i'h_i$ for some $h_i\in H_i$. Then we have $g(x_i,g_i)=(h_ix,g_i')$ for $i=1,2$, and so
\begin{align*}
\theta (g((x_1,g_1),(x_2,g_2)))
&=\theta ((h_1 x_1,g_1'),(h_2 x_2,g_2'))\\
&=(g_1'H_1,g_2'H_2, g_1'\circ \tilde{s_1}(h_1 x_1))\\
&=(g g_1 H_1,g g_2 H_2, g g_1\circ \tilde{s_1}(x_1))\\
&=g \theta ((x_1,g_1),(x_2,g_2)).
\end{align*}

Finally, for $(g_1 H_1,g_2 H_2,x)\in Z$ take $((x_1,g_1'),(x_2,g_2'))\in \hat{S_1}\times _{\tilde \gS}\hat{S_2}$ such that
\[ \theta ((x_1,g_1'),(x_2,g_2'))=(g_1 H_1,g_2 H_2,x).\]
Note that $\hat{s_1}\circ \phi_1=\hat{s_2}\circ \phi_2$ and $g_1'\circ \tilde{s_1}(x_1)=x$.
Then
\begin{align*}
\hat{s_i}\circ \phi_i \circ \theta^{-1} (g_1 H_1,g_2 H_2,x)
&=\hat{s_1}\circ \phi_1((x_1,g_1'),(x_2,g_2')) \\
&=g_1'\circ \tilde{s_1}(x_1)\\
&=x
\end{align*}
for $i=1,2$. This completes the proof.
\end{proof}

\begin{remark}\label{rem:double coset}
From the above proposition, we can identify $\hat{S_1}\times _{\tilde \gS}\hat{S_2}$ with $Z$ and we can see that the choice of $\Lambda_i$ does not influence the fiber product $\hat{S_1}\times _{\tilde \gS}\hat{S_2}$.

For $(g_1,g_2)\in \Lambda_1\times \Lambda_2$, if $(\tilde{S_1},g_1)\times_{\tilde{\gS}}(\tilde{S_2}, g_2)$ is not empty, then the stabilizer of $(\tilde{S_1},g_1)\times_{\tilde{\gS}}(\tilde{S_2}, g_2)$ is $g_1H_1g_1^{-1}\cap g_2H_2g_2^{-1}$.
Hence the quotient space $G\backslash \hat{S_1}\times _{\tilde \gS}\hat{S_2}$ is homeomorphic to the disjoint union of
\[ (g_1H_1g_1^{-1}\cap g_2H_2g_2^{-1})\backslash (g_1\circ\tilde{s_1}(\tilde{S_1})\cap g_2\circ \tilde{s_2}(\tilde{S_2}))\]
over $[g_1 H_1,g_2 H_2] \in G\backslash (G/H_1\times G/H_2)$, which is the quotient set associated with the diagonal action of $G$ on $G/H_1\times G/H_2$.
Actually, for any $g\in G$ and $(g_1, g_2)\in \gL_1 \times \gL_2$, there exists unique $(g_1' ,g_2')\in \gL_1\times \gL_2$ such that $(gg_1H_1,gg_2H_2)=(g_1'H_1, g_2'H_2)$, and then we have
\[ g\left( (\tilde{S_1},g_1)\times_{\tilde{\gS}}(\tilde{S_2}, g_2) \right) = (\tilde{S_1},g_1') \times_{\tilde{\gS}}(\tilde{S_2},g_2' )).\]
Hence each $[g_1 H_1,g_2 H_2] \in G\backslash (G/H_1\times G/H_2)$ corresponds to a connected component of $G\backslash \hat{S_1}\times _{\tilde \gS}\hat{S_2}$, which is possibly empty.
\end{remark}

\begin{lemma}\label{lem:hat s i circ phi i is a proper map}
The map $\hat{s_i}\circ \phi_i \: \hat{S_1}\times _{\tilde \gS}\hat{S_2}\rightarrow \tilde{\gS}$ is a proper map and $G$ acts on $\hat{S_1}\times_{\tilde{\gS}} \hat{S_2}$ freely and properly discontinuously.
\end{lemma}
\begin{proof}
Recall that $\tilde {s_i}\: \tilde{S_i}\rightarrow \tilde{\gS}$ is a proper map because $\tilde{s_i}(\tilde{S_i})$ is a closed subset of $\tilde{\gS}$ and $\tilde{s_i}$ is an embedding map.
Let $J$ be a compact subset of $\tilde{\gS}$. Recall the equation:
\[\hat{S_1}\times _{\tilde \gS}\hat{S_2}=\bigsqcup_{(g_1,g_2)\in \Lambda_1\times \Lambda_2}(\tilde{S_1},g_1)\times_{\tilde{\gS}}(\tilde{S_2}, g_2).\]
For each $(g_1,g_2)\in \Lambda_1\times \Lambda_2$ the intersection
\begin{align*}
&(\tilde{S_1},g_1)\times_{\tilde{\gS}}(\tilde{S_2}, g_2)\cap(\hat{s_i}\circ \phi_i )^{-1}(J)\\
=&\{ ((x_1,g_1),(x_2,g_2))\in (\tilde{S_1},g_1)\times (\tilde{S_2},g_2) \mid g_1\circ \tilde{s_1}(x_1)=g_2\circ \tilde{s_2}(x_2)\in J\}\\
=&\{ ((x_1,g_1),(x_2,g_2))\in ((g_1\circ \tilde{s_1})^{-1}(J))\times ((g_2\circ \tilde{s_2})^{-1}(J))\ | \\ 
&\qquad g_1\circ \tilde{s_1}(x_1)=g_2\circ \tilde{s_2}(x_2)\} \\
=& ((g_1\circ \tilde{s_1})^{-1}(J), g_1)\times_{\tilde \gS} ((g_2\circ \tilde{s_2})^{-1}(J), g_2)
\end{align*}
is compact since $(g_i\circ \tilde{s_i})^{-1}(J)=\tilde{s_i}^{-1}(g_i^{-1}J)$ is compact for $i=1,2$.

We prove that there are only finitely many $g_i\in \Lambda_i$ such that $(g_i\circ \tilde{s_i})^{-1}(J)$ is not empty, that is, $g_i\circ \tilde{s_i}(\tilde{S_i})\cap J\not= \emptyset$ for $i=1,2$.
In the case that $\gS$ is a cylinder or a torus, the fundamental group $G$ of $\gS$ acts on $\tilde{\gS}$ as parallel translations and our claim follows immediately.

In the case that $\gS$ is a compact hyperbolic surface, we apply Lemma \ref{lem:A(K) is relatively compact} to the counting subset current $\eta_{H_i}$ on $G$. Since $\eta_{H_i}(A(J))$ is finite, there are only finitely many $gH_i\in G/H_i$ such that $gCH_{H_i}$ intersects the compact subset $J$. Note that the Hausdorff distance between $gCH_{H_i}$ and $g\circ \tilde{s_i}(\tilde{S_i})$ is finite by Lemma \ref{lem:finite distance}. Hence there are only finitely many $g_i\in \Lambda_i$ such that $g_i\circ \tilde{s_i}(\tilde{S_i})\cap J\not= \emptyset$.

Therefore $(\hat{s_i}\circ \phi_i )^{-1}(J)$ is a finite union of compact subsets and so compact.
Since $\hat{s_i}\circ \phi_i$ is a $G$-equivariant map and $G$ acts on $\tilde{\gS}$ freely and properly discontinuously, $G$ also acts on $\hat{S_1}\times_{\tilde{\gS}} \hat{S_2}$ freely and properly discontinuously.
\end{proof}

From the above lemma, we can see that for any connected component $M$ of $\hat{S_1}\times_{\tilde \gS} \hat{S_2}$ if the stabilizer $\mathrm{Stab}(M)$ of $M$ is non-trivial, then the fundamental group of the quotient space $\Stab(M) \backslash M$ is isomorphic to $\Stab (M)$, which implies that $\Stab(M) \backslash M$ is not contractible. Since $G$ does not have a torsion, the stabilizer of a connected component $M$ of $\hat{S_1}\times _{\tilde \gS}\hat{S_2}$ is trivial if and only if $M$ is compact.

Since the maps $\hat{p_i}\circ \phi_i$ from $\hat{S_1}\times _{\tilde \gS}\hat{S_2}$ to $S_i$ satisfy the condition that
$s_1\circ (\hat{p_1}\circ \phi_1)=s_2\circ (\hat{p_1}\circ \phi_1)$, we can obtain a map $\Phi$ from $\hat{S_1}\times _{\tilde \gS}\hat{S_2}$ to $S_1\times_{\gS}S_2$ (see the following commutative diagram).
\[
\xymatrix{
\hat{S_1}\times _{\tilde \gS}\hat{S_2} \ar[ddr]_{\hat{p_1} \circ \phi_1} \ar[rrd]^{\hat{p_2} \circ \phi_2}\ar@{.>}[dr]_-\Phi & & \\
& S_1\times_{\gS}S_2 \ar[d]\ar[r] & S_2 \ar[d]^{s_2} \\
& S_1 \ar[r]_{s_1} & \gS}
\]
Explicitly, for $ (x_1,g_1),(x_2,g_2)\in \hat{S_1}\times _{\tilde \gS}\hat{S_2}$,
\[ \Phi((x_1,g_1),(x_2,g_2))=( \hat{p_1}(x_1,g_1),\hat{p_2} (x_2,g_2)) .\]

\begin{proposition}\label{prop:cubic diagram}
Let $\alpha=((x_1,g_1),(x_2,g_2)), \beta=((y_1,u_1),(y_2,u_2)) \in \hat{S_1}\times _{\tilde \gS}\hat{S_2}$. There exists $g\in G$ such that $ g(\alpha)=\beta$ if and only if $\Phi (\alpha)=\Phi (\beta)$. Therefore, the map $\Phi$ induces an injective continuous map $\Psi$ from the quotient space $G\backslash (\hat{S_1}\times _{\tilde \gS}\hat{S_2})$ to $S_1\times_{\gS}S_2$. Moreover, $\Psi$ is a homeomorphism. Then we obtain the following cubic commutative diagram.
\[
\xymatrix{
\hat{S_1}\times _{\tilde \gS}\hat{S_2} \ar[rrr]^{\phi_2}\ar[ddd]_{\Phi}\ar[ddr]_{\phi_1} & & &\text{\quad }\hat{S_2}\text{\quad }\ar[ddr]^{\hat{s_2}}\ar[ddd]_(.30){\hat{p_2}}|(.64)\hole &\\
&&&&\\
&\text{\quad }\hat{S_1}\text{\quad }\ar[rrr]^{\hat{s_1}}\ar[ddd]_{\hat{p_1}} &&&\text{\quad }\tilde{\gS}\text{\quad }\ar[ddd]^{\pi}\\
S_1\times_{\gS}S_2 \ar[ddr]\ar[rrr]|(.42)\hole &&&\text{\quad }S_2\text{\quad }\ar[ddr]_{s_2}&\\
&&&&\\
&\text{\quad }S_1\text{\quad } \ar[rrr]_{s_1}&&&\text{\quad }\gS \text{\quad }
}
\]
Every map from a space in the upper square to a space in the lower square is a canonical projection with respect to $G$-action, and every map in the upper square is $G$-equivariant.
\end{proposition}
\begin{proof}
Assume that there exists $g\in G$ such that $g(\alpha)=\beta$. Since $\phi_i$ is $G$-equivariant and $\hat{p_i}$ is a canonical projection with respect to the action of $G$ on $\hat{S_i}$,
\begin{align*}
\Phi (\beta )
&=\Phi (g\alpha ) =\Phi (g(x_1,g_1),g(x_2,g_2)) \\
&=( \hat{p_1}(g(x_1,g_1)),\hat{p_2} (g(x_2,g_2)))\\
&=(\hat{p_1}(x_1,g_1),\hat{p_2}(x_2,g_2) )\\
&=\Phi (\alpha ).
\end{align*}

Next, we assume that $\Phi (\alpha )=\Phi (\beta )$, that is,
\[ (\hat{p_1}(x_1,g_1 ),\hat{p_2}(x_2,g_2))=(\hat{p_1}(y_1,u_1) ,\hat{p_2}(y_2,u_2) ).\]
There exist $v_1,v_2\in G$ such that
\[ v_1(x_1,g_1) =(y_1,u_1), v_2( x_2,g_2)=(y_2,u_2).\]
It is sufficient to see that $v_1=v_2$, which implies that $v_1\alpha =\beta$.
Since $\alpha ,\beta$ belong to $\hat{S_1}\times _{\tilde \gS}\hat{S_2}$, we have
\[ \hat{s_1}(x_1,g_1)=\hat{s_2}(x_2,g_2),\ \hat{s_1}(y_1,u_1)=\hat{s_2}(y_2,u_2).\]
Therefore
\begin{align*}
v_1\hat{s_1}(x_1,g_1)
&=\hat{s_1}(v_1(x_1,g_1) )
=\hat{s_1}(y_1,u_1)
=\hat{s_2}(y_2,u_2)\\
&=\hat{s_2}(v_2(x_2,g_2))
=v_2\hat{s_2}(x_2,g_2)
=v_2\hat{s_1}(x_1,g_1). 
\end{align*}
This implies that $v_1=v_2$ since $G$ acts on $\tilde{\gS}$ freely.

To see the surjectivity of $\Psi$, we check that $\Phi$ is surjective. Take an arbitrary $(z_1,z_2)\in S_1\times_{\gS}S_2$.
Take $(x_i,g_i)\in \hat{S_i}$ such that $\hat{p_i}(x_i,g_i)=z_i$ for $i=1,2$.
Since $s_1(z_1)=s_2(z_2)$ and $s_i\circ \hat{p_i}=\pi \circ \hat{s_i}$, we can see that $\hat{s_1}(x_1,g_1),\hat{s_2}(x_2,g_2) \in \pi^{-1}(s_1(z_1))$. Hence there exists $g\in G$ such that $g \hat{s_1}(x_1,g_1)=\hat{s_2}(x_2,g_2)$, that is,
$(g(x_1,g_1),(x_2,g_2))\in \hat{S_1}\times _{\tilde \gS}\hat{S_2}$. Therefore we have
\[ \Phi (g(x_1,g_1),(x_2,g_2))=(\hat{p_1}(g(x_1,g_1)),\hat{p_2}(x_2,g_2))=(z_1,z_2).\]

From the above, $\Psi$ is a bijective continuous map. Hence it is sufficient to prove that the quotient space $G\backslash (\hat{S_1}\times _{\tilde \gS}\hat{S_2})$ is compact.
Since $\gS$ is compact, there exist a compact subset $K$ of $\tilde \gS$ such that $\pi (K)=\gS$, that is, $G(K)=\tilde{\gS}$.
Then $(\hat{s_i}\circ \phi_i)^{-1}(K)$ is also a compact subset of $\hat{S_1}\times _{\tilde \gS}\hat{S_2}$ by Lemma \ref{lem:hat s i circ phi i is a proper map}. Then we can see that 
\[ G((\hat{s_i}\circ \phi_i)^{-1}(K))=\hat{S_1}\times _{\tilde \gS}\hat{S_2}\]
since $\hat{s_i}\circ \phi_i$ is $G$-equivariant. Therefore the quotient space $G\backslash (\hat{S_1}\times _{\tilde \gS}\hat{S_2})$ is compact, which completes the proof.
\end{proof}

Let $(T_i,t_i)$ be a simple compact surface on $\gS$ homotopic to $(S_i,s_i)$ for $i=1,2$.
We identify $S_i$ with $T_i$ for simplicity of notation.
Recall that we have a lift $\tilde{t_i}\:\tilde{S}\rightarrow \tilde{\gS}$ of $t_i$ such that the Hausdorff distance between $\tilde{s_i}(\tilde{S_i})$ and $\tilde{t_i}(\tilde{S_i})$ is finite by Lemma \ref{lem:finite distance}.
Then, we can obtain the same diagram in Proposition \ref{prop:cubic diagram} for simple compact surfaces $(S_1,t_1)$, $(S_2,t_2)$ on $\gS$ and their lifts $(\tilde{S_1},\tilde{t_1})$, $(\tilde{S_2},\tilde{t_2})$.

\begin{proof}[Proof of Theorem \ref{thm:bigon criterion 3} (the Bigon criterion 3) ]\label{proof:proof of bigon criterion 3}
We classify our proof into several cases. 
We use Lemma \ref{lem:cubic diagram and immersed bigon} and consider $g_1\circ \tilde{s_1},g_2\circ \tilde{s_2}$ for $(g_1,g_2)\in \Lambda_1\times \Lambda_2$ instead of $s_1,s_2$. We will say that a boundary component $\tilde{B_1}$ of $g_1\circ \tilde{s_1}(\tilde{S_1})$ and a boundary component $\tilde{B_2}$ of $g_2\circ \tilde{s_2}(\tilde{S_2})$ form a bigon if there exist a boundary component $\tilde{B_1}'$ of $\tilde{S_1}$ and a boundary component $\tilde{B_2}'$ of $\tilde{S_2}$ such that $g_i\circ \tilde{s_i}(\tilde{B_i}')=\tilde{B_i}$ for $i=1,2$ and sub-arcs of $g_1\circ \tilde{s_1}|_{\tilde{B_1}'}$ and $g_2\circ \tilde{s_2}|_{\tilde{B_2}'}$ form a bigon.

\underline{Case 1:} The surface $\gS$ is a sphere.

See Proposition \ref{prop:simple compact surface on sphere}.
Note that if $\gS$ is not a sphere, $S_i$ can not be a sphere by Lemma \ref{lem:good compact surface}.

\underline{Case 2:} The surface $\gS$ is a cylinder.

The simple compact surface $S_i$ must be $S^1$ or a cylinder since $s_i$ induces an injective group homomorphism from the fundamental group of $S_i$ to that of $\gS$, which is isomorphic to $\mathbb{Z}$. 
Then we can see that $i([s_1],[s_2])=0$ for any simple compact surfaces $s_1,s_2$ on $\gS$ since we can deform $s_1$ and $s_2$ by homotopies such that their images do not intersect.
Now, we consider the intersection of $g_1\circ \tilde{s_1}(\tilde{S_1})$ and $g_2\circ \tilde{s_2}(\tilde{S_2})$ for $(g_1,g_2)\in \Lambda_1\times \Lambda_2$. The universal cover $\tilde{S_i}$ of $S_i$ is homeomorphic to $\RR$ or $\RR \times [0,1]$. Note that an infinite cyclic group $H_i$ acts on $\tilde{S_i}$ and the stabilizer of $g_1\circ \tilde{s_1}(\tilde{S_1})\cap g_2\circ \tilde{s_2}(\tilde{S_2})$ is $g_1H_1g_1^{-1}\cap g_2H_2g_2^{-1}=H_1\cap H_2$, which is also an infinite cyclic group.
If $g_1\circ \tilde{s_1}(\tilde{S_1})\cap g_2\circ \tilde{s_2}(\tilde{S_2})\not=\emptyset$, then $H_1\cap H_2$ acts on $g_1\circ \tilde{s_1}(\tilde{S_1})\cap g_2\circ \tilde{s_2}(\tilde{S_2})$, which implies that $g_1\circ \tilde{s_1}(\tilde{S_1})\cap g_2\circ \tilde{s_2}(\tilde{S_2})$ is non-compact connected, or an infinite union of compact connected components.

If $g_1\circ \tilde{s_1}(\tilde{S_1})\cap g_2\circ \tilde{s_2}(\tilde{S_2})$ is an infinite union of compact connected components, then we can see that $g_1\circ \tilde{s_1}$ and $g_2\circ \tilde{s_2}$ form a bigon.
Actually, any compact component of $g_1\circ \tilde{s_1}(\tilde{S_1})\cap g_2\circ \tilde{s_2}(\tilde{S_2})$ is surrounded by both $g_1\circ \tilde{s_1}(\partial \tilde{S_1})$ and $g_2\circ \tilde{s_2}(\partial \tilde{S_2})$, which implies that there exists a boundary component $\tilde{B}$ of $g_1\circ \tilde{s_1}(\partial \tilde{S_1})$ such that $\tilde{B}$ intersects a boundary component of $g_2\circ \tilde{s_2}(\partial \tilde{S_2})$ infinitely many times. Note that the restriction of $s_1$ and $s_2$ to any components of their boundaries are transverse. Therefore $g_1\circ \tilde{s_1}$ and $g_2\circ \tilde{s_2}$ form a bigon by Lemma \ref{lem:bigon easy}.

From the above, we can see that if $g_1\circ \tilde{s_1}$ and $g_2\circ \tilde{s_2}$ do not form a bigon, then $g_1\circ \tilde{s_1}(\tilde{S_1})\cap g_2\circ \tilde{s_2}(\tilde{S_2})$ is empty or non-compact connected.
By Lemma \ref{lem:cubic diagram and immersed bigon}, if $s_1$ and $s_2$ do not form an immersed bigon, then $S_1\times_\gS S_2$ does not have any contractible components, that is, $i(s_1,s_2)=0=i([s_1],[s_2])$.

The converse does not follow if $S_1,S_2$ are cylinders. For example, consider the case that
\[ \tilde{\gS}=\RR \times [-4,4],\ g_1\circ \tilde{s_1}( \tilde{S_1})=\RR \times [-2,2] \]
and
\[ g_2\circ \tilde{s_2}( \tilde{S_2})=\{ (x,y)\in \RR^2 \mid \sin x-2\leq y \leq \sin x +2\}. \]
Then $g_1\circ \tilde{s_1}$ and $g_2\circ \tilde{s_2}$ form a bigon but $g_1\circ \tilde{s_1}(\tilde{S_1})\cap g_2\circ \tilde{s_2}(\tilde{S_2})$ is non-compact connected.
If either $S_1$ or $S_2$ is $S^1$, then the converse follows immediately from the above argument.

\underline{Case 3:} The surface $\gS$ is a torus.

We assume that $\gS=\mathbb{Z}^2\backslash \RR^2$, which is the quotient space of $\RR^2$ by the natural action of $\mathbb{Z}^2$. Note that a non-trivial subgroup of $G=\mathbb{Z}^2$ is isomorphic to $\mathbb{Z}^2$ or $\mathbb{Z}$. First, we consider the case that $H_1$ is isomorphic to $\mathbb{Z}^2$. Then $H_1$ is a subgroup of $G$ of finite index, which implies that $S_1$ is a torus and $s_1$ is a covering map. Therefore $\tilde{s_1}(\tilde{S_1})=\tilde{\gS}=\RR^2$, and so 
$(\tilde{S_1},g_1)\times_{\tilde{\gS}}(\tilde{S_2}, g_2)$ does not include a compact component for any $(g_1,g_2)\in \Lambda_1\times \Lambda_2$. As a result, $i(s_1,s_2)=0$. 

Now, we assume that both $H_1$ and $H_2$ are isomorphic to $\mathbb{Z}$, which implies that $S_i$ is $S^1$ or a cylinder for $i=1,2$. 
If $H_1\cap H_2$ is not trivial, then we can apply the same argument in the case that $\gS$ is a cylinder to this case.
Therefore we consider the case that $H_1\cap H_2$ is trivial. Take $(a_i,b_i)\in H_i$ such that $(a_i,b_i)$ generates $H_i$. Then two vectors $(a_1,b_1)$ and $(a_2,b_2)$ are linearly independent over the ring $\mathbb{Z}$.

Note that the image $\tilde{s_i}(\tilde{S_i})$ divides $\tilde{\gS}$ into two regions since $H_i$ acts on $g_i\circ\tilde{s_i}(\tilde{S_i})$ for $i=1,2$. Hence $g_1\circ \tilde{s_1}(\tilde{S_1})$ intersects $g_2\circ \tilde{s_2}(\tilde{S_2})$ for any $(g_1,g_2)\in \Lambda_1\times \Lambda_2$ and the intersection includes at least one compact connected component of $\tilde{\gS}$. Moreover, we can see that if $g_1\circ \tilde{s_1}(\tilde{S_1})\cap g_2\circ \tilde{s_2}(\tilde{S_2})$ includes more than one compact components, then $g_1\circ \tilde{s_1}$ and $g_2\circ \tilde{s_2}$ form a bigon.
Actually, any boundary components of $g_1\circ \tilde{s_1}(\tilde{S_1})$ must go into $g_2\circ \tilde{s_2}(\tilde{S_2})$ and go out the opposite side at least once. If $g_1\circ \tilde{s_1}(\tilde{S_1})\cap g_2\circ \tilde{s_2}(\tilde{S_2})$ have more than one compact components, then a boundary component of $g_1\circ \tilde{s_1}(\tilde{S_1})$ must intersect a boundary component of $g_2\circ \tilde{s_2}(\tilde{S_2})$ more than once, and their sub-arcs form a bigon by Lemma \ref{lem:bigon easy}.

From the above, if $s_1$ and $s_2$ do not form an immersed bigon, then $s_1,s_2$ are in minimal position.
If either $S_1$ or $S_2$ is $S^1$, then the converse follows immediately from the above argument.

\underline{Case 4:} The surface $\gS$ is a compact hyperbolic surface.

In this case we thought of $\tilde{\gS}$ as a closed convex subspace of the hyperbolic plane $\HH$.
See the beginning part of Section \ref{sec:Volume functionals for Kleinian groups} for some definitions and notation related to hyperbolic geometry.

Take $(g_1,g_2)\in \Lambda_1\times \Lambda_2$. 
We prove that if $g_1\circ \tilde{s_1}$ and $g_2\circ \tilde{s_2}$ do not form a bigon, then the number of compact connected components of $g_1\circ \tilde{s_1}(\tilde{S_1})\cap g_2\circ \tilde{s_2}(\tilde{S_2})$ is minimum in the homotopy classes $[s_1]$ and $[s_2]$.
Note that the limit set $(g_i\circ \tilde{s_i}(\tilde{S_i}))(\infty )=g_i\gL(H_i)$ coincides with $(g_i\circ \tilde{t_i}(\tilde{S_i}))(\infty )$ from Lemma \ref{lem:finite distance}.
We classify our proof into several cases under the relation between $g_1\gL(H_1)$ and $g_2 \gL(H_2)$.
Since $H_1,H_2$ are finitely generated, we have 
\[ g_1\Lambda (H_1)\cap g_2\Lambda (H_2)=\Lambda (g_1H_1g_1^{-1}\cap g_2H_2g_2^{-1}).\]

\underline{Case 4-1:} The intersection $g_1\Lambda (H_1)\cap g_2\Lambda (H_2)$ is not empty.

In this case, $g_1H_1g_1^{-1}\cap g_2H_2g_2^{-1}$ is not trivial and acts on $g_1\circ \tilde{s_1}(\tilde{S_1})\cap g_2\circ \tilde{s_2}(\tilde{S_2})$. We prove that if $g_1\circ \tilde{s_1}(\tilde{S_1})\cap g_2\circ \tilde{s_2}(\tilde{S_2})$ includes a compact connected component $M$, then $g_1\circ \tilde{s_1}$ and $g_2\circ \tilde{s_2}$ form a bigon. In other words, if $g_1\circ \tilde{s_1}$ and $g_2\circ \tilde{s_2}$ do not form a bigon, then $g_1\circ \tilde{s_1}(\tilde{S_1})\cap g_2\circ \tilde{s_2}(\tilde{S_2})$ does not have a compact connected component.

Consider the case that $S_1$ is $S^1$, which implies that $H_1$ is an infinite cyclic group. Since $g_1H_1g_1^{-1}\cap g_2H_2g_2^{-1}$ is not trivial, $g_1H_1g_1^{-1}\cap g_2H_2g_2^{-1}$ is a finite index subgroup of $g_1H_1g_1^{-1}$. 
Assume that $g_1\circ \tilde{s_1}(\tilde{S_1})\cap g_2\circ \tilde{s_2}(\tilde{S_2})$ includes a compact connected component $M$.
Then the compact connected component $M$ must be a point or homeomorphic to a closed interval by the assumption on the simple compact surfaces $s_1$ and $s_2$ . 
If $M$ is a point, then $S_2$ is also $S^1$ and $g_1\circ \tilde{s_1}(\tilde{S_1})$ intersects $g_2\circ \tilde{s_2}(\tilde{S_2})$ transversely infinitely many times and their sub-arcs form a bigon by Lemma \ref{lem:bigon easy}.
Hence we consider the case that $M$ is homotopic to a closed interval.
Note that each endpoint of $M$ is the intersection point of $g_1\circ \tilde{s_1}(\tilde{S_1})$ with a boundary component of $g_2\circ \tilde{s_2}(\tilde{S_2})$. Since $g_1H_1g_1^{-1}\cap g_2H_2g_2^{-1}$ acts on $g_1\circ \tilde{s_1}(\tilde{S_1})\cap g_2\circ \tilde{s_2}(\tilde{S_2})$, $g_1\circ \tilde{s_1}(\tilde{S_1})$ intersects boundary components of $g_2\circ \tilde{s_2}(\tilde{S_2})$ infinitely many times. By giving an orientation to $g_1\circ \tilde{s_1}(\tilde{S_1})$ we can see that if $g_1\circ \tilde{s_1}(\tilde{S_1})$ goes out from a boundary component $\tilde{B}$ of $g_2\circ \tilde{s_2}(\tilde{S_2})$, then $g_1\circ \tilde{s_1}(\tilde{S_1})$ must go into $g_2\circ \tilde{s_2}(\tilde{S_2})$ through the same boundary component $\tilde{B}$. This implies that $g_1\circ \tilde{s_1}(\tilde{S_1})\cap g_2\circ \tilde{s_2}(\tilde{S_2})$ and $\tilde{B}$ form a bigon by Lemma \ref{lem:bigon easy}.

Next, consider the case that neither $S_1$ nor $S_2$ is $S^1$. Assume that $g_1\circ \tilde{s_1}(\tilde{S_1})\cap g_2\circ \tilde{s_2}(\tilde{S_2})$ includes a compact connected component $M$.
By Lemma \ref{lem:intersection is contractible}, a compact connected component $M$ of $g_1\circ \tilde{s_1}(\tilde{S_1})\cap g_2\circ \tilde{s_2}(\tilde{S_2})$ is a region surrounded by $g_1\circ \tilde{s_1}(\partial \tilde{S_1})$ and $g_2\circ \tilde{s_2}(\partial \tilde{S_2})$. Take a boundary component $\tilde{B}$ of $g_1\circ \tilde{s_1}(\tilde{S_1})$ and a non-trivial element $u\in g_1H_1g_1^{-1}$ such that $\tilde{B}$ form a side of $M$ and $\langle u\rangle$ acts on $\tilde{B}$. If $\tilde{B}(\infty )\cap g_2\Lambda (H_2)\not=\emptyset$, then there is $m\in \NN$ such that $u^m\in g_2H_2g_2^{-1}$ and $\tilde{B}(\infty )\subset g_2\Lambda (H_2)$ since $u$ is a hyperbolic element of the isometry group of $\HH$.
By applying the above argument in the case that $S_1=S^1$ to $\tilde{B}$ and $\langle u^m\rangle$, we can see that $\tilde{B}$ and a boundary component of $g_2\circ \tilde{s_2}(\tilde{S_2})$ form a bigon.

To obtain a contradiction, we assume that $g_1\circ \tilde{s_1}$ and $g_2\circ \tilde{s_2}$ do not form a bigon.
Then any boundary component of $g_1\circ \tilde{s_1}(\tilde{S_1})$ forming a side of $M$ goes into $g_2\circ \tilde{s_2}(\tilde{S_2})$ and goes out from $g_2\circ \tilde{s_2}(\tilde{S_2})$ only once.
Note that every non-trivial element of $G$ is a hyperbolic element in $\Isom (\HH)$ and for non-trivial $\gamma_1,\gamma_2 \in G$ either the intersection of $\gL(\langle g_1 \rangle )$ and $\gL (\langle g_2 \rangle )$ is empty or $\gL(\langle g_1 \rangle )=\gL (\langle g_2 \rangle )$. Hence if a boundary component $\tilde{B}$ of $g_1\circ \tilde{s_1}(\tilde{S_1})$ goes into $g_2\circ \tilde{s_2}(\tilde{S_2})$ and goes out from $g_2\circ \tilde{s_2}(\tilde{S_2})$ exactly once, then the limit set $\tilde{B}(\infty)$ of $\tilde{B}$ does not intersect $g_2 \gL (H_2)$ from the above argument.
Therefore the intersection of $g_1 \gL (H_1)$ and $g_2\gL (H_2)$ is empty since $M$ is compact.
This contradicts our assumption that $g_1\Lambda (H_1)\cap g_2\Lambda (H_2)$ is not empty.
Hence $g_1\circ \tilde{s_1}$ and $g_2\circ \tilde{s_2}$ form a bigon.

\begin{figure}[h]
\begin{center}
\includegraphics[width=10cm]{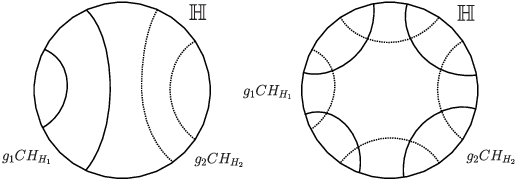}
\vspace{-0.3cm}
\caption{The left figure shows Case 4-2 and the right figure shows Case 4-3. In each figure, the region surrounded by solid lines is the convex hull $g_1CH_{H_1}$ of $g_1\gL(H_1)$ and the region surrounded by dotted lines is $g_2CH_{H_2}$.}\label{fig:intersection of convex hull}
\end{center}
\end{figure}

\underline{Case 4-2:} The intersection $g_1\Lambda (H_1)\cap g_2\Lambda (H_2)=\emptyset$ and there exist two closed intervals $I_1,I_2$ of $\partial \HH$ satisfying the condition that 
\[ I_1\cap I_2=\emptyset \text{ and }I_i\supset g_i\Lambda (H_i) \text{ for } i=1,2. \]

In this case, two convex hulls $CH(I_1),CH(I_2)$ do not intersect.
Take a boundary component $\tilde{B_i}$ of $g_i \circ \tilde{s_i}(\tilde{S_i})$ such that $CH(\tilde{B_i}(\infty))$ is closest to the geodesic line $CH(\partial I_i)$ for $i=1,2$. Then $\tilde{B_1}$ and $\tilde{B_2}$ form a bigon if and only if $g_1\circ \tilde{s_1}(\tilde{S_1})$ and $g_2\circ \tilde{s_2}(\tilde{S_2})$ intersects.
Therefore if $g_1\circ \tilde{s_1}(\tilde{S_1})\cap g_2\circ \tilde{s_2}(\tilde{S_2})$ includes a compact connected component $M$, then $\tilde{B_1}$ and $\tilde{B_2}$ form a bigon

\underline{Case 4-3:} The intersection $g_1\Lambda (H_1)\cap g_2\Lambda (H_2)=\emptyset$ and there do not exist two closed intervals $I_1,I_2$ of $\partial \HH$ satisfying the condition in Case 4-2.

This assumption implies that there exist a boundary component $\tilde{B}$ of $g_1\circ \tilde{s_1}(S_1)$ such that any interval of $\partial \HH$ connecting the two points in $\tilde{B}(\infty )$ must intersect $g_2\Lambda (H_2)$.
In this case $g_1\circ \tilde{s_1}(\tilde{S_1})$ must intersect $g_2\circ \tilde{s_2}(\tilde{S_2})$. Since $g_1\Lambda (H_1)\cap g_2\Lambda (H_2)=\emptyset$, the intersection $g_1\circ \tilde{s_1}(\tilde{S_1})\cap g_2\circ \tilde{s_2}(\tilde{S_2})$ is compact.
Therefore we prove that if $g_1\circ \tilde{s_1}$ and $g_2\circ \tilde{s_2}$ do not form a bigon, then $g_1\circ \tilde{s_1}(\tilde{S_1})\cap g_2\circ \tilde{s_2}(\tilde{S_2})$ includes exactly one compact connected component.

In the case that $S_1$ is $S^1$, if $g_1\circ \tilde{s_1}$ and $g_2\circ \tilde{s_2}$ do not form a bigon, then $g_1\circ \tilde{s_1}(\tilde{S_1})$ intersects $g_2 \circ \tilde{s_2}(\tilde{S_2})$ at a point, or goes into $g_2\circ \tilde{s_2}(\tilde{S_2})$ and goes out from $g_2\circ \tilde{s_2}(\tilde{S_2})$ exactly once, which implies that $g_1\circ \tilde{s_1}(\tilde{S_1})\cap g_2\circ \tilde{s_2}(\tilde{S_2})$ includes exactly one compact connected component.

Therefore, we assume that neither $S_1$ nor $S_2$ is $S^1$.
We also assume that $g_1\circ \tilde{s_1}(\tilde{S_1})$ and $g_2\circ \tilde{s_2}(\tilde{S_2})$ do not form a bigon. Then any boundary component $\tilde{B}$ of $g_1\circ \tilde{s_1}(\tilde{S_1})$ satisfies either one of the following two conditions
\begin{enumerate}
\item there exists an interval $I$ of $\partial \HH$ connecting the two points in $\tilde{B}(\infty) $ such that $I\cap g_2 \gL (H_2)=\emptyset$;
\item any interval $I$ of $\partial \HH$ connecting the two points in $\tilde{B}(\infty) $ must intersect $g_2\gL (H_2)$.
\end{enumerate}
If $\tilde{B}$ satisfies the condition (1), then $\tilde{B}$ does not intersect $g_2\circ \tilde{s_2}(\tilde{S_2})$ by the argument in the case that $S_1=S^1$.
If $\tilde{B}$ satisfies the condition (2), then $\tilde{B}$ goes into $g_2\circ \tilde{s_2}(\tilde{S_2})$ and goes out from $g_2\circ \tilde{s_2}(\tilde{S_2})$ exactly once, which divides $g_2\circ \tilde{s_2}(\tilde{S_2})$ into two connected components and one of the connected components contains $g_1\circ \tilde{s_1}(\tilde{S_1})\cap g_2\circ \tilde{s_2}(\tilde{S_2})$. Therefore, $g_1\circ \tilde{s_1}(\tilde{S_1})\cap g_2\circ \tilde{s_2}(\tilde{S_2})$ is connected, and our claim follows.

From Case 4-1, 4-2 and 4-3, we can see that if $g_1\circ \tilde{s_1}$ and $g_2\circ \tilde{s_2}$ do not form a bigon, then the number of compact connected components of $g_1\circ \tilde{s_1}(\tilde{S_1})\cap g_2\circ \tilde{s_2}(\tilde{S_2})$ is minimum in the homotopy classes $[s_1]$ and $[s_2]$. Moreover, $i([s_1],[s_2])$ equals the number of $[g_1H_1, g_2H_2]\in G\backslash G/H_1\times G/H_2$ satisfying the condition of Case 4-3 by Remark \ref{rem:double coset}. From Proposition \ref{prop:cubic diagram} and Lemma \ref{lem:cubic diagram and immersed bigon}, if $s_1$ and $s_2$ do not form an immersed bigon, then $s_1$ and $s_2$ are in minimal position.
If either $S_1$ or $S_2$ is $S^1$, then the converse follows by considering each case 4-1, 4-2 and 4-3.
\end{proof}

\begin{supply}
Let $\gS=\mathbb{Z}^2\backslash \RR^2$. Assume that $s_1,s_2$ do not form a bigon, both $H_1$ and $H_2$ are infinite cyclic groups and $H_1\cap H_2$ is trivial.
In this setting we calculate the intersection number $i(s_1,s_2)=i([s_1],[s_2])$.

We have proved that $g_1\circ \tilde{s_1}(\tilde{S_1})\cap g_2\circ \tilde{s_2}(\tilde{S_2})$ contains exactly one compact connected component for any $(g_1,g_2)\in \Lambda_1\times \Lambda_2$. From Remark \ref{rem:double coset} and Proposition \ref{prop:cubic diagram}, $S_1\times_\gS S_2$ is homeomorphic to the disjoint union of
\[ (g_1H_1g_1^{-1}\cap g_2H_2g_2^{-1})\backslash (g_1\circ\tilde{s_1}(\tilde{S_1})\cap g_2\circ \tilde{s_2}(\tilde{S_2}))\]
over $[g_1 H_1,g_2 H_2] \in G\backslash (G/H_1\times G/H_2)$. Therefore $i(s_1,s_2)$ equals the cardinality of $G\backslash (G/H_1\times G/H_2)$. 
Define a map $\tau \: G/\langle H_1\cup H_2\rangle \rightarrow G\backslash (G/H_1\times G/H_2)$ as
\[ \tau (g\langle H_1\cup H_2\rangle ) =[H_1,gH_2] \]
for $g\langle H_1\cup H_2\rangle \in G/\langle H_1\cup H_2\rangle$. The map $\tau $ is well-defined. Actually, since $G=\mathbb{Z}^2$ is commutative, for $(h_1,h_2)\in H_1 \times H_2$ we have
\[ [H_1,(gh_1h_2)H_2]=[H_1,(h_1g)H_2]=[H_1,gH_2].\]

We prove that $\tau $ is bijective. The surjectivity of $\tau$ follows immediately since $\tau $ is well-defined. We check the injectivity of $\tau$. For $g,g'\in G$, assume that 
\[ \tau (g\langle H_1\cup H_2\rangle )=\tau (g'\langle H_1\cup H_2\rangle ),\]
that is, $[H_1,gH_2]=[H_1,g'H_2]$. Then we can take $h_1\in H_1$ such that $h_1gH_2=g'H_2$, which implies that there exists $h_2\in H_2$ such that $h_1gh_2=g'$. Hence
\[ g' \langle H_1\cup H_2\rangle =gh_1h_2\langle H_1\cup H_2\rangle =g \langle H_1\cup H_2\rangle.\]

From the above, $i(s_1,s_2)$ equals the index $[G:\langle H_1 \cup H_2 \rangle ]$.
Note that $\langle H_1 \cup H_2 \rangle$ is a finite index subgroup of $G$.
Let $(a_i,b_i)$ be a generator of $H_i$. In order to calculate the index $[G:\langle H_1 \cup H_2 \rangle ]$ we consider the area of the covering space of $\gS$ corresponding to $\langle H_1 \cup H_2 \rangle$.
The area of the quotient space $\langle H_1\cup H_2\rangle \backslash \RR^2$ equals the area of the parallelogram formed by the two vectors $(a_1,b_1),(a_2,b_2)$, that is, $|a_1b_2-b_1a_2|$. Since the area of $\gS$ is $1$, $[G:\langle H_1 \cup H_2 \rangle ]= |a_1b_2-b_1a_2|$.
Therefore
\[ i(s_1,s_2)=i([s_1],[s_2])=|a_1b_2-b_1a_2|.\]
Even if $H_1,H_2$ are infinite cyclic and $H_1\cap H_2$ is not trivial, we have the same formula since $i(s_1,s_2)=0$ and the area of the parallelogram formed by the two vectors $(a_1,b_1),(a_2,b_2)$ equals $0$.

This result is well-known in the case that $s_1,s_2$ are simple closed curves on the torus $\gS =\mathbb{Z}^2\backslash \RR^2$ (see \cite[1.2.3 Intersection Numbers]{FM12}).
\end{supply}

\subsection{Continuous extension of intersection number}\label{conti ext of int number}

First, we recall several facts on geodesic currents on hyperbolic groups in \cite{Bon88b}.

Let $G$ be an infinite hyperbolic group.
Set
\[ \partial_2 G:= \{ S\in \H (\partial G)\mid \# S=2\}.\]
We endow $\partial_2 G$ with the subspace topology of $\H (\partial G)$, which coincides with the topology induced by the Hausdorff distance.
\begin{definition}[Geodesic currents on hyperbolic groups]
A \ti{geodesic current} on $G$ is a $G$-invariant locally finite Borel measure on $\partial_2 G$.
The space of geodesic currents on $G$ is denoted by $\GC (G)$, which is equipped with the weak-$\ast$ topology.
\end{definition}

Since $\partial_2G$ is a $G$-invariant closed subspace of $\H (\partial G)$, we can consider $\GC (G)$ as an $\RRR$-linear closed subspace of $\SC (G)$. A subset current on $G$ whose support is included in $\partial_2 G$ can be considered as a geodesic current on $G$. By restricting a subset current to $\partial _2 G$, we can obtain an $\RRR$-linear map from $\SC (G)$ to $\GC (G)$ but this map is not continuous in general (see Theorem \ref{thm:geodesic currents approximated by subset currents}). We will construct a continuous $\RRR$-linear projection $\B$ from $\SC (G)$ to $\GC (G)$ in the case that $G$ is the fundamental group of a compact hyperbolic surface (see Section \ref{sec: projection B}).

For $g\in G$ with infinite order, since its limit set $\gL (\langle g \rangle )$ belongs to $\partial_2 G$, the counting subset current $\eta_{\langle g\rangle }$ can be considered as a geodesic current on $G$. We will write $\eta_g$ in place of $\eta_{\langle g\rangle}$ and call $\eta_g$ the \ti{counting geodesic current} for $g\in G$. If $g\in G$ has a finite order, then we define $\eta_g$ to be the zero measure on $\partial_2 G$. A geodesic current $\mu$ is called \ti{rational} if there exist $g\in G$ and $r\in \RRR$ such that $\mu=c\eta_g$.

Bonahon \cite{Bon88b} proved the following theorem:

\begin{theorem}[See {\cite[Theorem 7]{Bon88b}}]\label{thm: Bonahon hyperbolic group geodesic currents dense}
For any infinite hyperbolic group $G$, the set of all rational geodesic currents on $G$ is a dense subset of $\GC (G)$.
\end{theorem}

In the case of subset currents, the same denseness property was proved for free groups of finite rank in \cite[Theorem 5.8]{KN13}.
In Subsection \ref{subsec:denseness property of surface groups} we will prove that surface groups have the denseness property of rational subset currents. 

If an infinite hyperbolic group $G$ is virtually cyclic, that is, $\# \partial G=2$, then $\SC (G)$ coincides with $\GC (G)$.
If $G$ is not virtually cyclic, which is called \ti{non-elementary}, then we have the following theorem, which means that any geodesic current can be approximated by a sequence of rational subset currents in $\SC(G)\setminus \GC(G)$.

\begin{theorem}\label{thm:geodesic currents approximated by subset currents}
Let $G$ be a non-elementary hyperbolic group. Then the (Gromov) boundary $\partial G$ includes uncountably many points. For any $\mu \in \GC (G)$ there exists a sequence $\{ H_n \}_{n\in \NN}$ of quasi-convex subgroups of $G$ and a sequence $\{ c_n \}_{n \in \NN}$ of $\RRR$ such that $H_n$ is non-cyclic and isomorphic to a free group of finite rank, and the sequence of rational subset currents $c_n \eta_{H_n}$ converges to $\mu$.
\end{theorem}
\begin{proof}
From Theorem \ref{thm: Bonahon hyperbolic group geodesic currents dense}, it is sufficient to prove the statement in the case that $\mu =\eta_g$ for $g\in G$ with infinite order.

Take $g\in G$ with infinite order. Take $h\in G$ with infinite order such that $\gL (\langle h \rangle )\cap \gL (\langle g \rangle )=\emptyset$.
By using the Ping-Pong Lemma, for a sufficiently large $m\in \NN$ the subgroup $H:=\langle g^m ,h^m \rangle $ is isomorphic to the free group of rank $2$ (see \cite[Part III, $\gG$, 3.20 Proposition]{FM12}). Moreover, we can see that if $m$ is sufficiently large, then $H$ is a quasi-convex subgroup of $G$.

Set $a:=g^m, b:=h^m$. Define a subgroup $H_n$ of $H$ by
\[ H_n :=\langle a^n ,b\rangle \]
for $n\in \NN$. Then we can see that the sequence of rational counting subset currents $\frac{1}{n} \eta_{H_n}^H$ on $H$ converges to the counting geodesic current $\eta_a^H$ on $H$ by using \cite[Proposition 3.7]{KN13} (see Proposition \ref{prop:convergence and subset cylinder} for detail).
By using the map $\iota_H$ in Section \ref{sec:relation between subgroups}, we see that $\frac{1}{n} \eta_{H_n}$ converges to $\eta_a$.
Note that
\[ \eta_a= \eta_{g^m} =m \eta_g \] 
by Proposition \ref{prop:property}.
Hence $\frac{1}{mn} \eta_{H_n}$ converges to $\eta_g$ as $n\rightarrow \infty$.
\end{proof}

\begin{notation}
Let $\gS$ be a compact hyperbolic surface (possibly with geodesic boundary) and $\tilde{\gS}$ the universal cover of $\gS$, which we considered as a convex subspace of $\HH$. Let $G$ be the fundamental group of $\gS$, which is isomorphic to a free group of finite rank or a surface group. 
Note that the limit set $\gL(G)$ with respect to the action of $G$ on $\tilde{\gS}$ coincides with the limit set $\tilde{\gS}(\infty)$ of $\tilde{\gS}$.
\emph{When we identify the boundary $\partial G$ of $G$ with the limit set $\gL(G)$ by using the action of $G$ on $\tilde{\gS}$, we will say subset currents on $\gS$ instead of subset currents on $G$. The term geodesic currents on $\gS$ is also used in the same way.}
The space of subset currents (or geodesic currents) on $\gS$ will be denoted by $\SC(\gS)$ (or $\GC(\gS)$, respectively).
\end{notation}

Recall that a non-trivial conjugacy class of $G$ is corresponding to a non-trivial free
homotopy class of an oriented closed curve on $G$, which contains a unique oriented closed geodesic.
Hence a non-trivial conjugacy class of $G$ is corresponding to an oriented closed geodesic on $G$.
In addition, for non-trivial $g\in G$ the conjugacy class of $\langle g \rangle$ is corresponding to an unoriented closed geodesic on $\gS$, which coincides with the convex core $C_{\langle g \rangle }$. The map $p_{\langle g \rangle}$ from $C_{\langle g \rangle }$ to $\gS$ is induced by the universal covering map.
We will write $C_g$ in place of $C_{\langle g \rangle}$ and call $C_g$ the (unoriented) closed geodesic corresponding to $g$.

Bonahon \cite{Bon86} proved the following theorem:

\begin{theorem}[See {\cite[Proposition 4.5]{Bon86}}]\label{thm:intersection number of geodesic currents}
Let $\gS$ be a compact hyperbolic surface. There exists a unique continuous symmetric $\RRR$-bilinear functional 
\[ i_{\GC}\: \GC (\gS)\times \GC (\gS)\rightarrow \RRR \]
such that for any non-trivial $g,h\in G$ we have
\[ i_{\GC}(\eta_g, \eta_h )=i (C_g,C_h).\]
\end{theorem}

Recall that a subgroup $H$ of $G=\pi_1(\gS)$ is quasi-convex if and only if $H$ is finitely generated.
For two non-trivial finitely generated subgroups $H$ and $K$ of $G$, we have the convex cores $(C_H, p_H)$ and $(C_K,p_K)$ of $H$ and $K$.
From Theorem \ref{thm:geodesic criterion}, $(C_H, p_H)$ and $(C_K,p_K)$ are simple compact surfaces on $\gS$ in minimal position.
We will prove the following theorem in this subsection, which is a generalization of Theorem \ref{thm:intersection number of geodesic currents}:

\begin{theorem}[Intersection number of subset currents]\label{thm:intersection number of subset currents}
Let $\gS$ be a compact hyperbolic surface. There exists a unique continuous symmetric $\RRR$-bilinear functional 
\[ i_{\SC}\: \SC (\gS)\times \SC (\gS)\rightarrow \RRR \]
such that for any non-trivial finitely generated subgroups $H$ and $K$ of $G$ we have
\[ i_{\SC}(\eta_H, \eta_K )=i (C_H,C_K).\]
\end{theorem}

\begin{remark}
In the case that $\gS$ has boundary, $G$ is a free group of finite rank. We remark that for a free group $F$ of finite rank a surface whose fundamental group is isomorphic to $F$ is not unique up to homeomorphism. Therefore the functional $i_\SC$ on $\SC (F)$ is not uniquely determined. 

However, if $G$ is a surface group, then a surface whose fundamental group is isomorphic to $G$ is unique up to homeomorphism.
Note that the property that $\partial G$ is homeomorphic to $S^1$ plays an essential role.
For two non-trivial finitely generated subgroup $H$ and $K$ of $G$ we can see that the intersection number of $C_H$ and $C_K$ equals the number of equivalence classes $[g_1H,g_2K]\in G \backslash G/H\times G/K$ satisfying the condition that $g_1CH_H\cap g_2CH_K$ is a compact convex polygon, which depends only on the ``positional relation'' between $g_1\gL(H)$ and $g_2\gL(K)$ (see Case 4 of the proof of Theorem \ref{thm:bigon criterion 3}).
Therefore, in the case that $G$ is a surface group, we can call $i_{\SC}$ the intersection number \it{on} $\SC (G)$.
\end{remark}

The strategy to prove Theorem \ref{thm:intersection number of subset currents} is almost the same as that for proving the existence of the volume functional in Section \ref{sec:Volume functionals for Kleinian groups}.
First, we construct an $\RRR$-bilinear functional on $\SC(\gS)$ such that the functional associates any pair of counting subset currents $(\eta_H, \eta_K)$ with $i(C_H ,C_K)$ for any non-trivial finitely generated subgroups $H$ and $K$ of $G$. Then we prove the continuity of the functional, which is the main part of the proof. The uniqueness of the functional follows by the denseness property of rational subset currents.

Note that by restricting $i_{\SC}$ to $\GC (\gS)\times \GC (\gS)$ we can obtain $i_{\GC}$. If we want to obtain only $i_{\GC}$, then by assuming that $H,K$ are cyclic and all $(S_1,S_2)\in \H (\partial G)\times \H (\partial G)$ belong to $\partial_2G\times \partial_2G$, several parts of the following argument will be shorter or obvious, and our argument will give a new proof to Theorem \ref{thm:intersection number of geodesic currents}.

Recall that for simple compact surfaces $(S_1,s_1),(S_2,s_2)$ on $\gS$ we constructed the pairs $(\hat{S_1},\hat{s_1}), (\hat{S_2},\hat{s_2})$ and the fiber product $\hat{S_1}\times_{\tilde{\gS}}\hat{S_2}$.
Let $H,K$ be non-trivial finitely generated subgroups of $G$.
From Remark \ref{rem:expression of hat S i} and Proposition \ref{prop:characterize fiber product}, we set
\[\hat{CH}_H:=\{ (gH,x)\in G/H\times \tilde{\gS}\mid x\in gCH_H\} \]
and set
\begin{align*}
\hat{CH}_H\times_{\tilde{\gS}}\hat{CH}_K:=\{ &(g_1H,g_2K,x)\in G/H\times G/K\times \tilde{\gS}\ |\\
	&x\in g_1CH_H\cap g_2CH_K\}. 
\end{align*}
Then $G$ acts on $\hat{CH}_H$ by
\[ u(gH,x):=(ugH,ux)\]
for $u\in G$ and $(gH,x)\in \hat{CH}_H$. Moreover, $G$ acts on $\hat{CH}_H\times_{\tilde{\gS}}\hat{CH}_K$ by
\[ u(g_1H,g_2K,x):=(ug_1H,ug_2K,ux)\]
for $u\in G$ and $(g_1H,g_2K,x)\in \hat{CH}_H\times_{\tilde{\gS}}\hat{CH}_K$.
By the same way as that for simple compact surfaces $(S_1,s_1),(S_2,s_2)$ on $\gS$ in Proposition \ref{prop:cubic diagram}, we can obtain the following cubic commutative diagram for $H$ and $K$:
\[
\xymatrix{
\CH_H \times _{\tilde \gS}\CH_K \ar[rrr]\ar[ddd]_{\Phi}\ar[ddr] & & &\text{\quad }\CH_K\text{\quad }\ar[ddr]\ar[ddd]|(.64)\hole &\\
&&&&\\
&\text{\quad }\CH_H\text{\quad }\ar[rrr]\ar[ddd] &&&\text{\quad }\tilde{\gS}\text{\quad }\ar[ddd]^{\pi}\\
C_H\times_{\gS}C_K \ar[ddr]\ar[rrr]|(.455)\hole &&&\text{\quad }C_K\text{\quad }\ar[ddr]_{p_K}&\\
&&&&\\
&\text{\quad }C_H\text{\quad } \ar[rrr]_{p_H}&&&\text{\quad }\gS \text{\quad }
}
\]

The map from $\CH_H$ to $\tilde{\gS}$ and the map from $\hat{CH}_H\times_{\tilde{\gS}}\hat{CH}_K$ to $\CH_H$ are the projections. Explicitly, $(gH_,x)\in \CH_H$ is mapped to $x\in \tilde{\gS}$ and $(g_1H,g_2K,x)\in \hat{CH}_H\times_{\tilde{\gS}}\hat{CH}_K$ is mapped to $(g_1H,x)\in \CH_H$.
The quotient space $G\backslash \CH_H$ is identified with $C_H$ and the quotient space $G\backslash \hat{CH}_H\times_{\tilde{\gS}}\hat{CH}_K$ is identified with $C_H\times_\gS C_K$ by Proposition \ref{prop:characterize fiber product}.

By the definition, $i(C_H,C_K)$ equals the number of contractible components of $C_H\times_\gS C_K$.
Each contractible component of $C_H\times_\gS C_K$ comes from the $G$-orbit of a compact connected component of $\hat{CH}_H\times_{\tilde{\gS}}\hat{CH}_K$.

In the following definition, we introduce the notion of the size of a compact connected component of $\hat{CH}_H\times_{\tilde{\gS}}\hat{CH}_K$ by using a fundamental domain $\F$ for the action of $G$ on $\tilde{\gS}$.

\begin{definition}[Size of a compact connected component]\label{def:size of cpt comp}
For $x\in \tilde{\gS}$ we take the Dirichlet domain $\F=\F_x$ centered at $x$, which is a compact convex polygon. By removing some edges and vertices of the boundary of $\F$ we can modify $\F$ such that $G(\F)=\tilde{\gS}$ and $g\F\cap \F=\emptyset$ for any non-trivial $g\in G$.
We define $\mathrm{Fin}(G)$ to be the family of all non-empty finite subset of $G$. Note that for any non-empty bounded subset $X$ of $\tilde{\gS}$ there exists a unique $G_0\in \mathrm{Fin}(G)$ such that $G_0(\F)$ covers $X$ precisely, that is, $X\subset G_0(\F)$ and $X\cap g\F \not=\emptyset$ for every $g\in G_0$.
Then we say that the \ti{size} of $X$ with respect to $\F$ is $G_0$.
For $G_0\in \mathrm{Fin}(G)$ we define $C_{\F}(G_0;H,K)$ to be the number of compact connected components of $\hat{CH}_H\times_{\tilde{\gS}}\hat{CH}_K$ whose size with respect to $\F$ are $G_0$.
\end{definition}

Now, we consider the natural action of $G$ on $\mathrm{Fin}(G)$ from left and take a complete system of representatives $\{ G_j\}_{j\in J}$ of $G\backslash \mathrm{Fin}(G)$.

\begin{lemma}
The following equality holds:
\[ i(C_H,C_K)=\sum_{j\in J}C_{\F}(G_j;H,K).\]
\end{lemma}
\begin{proof}
Since $i(C_H,C_K)$ is the number of contractible components of $G\backslash \hat{CH}_H\times_{\tilde{\gS}}\hat{CH}_K$, it is sufficient to see that for any compact connected component $M$ of $\hat{CH}_H\times_{\tilde{\gS}}\hat{CH}_K$ there exist unique $j\in J$ and $g\in G$ such that $M$ is precisely covered by $gG_j(\F)$.
Actually, we have a unique $G_0\in \mathrm{Fin}(G)$ such that $G_0(\F)$ cover $M$ precisely and there exists unique $j\in J$ and $g\in G$ such that $gG_j=G_0$. Hence our claim follows.
\end{proof}

For $G_0\in \mathrm{Fin}(G)$ set
\begin{align*}
C_{\F}(G_0):=\{ &(S_1,S_2)\in \H (\partial G)\times \H (\partial G)\ |\\ 
&CH(S_1)\cap CH(S_2) \text{ is precisely covered by }G_0(\F)\}.
\end{align*}
We can check that $C_\F(G_0)$ is a Borel subset of $\H (\partial G)\times \H (\partial G)$ from Lemma \ref{lem:taking convex hull is continuous}.
Then for the product measure $\eta_H\times \eta_K$ we have 
\[ \eta_H\times \eta_K (C_\F (G_0))=C_\F (G_0;H,K).\]
Actually,
\begin{align*}
\eta_H\times \eta_K 
&=\left( \sum_{gH\in G/H}\delta_{g\gL (H)}\right)\times \left( \sum_{gK\in G/K}\delta_{g\gL (K)}\right) \\
&= \sum_{(g_1H,g_2K) \in G/H\times G/K}\delta_{g_1\gL (H)}\times \delta_{g_2\gL (K)} \\
&=\sum_{(g_1H,g_2K) \in G/H\times G/K}\delta_{(g_1\gL (H),g_2\gL (K))},
\end{align*}
where $\delta_{(g_1\gL (H),g_2\gL (K))}$ is the Dirac measure at $(g_1\gL (H),g_2\gL (K))$ on $\H (\partial G )\times \H (\partial G)$.
In addition,
\[ \hat{CH}_H\times_{\tilde{\gS}}\hat{CH}_K\cong \bigsqcup_{(g_1H,g_2K) \in G/H\times G/K}g_1CH_H\cap g_2CH_K.\]
Hence
\begin{align*}
\eta_H\times \eta_K (C_\F (G_0))
=&\# \{ (g_1H,g_2K) \in G/H\times G/K \ |\\
 &\qquad g_1 CH_H\cap g_2 CH_K\text{ is precisely covered by }G_0(\F)\} \\
=&C_\F (G_0;H,K).
\end{align*}
As a result, we obtain the following equation:
\[ i(C_H,C_K)=\sum_{j\in J}\eta_H\times \eta_K (C_{\F}(G_j)).\]

Note that for $G_1,G_2\in \mathrm{Fin}(G)$ with $G_1\not=G_2$ the intersection $C_{\F}(G_1)\cap C_{\F}(G_2)$ is empty by the definition.

\begin{definition}
We define a map $i_{\SC}$ from $\SC (\gS)\times \SC (\gS)$ to $\RRR$ by
\[ i_{\SC}(\mu,\nu ):=\mu \times \nu \left( \bigsqcup_{j\in J}C_\F (G_j)\right) \]
for $\mu ,\nu \in \SC (\gS)$.
\end{definition}

By the definition of $i_{\SC}$ we can see that $i_{\SC}(\eta_H,\eta_K)=i(C_H,C_K)$ for any non-trivial finitely generated subgroups $H$ and $K$ of $G$. Moreover, $i_{\SC}$ is a symmetric $\RRR$-bilinear functional. The remaining problem is proving the continuity of $i_{\SC }$.

First, we check that the definition of $i_{\SC}$ is independent of the choice of the fundamental domain $\F$ and $\{ G_j \}$.
Set
\[ \mathcal{I}:=\{ (S_1,S_2)\in \H (\partial G)\times \H(\partial G)\mid CH(S_1)\cap CH(S_2)\text{ is $\not=\emptyset$ and bounded} \}. \]
Then $\mathcal{I}$ is a $G$-invariant open subset of $\H (\partial G)\times \H (\partial G)$ with respect to the diagonal action of $G$.
Moreover, $G$ acts on $\mathcal{I}$ freely.

\begin{lemma}
The set
\[ Q:=\bigsqcup_{j\in J}C_\F (G_j)\]
is a Borel fundamental domain for the action of $G$ on $\mathcal{I}$ satisfying the condition that $G(Q)=\mathcal{I}$ and $gQ\cap Q$ is empty for any non-trivial $g\in G$.
Therefore, the definition of $i_{\SC }$ is independent of the choice of $\F$ and $\{ G_j \}$.
\end{lemma}
\begin{proof}
First, we remark that the definition $\mathcal{I}$ is independent of the choice of $\F$ and $\{ G_j\}$. Moreover, in the case that $G$ is a surface group, the definition of $\mathcal{I}$ is independent of $\gS$.

For $(S_1,S_2)\in \mathcal{I}$ there exists a unique $G_0\in \mathrm{Fin}(G)$ such that $G_0(\F)$ cover $CH(S_1)\cap CH(S_2)$ precisely. Hence
\[ \mathcal{I}=\bigsqcup_{G_0\in \mathrm{Fin}(G)} C_\F(G_0) .\]
Then we can see that $G(Q)=\mathcal{I}$ and $gQ\cap Q$ is empty for any non-trivial $g\in G$, which implies that $Q$ is a Borel fundamental domain for the action of $G$ on $\mathcal{I}$.
By the same way as that for Lemma \ref{lem:volume of fundamental domain}, we can see that $i_{\SC}$ is independent of the choice of $\F$ and $\{ G_j \}$.
\end{proof}

The following proposition is known as the Portmanteau theorem for probability measures on a metric space (see \cite[Theorem 2.1]{Bil99}), which will be used later in order to prove the continuity of $i_{\SC}$. We will use the argument in this proof for proving the continuity of a certain functional in Section \ref{sec:an intersection functional}.

\begin{proposition}\label{prop:Portmanteau}
Let $\mu_n,\mu\in \SC (\gS)\ (n\in \NN)$. The following are equivalent:
\begin{enumerate}
\item $\mu_n$ converges to $\mu$ in the weak-$\ast$ topology;
\item $\limsup_{n\rightarrow \infty} \mu_n (K)\leq \mu(K)$ for any compact subset $K$ of $\H (\partial G)$, and\\
$\liminf_{n\rightarrow \infty} \mu_n (U)\geq \mu(U)$ for any relatively compact open subset $U$ of $\H(\partial G)$;
\item $\lim_{n\rightarrow \infty} \mu_n(E)=\mu (E)$ for any relatively compact Borel subset $E$ of $\H (\partial G)$ with $\mu (\partial E)=0$;
\item $\lim_{n\rightarrow \infty}\int f d\mu_n=\int f d\mu$ for any bounded function $f:\H (\partial G)\rightarrow \RRR$ with compact support which is continuous at $\mu$-$a.e$, that is, the set $\Delta (f)$ of non-continuous points of $f$ has measure zero with respect to $\mu$. Explicitly, 
\[ \Delta (f):=\{ S\in \H (\partial G)\mid f\ \text{is not continuous at } S\}.\] 
\end{enumerate}
For product measures $\mu_n \times \nu_n \ (\mu_n, \nu_n \in \SC (\gS ), n\in \NN)$ and $\mu \times \nu \ (\mu ,\nu \in \SC (\gS))$ the same result follows by the same proof.
\end{proposition}
\begin{proof}
Since $\H (\partial G)$ is a locally compact, separable and metrizable space, we can take a metric $d$ on $\H (\partial G)$ compatible with the topology such that $(\H (\partial G) ,d)$ is a proper metric space, that is, every closed ball with respect to $d$ is a compact subset of $\H (\partial G)$. We will use this property of $d$ in the proof of $(3)\Rightarrow (4)$.

$(4)\Rightarrow (1)$: Obvious.

$(1)\Rightarrow (2)$: For a compact subset $K$ of $\H(\partial G)$, set
\[ K_n:=\{ x\in \H (\partial G)\mid d(x,K)< \frac{1}{n}\}\]
for $n\in \NN$.
Then the characteristic function $\chi_{K_n}$ converges pointwise to $\chi_K$, which implies that
\[ \mu(K_n)=\int \chi_{K_n}d\mu \xrightarrow{n\rightarrow \infty} \int \chi_K d\mu=\mu (K) .\]
Fix $\varepsilon>0$. Then there exists $N\in \NN$ such that $\mu(K_N)\leq \mu(K)+\varepsilon$. By the Urysohn Lemma we have a continuous function $f\: \H (\partial G)\rightarrow \RR$ satisfying the condition that $f|_{K} \equiv 1$, $f|_{(K_N)^c}\equiv 0$ and $0\leq f(S)\leq 1$ for any $S\in \H(\partial G)$. Then we have
\[ \limsup_{n\rightarrow \infty} \mu_n(K)\leq \limsup_{n\rightarrow \infty} \int fd\mu_n= \int fd\mu \leq \mu (K_N)\leq \mu (K)+\varepsilon.\]
Since $\varepsilon>0$ is arbitrary, we have
\[ \limsup_{n\rightarrow \infty} \mu_n (K)\leq \mu(K).\]

For a relatively compact open subset $U$ of $\H(\partial G)$, set
\[ U_n:=\{ x\in U\mid d(x,U^c)\geq \frac{1}{n}\} \]
for $n\in \NN$. Then the characteristic function $\chi_{U_n}$ converges pointwise to $\chi_U$, which implies that
\[ \mu(U_n)=\int \chi_{U_n}d\mu \xrightarrow{n\rightarrow \infty} \int \chi_U d\mu=\mu(U) .\]
Fix $\varepsilon>0$. Then there exists $N\in \NN$ such that $\mu(U_N)\geq \mu(U)-\varepsilon$. By the Urysohn Lemma we have a continuous function $f\: \H (\partial G)\rightarrow \RR$ satisfying the condition that $f|_{U_N} \equiv 1$, $f|_{(U)^c}\equiv 0$ and $0\leq f(S)\leq 1$ for any $S\in \H(\partial G)$. Then we have
\[ \liminf_{n\rightarrow \infty} \mu_n(U)\geq \liminf_{n\rightarrow \infty}\int fd\mu_n = \int f d\mu \geq \mu (U_N)\geq \mu (U)-\varepsilon.\]
Since $\varepsilon>0$ is arbitrary, we have
\[ \liminf_{n\rightarrow \infty} \mu_n (U)\geq \mu(U).\]

$(2)\Rightarrow (3)$: Since $\mathrm{Int}(E)\subset E\subset \ol E$ and $\partial E=\ol E\setminus \mathrm{Int}(E)$, we have
\[ \mu (\mathrm{Int}(E))=\mu(E)=\mu (\ol E).\]
Therefore,
\begin{align*}
&\limsup_{n\rightarrow \infty}\mu_n(E)\leq \limsup_{n\rightarrow \infty }\mu_n (\ol E)\leq \mu(\ol E)=\mu(E) \\
=&\mu (\mathrm{Int}(E))\leq \liminf_{n\rightarrow \infty}\mu_n (\mathrm{Int}(E))\leq \liminf_{n\rightarrow \infty }\mu _n(E),
\end{align*}
and so
\[ \lim_{n\rightarrow \infty}\mu_n (E)=\mu (E).\]

$(3)\Rightarrow (4)$: This is the main part of this proof. We can assume that $f\geq 0$ without loss of generality. Let $\mathrm{supp}f$ denote the support of $f$. Set
\[ C:=\sup \{ f(x)\mid x\in \H(\partial G)\}\]
and
\[ A_t:=\{ x\in \H (\partial G)\mid f(x)\geq t\} \]
for each $t\in [0,C]$. Note that $\int f d\mu$ equals the area of
\[ U:=\{ (x,y)\in \H (\partial G)\times \RR \mid 0\leq y \leq f(x) \}\]
with respect to the product measure of $\mu \times m_\RR$, where $m_\RR$ is the Lebesgue measure on $\RR$.
Since
\[ U=\{ (x,y)\in \H (\partial G)\times \RR \mid y\in [0,C] , x\in A_y \}, \]
we have
\[ \int f d\mu =\int_0^C \mu (A_t) dm_\RR (t),\quad \int f d\mu_n =\int_0^C \mu_n (A_t) dm_\RR (t).\]
By the bounded convergence theorem, it is sufficient to prove that $\mu_n(A_t)\ (t\in (0,C])$ is uniformly bounded,
and $\mu_n(A_t)$ converges pointwise to $\mu(A_t)$ for $m_\RR$-a.e. $t\in [0,C]$.

First, we see that $\mu_n(A_t)\ (n\in \NN, t\in (0,C])$ is uniformly bounded. Note that $A_0=\H (\partial G)$.
For any $t\in (0,C]$, $A_t$ is included in $L:=\mathrm{supp} f$, which is compact.
Hence it is sufficient to see that $\mu_n(L)$ is bounded.
Since $(\H (\partial G), d)$ is a proper metric space, the closed $r$-neighborhood of $L$, denoted by $B(L,r)$, is also compact for $r\geq 0$.
Set $C(L,r):=\{ x\in \H (\partial G) \mid d(L,x)=r \}$ for $r>0$, which includes the boundary $\partial B(L,r)$.
Then we have
\[ B(L,1)=L\sqcup \bigsqcup_{0<r\leq 1} C(L,r).\]
Since the interval $(0,1]$ is an uncountable set, there exists $r_0\in (0,1]$ such that $C(L,r_0)$ has zero measure with respect to $\mu$ (see Lemma \ref{lem:uncountably family has measure zero set} for more general statement).
Then $\mu(\partial B(L,r_0))=0$, which implies that $\mu_n (B(L,r_0))$ converges to $\mu (B(L,r_0))$ by the assumption.
Hence there exists $M>0$ such that $\mu_n(B(L,r_0))\leq M$ for any $n\in \NN$, and so
\[ \mu_n(A_t)\leq \mu_n(L)\leq \mu_n(B(L,r_0))\leq M\]
for any $n\in \NN, t\in (0,C]$.

Next, we see that $\mu_n(A_t)$ converges pointwise to $\mu(A_t)$ for $m_\RR$-a.e. $t\in [0,C]$.
From the assumption (3), it is sufficient to see that for $m_\RR$-a.e. $t\in[0,C]$ we have $\mu(\partial A_t)=0$.
Set
\[ B_t:=\{ x\in \H (\partial G)\mid f(x)= t\} \]
for $t\in [0,C]$. We prove that $\partial A_t\subset B_t\cup \Delta (f)$ for each $t\in [0,C]$.
Take $x\in \partial A_t$ and assume that $f$ is continuous at $x$, which implies that $x\not\in \Delta (f)$. 
If $f(x)>t$, then there exists an open neighborhood $V$ of $x$ such that for any $x'\in V$ we have $f(x')>t$, which implies that $V\subset A_t$ and contradicts the assumption that $x\in \partial A_t$. Therefore $f(x)=t$ and $x\in B_t$.

Since $\mu (\Delta(f))=0$, it is sufficient to prove that for $m_\RR$-a.e. $t\in [0,C]$ we have $\mu (B_t)=0$.
Note that $\mu (A_t)$ is a decreasing function with respect to $t$. Therefore $\mu (A_t)$ has at most countably many non-continuous 
points. If $\mu(A_t)$ is continuous at $t_0>0$, then $B_{t_o}\subset (A_{t_0-\varepsilon}\setminus A_{t_0+\varepsilon})$ for any small $\varepsilon>0$ and
\[ 0\leq \mu (B_{t_0})\leq \lim_{\varepsilon \searrow 0}(\mu (A_{t_0-\varepsilon})-\mu (A_{t_0+\varepsilon}))=0.\]
Therefore $\mu (B_t)=0$ except countably many points of $[0,C]$. This completes the proof.
\end{proof}

In order to prove the continuity of $i_{\SC}$, we focus on the boundary of $C_\F (G_j)$ for $j\in J$.
We assume that $G_j$ contains $\id$ for every $j\in J$.

Since $CH(S_1)\cap CH(S_2)$ is a compact convex subset of $\HH$ surrounded by geodesics for $(S_1,S_2)\in \mathcal{I}$, $CH(S_1)\cap CH(S_2)$ can be considered as a polygon.
We define $B_\F$ to be a subset of $\mathcal{I}$ consisting of points $(S_1,S_2)$ satisfying one of the following conditions:
\begin{enumerate}
\item[$B_\F 1$)] a vertex of $CH(S_1)\cap CH(S_2)$ belongs to $\partial \F$;
\item[$B_\F 2$)] an edge of $CH(S_1)\cap CH(S_2)$ overlaps an edge of $\ol{\F}$;
\item[$B_\F 3$)] an edge of $CH(S_1)\cap CH(S_2)$ is tangent to a vertex of $\ol{\F}$.
\end{enumerate}
A geodesic $\ell$ in $\HH$ is said to be tangent to a vertex of a (convex) polygon $P$ of $\HH$ if the intersection of $\ell$ and $P$ is exactly the vertex. Note that $B_\F$ does not depend on edges and vertices removed from the Dirichlet domain $\F_x$. Hence for any $y\in \tilde{\gS}$ and the Dirichlet domain $\F_y$ centered at $y$ we can define $B_{\F_y}$ as above.
Set
\[
\Delta_\F :=\{ (S,S)\in \partial_2 G\times \partial _2 G\mid CH(S)\cap \ol{\F}\not= \emptyset \}.
\]
We see that $B_\F$ and $\gD_\F$ are closed subsets of $\H (\partial G)\times \H (\partial G)$.

\begin{lemma}\label{lem:boundary is included in the union}
For $\{ \id \}\in \mathrm{Fin}(G)$ the boundary $\partial C_\F (\{ \id \})$ of $C_\F (\{ \id \})$ in $\H (\partial G)\times \H (\partial G)$ is included in the union of $B_\F$ and $\Delta_\F$.
\end{lemma}
\begin{proof}
First of all, we remark that for $S\in \H(\partial G)$ with $\#S\geq 3$ for any interior point $z$ of $CH(S)$ there exists an open neighborhood $U$ of $S$ in $\H (\partial G)$ such that for any $S'\in U$ the convex hull $CH(S')$ also contains $z$ as an interior point from Lemma \ref{lem:taking convex hull is continuous}.

Let $(S_1,S_2)\in \partial C_\F (\{ \id \})$. By the definition, for any open neighborhood $O$ of $(S_1,S_2)$ both $O \cap C_\F (\{ \id \})$ and $O\cap C_\F (\{ \id \} )^c$ are non-empty.

\underline{Claim 1:} \emph{If $S_1=S_2=:S$, then $(S,S)\in \Delta_\F$.}

If $\#S\geq 3$, then the interior of $CH(S)$ is not bounded in $\tilde{\gS}$ and we can take $z\in \mathrm{Int}(CH(S))\setminus \ol\F$. Then take an open neighborhood $U$ of $S$ such that for any $S'\in U$ the convex hull $CH(S')$ also contains $z$ as an interior point. Now, we can see that $U\times U$ is an open neighborhood of $(S,S)$ and does not intersect $C_\F (\{ \id \})$, which contradicts the assumption that $(S,S)\in \partial C_\F (\{\id \})$. Hence $\#S=2$. If $CH(S)$ does not intersects $\ol \F$, then there exists a neighborhood $U$ of $S$ such that for any $S'\in U$ the convex hull $CH(S')$ does not intersects $\ol \F$, which contradicts the assumption that $(S,S)\in \partial C_\F(\{ \id \})$. Therefore $(S,S)\in \Delta_\F$. This argument will be used frequently in this proof but we will not remark it from now on.

\underline{Claim 2:} \emph{If $S_1\not=S_2$, then $S_1\cap S_2=\emptyset$.}

To obtain a contradiction, suppose that $S_1\not=S_2$ and $S_1\cap S_2\not=\emptyset$. From the proof of Claim 1, $\#(S_1\cap S_2)$ must be smaller than $3$ and the interior of $CH(S_1)\cap CH(S_2)$ must be included in $\ol{\F}$. Since $S_1\not= S_2$, we can assume that $S_1\geq 3$ from Claim 1.
If $\#S_2=2$, then $CH(S_2)$ is a boundary component of $CH(S_1)$ or included in the interior of $CH(S_1)$. In both cases, there exists an open neighborhood $U$ of $(S_1,S_2)$ such that $U\subset C_\F (\{ \id \} )^c$, a contradiction.

Now, we can assume that $\#S_1,\#S_2\geq 3$. If $\#(S_1\cap S_2)=1$, then $CH(S_1)\cap CH(S_2)$ must be empty and there exists an open neighborhood $U$ of $(S_1,S_2)$ such that $U\subset C_\F (\{ \id \} )^c$, a contradiction.
If $\#(S_1\cap S_2)=2$, then $CH(S_1)$ and $CH(S_2)$ have one common boundary component and $\mathrm{Int}(CH(S_1))\cap \mathrm{Int}(CH(S_2))$ is empty. Even in this case, there exists an open neighborhood $U$ of $(S_1,S_2)$ such that
$U\subset C_\F (\{ \id \} )^c$.
Therefore in any cases we can obtain a contradiction.

\underline{Claim 3:} \emph{If $S_1\not=S_2$, then $(S_1,S_2)\in B_\F$.}

Since $S_1\cap S_2=\emptyset$, the intersection $CH(S_1)\cap CH(S_2)$ should be non-empty and bounded. If $CH(S_1)\cap CH(S_2)$ contains an exterior point of $\ol \F$, then $(S_1,S_2)\not\in \partial C_\F(\{ \id \})$ from the proof of Claim 1. Hence $CH(S_1)\cap CH(S_2)$ is included in $\ol \F$. If $CH(S_1)\cap CH(S_2)$ is included in the interior of $\F$, then for $(S_1',S_2')$ sufficiently close to $(S_1,S_2)$ the intersection $CH(S'_1)\cap CH(S'_2)$ is also included in the interior of $\F$. Therefore, $CH(S_1)\cap CH(S_2)$ is not included in the interior of $\F$, which implies that $(S_1,S_2)$ satisfies the condition $(B_\F 1)$ or $(B_\F 2)$.
\end{proof}

\begin{lemma}
For $G_j\in \mathrm{Fin}(G)$ the boundary $\partial C_\F (G_j)$ is included in $G_j(B_\F\sqcup \Delta_\F)$.
\end{lemma}
\begin{proof}
Let $(S_1,S_2)\in \partial C_\F (G_j)$.
By the same way for Claim 1 in the above lemma, we can see that if $S_1=S_2=:S$, then $(S,S)\in G_j(\Delta_\F)$.
Note that
\[
G_j(\Delta_\F ) =\{ (S,S)\in \partial_2 G\times \partial _2 G\mid CH(S)\cap G_j(\ol{\F})\not=\emptyset \}
\]
and that such $(S,S)$ may not exist since in this case $\# S=2$ and the convex hull $CH(S)$ needs to intersect $g\ol \F$ for every $g\in G_j$. 

By the same way for Claim 2 in the above lemma, we can see that if $S_1\not=S_2$, then $S_1\cap S_2=\emptyset$.
Now, we prove that if $S_1\cap S_2=\emptyset$, then $(S_1,S_2)\in G_j(B_\F)$.
From the assumption, the intersection $CH(S_1)\cap CH(S_2)$ must be included in $G_j(\ol \F)$. Since $(S_1,S_2)\in \partial C_\F (G_j)$, for every $\varepsilon>0$ there exists a polygon $P$ such that the Hausdorff distance between $P$ and $CH(S_1)\cap CH(S_2)$ is smaller than $\varepsilon$, and $P$ is not precisely covered by $G_j(\F)$, which implies that $P$ is not included in $G_j(\F)$, or $P$ does not intersect $g(\F)$ for some $g\in G_j$.

If for every $\varepsilon>0$ the $\varepsilon$-neighborhood of $CH(S_1)\cap CH(S_2)$ is not included in $\ol{G_j( \F)}$, then a vertex of $CH(S_1)\cap CH(S_2)$ belongs to $\partial G_j(\F)$ or an edge of $CH(S_1)\cap CH(S_2)$ overlaps an edge of $\partial G_j(\F )$, which implies that for some $g\in G_j$, $g^{-1}(S_1,S_2)$ satisfies the condition $(B_\F 1)$ or $(B_\F 2)$ and belongs to $B_\F$.

If there exists $\varepsilon>0$ such that the $\varepsilon$-neighborhood of $CH(S_1)\cap CH(S_2)$ is included in $\ol{G_j( \F )}$, then there exists $g_0\in G_j$ such that $CH(S_1)\cap CH(S_2)$ does not contain an interior point of $g_0(\F)$. Since $CH(S_1)\cap CH(S_2)$ intersects $g_0(\ol \F)$ and both $CH(S_1)\cap CH(S_2)$ and $\F$ are polygons, $g_0^{-1}(S_1,S_2)$ satisfies at least one of the conditions to belong to $B_\F$. In this case we need the condition $(B_\F 3)$. Therefore in any cases $(S_1,S_2)\in G_j(B_\F)$.
\end{proof}

Our immediate goal is to prove Lemma \ref{lem:measure zero fundamental domain}, which says that for any $\mu,\nu \in \SC (\gS)$ there exists a Dirichlet domain $\F$ such that
\[\mu \times \nu (B_\F )=0.\]
By taking a path $c :[0,1]\rightarrow \tilde{\gS}$ starting from $x$ we can obtain a family of Dirichlet domains $\{ \F_{c (t)} \}_{t\in [0,1]}$.
We investigate how $\partial \F_x$ changes when $x$ moves along $c$. Recall that each edge of the Dirichlet domain $\F_x$ is a sub-arc of the perpendicular bisector of the geodesic joining $x$ to $g(x)$, denoted by $[x, g(x)]$, for $g\in G$.
We say that such perpendicular bisectors \ti{surround} $\F_x$.
Since $G$ acts on $\tilde{\gS}$ cocompactly and properly discontinuously, there are only finitely many perpendicular bisectors surrounding $\F_y$ for any $y\in \tilde{\gS}$. 
Fix $g\in G$ and consider how the perpendicular bisector of $[x,g(x)]$ moves when $x$ moves along $c$.
From now on, we consider the Poincar\'e disk model of $\HH$ and we will use the Euclidean geometry for considering geodesics of $\HH$.

\begin{lemma}\label{lem:perpendicular bisector}
Let $\ell $ be a geodesic line of $\HH$. Take $y_1,y_2\in \HH$ such that $y_1,y_2$ belong the same connected component of $\HH \setminus \ell$. Let $y_i'$ be the foot of the perpendicular line from $y_i$ to $\ell$ for $i=1,2$. If $d_\HH (y_1,y_1')=d_\HH (y_2,y_2')$ and $b$ is the perpendicular bisector of $[y_1,y_2]$, then $b$ is also the perpendicular bisector of $[y_1',y_2']\subset \ell$.
\end{lemma}
\begin{proof}
Take an isometry $\phi$ such that $\phi$ maps the midpoint between $y_1'$ and $y_2'$ to $0\in \HH$. Now, from the Euclidean geometry it is easy to see that the perpendicular bisector of $[\phi (y_1), \phi (y_2)]$ is also the perpendicular bisector of $[\phi (y_1'), \phi (y_2')]\subset \phi (\ell )$. Since $\phi $ is an isometry of $\HH$, this completes the proof.
\end{proof}

Fix non-trivial $g\in G$. For $y\in \tilde{\gS}$ we define $b_g(y)$ to be the perpendicular bisector of $[y,g(y)]$.
Let $x_0,y_0$ be the feet of the perpendicular lines from $x,y\in \tilde{\gS}$ to the axis $\mathrm{Ax}(g)$ of $g$, respectively.
For any $z\in \HH$ the hyperbolic distance from $z$ to $\mathrm{Ax}(g)$ coincides with that from $g(z)$ to $\mathrm{Ax}(g)$. Hence, we have $b_g(x)=b_g(x_0)$ and $b_g(y)=b_g(y_0)$ from the above lemma.
Therefore, the bisector $b_g(x)$ coincides with $b_g(y)$ if and only if $x_0=y_0$. Moreover, if $b_g(x)$ does not coincides with $b_g(y)$, then $b_g(x)$ does not intersect $b_g(y)$.

Recall that the translation length of $g$ is the hyperbolic distance between a point $z\in \mathrm{Ax}(g)$ and $g(z)$.
Take an isometry $\phi$ of $\Isom (\HH)$ such that $\phi$ fixes the axis of $g$ and $\phi ^2=g$. Then the translation length of $\phi$ is a half of that of $g$ and $b_g(y)$ equals $\phi (\ell _y)$ for the perpendicular line $\ell_y$ from $y$ to $Ax(g)$.

Now, we consider how the vertices of $\F_x$ moves when $x$ moves along $c$. Since a vertex of $\F_x$ is the intersection of two bisectors $b_{g_1}(x)$ and $b_{g_2}(x)$ for some $g_1,g_2\in G$, we have a map $\Phi_{g_1,g_2}$ from an open neighborhood of $x$ to a neighborhood of $b_{g_1}(x)\cap b_{g_2}(x)$.
Note that if $b_{g_1}(x)$ and $b_{g_2}(x)$ intersects at a point, then there exists an open neighborhood $U$ of $x$ such that $b_{g_1}(y)$ and $b_{g_2}(y)$ also intersects at a point for any $y\in U$.
From the above construction of $b_{g_i}(y)$ for $y\in \tilde{\gS}$, we can see that $\Phi_{g_1,g_2}$ is a $C^\infty$-map on $U$.
Therefore we have the following lemma:

\begin{lemma}\label{lem:construct a non-geodesic path}
Let $g_1,g_2$ be non-trivial elements of $G$. Assume that $b_{g_1}(x)$ and $b_{g_2}(x)$ intersects at a point for $x\in \tilde{\gS}$. Then there exists an open neighborhood $U$ of $x$ and an injective $C^\infty$-map $\Phi_{g_1,g_2}$ from $U$ to $\tilde{\gS}$ which maps $y\in U$ to the intersection point of $b_{g_1}(y)$ and $b_{g_2}(y)$.
Since $\Phi_{g_1,g_2}$ is injective, a subset of $U$ consisting of points $y$ satisfying the condition that the Jacobian of $\Phi_{g_1,g_2}$ at $y$ equals $0$ is a closed subset of $U$ without interior points.
\end{lemma}
\begin{proof}
We check only the injectivity of $\Phi_{g_1,g_2}$.
For any $y\in U$, the perpendicular line from $y$ to $\mathrm{Ax}(g_1)$ and that to $\mathrm{Ax}(g_2)$ intersects at $y$ and $b_{g_1}(y)$ and $b_{g_2}(y)$ intersects at a point.
Assume that $\Phi_{g_1,g_2}(y)=\Phi_{g_1,g_2}(z)$ for $y,z\in U$. Then $b_{g_1}(y)=b_{g_1}(z)$ and $b_{g_2}(y)=b_{g_2}(z)$.
Therefore the foot of the perpendicular line from $y$ to $\mathrm{Ax}(g_i)$ coincides with that from $z$ for $i=1,2$, which implies that $y=z$.
\end{proof}

Note that for any $x\in \tilde{\gS}$ and any non-trivial $g_1,g_2\in G$ with $g_1\not=g_2$, $b_{g_1}(x)$ never coincide with $b_{g_2}(x)$ since $g_1(x)\not=g_2(x)$. 

The following measure-theoretic lemma will plays an essential role in proving Lemma \ref{lem:measure zero fundamental domain}.

\begin{lemma}\label{lem:uncountably family has measure zero set}
Let $(X, \mu)$ be a measurable space, where $\mu$ is a measure on $X$.
Let $\{ A_\lambda \}_{\Lambda \in \gL}$ be an uncountable family of measurable subsets of $X$.
Let $B$ be a measurable subset of $X$ such that $B$ includes $\bigcup_{\lambda \in \gL}A_\lambda$.
Assume that $\mu (B)<\infty$ and there exists $M>0$ such that for any $x\in X$ we have
\[ \# \{ \lambda \in \Lambda \mid A_\lambda \ni x\} \leq M.\]
Such a family $\{ A_\lambda \}$ is said to be $M$-essentially disjoint. Then the subset
\[ \gL_{>0}:=\{ \lambda \in \Lambda \mid \mu (A_\lambda )>0 \} \]
is countable.
\end{lemma}
\begin{proof}
To obtain a contradiction, suppose that $\Lambda _{>0}$ is uncountable.
For each $n\in \NN$ consider a subset
\[ \gL_n:=\{ \lambda \in \gL \mid \frac{1}{n}\leq \mu (A_\lambda )< \frac{1}{n-1} \} ,\]
where if $n=1$, then $1/(n-1)$ means $\infty$. Since $\mu (A_\lambda )\leq \mu (B)<\infty$ for any $\lambda \in \gL$, we have
\[ \gL_{>0}=\bigsqcup_{n\in \NN }\gL_n.\]
Then we can see that there exists $n_0\in \NN$ such that $\gL_{n_0}$ is uncountable.
Since $\{ A_\lambda \}$ is $M$-essentially disjoint, for any finitely many $\lambda _1,\dots ,\lambda_k\in \gL_{n_0}$ we have
\[ \mu (\bigcup_ {i=1}^k A_{\lambda _k})\geq \frac{1}{M}\sum_{i=1}^k\mu (A_{\lambda_i})\geq \frac{1}{M}\cdot k\cdot \frac{1}{n_0}.\]
Therefore for a countably infinite subset $L\subset \gL_{n_0}$ we have
\[ \mu (\bigcup_{\lambda \in L }A_\lambda )\geq \frac{k}{Mn_0}\]
for any $k\in \NN$. Hence
\[ \mu (\bigcup_{\lambda \in L }A_\lambda )=\infty, \]
which contradicts our assumption that $\mu (B)<\infty $.
\end{proof}


\begin{lemma}\label{lem:measure zero fundamental domain}
There exists a smooth curve $c\:[0,1]\rightarrow \tilde{\gS}$ such that for any $\mu,\nu \in \SC (\gS)$, the set
\[ \{ t\in [0,1]\mid \mu \times \nu (B_{\F_{c(t)}} )>0 \} \]
is countable. Therefore for almost all $t\in [0,1]$ we have $\mu \times \nu (B_{\F_{c(t)}} )=0$.
\end{lemma}
\begin{proof}
Take a relatively compact open subset $U$ of $\tilde{\gS}$ and a compact subset $K$ of $\tilde{\gS}$ such that $K$ includes the Dirichlet domain $\F_{y}$ for any $y\in U$. Since $G$ acts on $\tilde{\gS}$ properly discontinuously, there exists a constant $M_1>0$ such that
\[ \# \{ g\in G\mid b_g(y)\cap K \not =\emptyset \ \text{for some }y\in U \} <M_1.\]
Note that if $b_g(y)\cap K \not =\emptyset$, then the hyperbolic distance from $y$ to $g(y)$ is smaller than or equal to twice the diameter of $K$.
Take all $g_1,\dots ,g_m\in G\setminus \{ \id \}$ such that $b_{g_i}(y)\cap K\not=\emptyset $ for some $y\in U$.
Then $m< M_1$, which implies that the number of edges of the Dirichlet domain $\F_y$ for any $y\in U$ is less than $M_1$.

From Lemma \ref{lem:perpendicular bisector} and \ref{lem:construct a non-geodesic path}, we can take a smooth curve $c\:[0,1]\rightarrow U$ satisfying the following condition:
\begin{enumerate}
\item[($\ast$)] for any $t_1,t_2\in [0,1]$ with $t_1\not=t_2$ the foot of the perpendicular line from $c(t_1)$ to $\mathrm{Ax}(g_i)$ is different from that from $c(t_2)$ for any $i=1,\dots ,m$.
\end{enumerate}
Then for any $t_1,t_2\in [0,1]$ with $t_1\not= t_2$ and $g_i$, the bisector $b_{g_i}(c(t_1))$ and $b_{g_i}(c(t_2))$ are disjoint.
We will modify $c$ later. 

In order to apply Lemma \ref{lem:uncountably family has measure zero set} to the family $\{ B_{\F_{c(t)}}\} _{t\in [0,1]}$, we prove that for any $(S_1,S_2)\in \mathcal{I}$, the cardinality of $\{ t\in [0,1]\mid (S_1,S_2)\in B_{\F_{c(t)}} \}$ is uniformly bounded.
Since $K$ is compact, there exists $M_2>0$ such that the number of boundary components of $CH(S)$ intersecting $K$ is less than $M_2$ for any $S\in \H (\partial G)$, which implies that for $(S_1,S_2)\in \mathcal{I}$ the number of edges of the polygon $CH(S_1)\cap CH(S_2)$ intersecting $K$ is less than $2M_2$.

For $(S_1,S_2)\in \mathcal{I}$ and each vertex $v$ of $CH(S_1)\cap CH(S_2)$, $v$ belongs to $b_{g_i}(c(t))$ at most once for $t\in [0,1]$ for each $g_i$, that is, the number of $t\in [0,1]$ such that $v\in \partial \F_{c(t)}$ is less than $M_1$. This corresponds to the condition $(B_\F1)$.
By the same way we can see that for each edge $e$ of $CH(S_1)\cap CH(S_2)$ the number of $t\in [0,1]$ such that $e$ overlaps an edge of $\F_{c(t)}$ is less than $M_1$. This corresponds to the condition $(B_\F 2)$.

Now, we want to see that for each edge $e$ of $CH(S_1)\cap CH(S_2)$ the number of $t\in [0,1]$ such that $e$ is tangent to a vertex of $\F_{c(t)}$ is uniformly bounded.
For any pair of $g_i,g_j$ such that $b_{g_i}(c(0))$ and $b_{g_j}(c(0))$ intersect at a point belonging to $K$, we can assume that $U$ is sufficiently small and we can define the map $\Phi_{g_i,g_j}$ on $U$.
We can also assume that if $b_{g_i}(c(0))$ and $b_{g_j}(c(0))$ intersect at a point belonging to the complement $K^c$, then $b_{g_i}(x)$ and $b_{g_j}(x)$ do not intersect at a point belonging to $K$ for any $x\in U$.
If $\Phi_{g_i,g_j}\circ c$ is a geodesic and $\Phi_{g_i,g_j}(c(t))$ is a vertex of $\F_{c(t)}$ for every $t\in [0,1]$, then an edge $e$ of $CH(S_1)\cap CH(S_2)$ can be tangent to $\Phi_{g_i,g_j}(c(t))$ for every $t\in [0,1]$. This is an undesirable case.

We modify $c$ such that $c$ satisfies the above condition $( \ast )$ and the condition that any geodesic meets the image of $\Phi_{g_i,g_j}\circ c$ at most 2 times for any pair of $g_i,g_j$.
From Lemma \ref{lem:construct a non-geodesic path} we can assume that the Jacobian of $\Phi_{g_i,g_j}$ at $y$ is not $0$ for every $y\in U$ and every pair of $g_i,g_j$.

We use the Euclidean geometry on the Poincar\'e disk model of $\HH$ in oder to modify $c$.
Since $K$ is bounded in $\HH$, there exists a constant $R_0>0$ such that any geodesic in $\HH$ passing through $K$ is a sub-arc of a line or a circle with radius larger than $R_0$ in the Euclidean plane containing $\HH$, whose absolute value of curvature is less than $1/R_0$.
If the absolute value of the curvature of a smooth curve $\gamma $ is larger than $1/R_0$ and the length of $\gamma $ is small enough, then $\gamma $ is approximated by a sub-arc of a circle with radius smaller than $R_0$ and any line or a circle with radius larger than $R_0$ in the Euclidean plane meets $\gamma $ at most twice.
Note that if the absolute value of the curvature of $\gamma$ is larger than $1/R_0$ and smaller than some constant $L>0$, then we take $\gamma$ such that the length of $\gamma$ is smaller than $\pi/L$, which is the length of a half-circle with radius $1/L$.
Now, we want to prove the following claim:

\underline{Claim:} \emph{We can modify $c$ so that $c$ satisfies the condition $(\ast )$, and the absolute of the curvature of $\Phi_{g_i,g_j}\circ c$ is larger than $1/R_0$ for any pair of $g_i,g_j$.}

We suppose for a moment that the above claim is true and prove the statement of the lemma. First, we can see that for each edge $e$ of $CH(S_1)\cap CH(S_2)$ the number of $t\in [0,1]$ such that $e$ is tangent to a vertex of $\F_{c(t)}$ is less than $2M_1$ since the number of vertices of $\F_{c(t)}$ is less than $2M_1$.
Recall that the number of edges of $CH(S_1)\cap CH(S_2)$ intersecting $K$ is at most $2M_2$.
Therefore for each $(S_1,S_2)\in \mathcal{I}$ the number of $t\in [0,1]$ such that $B_{\F_{c(t)}}$ containing $(S_1,S_2)$ is at most $2M_2(M_1+M_1+2M_1)$. Note that the union of $B_{\F_{c(t)}}$ over $t\in[0,1]$ is included in
\[ \{ (S_1,S_2)\in \H (\partial G)\times \H (\partial G)\mid CH(S_1)\cap CH(S_2)\cap K\not=\emptyset \} ,\]
which is compact.
Hence by applying Lemma \ref{lem:uncountably family has measure zero set} to $\mu\times \nu$ and the family $\{ B_{\F_{c(t)}}\} _{t\in [0,1]}$, the set
\[ \{ t\in [0,1]\mid \mu \times \nu (B_{\F_{c(t)}} )>0 \} \]
is countable.

Now we prove Claim in the above. Set $c(t)=(u(t),v(t))$ for $t\in [0,1]$ and set $\Phi(u,v):= \Phi_{g_i,g_j}(u,v)=(\alpha (u,v), \beta (u,v))$ for $(u,v)\in U$. Let $c'$ denote the derivative of $c$. We denote by $\alpha _u$ the partial derivative of $\alpha $ with respect to $u$ at $c(t)$ for some $t\in [0,1]$.
Recall that the curvature $\kappa_c $ of $c$ is
\[ \kappa _c=\frac{u'v''-v'u''}{(u'^2+v'^2)^{\frac{3}{2}}}.\]
We have
\[ (\alpha \circ c)'=\alpha_u u'+\alpha_v v',\]
\[ (\alpha \circ c)''=\alpha_{uu} u'^2+2\alpha_{uv}u'v'+\alpha_{vv}v'^2+\alpha_u u''+\alpha_v v'',\]
and
\[(\alpha \circ c)'(\beta \circ c)''-(\beta \circ c)'(\alpha \circ c)''= \phi +\psi ,\]
where 
\[ \phi =(\alpha_u \beta_v -\beta_u \alpha_v ) (u'v''-v'u'')\]
and
\begin{align*}
 \psi =& (\alpha_u u'+\alpha_v v')(\beta_{uu} u'^2+2\beta_{uv}u'v'+\beta_{vv}v'^2)\\
&-(\beta_u u'+\beta_v v')(\alpha_{uu} u'^2+2\alpha_{uv}u'v'+\alpha_{vv}v'^2).
\end{align*}
Then
\[ \kappa_{\Phi \circ c}=\frac{ \phi +\psi}{ ((\alpha \circ c)'^2+(\beta \circ c)'^2 )^{\frac{3}{2}}}.\]
Since $\Phi_{g_i,g_j}$ is given for any pair $g_i,g_j$, we can regard the partial derivatives $\alpha_u ,\beta_u ,\dots ,\beta_{vv}$ appeared in $\kappa_{\Phi \circ c}$ as almost constant.
Note that the Jacobian of $\Phi$, which equals $\alpha_u\beta_v-\beta_u \alpha_v$, is not $0$. We modify the second derivatives $u'' ,v''$ so that $(u'v''-v'u'')>0$ is large. Then $\kappa_{\Phi \circ c}(t)$ is larger than $1/R_0$. Note that $u'$ and $v'$ do not have to change so much if we restrict $c$ to a short interval $[0,r]$ for some small $r>0$.

For example, consider a function $f(t)=(t+1)^{a}-at-1$ around $0$ for a large $a\in \NN$. Then we have $f'(t)=a(t+1)^{a-1}-a$, $f''(t)=a(a-1)(t+1)^{a-2}$.
Consider the case that $u'>0$. Set $\hat{c}(t):=(u(t),v(t)+f(t))$ for $t\in [0,r]$ for a sufficiently small $r>0$. Then $\hat c$ is close to $c$, $\hat{c}'$ is close to $c'$, and $(v(t)+f(t))''$ is sufficiently large for $t\in [0,r]$. Since $u',v',u'',v''$ is bounded in $U$, $u'(v''+f'')-(v'+f')u''$ is sufficiently large, which implies that the absolute value of the curvature of $\Phi \circ \hat{c}$ is sufficiently large.
Note that if $\hat c$ is close to $c$ and $\hat{c}'$ is close to $c'$ on $[0,r]$, then $\hat{c}$ also satisfies the condition $(\ast )$.
This completes the proof.
\end{proof}

\begin{remark}\label{rem:boundary of A K has zero measure}
For a subset $K$ of $\tilde{\gS}$ set
\[A(K):=\{ S\in \H (\partial G)\mid CH(S)\cap K\not=\emptyset \} .\]
If $K$ is open or compact, then so is $A(K)$ respectively from Lemma \ref{lem:taking convex hull is continuous}.
By using the curve $c$ in Lemma \ref{lem:measure zero fundamental domain}, we can see that for any $\mu\in \SC (\gS)$, a set
\[ \{ t\in [0,1]\mid \mu (\partial A (\F_{c(t)}))>0\}\]
is countable. In fact, the boundary of $A(\F_{c(t)})$,
\begin{align*}
\partial A(\F_{c(t)})=\{ &S\in \H (\partial G)\ |\\
&CH(S)\cap \mathrm{Int}(\F_{c(t)})=\emptyset \text{ and } CH(S)\cap \partial \F_{c(t)}\not=\emptyset \}.
\end{align*}
Hence for $S\in \H (\partial G)$, if $S\in \partial A(\F_{c(t)})$, then there exists a boundary component $B$ of $CH(S)$ such that $B$ is tangent to a vertex of $\F_{c(t)}$ or overlaps an edge of $\F_{c(t)}$.
\end{remark}

\begin{proof}[Proof of Theorem \ref{thm:intersection number of subset currents}]
Take $(\mu_n,\nu_n ),(\mu ,\nu)\in \SC (\gS)\times \SC (\gS)\ (n\in \NN)$ such that $(\mu_n,\nu_n)$ converges to $(\mu,\nu)$.
Then the product measure $\mu_n\times \nu_n$ converges to $\mu\times \nu$ in the weak-$\ast$ topology in $\H (\partial G)\times \H (\partial G)$ by general theory (see \cite[Theorem 2.8]{Bil99} for the case of probability measures).
From Lemma \ref{lem:measure zero fundamental domain} and Remark \ref{rem:boundary of A K has zero measure} there exists $x\in \tilde{\gS}$ such that
\[ \mu (\partial A(\F_x ))=\nu (\partial A(\F_x )) =\mu \times \nu (B_{\F_x} )=0.\]
Then $\mu_n (A(\F)), \nu_n (A(\F))$ converges to $\mu (A (\F)),\nu (A(\F))$ respectively by Proposition \ref{prop:Portmanteau}. Set
\[ M:=\sup \{ \mu_n (A(\F )),\nu_n (A(\F ))\mid n\in \NN \} .\]
Remove some vertices and edges from $\F=\F_x$ such that $G (\F )=\tilde{\gS}$ and $g\F \cap \F =\emptyset$ for any non-trivial $g\in G$.
We are going to prove the following claim.

\underline{Claim:} $\mu_n\times \nu_n (C_\F (G_j))$ converges to $\mu \times \nu (C_\F(G_j))$ for any $j\in J$.

Suppose for a moment that Claim is true and prove that $i_{\SC}(\mu_n ,\nu_n )$ converges to $i_{\SC}(\mu ,\nu )$.
Recall that $G_j$ contains $\id $ for every $j\in J$. Hence
\[ \bigsqcup_{j\in J}C_\F (G_j)\subset A(\F )\times A(\F ),\]
which implies that
\[ \sum_{j\in J}\mu_n \times \nu_n (C_\F (G_j))\leq M^2\]
for any $n\in \NN$. Therefore
\begin{align*}
\lim_{n\rightarrow \infty} i_{\SC}(\mu_n ,\nu_n )&=\lim_{n\rightarrow \infty }\sum_{j\in J}\mu_n\times \nu_n (C_\F (G_j)) \\
&=\sum_{j\in J}\mu \times \nu (C_\F (G_j))=i_{\SC}(\mu ,\nu ),
\end{align*}
which proves the theorem.

Now, we prove Claim in the above.
Fix $j\in J$ and $\varepsilon>0$.
From Proposition \ref{prop:Portmanteau}, we consider the boundary $\partial C_\F (G_j)$.
Recall that $\partial C_\F(G_j)\subset G_j(B_\F \sqcup \gD_\F)$ for $j\in J$ and we have
\[ G_j(\gD_\F ) =\{ (S,S)\in \H (\partial G)\times \H (\partial G) \mid \# S=2, CH(S)\cap G_j(\ol{\F} )\not =\emptyset \}, \]
which is included in the compact set $A(G_j(\ol{\F}))$.
Note that $\mu \times \nu ( G_j (B_\F )) =0$ since $\mu \times \nu (B_\F) =0$.
Hence, if $\mu \times \nu ( G_j (\gD_\F)) =0$, then immediately we can see that $\mu \times \nu (\partial C_\F (G_j ))=0$, which implies that $\mu_n\times \nu_n (C_\F (G_j))$ converges to $\mu \times \nu (C_\F(G_j))$ by Proposition \ref{prop:Portmanteau}.
From now on, we assume that $\mu \times \nu (G_j (\gD_\F ))>0$.

From the Fubini's theorem we have
\begin{align*}
\mu \times \nu (G_j(\gD_\F) ) 
&=\int \chi_{G_j(\gD_\F)}(S_1,S_2) d\mu \times \nu \\
&=\int \left( \int \chi_{G_j(\gD_\F)}(S_1,S_2) d\mu (S_1)\right) d\nu (S_2)\\
&=\int_{\partial_2 G\cap A(G_j(\ol \F ))} \mu (\{ S_2\} ) d\nu (S_2)\\
&=\sum_{S\in \partial_2 G\cap A (G_j(\ol \F ))\text{: common atom of }\mu,\nu }\mu (\{ S\} ) \nu (\{ S\} ),
\end{align*}
where $\chi_{G_j(\gD_\F)}$ is the characteristic function of $G_j(\gD_\F)$ on $\H (\partial G)\times \H(\partial G)$.
Recall that a point $S\in \H (\partial G)$ is called an atom of $\mu$ if $\mu (\{ S\})>0$. Since $\mu, \nu$ are locally finite measures, they have at most countably many atoms.
Therefore there exist finite common atoms $S_1,\dots ,S_m\in \partial_2 G\cap A (G_j(\ol \F ))$ of $\mu ,\nu$ such that
\begin{equation}
\mu \times \nu (G_j(\gD_\F ))<\sum_{k=1}^{m}\mu \times \nu (\{ (S_k,S_k)\} )+\varepsilon. \tag{$\ast$}
\end{equation}
We will construct an open neighborhood $V$ of $\{ (S_1,S_1),\dots ,(S_m,S_m)\}$ such that $\mu_n \times \nu_n (V\cap C_\F (G_j))<\varepsilon$ for any $n\in \NN$.

Since $\#S_k=2$ and $S_k$ is an atom of a subset current, there exists $g_k\in G$ such that
$S_k =\gL (\langle g_k \rangle )$ from Lemma \ref{thm:counting sc}.
Hence $g_k (S_k)=S_k$.
Since $\mu (\partial A(\F ))=\nu (\partial A(\F) )=0$, we have $\mu (\partial A(g\F ))=\nu (\partial A(g\F) )=0$ for any $g\in G_j$, which implies that $CH(S_k)$ passes through the interior of $g\F$ for any $g\in G_j$.
Hence $S_k\in \mathrm{Int}(A(G_j(\F )))$.
Then we can take an open neighborhood $O_k\subset \mathrm{Int}(A(G_j(\F)))$ of $S_k$.
Take an arbitrary $L\in \NN$ and set
\[ U_k:=\bigcap_{l=1}^{L} (g_k)^{-l}(O_k).\]
Then $U_k$ is also an open neighborhood of $S_k$ and
\[ g_k (U_k), \dots , (g_k)^L (U_k)\subset O_k \subset \mathrm{Int}(A(G_j(\F ) )).\] 

Now, we consider the intersection of $U_k\times U_k$ and $C_\F (G_j)$. Note that $gC_\F(G_j)\cap C_\F (G_j)=\emptyset $ for any non-trivial $g\in G$. Therefore 
\[ g_k\left( U_k\times U_k\cap C_\F (G_j) \right) ,\dots ,(g_k)^L \left( U_k\times U_k \cap C_\F (G_j)\right) \] 
are pairwise disjoint, and for any $n\in \NN$ we have
\begin{align*}
&\mu _n\times \nu _n(U_k\times U_k\cap C_\F (G_j) )\\
=&\frac{1}{L}\mu_n \times \nu_n \left( \bigsqcup_{l=1}^L (g_k)^l (U_k\times U_k\cap C_\F (G_j) )\right) \\
\leq &\frac{1}{L}\mu_n \times \nu_n (A (G_j(\F ) )\times A(G_j (\F ) ))\\
\leq &\frac{1}{L}\sum_{g_1,g_2\in G_j} \mu_n \times \nu_n ((g_1A (\F ) )\times (g_2A(\F )))\\
\leq &\frac{(\# G_j M)^2}{L}.
\end{align*}
Set $V:= (U_1\times U_1) \cup \cdots \cup (U_m \times U_m)$. Then we have
\[ \mu_n \times \nu_n (V\cap C_\F (G_j ))\leq \sum_{k=1}^m \frac{(\# G_jM)^2}{L}\leq \frac{m(\# G_j M)^2}{L}.\]
Hence if $L$ is sufficiently large, then we have
\[ \mu_n \times \nu_n (V\cap C_\F (G_j))<\varepsilon.\]
Note that $V$ contains all of $(S_1,S_1),\dots ,(S_m,S_m)$.

Since $C_\F(G_j)\cap G_j (\gD_\F )=\emptyset $, we can see that
\[ \mathrm{Int}(C_\F (G_j))=C_\F (G_j)\setminus G_j(B_\F). \]
Then from Proposition \ref{prop:Portmanteau} and Equation ($\ast$), we have
\begin{align*}
&\mu\times \nu (C_\F (G_j))=\mu \times \nu (\mathrm{Int}(C_\F (G_j)))\\
\leq &\liminf_{n\rightarrow \infty} \mu_n\times \nu_n (\mathrm{Int}(C_\F (G_j)))\\
\leq &\liminf_{n\rightarrow \infty} \mu_n\times \nu_n (C_\F (G_j))\\
\leq &\limsup_{n\rightarrow \infty} \mu_n\times \nu_n (C_\F (G_j))\\
\leq &\limsup_{n\rightarrow \infty} \mu_n\times \nu_n (C_\F (G_j)\setminus V)\\
		&+\limsup_{n\rightarrow \infty} \mu_n\times \nu_n (C_\F (G_j)\cap V)\\
\leq &\limsup_{n\rightarrow \infty} \mu_n\times \nu_n (\ol{C_\F (G_j)}\setminus V)+\varepsilon \\
\leq &\mu \times \nu (\ol{C_\F (G_j)}\setminus V)+\varepsilon \\
\leq &\mu \times \nu (C_\F (G_j))+\mu \times \nu (G_j(\gD_\F )\setminus V)+\varepsilon \\
< &\mu \times \nu (C_\F (G_j))+2\varepsilon .
\end{align*}
Since $\varepsilon >0$ is arbitrary,
\[ \lim_{n\rightarrow \infty }\mu_n \times \nu_n (C_\F (G_j))=\mu \times \nu (C_\F (G_j)).\]
This completes the proof.
\end{proof}

\if{0}
\begin{proposition}
Let $G$ be a free group of finite rank or a surface group.
For the reduced rank functional $\rk$ on $\SC (\gS)$ we have
\[ \rk^{-1}(0)=\GC (\gS).\]
\end{proposition}
\fi

\section{Intersection functional on subset currents}\label{sec:an intersection functional}

Let $\gS $ be a compact hyperbolic surface (possibly with geodesic boundary) and $G$ the fundamental group of $\gS$. The notation in this section is based on that in Section \ref{sec:intersection number of subset currents} and we consider the universal cover $\tilde{\gS}$ of $\gS$ as a subspace of the hyperbolic plane $\HH$. We identify $\partial G$ with the limit set $\gL (G) \subset \partial \HH $.

Recall that for two non-trivial finitely generated subgroups $H$ and $K$ of $G$ we have considered the fiber product $C_H\times_\gS C_K$ corresponding to the convex cores $(C_H, p_H)$ and $(C_K,p_K)$.
Now, instead of contractible components of $C_H\times_\gS C_K$ we focus on the \ti{non-contractible} components of $C_H\times_\gS C_K$.
Note that $C_H\times_\gS C_K$ can be considered as the quotient space $G\setminus \hat{CH}_H\times_{\tilde{\gS}} \hat{CH}_K$ and every non-contractible component of $C_H\times_\gS C_K$ is corresponding to
\[ (g_1Hg_1^{-1}\cap g_2Kg_2^{-1})\backslash (g_1CH_H\cap g_2CH_K) \]
for $[g_1H,g_2K]\in G\backslash (G/H\times G/K)$ with $g_1Hg_1^{-1}\cap g_2Kg_2^{-1}\not=\emptyset$. Note that if $g_1CH_H\cap g_2CH_K=\emptyset $ for $[g_1H,g_2K]\in G\backslash (G/H\times G/K)$, then $g_1Hg_1^{-1}\cap g_2Kg_2^{-1}$ is trivial.
Recall that the reduced rank of a non-trivial finitely generated subgroup $J$ of $G$ equals $-\chi(C_J)$ and that of the trivial subgroup $\{ \id \}$ equals $0$ (see the above of Corollary \ref{cor:reduced rank functional} for detail).

\begin{definition}[Product $\N$]
We define the product $\N$ of two finitely generated subgroups $H$ and $K$ of $G$ by
\[ \N (H,K):=\sum_{[g_1H,g_2K]\in G\backslash (G/H\times G/K)}\rk (g_1Hg_1^{-1}\cap g_2Kg_2^{-1}).\]
\end{definition}

See \cite{Sas15} for the background of the product $\N$ in the case that $G$ is a free group of finite rank.

\begin{remark}
Let $H,K$ be finitely generated subgroups of $G$.
From \cite[Theorem 4.1]{Sas15}, we have a bijective map from $G\backslash (G/H\times G/K)$ to the set of all double cosets $H\backslash G/K$, which maps $[g_1H,g_2K]$ to $Hg_1^{-1}g_2K$. Since $\rk$ is invariant up to conjugacy, we obtain
\[ \N (H,K)=\sum_{HgK\in H\backslash G/K}\rk (H\cap gKg^{-1}).\]
In the case that $G$ is a free group of finite rank, this expression of the product $\N$ is often used for stating the Strengthened Hanna Neumann Conjecture, which can be written as follows: for any finitely generated subgroups $H$ and $K$ of $G$ the inequality
\[ \N (H , K )\leq \rk (H ) \rk (K )\]
holds. This conjecture was individually proved by by Friedman \cite{Fri15} and Mineyev \cite{Min12}. As far as the author knows, the surface group version of the Strengthened Hanna Neumann Conjecture is still an open problem.

Next, we consider a geometrical expression of the product $\N$.
For each $[g_1H,g_2K]\in G\backslash (G/H\times G/K)$, if $g_1Hg_1^{-1}\cap g_2Kg_2^{-1}\not=\{\id \}$, then $g_1CH_H\cap g_2CH_K$ is non-empty and there exists a corresponding connected component of $C_H\times_\gS C_K$ whose fundamental group is isomorphic to $g_1Hg_1^{-1}\cap g_2Kg_2^{-1}$. We define the reduced rank $\rk(M)$ of a non-contractible compact surface or a circle $M$ to be $-\chi (M)$ and the reduced rank of a contractible space $M$ to be $0$.
Then we can see that
\[ \N (H,K)=\sum_{M\subset C_H \times_\gS C_K} \rk (M),\]
where the sum is taken over all connected components of $C_H\times_\gS C_K$.
Note that a non-contractible connected component of $C_H \times_\gS C_K$ can be a circle or a compact surface each of whose boundary components is piecewise geodesic.
\end{remark}

Our goal in this section is to prove the following theorem. In the case that $G$ is a free group of finite rank, this theorem is proved in \cite[Theorem 3.2]{Sas15}.

\begin{theorem}\label{thm:intersection functional}
Let $\gS$ be a compact hyperbolic surface.
There exists a unique symmetric continuous $\RRR$-bilinear functional 
\[ \N \: \SC (\gS)\times \SC (\gS)\rightarrow \RRR \]
such that for any non-trivial finitely generated subgroups $H$ and $K$ of $G$ we have
\[ \N (\eta_H,\eta_K )=\N (H,K).\]
\end{theorem}

In the case that $G$ is a free group $F$ of finite rank, from the above theorem, we can see that the inequality
\[ \N (\mu , \nu )\leq \rk (\mu )\rk (\nu )\]
holds for any $\mu ,\nu \in \SC (F)$, which is a direct corollary to the Strengthened Hanna Neumann Conjecture.

Note that for any finitely generated subgroup $H$ of $G$ we have $\N (G , H )=\rk (H)$ since $C_H\times_\gS \gS$ is naturally homeomorphic to $C_H$. Hence $\N (\eta_G ,\cdot )$ coincides with the reduced rank functional $\rk$ by the denseness property of rational subset currents for $G$.

The guidelines for proving Theorem \ref{thm:intersection functional} is almost the same as that in \cite{Sas15}.
The main objects considered in \cite{Sas15} are graphs and trees but our main objects here are surfaces and circles.
One of the keys for proving Theorem \ref{thm:intersection functional} is the Gauss-Bonnet Theorem.
Note that a boundary component of a 2-dimensional connected component of $C_H\times_\gS C_K$ is not totally geodesic in general but piecewise geodesic. Moreover, $C_H \times_\gS C_K$ contains a 1 or 0-dimensional object if either $H$ or $K$ is cyclic.

In order to apply the Gauss-Bonnet Theorem to $C_H\times_\gS C_K$, we introduce some notation.
For a corner $v$ of a piecewise geodesic, which is called a vertex, we define $\An(v)$ to be the exterior angle of $v$.
If a 1-dimensional connected component $M$ of $C_H\times_\gS C_K$ has a boundary, then $M$ is a geodesic segment.
For an end-point $v$ of the geodesic segment, which is also called a vertex, we define $\An(v)$ to be $\pi$.
If $M=\{ v \}$ is a connected component of $C_H\times_\gS C_K$, then we also call $v$ a vertex of $C_H\times_\gS C_K$ and define $\An(v)$ to be $2\pi$.
By applying the Gauss-Bonnet Theorem to each connected component $M$ of $C_H\times_\gS C_K$, we have the following formula
\[ 2\pi \chi (M)=-\Area (M)+\sum_{\text{$v$: vertex of $M$}}\An (v),\]
and so
\begin{equation}
2\pi \chi (C_H\times_\gS C_K )=-\Area (C_H\times_\gS C_K)+\sum_{\text{$v$: vertex of }C_H\times_\gS C_K}\An (v), \tag{GB} \label{eq:Gauss-Bonnet}
\end{equation}
where $\chi (C_H\times_\gS C_K)$ (or $\Area (C_H\times_\gS C_K)$) is the sum of the Euler characteristic (or the area, respectively) of each connected component of $C_H\times_\gS C_K$. If $M$ is a 1-dimensional or 0-dimensional connected component of $C_H\times_\gS C_K$, then the area of $M$ is 0.

Since the Euler characteristic of a contractible component is 1, we have the following equation:
\[ \N (H,K)=-\chi (C_H\times_\gS C_K)+i(C_H,C_K).\]
We will extend $\chi$ to a symmetric continuous $\RRR$-bilinear functional from $\SC (\gS)\times \SC (\gS)$ to $\RR$ by using Formula (\ref{eq:Gauss-Bonnet}). For that purpose we will extend both the ``area term'' and the ``angle term'' in Formula (\ref{eq:Gauss-Bonnet}) to symmetric continuous $\RRR$-bilinear functionals from $\SC (\gS)\times \SC (\gS)$ to $\RRR$.

First, we extend the ``area term'' by using the same method of Theorem \ref{thm:volume functional}.
Take a Dirichlet domain $\F$ for the action of $G$ on $\tilde{\gS}$. Recall that $m_\HH $ is the measure on $\HH$ induced by the Riemannian metric on $\HH$.
We define a function $f$ from $\H (\partial G)\times \H (\partial G)$ to $\RR$ by
\[ f(S_1,S_2):=m_\HH (CH(S_1)\cap CH(S_2)\cap \F )\]
for $(S_1,S_2)\in \H (\partial G)\times \H (\partial G)$.
\begin{proposition}\label{prop:area bilinear functional}
The function $f$ is a continuous function with compact support. The functional $f^\ast$ from $\SC (\gS)\times \SC (\gS)$ to $\RRR$ defined by
\[ f^\ast (\mu ,\nu ):= \int f d \mu\times \nu \quad (\mu ,\nu \in \SC (\gS))\]
is a symmetric continuous $\RRR$-bilinear functional satisfying the condition that
for any non-trivial finitely generated subgroups $H$ and $K$ of $G$ we have
\[ f^\ast (\eta_H ,\eta_K )=\Area (C_H\times_\gS C_K).\]
\end{proposition}
\begin{proof}
Recall that $A(\ol{\F})=\{ S\in \H (\partial G)\mid CH(S)\cap \ol{F}\not= \emptyset\}$ is a compact subset of $\H (\partial G)$.
For any $(S_1,S_2)\in \H (\partial G)\times \H (\partial G)$ if either $CH(S_1)\cap \F$ or $CH (S_2)\cap \F$ is empty, we have $f(S_1,S_2)=0$. This implies that the support of $f$ is included in $A(\ol{\F })\times A (\ol{\F})$. Hence $f$ has a compact support.
We can prove that $f$ is continuous by the same way as the proof of Proposition \ref{prop:f is continuous}.

Now, we check that $f^\ast(\eta_H, \eta_K)=\Area (C_H\times_\gS C_K)$ for any non-trivial finitely generated subgroups $H$ and $K$ of $G$.
First we have
\begin{align*}
f^\ast (\eta_H ,\eta_K )
&=\int f d\eta_H\times \eta_K \\
&=\sum_{(g_1H,g_2K)\in G/H\times G/K}f(g_1\gL(H), g_2\gL (K))\\
&=\sum_{(g_1H,g_2K)\in G/H\times G/K}m_{\HH}(g_1CH_H\cap g_2 CH_K\cap \F ).
\end{align*}
Set
\[ P:= \{ (g_1 H, g_2 K, x) \in G/H \times G/K \times \tilde{\gS} \mid x \in g_1 CH_H \cap g_2CH_K \cap \F \}.\]
We can extend the measure $m_{\HH}$ to the measure on $G/H \times G/K \times \tilde{\gS}$ naturally since $G/H \times G/K$ is a countable discrete space. Then we have $m_{\HH} (P )=f^\ast (\eta_H ,\eta_K )$.
From the proof of Lemma \ref{lem:volume of fundamental domain}, it is sufficient to see that $P$ is a measure-theoretic fundamental domain for the action of $G$ on $\hat{CH}_H\times_{\tilde{\gS}} \hat{CH}_K$, that is, $G(P)=\hat{CH}_H\times_{\tilde{\gS}} \hat{CH}_K$ and $gP\cap P$ has measure zero for any non-trivial $g\in G$.

For any $(g_1H ,g_2K ,x)\in \hat{CH}_H\times_{\tilde{\gS}} \hat{CH}_K$ there exists $g\in G$ such that $gx \in \F$.
Then $g(g_1H,g_2K,x)\in P$. Hence $G(P)=\hat{CH}_H\times_{\tilde{\gS}} \hat{CH}_K$.
For any $g\in G$ the projection of the intersection $gP\cap P$ onto $\tilde{\gS}$ equals $g\F \cap \F$. Hence $gP \cap P$ has measure zero.

From now on, we give another proof of the equality $f^\ast (\eta_H ,\eta_K )=\Area (C_H\times_\gS C_K)$ by considering each connected component of $C_H \times_{\gS} C_K$ in order to give its another description.
Recall that the fiber product $C_H\times_\gS C_K$ is the disjoint union of
\[ M_{g_1H,g_2K}:=(g_1Hg_1^{-1}\cap g_2Kg_2^{-1})\backslash (g_1CH_H\cap g_2 CH_K)\]
over all $[g_1H,g_2K]\in G\backslash (G/H\times G/K)$.
Fix $g_1,g_2\in G$ and set $J:=g_1Hg_1^{-1}\cap g_2Kg_2^{-1}$, which is the stabilizer of $g_1CH_H\cap g_2CH_K$ in $\CH_H\times_{\tilde{\gS}}\CH_K$. 
Note that $M_{g_1H,g_2K}$ can be empty and $J$ can be $\{ \id \}$.
The preimage of $M_{g_1H,g_2K}$ with respect to the quotient map $\Phi$ from $\CH_H\times_{\tilde \gS}\CH_K$ to $C_H\times_\gS C_K$ coincides with
\begin{align*}
&\{ (gg_1H,gg_2K,x)\in G/H\times G/K\times \HH \mid g\in G, x\in g(g_1CH_H\cap g_2CH_K)\} \\
\cong &\bigsqcup_{gJ\in G/J}g(g_1CH_H\cap g_2CH_K).
\end{align*}
Take a complete system of representatives $R$ of $G/J$. Now, we prove that
\[ A:=(g_1CH_H\cap g_2 CH_K)\cap \bigcup_{g\in R}g^{-1}\F \]
is a measure-theoretic fundamental domain for the action of $J$ on $g_1CH_H\cap g_2CH_K$,
Note that $R^{-1}$ is a complete system of representatives of $J\backslash G$, which implies that
\[ J\left(\bigcup_{g\in R} g^{-1}\F \right) =G(\F)\]
and $u_1g^{-1}\not= u_2g^{-1}$ for any $g\in R$ and $u_1,u_2\in J$ with $u_1\not=u_2$.
Hence $J(A)=g_1CH_H\cap g_2CH_K$ and $u(A)\cap A$ has measure zero for any non-trivial $u\in J$.
From the proof of Lemma \ref{lem:volume of fundamental domain}, we can see that $m_\HH(A)$ equals the area of $M_{g_1H,g_2K}$.

Now, we prove that $\Area (C_H\times_\gS C_K)=f^\ast (\eta_H ,\eta_K )$.
We have a natural bijective map from $G/J$ to the equivalence class $[g_1H, g_2K]$ which maps $gJ\in G/J$ to $(gg_1H,gg_2K)\in [g_1H,g_2K]$.
Since $m_\HH $ is a $G$-invariant measure, we have
\begin{align*}
m_\HH (A) &=\sum_{g\in R}m_\HH ((g_1CH_H\cap g_2 CH_K)\cap g^{-1}\F ) \\
&=\sum_{(g_1'H,g_2'K)\in [g_1H,g_2K]}m_\HH (g_1'CH_H\cap g_2' CH_K\cap \F ).
\end{align*}
Note that $G/H\times G/K$ is the disjoint union of $[g_1H,g_2K]\in G\backslash (G/H\times G/K)$.
Hence
\begin{align*}
&f^\ast (\eta_H ,\eta_K )\\
=&\sum_{(g_1H,g_2K)\in G/H\times G/K}m_{\HH}(g_1CH_H\cap g_2 CH_K\cap \F )\\
=&\sum_{[g_1H,g_2K]\in G\backslash (G/H\times G/K)}\sum_{(g_1'H,g_2'K)\in [g_1H,g_2K]}m_\HH (g_1'CH_H\cap g_2' CH_K\cap \F )\\
=&\sum_{[g_1H,g_2K]\in G\backslash (G/H\times G/K)}\Area (M_{g_1H,g_2K})\\
=&\Area (C_H\times_\gS C_K).
\end{align*}
This completes the proof.
\end{proof}

From now on, we are going to extend the ``angle term'' to a symmetric continuous $\RRR$-bilinear functional on $\SC (\gS)\times \SC (\gS)$ by using the method of proving the continuous extension of the intersection number.
Let $\F=\F_x$ be the Dirichlet domain centered at $x\in \tilde{\gS}$. We remove some edges and vertices of $\F$ such that $G(\F)=\tilde{\gS}$ and $g\F\cap \F=\emptyset$ for any non-trivial $g\in G$.
For $(S_1,S_2)\in \H (\partial G )\times \H (\partial G)$ with $CH(S_1)\cap CH(S_2)\not= \emptyset$, a vertex of $CH(S_1)\cap CH(S_2)$ is the intersection point of a boundary component of $CH(S_1)$ and that of $CH(S_2)$. We define the angle $\An (v)$ at $v$ to be the exterior angle at $v$.
Define a function $\phi$ from $\H(\partial G)\times \H (\partial G)$ to $\RR$ by
\[ \phi _\F(S_1,S_2):=\sum_{\text{$v$: vertex of }CH(S_1)\cap CH(S_2)\text{ in }\F}\An (v)\]
for $(S_1,S_2)\in \H(\partial G)\times \H (\partial G)$.
From the proof of Proposition \ref{prop:area bilinear functional}, we can see that for non-trivial finitely generated subgroups $H$ and $K$ of $G$ the restriction of the quotient map $\Phi$ to
\[ \{ (g_1H,g_2K,x)\in G/H\times G/K\times \tilde{\gS} \mid x\in g_1CH_H\cap g_2CH_K\cap \F \} \]
is a bijection onto $C_H\times_\gS C_K$. Therefore we obtain
\begin{align*}
&\int \phi_\F d\eta_H\times \eta_K \\
=&\sum_{(g_1H,g_2K)\in G/H\times G/K}\phi_\F (g_1\gL(H) , g_2\gL(K) )\\
=&\sum_{(g_1H,g_2K)\in G/H\times G/K}\sum_{\text{$v$: vertex of }g_1CH_H\cap g_2CH_K\text{ in }\F}\An (v)\\
=&\sum_{\text{$v$: vertex of }C_H\times_\gS C_K}\An (v).
\end{align*}

We define the symmetric $\RRR$-bilinear functional $\phi_\F^\ast$ from $\SC(\gS)\times \SC(\gS)$ to $\RRR$ by
\[ \phi_\F^\ast (\mu, \nu )=\int \phi_\F d\mu \times \nu \]
for $\mu ,\nu \in \SC (\gS)$.
Recall that $\mathrm{Span}(\SC_r(\gS))$ is the $\RRR$-linear subspace of $\SC(\gS)$ generated by the set $\SC_r(\gS)$ of rational subset currents on $\gS$, which is a dense subset of $\SC(\gS)$.
Because of the technical problem of the proof we do not prove the continuity of $\phi_\F^\ast$ but prove the continuity of the restriction of $\phi_\F^\ast $ to the $\mathrm{Span} (\SC_r(\gS))\times \mathrm{Span} (\SC_r(\gS))$,
which is uniquely extended to a continuous $\RRR$-bilinear functional on $\SC(\gS)\times \SC(\gS)$ by using the Hahn-Banach extension theorem.

The point is that $\phi_\F$ is not continuous and we need to understand the set $\gD(\phi_\F )$ of non-continuous points of $\phi_\F$ from Proposition \ref{prop:Portmanteau} and the proof of the continuity of $i_\SC$.
For any $S\in \H (\partial G)$, the number of boundary components of $CH(S)$ intersecting a bounded subset of $\HH$ is bounded by a constant independent of $S$. Hence it is sufficient to consider a finite number of boundary components of $CH(S)$ intersecting a neighborhood of $\F$ for $S\in \H (\partial G)$ when we see how the value of $\phi_\F$ changes.

Let $S\in \H (\partial G)$ and $B_1,\dots ,B_k$ the boundary components of $CH(S)$ intersecting a neighborhood of $\F$.
Assume that $\#S\geq 3$.
For a sufficiently small neighborhood $U$ of $S$ we can see that for any $S'\in U$ there exist boundary components $B_1',\dots ,B_k'$ of $CH(S')$ such that $B_1',\dots ,B_k'$ are the boundary components of $CH(S')$ intersecting the neighborhood of $\F$ and the Hausdorff distance between $B_i$ and $B_i'$, which is induced by the Euclidean metric, is small for every $i=1,\dots ,k$ from Lemma \ref{lem:taking convex hull is continuous}.
Moving the boundary component $B_1$ of $CH(S)$ in $U$ means taking a path from $S$ to a point $S''\in U$ such that for every point $S'$ in the path $B_i=B_i'$ for $i=2,\dots ,k$.
Moving the boundary component $B_1$ of $CH(S)$ ``a little'' means taking a (sufficiently) small open neighborhood $U$ of $S$ and moving $B_1$ of $CH(S)$ in $U$.

Let $(S_1,S_2)\in \H (\partial G)\times \H (\partial G)$. Assume that a boundary component $B_1$ of $CH(S_1)$ and a boundary component $B_2$ of $CH(S_2)$ intersect at a point $v$. If we move the boundary components $B_1$ and $B_2$ a little , then the intersection point and the exterior angle at the point change continuously.

Now, we define $C_{\F}$ to be a subset of $\H(\partial G)\times \H (\partial G)$ consisting of points $(S_1,S_2)$ satisfying the condition that there exists a vertex of $CH(S_1)\cap CH(S_2)$ belonging to $\partial \F$. We can see that $C_\F$ is included in $\gD(\phi_\F)$ by Lemma \ref{lem:taking convex hull is continuous}.
Moreover, $C_{\F}$ is a closed subset of $\H (\partial G)\times \H (\partial G)$ by Lemma \ref{lem:taking convex hull is continuous}.
Note that for $(S_1,S_2)\in C_\F$ the intersection $CH(S_1)\cap CH(S_2)$ is not necessarily bounded.

Next, we define $D_\F$ to be a subset of $\H (\partial G) \times \H (\partial G)$ consisting of points $(S_1,S_2)$ satisfying the condition that $CH(S_1)$ and $CH(S_2)$ share one boundary component intersecting $\ol{\F}$.
In other words, an element $(S_1,S_2)\in D_\F$ satisfies $S:=S_1\cap S_2 \in \partial_2 G$, $CH(S_1)\cap CH(S_2)=CH(S)$ and $CH(S)\cap \ol{\F}\not= \emptyset$.
Note that for $(S_1,S_2)\in D_\F$ we have $\phi_\F(S_1,S_2)=0$ and the cardinality of $S_i$ can be $2$.
Let $(S_1,S_2)\in D_\F$ and $S=S_1\cap S_2$. Assume that $\#S_1, \#S_2 \geq 3$ and $CH(S)$ passes through the interior $\mathrm{Int}(\F)$ of $\F$. Then we see that there exists $S'\in \partial_2 G$ close to $S$ such that $CH(S)$ and $CH(S')$ intersect at a point in $\mathrm{Int}(\F)$.
Hence there exists $S_1'\in \H (\partial G)$ close to $S_1$ such that $\phi_\F (S_1', S_2)$ is close to $\pi$, which implies that $\phi_\F$ is not continuous at $(S_1,S_2)$.
We see that $D_\F$ is a closed subset of $\H (\partial G)\times \H (\partial G)$ by Lemma \ref{lem:taking convex hull is continuous}.

For $S\in \partial_2 G$ with $CH(S)\cap \ol{\F}\not=\emptyset$, we see that $(S,S) \in D_\F$ and $\phi_\F$ is not continuous at $(S,S)$.
Remark that $D_\F$ includes the subset
\[ \gD_\F =\{ (S,S)\in \partial_2 G\times \partial _2 G\mid CH(S)\cap \ol{\F}\not= \emptyset \}, \]
which we used for proving Theorem \ref{thm:intersection number of subset currents}. During the proof of the continuity of $\phi_\F^\ast$, $D_\F$ will plays the same role as $\gD_\F$ in the proof of Theorem \ref{thm:intersection number of subset currents}.

\begin{lemma}\label{lem:non-continuous points of phi}
The set $\gD(\phi_\F )$ of non-continuous points of $\phi_\F$ is included in $C_{\F}\sqcup D_\F$.
\end{lemma}
\begin{proof}
Take any $(S_1,S_2)\in \H(\partial G)\times \H (\partial G)\setminus (C_{\F}\sqcup D_\F)$. It is sufficient to prove that $\phi_\F$ is continuous at $(S_1,S_2)$.
Since $C_{\F}\sqcup D_\F$ is a cosed subset of $\H (\partial G)\times \H (\partial G)$ we can take an open neighborhood $U$ of $(S_1,S_2)$ such that $U\cap (C_\F\sqcup D_\F )=\emptyset$. Since $(S_1,S_2)\not \in C_\F$, any vertex of $CH(S_1)\cap CH(S_2)$ is contained in the interior of $\F$ or the exterior of $\F$.
We divides the proof into several cases.
We assume that $U$ is sufficiently small in each case.

\underline{Case 1:} $\F$ does not contain any vertex of $CH(S_1)\cap CH(S_2)$.

If there exists no $S\in \partial _2 G$ such that $CH(S)$ is a common boundary component of $CH(S_1)$ and $CH(S_2)$, and $CH(S)\cap \ol{\F}\not=\emptyset$, then we can take a sufficiently small open neighborhood $U$ of $(S_1,S_2)$ such that $\F$ does not contain any vertex of $CH(S_1')\cap CH(S_2')$ for $(S_1',S_2')\in U$, which implies that $\phi_\F(S_1',S_2')=0$ and $\phi_\F$ is continuous at $(S_1,S_2)$.
Now, we assume that there exists $S\in \partial _2 G$ such that $CH(S)$ is a common boundary component of $CH(S_1)$ and $CH(S_2)$, and $CH(S)\cap \ol{\F}\not=\emptyset$. Since $(S_1,S_2)\not\in D_\F$, $(CH(S_1)\cap CH(S_2))\setminus CH(S)$ is not empty.
Hence even if $CH(S_1')\cap CH(S_2')$ has a vertex contained in $\F$ for $(S_1',S_2')\in U$, the exterior angle at the vertex is small.
Therefore, $\phi_\F$ is continuous at $(S_1,S_2)$.

From now on, we assume that $\F$ contains at least one vertex of $CH(S_1)\cap CH(S_2)$.

\underline{Case 2:} Both $S_1$ and $S_2$ belong to $\partial_2G$, that is, $CH(S_1)\cap CH(S_2)$ is a point contained in $\F$.

Since $(S_1,S_2)\not\in C_\F$, $CH(S_1)\cap CH(S_2)$ is an interior point of $\F$.
Then we can take an small open neighborhood $V$ of $CH(S_1)\cap CH(S_2)$ included in $\F$ such that if $U$ is sufficiently small, then for any $(S_1',S_2')\in U$ we have $CH(S_1')\cap CH(S_2')\subset V$. Hence the area of $CH(S_1')\cap CH(S_2')$ is smaller than that of $V$ for any $(S_1',S_2')\in U$. From the Gauss-Bonnet Theorem, we have
\[ 2\pi \leq \sum_{\text{$v$: vertex of }CH(S_1')\cap CH(S_2')}\An (v) \leq \Area (V)+2\pi .\]
Since $\phi_\F (S_1,S_2)=2\pi$, $\phi$ is continuous at $(S_1,S_2)$.

\underline{Case 3:} Only one of $S_1$ and $S_2$ belongs to $\partial_2G$.

In this case $CH(S_1)\cap CH(S_2)$ is a geodesic segment or a geodesic half-line.
We assume that $\#S_1=2$ and $\#S_2\geq 3$. Let $v$ be a vertex of $CH(S_1)\cap CH(S_2)$ contained in the interior of $\F$.
Note that the geodesic line $CH(S_1)$ meets a boundary component $B$ of $CH(S_2)$ at $v$. Take $(S_1',S_2')\in U$ and assume that $U$ is sufficiently small. If $\#S_1'=2$, then $CH(S_1')$ meets $CH(S_2')$ at a point $v'$ close to $v$, which is also contained in $\F$.
Hence $\An (v)=\pi =\An (v')$.
If $\#S_1'>2$, then $CH(S_1')$ has two boundary components $B_1,B_2$ meeting a boundary component $B'$ of $CH(S_2')$, which is close to $B$, at $w_1,w_2$ respectively, which are contained in $\F$.
The vertices $w_1,w_2$ are contained in a small open neighborhood of $v$. Then the interior angle at $w_1$ is close to the exterior angle at $w_2$, which implies that the sum $\An (w_1)+\An(w_2)$ is close to $\pi=\An (v)$.
Therefore $\phi_\F$ is continuous at $(S_1,S_2)$.

\underline{Case 4:} Both $\#S_1$ and $\#S_2$ are larger than $2$.

Recall that at most finitely many boundary components of $CH(S_1), CH(S_2)$ intersect a neighborhood of $\F$, which implies that $\F$ includes at most finitely many vertices of $CH(S_1)\cap CH(S_2)$.
Hence we can see that $\phi_\F$ is continuous at $(S_1,S_2)$ by considering the movement of boundary components of $CH(S_1)$ and $CH(S_2)$ in $U$.
\end{proof}

From the argument in the above proof, we can prove that $\phi_\F$ is a Borel function. Moreover, the support of $\phi_{\F}$ is included in the compact subset $A(\ol{\F})\times A(\ol{\F})$ since $\ol{\F}$ is compact.
Recall that the number of vertices of $CH(S_1)\cap CH(S_2)$ in $\F$ is uniformly bounded for any $(S_1,S_2)\in \H (\partial G)\times \H (\partial G)$.
Hence $\phi_{\F}$ is a bounded Borel function with compact support. 

For the Dirichlet domain $\F_x$ centered at $x\in \tilde{\gS}$ for the action of $G$ on $\tilde{\gS}$ we can define $C_{\F_x}$ by the same way as $C_\F$.
From the proof of Lemma \ref{lem:measure zero fundamental domain}, there exists a smooth curve $c\:[0,1]\rightarrow \tilde{\gS}$ such that for any $\mu,\nu \in \SC (\gS)$, the set
\[ \{ t\in [0,1]\mid \mu \times \nu (C_{\F_c (t)} )>0 \}\]
is countable. In order to apply the same method of proving the continuous extension of the intersection number on $\SC (\gS)$, we prove the following lemma:

\begin{lemma}\label{lem:angle is independent of fundamental domain}
Let $\F'$ be a Dirichlet domain for the action of $G$ on $\tilde{\gS}$.
By removing some edges and vertices of $\F'$, we assume that $G(\F')=\tilde{\gS}$ and $g\F'\cap \F'=\emptyset $ for any non-trivial $g\in G$. Then for any $\mu, \nu\in \SC (\gS)$ we have
\[ \int \phi_\F d\mu\times\nu =\int \phi_{\F'} d\mu\times\nu.\]
\end{lemma}
\begin{proof}
For a subset $U$ of $\tilde{\gS}$ we define a function $\phi_U$ by
\[ \phi_U (S_1,S_2):=\sum_{\text{$v$: vertex of }CH(S_1)\cap CH(S_2) \text{ in } U}\An (v)\]
for $(S_1,S_2)\in \H (\partial G)\times \H(\partial G)$. Then we can see that for any countable family of pairwise disjoint subsets $\{U_\lambda \}$ of $\tilde{\gS}$ we have
\[ \phi_{\sqcup_\lambda U_\lambda}=\sum_{\lambda}\phi_{U_\lambda}.\]
For a subset $U$ of $\tilde{\gS}$ and $g\in G$ we have $\phi_{gU}(S_1,S_2)=\phi_U(g^{-1}S_1,g^{-1}S_2) $ for any $(S_1,S_2)\in \H (\partial G)\times \H(\partial G)$.
Note that $\phi_{g_1\F \cap g_2\F'}$ is a Borel function for any $g_1,g_2\in G$.
Since $\mu\times \nu$ is $G$-invariant, we obtain
\begin{align*}
\int \phi_\F d\mu\times\nu 
&=\int \phi_{\sqcup_{g\in G}\F \cap g\F'} d\mu\times \nu \\
&=\sum_{g\in G} \int \phi_{\F \cap g\F'} d\mu\times \nu \\
&=\sum_{g\in G} \int \phi_{g^{-1}\F \cap \F'} d\mu\times \nu \\
&=\int \phi_{\F'} d\mu\times\nu.
\end{align*} 
This completes the proof.
\end{proof}

The following proposition is the main part of our proof of Theorem \ref{thm:intersection functional}.

\begin{proposition}\label{prop: angle bilinear functional}
There exists a unique symmetric continuous $\RRR$-bilinear functional
\[ \psi \: \SC(\gS)\times \SC(\gS) \rightarrow \RRR \] 
such that for any non-trivial finitely generated subgroups $H$ and $K$ of $G$ we have
\[ \psi (\eta_H,\eta_K )=\sum_{\text{$v$: vertex of }C_H\times_\gS C_K}\An (v).\]
\end{proposition}
\begin{proof}
It is sufficient to prove that the restriction of $\phi_\F^\ast$ to $\mathrm{Span}(\SC_r(\gS))\times \mathrm{Span}(\SC_r(\gS))$ is continuous.
Take $(\mu_n ,\nu_n)\in \mathrm{Span}(\SC_r(\gS))\times \mathrm{Span}(\SC_r(\gS))\ (n\in \NN)$ converging to $(\mu, \nu )\in \mathrm{Span}(\SC_r(\gS))\times \mathrm{Span}(\SC_r(\gS))$.
We prove that $\phi_\F^\ast (\mu_n,\nu_n)$ converges to $\phi_\F^\ast (\mu,\nu)$ partially following the proof of $(3)\Rightarrow (4)$ in Proposition \ref{prop:Portmanteau}.
We will also use the method that we used in the proof of Theorem \ref{thm:intersection number of subset currents}.

Fix $\varepsilon >0$. By moving the base point of the Dirichlet domain, we can assume that a Dirichlet $\F$ satisfies the condition that
\[ \mu (\partial A(\F ))=\nu (\partial A(\F )) =\mu\times \nu (C_{\F})=0.\]
Set
\[ M:=\sup \{ \mu_n (A(\F )),\nu_n (A(\F ))\mid n\in \NN \} ,\]
\[ C:=\sup \{ \phi_\F (S_1,S_2)\mid (S_1,S_2)\in \H (\partial G)\times \H (\partial G)\}, \]
and set
\begin{align*}
A_t&:=\{ (S_1,S_2)\in \H (\partial G)\times \H (\partial G)\mid \phi_\F(S_1,S_2) \geq t\}, \\
B_t&:=\{ (S_1,S_2)\in \H (\partial G)\times \H (\partial G)\mid \phi_\F(S_1,S_2) = t\}
\end{align*}
for $t\in [0,C]$. Then we have
\[ \int \phi_\F d\mu\times \nu=\int_0^C \mu \times \nu (A_t) dm_\RR(t) \]
and 
\[ \int \phi_\F d\mu_n\times \nu_n=\int_0^C \mu_n \times \nu_n (A_t) dm_\RR(t). \]
Now, it is sufficient to prove that $\mu_n\times \nu_n(A_t)$ converges pointwise to $\mu\times \nu(A_t)$ for $m_\RR$-a.e. $t\in [0,C]$.
Note that $A_t\subset A(\F)$ for any $t>0$ and so
\[\mu_n\times \nu_n(A_t),\mu\times \nu(A_t)\leq M^2.\]

We know that $\partial A_t \subset B_t \cup \gD (\phi_\F )$ and $\mu \times \nu (B_t)=0$ for $m_\RR$-a.e. $t\in [0,C]$.
From the proof of Proposition \ref{prop:Portmanteau}, if $\mu\times \nu (\gD (\phi_\F ))=0$, then $\mu_n\times \nu_n(A_t)$ would converge pointwise to $\mu\times \nu(A_t)$ for $m_\RR$-a.e. $t\in [0,C]$. However, we have $\gD (\phi_\F)\subset C_\F \sqcup D_\F$ from Lemma \ref{lem:non-continuous points of phi}, and $\mu\times \nu (D_\F )$ is not necessarily equal to zero.
Hence, we need to evaluate the influence of $\mu\times \nu (D_\F )$.

From now on, we assume that $\mu\times \nu(D_\F)>0$.
Note that for any $(S_1,S_2)\in D_\F$ we have $\phi_\F (S_1,S_2)=0$. Therefore $A_t\cap D_\F=\emptyset$ for every $t>0$.
Moreover, if $\mu \times \nu (B_t)=0$, then $\mu \times \nu (A_t)=\mu \times \nu (\mathrm{Int}(A_t))$, which will be used later.
Fix $\delta>0$ such that $M^2\delta<\varepsilon$. Then we have
\[ \int_0^{\delta} \mu \times \nu (A_t) dm_\RR(t), \int_0^{\delta} \mu_n \times \nu_n (A_t) dm_\RR(t) <\varepsilon.\]

Similarly to the proof of Theorem \ref{thm:intersection number of subset currents}, we construct an open subset $V$ of $\H (\partial G )\times \H (\partial G)$ such that $\mu_n \times \nu_n (A_t \cap V)\leq \varepsilon $ for any $n\in \NN, t\in [\delta,C]$, and $\mu \times \nu (D_\F \setminus V)=0$. 

Since $\mu ,\nu \in \mathrm{Span}(\SC_r(\gS))$ and $D_\F$ is compact, there exists $(S_1^k, S_2^k)\in D_\F\ (k=1,\dots, m)$ such that $(S_1^k,S_2^k)$ is an atom of $\mu \times \nu$ for every $k$ and
\[ \mu \times \nu (D_\F )= \sum_{k=1}^m \mu \times \nu (\{ (S_1^k,S_2^k)\} ).\]
In order to obtain this equation, we have restricted $\phi_\F^\ast$ to $\mathrm{Span}(\SC_r(\gS))\times \mathrm{Span}(\SC_r(\gS))$.

Let $(S_1,S_2)\in \{ (S_1^k, S_2^k)\}_{k=1,\dots ,m}$.
Let $B:=CH(S_1)\cap CH(S_2)$.
Since $S_i$ is the limit set of a finitely generated subgroup of $G$ for $i=1,2$, there exists $g\in G$ such that
$\gL(\langle g \rangle )=B(\infty )$ and $g(S_1,S_2)=(S_1,S_2)$.
Since $\mu (\partial A(\F) ) = \nu (\partial A(\F))=0$, $(S_1,S_2)$ belongs to $\mathrm{Int}(A(\F))\times \mathrm{Int}(A(\F))$, that is, $B$ passes through $\mathrm{Int}(\F)$.

Since $g$ can be considered as a self-homeomorphism of $\H (\partial G)\times \H (\partial G)$ fixing $(S_1,S_2)$, for any $L\in \NN$ we can take an open neighborhood $U$ of $(S_1,S_2)$ such that
\[ g(U),\dots ,g^L (U)\subset \mathrm{Int}(A(\F ))\times \mathrm{Int}(A(\F )).\]
Take a compact convex polygon $O$ of $\HH$ such that $O$ includes $g(\F), \dots ,g^L(\F)$.
We can also assume that $U$ is sufficiently small such that 
\[ \Area (CH(T_1)\cap CH (T_2)\cap O )<1\]
for any $\ell=1,\dots ,L$ and $(T_1,T_2)\in g^\ell(U)$.

Now, we consider $W_\ell:=g^{\ell}(U\cap A_t)$ for $t\in [\delta, C]$ and $\ell=1,\dots ,L$. Take $\alpha \in \NN$ such that $\alpha \delta > 2\pi +1$.
Note that $\alpha $ is independent of $L$.
We prove that $W_1,\dots ,W_L$ are $(\alpha-1)$-essentially disjoint, that is, for any $(T_1,T_2)\in \H (\partial G)\times \H (\partial G)$ we have
\[ \# \{ \ell \mid W_\ell \ni (T_1,T_2) \} \leq \alpha -1.\]
To obtain a contradiction, suppose that there exist $1\leq \ell_1<\ell_2<\cdots <\ell_\alpha \leq L$ such that
\[ W:=\bigcap_{s=1}^\alpha W_{\ell_s}\]
is not empty.
Take $(T_1,T_2)\in W$. Since $(T_1,T_2)\in W_{\ell_s}$, we have $\phi_\F (g^{-\ell_s}T_1,g^{-\ell_s}T_2)\geq \delta$, which implies that $\phi_{g^{\ell_s}\F}(T_1,T_2)\geq \delta$. Note that $\phi_{g^{\ell_s}\F}(T_1,T_2)$ equals the sum of the exterior angle of vertices of $CH(T_1)\cap CH(T_2)$ in $g^{\ell_s}\F$. Hence the sum of the exterior angle of vertices of $CH(T_1)\cap CH(T_2)$ in $O$ is larger than or equal to $\alpha \delta$.
Note that $CH(T_1)\cap CH(T_2)\cap O$ is a convex polygon.
From the Gauss-Bonnet Theorem, we have
\[ \Area (CH(T_1)\cap CH(T_2)\cap O)\geq \alpha \delta -2\pi >1,\]
a contradiction.

Hence $W_1,\dots ,W_L$ are in particular $\alpha$-essentially disjoint and
\begin{align*}
\mu_n\times \nu_n \left( \bigcup_{\ell =1}^LW_\ell \right) 
&\geq \frac{1}{\alpha }\sum_{\ell=1}^L \mu_n\times \nu_n ( W_\ell )\\
&= \frac{1}{\alpha }\sum_{\ell=1 }^L \mu_n \times \nu_n (U \cap A_t)\\
&=\frac{L}{\alpha }\mu_n \times \nu_n (U \cap A_t).
\end{align*}
Since $W_\ell$ is included in $A(\F)\times A (\F)$ for every $\ell=1,\dots , L$, we have
\[ \mu_n\times \nu_n (A_t \cap U)\leq \frac{\alpha M^2}{L}.\]

From the above, we can take an open neighborhood $U_k$ of $(S_1^k,S_2^k)$ such that
\[ \mu_n\times \nu_n (A_t \cap U_k)\leq \frac{\alpha M^2}{L} \]
for every $k=1,\dots , m$.
Set $V:= U_1 \cup \cdots \cup U_m$. Then
\[ \mu_n \times \nu_n (A_t \cap V)\leq \sum_{k=1}^m \frac{\alpha M^2}{L}\leq \frac{m\alpha M^2}{L}.\]
By taking a sufficiently large $L$, we have
\[ \mu_n \times \nu_n (A_t \cap V)\leq \varepsilon \]
for any $n\in \NN$ and $t\in [\delta, C]$. Moreover, $\mu\times \nu (D_\F \setminus V)=0$.

From Proposition \ref{prop:Portmanteau}, for any $t\in [\delta, C]$ with $\mu \times \nu (B_t)=0$ we have
\begin{align*}
&\mu \times \nu (A_t)=\mu \times \nu (\mathrm{Int}(A_t))\\
\leq &\liminf_{n\rightarrow \infty} \mu_n\times \nu_n (\mathrm{Int}(A_t))
\leq \liminf_{n\rightarrow \infty} \mu_n\times \nu_n (A_t)\\
\leq &\limsup_{n\rightarrow \infty}\mu_n\times \nu_n (A_t)\\
\leq &\limsup_{n\rightarrow \infty}\mu_n\times \nu_n (A_t\setminus V)+\limsup_{n\rightarrow \infty}\mu_n\times \nu_n (A_t\cap V)\\
\leq &\limsup_{n\rightarrow \infty}\mu_n\times \nu_n (\ol{A_t}\setminus V)+\varepsilon\\
= &\mu \times \nu (\ol{A_t}\setminus V)+\varepsilon\\
\leq &\mu \times \nu (A_t)+\mu \times \nu (\partial A_t\setminus V)+\varepsilon\\
\leq &\mu \times \nu (A_t)+\mu \times \nu (D_\F \setminus V)+\varepsilon\\
\leq &\mu \times \nu (A_t)+\varepsilon.
\end{align*}
Since $\varepsilon>0$ is arbitrary, for $m_\RR$-a.e. $t\in [\delta, C]$,
\[ \lim_{n\rightarrow \infty }\mu_n \times \nu_n (A_t)=\mu \times \nu (A_t).\]
Then
\begin{align*}
&\left| \int_0^C \mu_n \times \nu_n (A_t) dm_\RR(t) -\int_0^C \mu \times \nu (A_t) dm_\RR(t) \right| \\
\leq 	&\int _0^\delta | \mu_n\times \nu_n (A_t)-\mu \times \nu (A_t)| dm_\RR (t) \\
		&+\int _\delta^C | \mu_n\times \nu_n (A_t)-\mu \times \nu (A_t)| dm_\RR (t)\\
\leq 	&2M^2\delta+\int _\delta^C | \mu_n\times \nu_n (A_t)-\mu \times \nu (A_t)| dm_\RR (t)\\
\leq 	&2\varepsilon +\int _\delta^C | \mu_n\times \nu_n (A_t)-\mu \times \nu (A_t)| dm_\RR (t).
\end{align*}
Note that the last term
\[ \int _\delta^C | \mu_n\times \nu_n (A_t)-\mu \times \nu (A_t)| dm_\RR (t)\]
converges to $0$ when $n\rightarrow \infty$. Since $\varepsilon>0$ is arbitrary,
\[\int_0^C \mu_n \times \nu_n (A_t) dm_\RR(t) \underset{n\rightarrow \infty}{\rightarrow} \int_0^C \mu \times \nu (A_t) dm_\RR(t) .\]
This completes the proof.
\end{proof}

\begin{proof}[Proof of Theorem \ref{thm:intersection functional}]
Recall that by the Gauss-Bonnet Theorem for non-trivial finitely generated subgroups $H$ and $K$ of $G$ we have
\[ 2\pi \chi (C_H\times_\gS C_K )=-\Area (C_H\times_\gS C_K)+\sum_{\text{$v$: vertex of }C_H\times_\gS C_K}\An (v).\]
From Propositions \ref{prop:area bilinear functional} and \ref{prop: angle bilinear functional} we define a functional $\hat{\chi }$ to be
\[ \frac{1}{2\pi}(-f^\ast +\psi ),\]
which is a continuous $\RRR$-bilinear functional from $\SC (\gS)\times \SC (\gS)$ to $\RR$ sending $(\eta_ H, \eta_K)$ to $\chi (C_H \times_\gS C_K)$. Since $f^\ast$ and $\psi$ are symmetric, so is $\hat{\chi}$.

Recall that
\[ \N (H,K)=-\chi (C_H\times_\gS C_K)+i(C_H,C_K).\]
Hence we define a functional $\N$ to be $ -\hat{\chi}+i_{\SC}$.
Then $\N$ is a symmetric continuous $\RRR$-bilinear functional from $\SC (\gS)\times \SC (\gS)$ to $\RR$ sending $(\eta_H ,\eta_K )$ to $\N (H ,K )$.
Since $\N(H,K)\geq 0$ for any finitely generated subgroups $H$ and $K$ of $G$, we have $\N (\mu,\nu )\geq 0$ for any $\mu,\nu \in \SC (\gS)$ from the denseness property of rational subset currents for $G=\pi_1(\gS)$.
The uniqueness of $\N$ also follows by the denseness property of rational subset currents.
\end{proof}

\section{Projection from subset currents onto geodesic currents}\label{sec: projection B}

Let $\gS $ be a compact hyperbolic surface (possibly with geodesic boundary) and $G$ the fundamental group of $\gS$.
We consider the universal cover $\tilde{\gS}$ of $\gS$ as a subspace of the hyperbolic plane $\HH$ and identify $\partial G$ with the limit set $\gL (G)=\tilde{\gS}(\infty) \subset \partial \HH $.
The notation in this section is based on that in Sections \ref{sec:intersection number of subset currents} and \ref{sec:an intersection functional}.

Take a non-trivial finitely generated subgroup $H$ of $G$. We consider the case that convex core $C_H$ has a boundary. The restriction of the map $p_H\:C_H\rightarrow \gS$ to each boundary component of $C_H$ can be considered as a closed geodesic on $\gS$.
We denote by $\partial C_H$ the set of all boundary components of $C_H$.
In the case that $H$ is an infinite cyclic group, that is, $C_H$ itself is a closed geodesic on $\gS$, we consider a copy of $C_H$, denoted by $C_H'$.
Then we set $\partial C_H:=\{ C_H, C_H' \}$, which can be interpreted as circles with different orientation. If $C_H$ does not have a boundary, that is, $H$ is a surface group, then $\partial C_H$ is an empty set.

For a closed geodesic $c$ on $\gS$ we define a counting geodesic current $\eta_c$ to be $\eta_g$ for $g\in G$ satisfying the condition that a representative of $g$ is freely homotopic to $c$. If $c$ is a boundary component of $C_H$, then we can take $h\in H$ such that $\eta_c=\eta_h$.
The following theorem is the main theorem in this section:

\begin{theorem}\label{thm:projection from sc to gc}
Let $\gS$ be a compact hyperbolic surface.
There exists a unique continuous $\RRR$-linear map
\[ \B\: \SC (\gS)\rightarrow \GC (\gS)\]
such that for any non-trivial finitely generated subgroup $H$ of $G$ we have
\[ \B (\eta_H )=\frac{1}{2}\sum_{c\in \partial C_H} \eta_c.\]
Especially, the restriction of $\B$ to $\GC (\gS)$ is the identity map.
\end{theorem}

Note that if $\partial C_H$ is empty, then $\B (\eta_H)$ is the zero measure in the above theorem.

\subsection{Construction of projection}
Take a non-trivial finitely generated subgroup $H$ of $G$ with $\partial C_H\not=\emptyset$.
For a boundary component $c$ of $C_H$ we regard $c$ as an element of $H$ such that $\eta_c=\eta_{\langle c \rangle}$.
Note that an element $h\in H$ satisfying the condition that $\eta_c =\eta_{\langle h \rangle }$ is not unique.
Recall that we have the continuous $\RRR$-linear map $\iota_H$ from $\SC(H)$ to $\SC (\gS)$ (see Section \ref{sec:relation between subgroups}).
Then we have
\[\sum_{c\in \partial C_H} \eta_c 
=\sum_{c\in \partial C_H}\iota_H (\eta_{\langle c \rangle }^H)
=\iota_H \left (\sum_{c\in \partial C_H} \sum_{h\langle c\rangle \in H/\langle c \rangle }\delta_{h\gL (\langle c \rangle )}\right).
\]

For $S\in \H (\partial G)$ we define $b(S)$ to be the set of all connected components of $\partial \HH\setminus S$.
Since $\partial \HH$ is homeomorphic to $S^1$, $b(S)$ consists of at most countably many open intervals.
For $S\in \H (\partial G)$ and $\alpha \in b(S)$, the boundary $\partial \alpha$ belongs to $\partial _2 G$.
Now, we prove the following lemma:
\begin{lemma}\label{lem:relation betwee C_H and b}
The following equality holds:
\[ \sum_{c\in \partial C_H} \sum_{h\langle c\rangle \in H/\langle c \rangle }\delta_{h\gL (\langle c \rangle )} =\sum_{\alpha \in b(\gL (H))} \delta_{\partial \alpha }.\]
\end{lemma}
\begin{proof}
First, we consider the case that $H$ is an infinite cyclic group. Then $c$ is a generator of $H$ and the left hand side equals $2\delta_{\gL (H)}$, which coincides the right hand side. Actually, this is the reason of the definition of $\partial C_H$.

From now on, we assume that $H$ is not an infinite cyclic group and $\partial C_H$ is not empty.
We define a map $\psi$ from $\partial C_H \times H/\langle c \rangle $ to $b(\gL (H))$ as follows.
For each $c\in \partial C_H$ we have a cyclic subgroup $\langle c \rangle $ of $H$, and the convex hull $CH_{\langle c \rangle }$ of the limit set $\gL (\langle c \rangle )$ is a boundary component of $CH_H$.
For $h\langle c \rangle \in H/\langle c \rangle$, we define $\psi (c , h\langle c \rangle )$ to be the open interval connecting the two points of $h\gL (\langle c \rangle )$ and not intersecting $\gL (H)$, which implies that $\psi (c ,h \langle c\rangle )\in b (\gL (H))$.
Then $\partial \psi (c, h \langle c \rangle )=h \gL (\langle c \rangle )$.
Hence, it is sufficient to see that $\psi $ is a bijective map.

First, we see that $\psi$ is surjective.
Take $\alpha \in b(\gL(H))$. Then there exists a boundary component $B$ of $CH_H$ such that $B(\infty )=\partial \alpha$. Take $c\in \partial C_H$ corresponding to $B$ and $h\in H$ such that $h\gL(\langle c \rangle) =B(\infty )=\partial \alpha$. Hence $\psi (c, h \langle c \rangle )=\alpha$.

Next, we see that $\psi$ is injective.
Take $c_1,c_2\in \partial C_H$ and $h_1\langle c_1 \rangle \in H/\langle c_1\rangle, h_2\langle c_2\rangle\in H/\langle c_2\rangle$.
It is sufficient to see that if $h_1\gL (\langle c_1\rangle )=h_2\gL (\langle c_2\rangle )$, then $c_1=c_2$ and $h_1\langle c_1 \rangle= h_2\langle c_2 \rangle$. Since $h_2^{-1}h_1CH_{\langle c_1\rangle }=CH_{\langle c_2\rangle }$, we can see that $c_1=c_2$.
Set $h:=h_2^{-1}h_1$, which fixes $\gL (\langle c_1 \rangle )$. Since $c_1$ is a simple closed geodesic on $C_H$, there exists no element $h_0$ of $H$ such that $c_1=h_0^k$ for some $k\geq 2$.
Therefore $h=c_1^k$ for some $k\in \Z$, which implies that $h_1\langle c_1\rangle =h_2\langle c_2 \rangle$.
This completes the proof.
\end{proof}

From the above lemma, we have
\[ \sum_{c\in \partial C_H} \eta_c = \iota_H \left( \sum_{\alpha \in b(\gL (H))} \delta_{\partial \alpha } \right) .\]

The strategy to prove Theorem \ref{thm:projection from sc to gc} is as follows.
First, we construct a measure $\B (\mu ) $ on $\partial _2 G$ for $\mu \in \SC (\gS )$. Next, we check that $\B ( \eta_H )$ equals $1/2\sum_{c \in \partial C_H } \eta_c$ for any non-trivial finitely generated subgroup $H$ of $G$.
Then we prove that $\B (\mu )$ is a geodesic current on $\gS$ for any $\mu \in \SC (\gS )$ and $\B$ is an $\RRR$-linear map from $\SC (\gS )$ to $\GC (\gS )$.
Finally, we prove that $\B$ is continuous. The uniqueness of $\B$ follows by the denseness property of rational subset currents for $G$.

We will denote by $\mathcal{O}$ the set of all open intervals of $\partial \HH$. We endow $\mathcal{O}$ with the topology induced by the Hausdorff distance. A set $b(S)$ is a subset of $\mathcal{O}$ for $S\in \H (\partial G)$. Define a function $\varphi$ from $\H (\partial G)\times \mathcal{O}$ to $\RR$ by
\[ \varphi (S,\alpha ):=\chi _{b(S)}(\alpha )\quad ((S,\alpha )\in \H (\partial G)\times \mathcal{O}) ,\]
that is, if $\alpha \in b(S)$, then $\varphi (S,\alpha )=1$; if $\alpha \not\in b(S)$, then $\varphi(S,\alpha)=0$.
For $\alpha \in \mathcal{O}$ we have a Dirac measure $\delta_\alpha $ on $\mathcal {O}$. Then $\varphi(S,\alpha )=\delta_\alpha (b(S))$ for $(S,\alpha )\in \H (\partial G)\times \mathcal{O}$. We have $\varphi(S,\alpha)=1$ if and only if $CH(\partial \alpha )$ is a boundary component of $CH(S)$.
We denote by $\mathcal{M}$ the counting measure on $\mathcal{O}$, that is, for any subset $U$ of $\mathcal {O}$, $\mathcal{M}(U)$ is the cardinality of $U$. For a Borel subset $E$ of $\partial_2 G$, set
\[ b(E):=\bigcup_{S\in E} b(S)\subset \mathcal{O}.\]
Then for any $\alpha \in \mathcal{O}$, $\alpha$ belongs to $b(E)$ if and only if $\partial \alpha $ belongs to $E$.

Now, for $\mu \in \SC (\gS)$ we define a measure $\B(\mu)$ on $\partial_2 G$ by
\[ \mathcal{B}(\mu )(E):=\frac{1}{2}\int_{b(E)} \left( \int \varphi (S ,\alpha ) d\mu (S)\right) d\mathcal{M}(\alpha )\]
for a Borel subset $E$ of $\partial_2 G$.
We can see that the preimage $\varphi ^{-1}(0)$ is an open subset of $\H (\partial G)\times \mathcal{O}$, which implies that $\varphi$ is a Borel function on $\H (\partial G )\times \mathcal{O}$.
Actually, $(S,\alpha )\in \varphi^{-1}(0)$ implies that $\partial \alpha $ is not a boundary component of $CH(S)$. It is easy to see that this is an ``open condition'' from Lemma \ref{lem:taking convex hull is continuous}.

Take a non-trivial finitely generated subgroup $H$ of $G$. Note that the action of $G$ on $\HH$ induces the action of $G$ on $\mathcal{O}$. Then for any Borel subset $E$ of $\partial_2G$ we have
\begin{align*}
&2\B (\eta_H )(E)\\
=&\int_{b(E)} \left( \int \varphi (S ,\alpha ) d\eta_H (S)\right) d\mathcal{M}(\alpha )\\
=&\int \left( \int_{b(E)} \varphi (S ,\alpha ) d\mathcal{M}(\alpha ) \right) d\eta_H (S)\\
=&\sum_{gH\in G/H} \int_{b(E)} \varphi (g\gL(H) ,\alpha ) d\mathcal{M}(\alpha )
=\sum_{gH\in G/H} \int_{b(E)} \delta_{\alpha }(b(g\gL( H))) d\mathcal{M}(\alpha ) \\
=&\sum_{gH\in G/H} \int_{b(g \gL (H))} \delta_{\alpha }(b(E)) d\mathcal{M}(\alpha )
=\sum_{gH\in G/H} \sum_{\alpha \in b(\gL (H))}\delta_{g(\alpha ) }(b(E)) \\
=&\sum_{gH\in G/H} \sum_{\alpha \in b(\gL (H) )}g_\ast (\delta_{ \alpha})(b(E))
=\iota_H \left( \sum_{\alpha \in b(\gL (H))}\delta_{\alpha}\right) (b(E))\\
=&\iota_H \left( \sum_{\alpha \in b(\gL (H))}\delta_{\partial \alpha}\right) (E)
=\sum_{c\in \partial C_H} \eta_c (E).
\end{align*}
Hence we see that
\[ \B (\eta _H )(E)=\frac{1}{2}\sum_{c\in \partial C_H} \eta_c (E). \]

\begin{lemma}\label{lem:B(mu) is a geodesic current}
For any $\mu \in \SC (\gS)$ the measure $\B(\mu)$ on $\partial_2G$ is a geodesic current on $\gS$.
\end{lemma}
\begin{proof}
First, we check that $\B (\mu)$ is $G$-invariant.
Take a Borel subset $E$ of $\partial_2 G$ and $g\in G$. Since $\mu$ is $G$-invariant, we have
\begin{align*}
2\B (gE)
&=\int_{b(gE)}\left( \int \varphi (S ,\alpha ) d\mu (S)\right) d\mathcal{M}(\alpha )\\
&=\int_{gb(E)}\left( \int \varphi (S ,\alpha ) d\mu (S)\right) d\mathcal{M}(\alpha )\\
&=\int_{b(E)}\left( \int \varphi (S ,g\alpha ) d\mu (S)\right) d\mathcal{M}(\alpha )\\
&=\int_{b(E)}\left( \int \varphi (g^{-1}S ,\alpha ) d\mu (S)\right) d\mathcal{M}(\alpha )\\
&=\int_{b(E)}\left( \int \varphi (S ,\alpha ) d\mu (S)\right) d\mathcal{M}(\alpha )\\
&=2\B (E).
\end{align*}

Next, we check that $\B(\mu )$ is a locally finite measure. Take a compact subset $K$ of $\tilde{\gS}$. From Lemma \ref{lem:A(K) is relatively compact}, it is sufficient to see that $\B(\mu)(A_2(K))<\infty$ for
\[ A_2(K)=\{ S\in \partial_2 G\mid CH(S)\cap K\not =\emptyset \}.\]
From the Fubini Theorem we have
\begin{align*}
2\B(\mu )(A_2(K))
&=\int_{b(A_2(K))}\left( \int \varphi (S ,\alpha ) d\mu (S)\right) d\mathcal{M}(\alpha )\\
&=\int \left( \int_{b(A_2(K))} \varphi (S ,\alpha ) d\mathcal{M}(\alpha ) \right) d\mu (S).
\end{align*}
Set
\[ \hat{\varphi}(S):=\int_{b(A_2(K))} \varphi (S ,\alpha ) d\mathcal{M}(\alpha ) \]
for $S\in \H (\partial G)$. It is sufficient to prove that $\hat{\varphi}$ is a bounded function with compact support.
Take $S\in \H (\partial G)$. We can see that $\hat{\varphi}(S)$ equals the number of boundary components of $S$ passing through $K$, which is uniformly bounded since $K$ is bounded.
Finally, we see that the support of $\hat{\varphi}$ is included in $A(K)$. Take $S\in \H (\partial G)\setminus A(K)$. Then 
$CH(S)\cap K=\emptyset $, which implies that $\hat{\varphi} (S)=0$.
This completes the proof.
\end{proof}

\begin{proof}[Proof of Theorem \ref{thm:projection from sc to gc}]
From the above lemma, we can see that $\B$ is an $\RRR$-linear map from $\SC (\gS)$ to $\GC(\gS)$.
It is sufficient to prove that $\B$ is continuous.
Take $\mu_n \in \SC (\gS)\ (n\in \NN)$ converging to $\mu \in \SC (\gS)$.
From Proposition \ref{prop:Portmanteau}, it is sufficient to prove that for any relatively compact Borel subset $E$ of $\partial_2G$ with $\B(\mu)(\partial E)=0$ the sequence $\B (\mu_n )(E)$ converges to $\B (\mu )(E)$.

Take a relatively compact Borel subset $E$ of $\partial_2G$ with $\B(\mu)(\partial E)=0$.
Define a map $\hat{\varphi }\: \H (\partial G )\rightarrow \RRR$ by
\[ \hat{\varphi}(S):=\int_{b(E)}\varphi (S,\alpha )d \mathcal{M}(\alpha )\]
for $S\in \H (\partial G)$. Then we have
\[ 2 \B (\mu_n)(E) =\int_{b(E)}\left( \int \varphi (S ,\alpha ) d\mu_n (S)\right) d\mathcal{M}(\alpha )= \int \hat{\varphi} d\mu_n \]
and
\[ 2\B (\mu )(E)=\int \hat{\varphi} d\mu.\]
From the proof of Lemma \ref{lem:B(mu) is a geodesic current}, $\hat{\varphi}$ is a bounded function with compact support.
It is sufficient to prove that the set $\gD (\hat{\varphi})$ of non-continuous points of $\hat{\varphi}$ has measure zero with respect to $\mu$ from Proposition \ref{prop:Portmanteau}.

Since $\B(\mu) (\partial E)=0$, we obtain
\[ 0=2\B(\mu )(\partial E )=\int \left( \int_{b(\partial E)} \varphi (S,\alpha )d \mathcal{M}(\alpha ) \right) d \mu (S) .\]
Note that for $S_1,S_2\in \partial_2G$, if $b(S_1)\cap b(S_2)\not=\emptyset$, then $S_1=S_2$. We can see that for $S\in \H (\partial G)$
\[ \int_{b(\partial E)}\varphi (S,\alpha )d\mathcal{M}(\alpha )=\# (b(S)\cap b(\partial E )).\]
Set
\[ U:=\{ S\in \H (\partial G)\mid b(S)\cap b (\partial E)\not=\emptyset \}.\]
Then for the characteristic function $\chi_U$ of $U$ on $\H (\partial G)$ we have
\[ \chi_U(S) \leq \int_{b(\partial E)}\varphi (S,\alpha )d\mathcal{M}(\alpha ) \]
for $S\in \H (\partial G)$, which implies that
\[ \mu(U)=\int \chi_U d\mu \leq \int \left( \int_{b(\partial E)} \varphi (S,\alpha )d \mathcal{M}(\alpha ) \right) d\mu (S)= 0.\]
Therefore $\mu (U)=0$. 

Now, we prove that $\gD (\hat{\varphi})\subset U$. Take $S\in \H (\partial G )\setminus U$, which implies that $b(S)\cap b(\partial E)=\emptyset$. First, we see that
\[ b(S)\subset b(\mathrm{Int}(E)\sqcup \ol{E}^c ).\]
Hence, $\hat{\varphi}(S)=\# (b(S)\cap b(\mathrm{Int}(E) ))$.
Since $E$ is relatively compact, we can take a compact subset $K$ of $\tilde{\gS}$ such that $E\subset A_2(K)$ by Lemma \ref{lem:A(K) is relatively compact}.
Note that there are only finitely many $\alpha_1,\dots \alpha_m\in b(S)$ satisfying the condition that $CH(\partial \alpha _i)\cap K\not=\emptyset$. Hence we do not need to consider $\alpha \in b(S) \setminus \{ \alpha _1, \dots ,\alpha_m \}$. Since each $\alpha_i$ belongs to $b(\mathrm{Int}(E))$ or $b(\ol{E}^c)$, we can take an open neighborhood $V$ of $S$ in $\H (\partial G)$ such that for any $S'\in V$ we have $b(S')\cap b(\partial E)=\emptyset$ and
\[ \# (b(S')\cap b(\mathrm{Int}(E)) )=\# (b(S)\cap b(\mathrm{Int}(E)) ).\]
This implies that $\hat{\varphi}$ is constant on $V$. Hence $\hat{\varphi}$ is continuous at $S$.
\end{proof}

\subsection{Application of projection}
In this subsection, we consider the application of the projection $\B$.
The following theorem relates the intersection number $i_\SC $ on $\SC (\gS )$ to the intersection number $i_\GC$ on $\GC (\gS )$.

\begin{theorem}\label{thm:projection and inequality}
For any subset currents $\mu,\nu \in \SC (\gS)$ the following inequality holds:
\[ i_{\SC }(\mu,\nu )\leq i_{\GC}(\B (\mu ),\B (\nu )).\]
If either $\mu$ or $\nu$ belongs to $\GC (\gS)$, then the equality holds.
\end{theorem}
\begin{proof}
From the denseness property of rational subset currents and rational geodesic currents on $\gS$, it is sufficient to prove that the inequality and the equality holds for $\eta_H$ and $\eta_K$ for non-trivial finitely generated subgroups $H$ and $K$ of $G$.
Recall that
\[ \CH_H =\{ (gH,x)\in G/H\times \tilde{\gS}\mid x\in gCH_H\}.\]
Set
\[ \partial \CH_H :=\{ (gH,x)\in G/H\times \tilde{\gS}\mid x\in g(\partial CH_H)\} \subset \CH_H.\]

First we consider the case that neither $H$ nor $K$ is cyclic and $C_H$ and $C_K$ have boundaries.
Note that if $C_H$ has no boundary, then the equality holds immediately since $i(C_H,C_K)=0$ (see Example \ref{example:intersection number 0}).
Recall that from Lemma \ref{lem:relation betwee C_H and b}, we have
\[ \sum_{c\in \partial C_H} \eta_c = \iota_H \left( \sum_{\alpha \in b(\gL (H))} \delta_{\partial \alpha } \right) .\]
By considering the correspondence between the Dirac measures in the equality, we can identify $\partial \CH_H$ with $\bigsqcup_{c\in \partial C_H} \CH_{\langle c \rangle }$.
Moreover, we obtain a natural inclusion map
\[ \iota \: \bigsqcup_{(c,d) \in \partial C_H \times \partial C_K}\CH_{\langle c\rangle}\times_{\tilde{\gS}}\CH_{\langle d\rangle } \hookrightarrow 
\CH_H \times_{\tilde{\gS}} \CH_K.\]
Since the inclusion map $\iota$ is $G$-equivariant, $\iota$ induces an inclusion map
\[ \bigsqcup_{(c,d) \in \partial C_H \times \partial C_K} C_{\langle c \rangle}\times_{\gS}C_{\langle d \rangle } \hookrightarrow 
C_H \times_{\gS} C_K.\]
Since $CH_H$ and $CH_K$ are surfaces with boundaries, we can see that any contractible component of $C_H\times_{\gS}C_K$ is a polygon with $2\ell$-edges for $\ell \geq 2$, each of whose vertices is the intersection point of a boundary component of $C_H$ and that of $C_K$. Therefore we have
\[ i(C_H,C_K)\leq \frac{1}{4} \sum_{c\in \partial C_H }\sum_{d\in \partial C_K } i(c ,d),\]
that is,
\[ i_{\SC}(\eta_H,\eta_K )\leq i_{\GC}(\B (\eta_H ), \B (\eta_K )).\]

In the case that both $H$ and $K$ are infinite cyclic groups, the equality is obvious.
Assume that $H$ is an infinite cyclic group, $K$ is not cyclic and $C_K$ has a boundary.
By the same way as the above, we have an inclusion map 
\[ \bigsqcup_{d\in \partial C_K} C_H\times_{\gS}C_{\langle d \rangle }\hookrightarrow C_H\times_{\gS} C_K.\]
We can see that any contractible component of $C_H\times_{\gS}C_K$ is a geodesic segment, each of whose endpoints is the intersection point of $C_H$ and a boundary component of $C_K$. Therefore we have
\[ i(C_H,C_K)= \frac{1}{2} \sum_{d\in \partial C_K} i(C_H,d),\]
that is,
\[ i_{\SC}(\eta_H,\eta_K )= i_{\GC}(\B (\eta_H ), \B (\eta_K )).\]
This completes the proof.
\end{proof}

\begin{supply}\label{supply:intersection number inequality}
For two transverse simple compact surfaces $(S_1,s_1), (S_2, s_2)$ on $\gS$ not forming a bigon, the same inequality also follows by the same proof as above, that is, we have
\[ i (s_1, s_2) \leq \frac{1}{4} \sum_{(c_1,c_2)\in \partial S_1 \times \partial S_2 }i (c_1, c_2) \]
if $S_1$ and $S_2$ are not $S^1$, where $\partial S_i$ is the set of boundary components of $S_i$. We also have
\[ i (s_1, s_2) = \frac{1}{2} \sum_{c\in \partial S_2 }i (s_1, c), \]
if $S_1=S^1$. Roughly speaking, this equality holds since the ``intersection'' of a closed curve $S_1$ and a simple compact surface $S_2$ consists of arcs whose endpoints are the intersection of $S_1$ and the $\partial S_2$.
\end{supply}

Consider the case that $\gS$ is a closed hyperbolic surface.
Bonahon \cite{Bon88} proved that there exists an embedding $L$ from the Teichm\"uller space $\mathcal{T}(\gS)$ of $\gS$ to $\GC(\gS)$, and for $m\in \mathcal{T}(\gS)$ and a non-trivial $g\in G$ the intersection number $i_{\GC}(L(m), \eta_g)$ coincides with the $m$-length of the (unoriented) geodesic corresponding to $g$, which is denoted by $\ell_m(g)$ and called the $m$-length of $g$.
We call $\ell_m:=i_{\GC}(L(m),\cdot )$ the \ti{$m$-length functional} on $\GC (\gS )$, which sends $\eta_g$ to $\ell_m (g)$ for every non-trivial element $g\in G$.


From Theorem \ref{thm:projection and inequality}, we obtain the following corollary, which generalizes the $m$-length functional $\ell_m$ on $\GC (\gS)$ to the $m$-length functional on $\SC (\gS)$:

\begin{corollary}\label{cor: m-length functinal on SC}
Let $\gS$ be a closed hyperbolic surface.
For $m\in \mathcal{T}(\gS)$ define the functional $\ell_m \: \SC(\gS )\rightarrow \RRR$ as 
\[ \ell_m(\mu ):=i_{\SC}( L(m), \mu )\]
for $\mu \in \SC (\gS)$. Then for every non-trivial finitely generated subgroup $H$ of $G$ we have
\[ \ell_m (\eta_H )=\frac{1}{2}\sum _{c \in \partial C_H} \ell_m (c),\]
where $\ell_m (c)$ is the $m$-length of $c$.
\end{corollary}
\begin{proof}
Since $L(m)\in \GC (\gS)$, by applying Theorem \ref{thm:projection and inequality} we have the equality
\[ i_{\SC} (L (m),\mu )=i_{\GC} (L(m), \B (\mu ))\]
for any $\mu \in \SC (\gS)$. Hence for a non-trivial finitely generated subgroup $H$ of $G$ we have 
\[\ell_m (\eta_H )=i_{\SC}( L(m), \eta_H )=i_{\GC} (L(m), \B (\mu ))=\frac{1}{2}\sum _{c \in \partial C_H} \ell_m (c),\]
as required.
\end{proof}

\begin{remark}[Intersection number and more general $m$-length functionals]
In the case that $\gS$ has no boundary, the above Bonahon's result was extended to all negatively curved Riemannian metrics by Otal in \cite{Ota90}, to negatively curved cone metrics by Hersonsky and Paulin in \cite{HP97}, to singular flat metrics (or nonpositively
curved Euclidean cone metrics) coming from quadratic differentials by Duchin-Leininger-Rafi in \cite{DLR10}, and to all singular flat metrics by Bankovic-Leininger in \cite{BL18}.
For any such metric $m$ on $\gS$, we can obtain an associated geodesic current $L_m\in \GC (\gS)$, and for non-trivial $g\in \pi_1(\gS)$, the intersection number $i_{\GC}(L_m ,\eta_g)$ equals the $m$-length of $g$.

Even if $\gS$ has a boundary, we can consider the doubled surface $D\gS$ of $\gS$ (see Example \ref{exa:doubled surface inclusion}).
Let $D=\pi_1(D\gS)$. Then the map $\iota_{G}^D \: \SC (\gS )\rightarrow \SC (D\gS)$ is injective and 
\[ \iota_G^D(\eta_H^G) =\eta_H^D\]
for any finitely generated subgroup $H$ of $G=\pi_1(\gS)$.
Any hyperbolic structure $m$ on $\gS$ whose boundary components are closed geodesics induces the hyperbolic structure $m_d$ on $D\gS$.
Then we define $\ell_m\: \SC(\gS )\rightarrow \RRR$ as
\[ \ell_m (\mu ):= i_{\SC }(L(m_d),\iota_G^D( \mu) )\]
for $\mu \in \SC (\gS)$. Then we see that for any non-trivial finitely generated subgroup $H$ of $G$,
\[ \ell_m (\eta_H^G )=i_{\SC }(L(m_d), \eta_H^D )= \frac{1}{2}\sum _{c \in \partial C_H}\ell_{m_d}(c)=\frac{1}{2}\sum _{c \in \partial C_H}\ell_{m}(c).\]
Note that $c\in \partial C_H$ can be considered as a closed geodesic on $\gS$ since $H$ is a subgroup of $G$.
\end{remark}

\begin{supply}
We can construct the functional $\ell_m$ on $\SC (\gS)$ more directly in the case that $m$ is a hyperbolic metric on $\gS$.
We can apply the method which we have used for the construction of the volume functional and the intersection number on $\SC (\gS )$.

Assume that $m$ coincides with the given hyperbolic metric on $\gS$.
Take the Dirichlet domain $\F=\F_x$ centered at $x \in \tilde{\gS}$ with respect to the action of $G$ on $\tilde{\gS}$, and modify $\F$ by removing some edges and vertexes from $\F$ such that $G(\F )=\tilde{\gS}$ and $g\F \cap \F =\emptyset $ for any non-trivial $g\in G$.
For $S\in \H (\partial G)\setminus \partial_2 G$ we define $\lambda_\F (S)$ to be the half of the sum of the length of each component of $\F\cap \partial CH (S)$ and for $S\in \partial_2 G$ we define $\lambda_\F (S)$ to be the length of $\F\cap \CH(S)$.
Then $\lambda_\F\: \H (\partial G)\rightarrow \RRR$ is a non-continuous bounded Borel function with compact support.
We can see that the $\RRR$-linear functional $\lambda_\F^\ast$ defined by
\[ \lambda_\F^\ast (\mu ):= \int \lambda_\F d \mu \quad (\mu \in \SC (\gS ) )\]
associates a counting subset currents $\eta_H$ with $\ell_m (\eta_H)$ for any non-trivial finitely generated subgroup $H$ of $G$ by the same way as that for the volume functional in Section \ref{sec:Volume functionals for Kleinian groups}.
Note that for $S\in \partial _2 G$ such that $CH(S)$ passes through the interior of $\F$, $\lambda_\F$ is continuous at $S$.
Hence the set $\Delta (\lambda_\F)$ of non-continuous points of $\lambda_\F$ consists of $S\in \H (\partial G)$ satisfying the condition that a boundary component of $CH(S)$ partially coincides with an edge of $\ol{\F}$.

We can prove the continuity of $\lambda_\F^\ast$ by using the technique of moving the center of the Dirichlet domain $\F$ in Lemma \ref{lem:measure zero fundamental domain}.
Actually, we can see that $\lambda_\F^\ast $ does not depend on $\F$ by the same way as Lemma \ref{lem:angle is independent of fundamental domain}.
For any $\mu \in \SC (\gS)$ there exists $x \in \tilde{\gS}$ such that $\mu (\gD (\lambda_{\F_x}))=0$.
Hence if a sequence $\{ \mu _n\}_{n\in \NN}$ of $\SC (\gS)$ converges to $\mu$, then $\lambda_{\F_x}^\ast (\mu_n)$ converges to $\lambda_{\F_x}^\ast (\mu)$ by Proposition \ref{prop:Portmanteau}. Therefore $\lambda_\F^\ast$ is continuous.
Moreover, by the denseness property of rational subset currents the functional $\lambda_\F^\ast$ coincides with $\ell_m$.
\end{supply}

Consider the case that $\gS$ is a closed hyperbolic surface.
We are interested in a ``fiber'' $\B^{-1}(\mu )$ for $\mu \in \GC (\gS)$, especially in the case that $\mu =\eta_c$ for some closed geodesic $c$ on $\gS$.
First, we consider the case that $c$ is a simple closed geodesic on $\gS$.
By cutting $\gS$ along $c$ and regarding the cut end as the boundary, we can obtain a compact hyperbolic surface or a pair of compact hyperbolic surfaces $\gS - c$.
Moreover, the inclusion map induces a locally injective continuous map $s$ from $\gS-c$ to $\gS$, which is a simple compact surface on $\gS$ or a pair of simple compact surfaces on $\gS$. Then we can obtain a finitely generated subgroup $H$ or a pair of finitely generated subgroups $H_1, H_2$ of $G$ corresponding to $\gS-c$. Set $ \eta (\gS-c) := \eta_H$ or $\eta_{H_1} +\eta_{H_2}$ respectively.
Then we have
\[ \B (\eta (\gS -c) )=\eta_c .\]
Hence the above construction of $\eta (\gS-c )$ can be regarded as a section of the projection $\B$.

However, in the case that $c$ has self-intersection, we can not perform the same construction.
Nevertheless, from the Scott theorem in \cite{Sco78, Sco85}, $c$ is geometric in a finite covering space of $\gS$, that is, there exists a finite index subgroup $G_1$ of $G$ such that $G_1$ contains an element corresponding to $c$ and $c$ lifts to a simple closed geodesic $c_1$ on the convex core $C_{G_1}$. Then we obtain a subset current $\eta( C_{G_1}-c_1)$ on $G_1$.
Moreover, we have the projection $\B_{G_1}$ from $\SC (G_1)=\SC (C_{G_1})$ to $\GC (G_1)$ and
\[ \B_{G_1}( \eta( C_{G_1}-c_1) ) =\eta_{c_1} ,\]
which is a counting geodesic current on $G_1$ corresponding to $c_1$.

Recall that for any non-trivial finitely generated subgroup $H$ of $G$ we have the map $\iota_H$ from $\SC (H)$ to $\SC (G)=\SC(\gS)$.
Then $\iota_{G_1}(\eta_{c_1} )=\eta_c$, and $\iota_{G_1}(\eta (C_{G_1}-c_1) )$ is a subset current on $G$. By Theorem \ref{thm: B and iota is commutative}, we see that
\[ \B (\iota_{G_1}(\eta (C_{G_1}-c_1) ) )= \iota_{G_1}(\B_{G_1}(\eta (C_{G_1}-c_1) ))=\iota_{G_1}(\eta_{c_1})=\eta_c .\]
Hence $\iota_{G_1}(\eta (C_{G_1}-c_1))$ is a required subset current on $G$, which is a counting subset current on $G$ or a sum of two counting subset currents on $G$.
Note that $\iota_{G_1}(\eta (C_{G_1}-c_1))$ depends on the choice of $G_1$. 

From now on, we do not assume that $\gS$ is a closed surface.
Let $H$ be a finitely generated subgroup of $G$. We mainly consider the case that $H$ is non-cyclic.
Then we have the projection $\B_H$ from $\SC (H)$ to $\GC (H)$ by considering $H$ as the fundamental group of $C_H$.
We will write $\B_G$ in place of $\B$ from now on.
Note that $\iota_H$ maps a geodesic current on $H$ to a geodesic current on $G$.

\begin{theorem}\label{thm: B and iota is commutative}
For any non-trivial finitely generated subgroup $H$ of $G$ we have the following commutative diagram:
\[
\xymatrix{
\SC (H )\ar[r]^{\B_H}\ar[d]_{\iota_H} &\GC (H) \ar[d]^{\iota_H |_{\GC(H)}}\\
\SC (G) \ar[r]_{\B_G} &\GC (G) \ar@{}[lu]|{\circlearrowright}.}
\]
\end{theorem}
\begin{proof}
In the case that $H$ is cyclic, then $\SC (H)$ coincides with $\GC (H)$ and the claim is trivial.
Hence we consider the case that $H$ is non-cyclic.

We can see that for any non-trivial finitely generated subgroup $K$ of $H$ we have
\[ \B_G\circ \iota_H (\eta_K^H) =\frac{1}{2}\sum_{c\in \partial C_K}\eta_c = \iota_H \circ \B_H (\eta_K^H ) \]
since the convex core $C_K$ and its boundary do not depend on $H$.
By the denseness property of rational subset currents we have the equality
\[ \B_G\circ \iota_H (\mu ) = \iota_H \circ \B_H (\mu )\]
for any $\mu \in \SC (H)$.

We also give a direct proof. 
Take a complete system of representatives $R$ of $G/H$.
For any $\mu \in \SC (H)$ and any Borel subset $E\subset \partial_2G$ we have
\begin{align*}
2\B_G(\iota_H(\mu) )(E)
=&\int_{b(E)} \int_{\H (\partial G)} \varphi (S,\alpha ) d\left(\sum_{gH\in G/H}g_\ast( \mu ) (S)\right) d\mathcal{M}(\alpha )\\
=&\sum_{g\in R} \int_{b(E)} \int_{\H (\partial H)} \varphi (gS,\alpha ) d \mu (S)d\mathcal{M}(\alpha )\\
=&\sum_{g\in R} \int_{b(E)} \int_{\H (\partial H)} \varphi (S,g^{-1}\alpha ) d \mu (S)d\mathcal{M}(\alpha )\\
=&\sum_{g\in R} \int_{g^{-1}(b(E))} \int_{\H (\partial H)} \varphi (S,\alpha ) d \mu (S)d\mathcal{M}(\alpha )\\
=&\sum_{g\in R} 2 \B_H (\mu )(g^{-1}E)=2 \iota_H\circ \B_H (\mu )(E),
\end{align*}
which is the required equality.
\end{proof}

\section{Denseness property of rational subset currents}\label{sec:denseness property of rational subset currents}

Recall that for an infinite hyperbolic group $G$ a subset current $\mu$ on $\SC (G)$ is called rational if there exist $c\in \RRR$ and a quasi-convex subgroup $H$ of $G$ such that $\mu= c\eta_H$ (see Subsection \ref{subsec:space of subset currents}). We denote by $\SC_r(G)$ the set of all rational subset currents on $G$.
We say that $G$ has the denseness property (of rational subset currents) if $\SC_r (G)$ is a dense subset of $\SC (G)$.
In this section, our goal is to prove that surface groups have the denseness property.

In Subsection \ref{subsec: denseness property of free groups}, we give a proof of the denseness property for a free group $F$ of finite rank.
We assume that $\SC_r(F)$ is a dense subset of the subspace $\mathrm{Span}(\SC_r(F))$ of $\SC(F)$ generated by $\SC_r(F)$(see \cite[Proposition 5.2]{KN13}) and prove that $\mathrm{Span}(\SC_r(F))$ is a dense subset of $\SC (F)$.
Our proof is based on that in \cite{Kap17} but we introduce the notion of an SC-graph on $F$, which corresponds to an element of $\mathrm{Span}(\SC_r(F))$ and will play a fundamental role in proving the denseness property for a surface group.

In Subsection \ref{subsec:approximation by a sequence of sugroups}, we consider a certain sequence of finitely generated subgroups $H_n$ of a free group $F$ of rank $2$ and we see that the sequence $\SC(H_n)$ approximates $\SC(F)$ (see Theorem \ref{thm:approximating theorem for free group}). 

In Subsection \ref{subsec:denseness property of surface groups}, we prove the denseness property of rational subset currents for a surface group $G$ by applying the method in the proof of Theorem \ref{thm:approximating theorem for free group} in Subsection \ref{subsec:approximation by a sequence of sugroups}. A certain sequence of finitely generated subgroups of $G$ that are isomorphic to a free group will play an essential role in the proof.

\subsection{Denseness property of free groups}\label{subsec: denseness property of free groups}

For a free group $F$ of finite rank, the denseness property for $F$ was first proved by Kapovich and Nagnibeda in \cite{KN13} (see \ref{thm:F_N dense}). Kapovich in \cite{Kap17} gave another self-contained proof to the denseness property for $F$.
We change some parts of the proof in \cite{Kap17} such that our method can apply to the proof of the denseness property for a surface group.
Our method of proving the denseness is constructing a sequence $\mu_n$ of $\mathrm{Span}(\SC_r(F))$ converging to a given $\mu \in \SC(F)$. 

Fix a free group $F$ of rank $N\geq 2$. Fix a free basis $B$ of $F$.
We denote by $X$ the Cayley graph of $F$ with respect to $B$.
The set of vertices of $X$ is denoted by $V(X)$, which is identified with $F$.
We give a path metric $d=d_X$ to $X$ such that each edge of $X$ has length one.
We identify $\partial F$ with $\partial X$.
The quotient space $F\backslash X$ is a graph consisting of one vertex attached $N$ loops.
For a closed subset $S$ of $\partial F=\partial X$ with $\#S \geq 2$ the convex hull $CH(S)$ of $S$ in $X$ is a union of all geodesic lines connecting two points of $S$.
We denote by $\H (\partial F)$ the space of closed subsets of $\partial F$ containing at least $2$ points and endow $\H (\partial F)$ with the Hausdorff distance $d_{\mathrm{Haus}}$ induced by a metric on $\partial F$ compatible with the topology.
The limit set $Y(\infty )$ of a subset $Y\subset X$ is the set of accumulation points of $Y$ in $\partial X$.

Recall that we have constructed $\CH_H$ for a non-trivial finitely generated subgroup $H$ of the fundamental group of a compact hyperbolic surface $\gS$.
Now, we define a similar space $\CH_H$ on $X$ for a non-trivial finitely generated subgroup $H$ of $F$.
For the convex hull $CH_H:=CH(\gL (H))\subset X$ of the limit set $\gL (H)$ we define
\[ \CH_H:=\{ (gH,x )\in F/H \times X\mid x\in gCH_H\} .\]
We have the projection map from $\CH_H$ to $X$.

We can consider $\CH_H$ as a geometric realization of the counting subset current $\eta_H$. Actually, for $gH\in F/H$ each connected component $gCH_H$ of $\CH_H$ corresponds to the Dirac measure at $g\gL (H)$.

\begin{definition}[SC-graph]
Let $Y$ be a graph, which is not necessarily connected, and $f$ a graph morphism from $Y$ to $X$, which is a continuous map sending vertices of $Y$ to vertices of $X$ and edges of $Y$ to edges of $X$.
We call the pair $(Y,f)$ a graph on $X$.
Now, we assume that $F$ acts on $Y$. When we consider a group action on a graph, we always assume that each element of the group acts as a graph isomorphism.
We call $(Y,f)$ a \ti{SC-graph} on $(F,X)$ (or simply $F$) if $(Y,f)$ satisfies the following conditions:
\begin{enumerate}
\item[SC1)] $f$ is an $F$-equivariant map;
\item[SC2)] the restriction of $f$ to each connected component $Y_0$ of $Y$ is injective and the image $f(Y_0)$ coincides with $CH (f(Y_0)(\infty ))$;
\item[SC3)] $\# f^{-1}(\id )<\infty $.
\end{enumerate}
We denote by $\Comp (Y)$ the set of all connected components of $Y$.
Since each $Y_0\in \Comp (Y)$ can be identified with $f(Y_0)\subset X$, we will write $f(Y_0)$ simply $Y_0$ when no confusion can arise.
Moreover, we often omit the projection $f$ when we consider an SC-graph on $F$.
\end{definition}

We note that the graph $\CH_H$ for a non-trivial finitely generated subgroup $H$ of $F$ is an SC-graph on $F$.
In this case, we usually omit the canonical projection from $\CH_H$ onto $X$.

For an SC-graph $(Y,f)$ on $F$ we can define a subset current $\eta(Y)$ on $F$ by
\[ \eta (Y):= \sum_{Y_0\in \Comp (Y)}\delta_{f(Y_0)(\infty )}.\]
We check that the measure $\eta(Y)$ is a subset current on $F$. Since $f$ is an $F$-equivariant map, $F$ acts on the set $\Comp (Y)$ of connected components of $Y$. Hence $\eta(Y)$ is an $F$-invariant measure. Explicitly,
for $g\in G$ and a Borel subset $E$ of $\H (\partial F)$ we have
\begin{align*}
\eta(Y) (g^{-1}E)
&=\sum_{Y_0\in \Comp (Y)} \delta_{Y_0(\infty )}(g^{-1}(E))\\
&=\# \{ Y_0\in \Comp (Y) \mid (gf(Y_0))(\infty ) \in E\} \\
&=\# \{ Y_0\in \Comp (Y) \mid (f(gY_0))(\infty ) \in E\} .
\end{align*}
Now we check that $\eta(Y)$ is locally finite. Recall that for $g\in F=V(X)$
\[ A_g=\{ S\in \H (\partial F) \mid CH(S)\ni g\} \]
is a compact subset of $\H (\partial F)$ and it is sufficient to see that $\eta(Y) (A_\id )<\infty$ from the proof of Lemma \ref{lem:continug subset currents}.
By the definition of an SC-graph on $F$,
\[\eta(Y) (A_{\id})=\# \{ Y_0\in \Comp(Y) \mid f(Y_0)\ni \id \} =\# f^{-1}( \id )<\infty. \]

\begin{remark}
If $Y_1,\dots ,Y_m$ are SC-graphs on $F$, then the disjoint union $\bigsqcup_k Y_k$ is also an SC-graph on $F$.
We can see that
\[ \eta (\bigsqcup_{k=1}^m Y_k)=\sum_{k=1}^m \eta (Y_k).\]

From Theorem \ref{thm:counting sc} and the condition (SC2), for an SC-graph $Y$ on $F$ there exist finitely generated subgroups $H_1,\dots H_m$ of $F$ such that $Y$ is isomorphic to
\[ \bigsqcup_{k=1}^m \CH_{H_k}\]
and we have
\[ \eta(Y)=\sum_{k=1}^m \eta_{H_k}.\]
Actually, for each connected component $Y_0\in \Comp(Y)$ and for the stabilizer $H=\Stab(Y_0)$ we have $Y_0=CH_H$. If $Y\setminus F(Y_0)$ is not empty, then $Y\setminus F(Y_0)$ can be considered as an SC-graph on $F$ and we can see that
\[ \eta(Y)=\eta (Y\setminus F(Y_0))+\eta_H.\]
Hence an SC-graph on $F$ corresponds to a finite sum of counting subset currents on $F$.
\end{remark}

Fix $\mu \in \SC (F)$. Assume that we have $\nu \in \mathrm{Span}(\SC_r(F))$ close to $\mu$.
Then $\nu $ can be represented by a finite sum of the rational subset currents, that is,
\[ \nu =\sum_{k=1}^m a_k \eta_{H_k}\]
for $a_k>0$ and non-trivial finitely generated subgroups $H_k$ of $F$ for $k=1,\dots ,m$. We can assume that $a_k$ is a rational number for $k=1,\dots , m$.
Then we can take $M\in \NN$ such that $b_k:=Ma_k$ is a positive integer for any $k$.
Therefore we can see that $M\mu$ is approximated by
\[ \sum_{k=1}^m b_k\eta_{H_k} =\eta (\bigsqcup_{k=1}^m \bigsqcup_{b_k} \CH_{H_k}) ,\]
where $\sqcup_{b_k}\CH_{H_k}$ means the $b_k$ copies of $\CH_{H_k}$. This consideration is useful when we prove that $\mathrm{Span}(\SC_r(F))$ is a dense subset of $\SC (F)$.

Now, we introduce the notion of a round-graph and the subset cylinder with respect to it, which was originally introduced in \cite{KN13,Kap17}.
We will introduce a generalized round-graph in Subsection \ref{subsec:approximation by a sequence of sugroups}.

\begin{definition}[Round-graph, see {\cite[Definition 3.3]{Kap17}}]
We denote by $\NN$ the set of positive integers.
Let $r\in \NN$. For $v\in V(X)$ we denote by $B(v,r)$ the closed ball centered at $v$ with radius $r$.
A subgraph $T$ of $B(v,r)$ is called a \ti{round-graph} centered at $v$ with radius $r$ if $T\ni v$ and there exists $S\in \H (\partial F)$ such that
\[ T=CH(S)\cap B(v,r).\]
We denote by $\R_r(v)$ the set of all round-graphs centered at $v$ with radius $r$.
For $T\in \R_r(v)$ we define the \ti{subset cylinder} $\SCyl (T)$ with respect to $T$ by
\[ \SCyl (T):=\{ S\in \H (\partial F) \mid CH(S)\cap B(v,r)=T\}.\]
We denote by $\R_r$ the union of $\R_r(v)$ over all $v\in V(X)$.
\end{definition}

\begin{remark}[Property of subset cylinders]\label{rem:property of subset cylinder on tree}
A subset cylinder $\SCyl(T)$ is an open and closed subset of $\H (\partial F)$ for any $T\in \R_r(v)$, which implies that if a sequence $\mu_n\in \SC(F)\ (n\in \NN)$ converges to $\mu \in \SC(F)$, then $\mu_n(\SCyl (T))$ converges to $\mu (\SCyl (T))$ by Proposition \ref{prop:Portmanteau}.
Moreover, for any $S\in \H (\partial F)$ and $v\in CH(S)\cap V(X)$ we have a sequence of round-graphs
\[\{ CH(S)\cap B(v,n) \}_{n\in \NN}, \]
and the family of $\SCyl (CH(S)\cap B(v,n) )$ for $n\in \NN$ forms a fundamental system of open neighborhoods of $S$.

For $T\in \R_r(v)$ and $g\in F$ we can see that $gT$ is a round-graph centered at $gv$ with radius $r$ and $\SCyl (gT)=g\SCyl (T)$.
This implies that $F$ acts on $\R_r$.
Since a subset current $\mu \in \SC (F)$ is $F$-invariant, $\mu (\SCyl (gT))=\mu (\SCyl (T) )$ for any $T\in \R_r(v)$ and $g\in F$.
Hence we often deal with round-graphs centered at $\id\in V(X)$.

For $T_1,T_2\in \R_r(v)$ if $T_1\not=T_2$, then $\SCyl (T_1)\cap \SCyl (T_2)=\emptyset$. Note that $\# \R_r(v)$ is finite for any $r\in \NN$ and $v\in V(X)$ since $X$ is a locally finite graph. Moreover, for any $r\in \NN$ we have
\[ \bigsqcup_{T\in \R_r (v )} \SCyl (T)=\{ S\in \H (\partial F) \mid CH(S)\ni v \}=A_{v}.\]

For $v_1,v_2\in V(X)$ and $T_1\in \R_r(v_1), T_2 \in \R_r (v_2)$, if $\SCyl (T_1 )\cap \SCyl (T_2)\not=\emptyset$ and $B(v_1,r)\cap B(v_2,r)\not=\emptyset$, then for $S\in \SCyl (T_1 )\cap \SCyl (T_2)$ we have
\[ CH(S)\cap B(v_1,r)=T_1 \text{ and }CH(S)\cap B(v_2,r)=T_2,\]
and so
\[ T_1\cap B(v_1,r)\cap B(v_2,r)=T_2 \cap B(v_1,r)\cap B(v_2,r). \]
\end{remark}

\begin{lemma}
Let $v\in V(X)$ and $r_1,r_2\in \NN$ with $r_1\leq r_2$. For any $T\in \R_{r_1}(v)$ we have the following equality:
\[ \SCyl (T)=\bigsqcup_{\substack{T'\in \R_{r_2}(v)\\[1pt] T'\cap B(v,r_1)=T}} \SCyl (T').\]
\end{lemma}
\begin{proof}
Let $S$ belong to the left side. Then $CH(S)\cap B(v,r_1)=T$ and $CH(S)\ni v$.
Set $T':= CH(S)\cap B(v,r_2)$. Then we can see that $T'\cap B(v,r_1)=T$ and $S\in \SCyl (T')$.

Let $S$ belong to the right side. There exists $T'\in \R_{r_2}(v)$ such that $T'\cap B(v,r_1)=T$ and $S\in \SCyl (T')$, which implies that $CH(S)\cap B(v,r_1)=T'\cap B(v,r_1)=T$, and so $S\in \SCyl (T)$. 
\end{proof}

From the above lemma, for $\mu \in \SC (F)$ if we know $\mu (\SCyl (T) )$ for every $T\in \R_r (\id)$, then we can calculate $\mu (\SCyl (T'))$ for every $r' \in \NN$ with $r' \leq r$ and every $T' \in \R_{r'}$.

The following proposition is useful for seeing that a sequence of subset currents on $F$ converges to a subset current on $F$:

\begin{proposition}[See {\cite[Proposition 3.7]{KN13}}]\label{prop:convergence and subset cylinder}
Let $\mu , \mu_n\in \SC (F)\ (n\in \NN)$. Then $\mu_n$ converges to $\mu$ if and only if for any $r\in \NN$ and any $T\in \R_r(\id )$ we have
\[ \lim_{n\rightarrow \infty }\mu_n (\SCyl (T))=\mu (\SCyl (T)).\]
\end{proposition}
\begin{proof}
The ``only if'' part follows immediately by Remark \ref{rem:property of subset cylinder on tree}. We prove the ``if'' part. 
Note that for any $r\in \NN$ and $T\in \R_r$ we have
\[ \lim_{n\rightarrow \infty }\mu_n (\SCyl (T))=\mu (\SCyl (T))\]
from the assumption.
Let $f$ be a continuous function from $\H (\partial F)$ to $\RR$ with compact support.
Fix $\varepsilon>0$. We construct a step function approximating $f$ by using subset cylinders.
From Lemma \ref{lem:A(K) is relatively compact}, since the support $\mathrm{supp}f$ of $f$ is compact, we can take $g_1,\dots ,g_m\in F$ such that
\[ \mathrm{supp}f\subset \bigcup_{i=1}^m A_{g_i}.\]
We can take $M>0$ such that 
\[ M>\sup_{n\in \NN} \left\{\mu_n \left( \bigcup_{i=1}^m A_{g_i}\right) \right\} ,\mu \left( \bigcup_{i=1}^m A_{g_i}\right) .\]
Take $\delta >0$ such that for any $S_1,S_2\in \H (\partial F)$, if the Hausdorff distance $d_{\mathrm{Haus}}(S_1,S_2)$ is smaller than $\delta$, then $|f(S_1)-f(S_2)|<\varepsilon/M$.
Take sufficiently large $r\in \NN$ such that for every $i=1,\dots ,m$ and any $T\in \R_r(g_i)$ the diameter of $\SCyl (T)$ is smaller than $\delta $ and $B(g_i ,r )$ contains $g_1, \dots ,g_m$.

Now, we prove that there exist $T_1,\dots ,T_L\in \R_{r}(g_1)\sqcup \cdots \sqcup\R_{r}(g_m)$ such that
\[ \bigcup_{i=1}^m A_{g_i}=\bigsqcup_{j=1}^L\SCyl (T).\]
Set $O:=\R_r(g_1)\sqcup \cdots \sqcup\R_r(g_m)$. If $\SCyl (T_1)\cap \SCyl (T_2)\not=\emptyset$ for $T_1\in \R_r(g_{i_1}), T_2\in \R_r(g_{i_2})$ and $i_1<i_2$, then we remove $T_2$ from $O$. We continue this operation for each pair of $T_1,T_2\in O$ one by one.
Finally, we can obtain $O$ satisfying the condition that for any $T_1,T_2\in O$, if $T_1\not=T_2$, then $\SCyl (T_1)\cap \SCyl (T_2)=\emptyset$.

Take any $S\in \cup_i A_{g_i}$ and take the smallest $i_0$ such that $S\in A_{g_{i_0}}$. Then there exists $T\in \R_r(g_{i_0})$ such that $S\in \SCyl (T)$. Since $CH(S)$ does not contain $g_1,\dots , g_{i_0 -1}$, $T=CH(S)\cap B(g_{i_0} ,r)$ also does not contain $g_1,\dots , g_{i_0 -1}$. 
Note that $B(g_{i_0} ,r)$ contains $g_1 ,\dots , g_{i_0 -1}$, which implies that $\SCyl (T)\cap \SCyl (T')=\emptyset$ for any $T' \in \R_r(g_1) \sqcup \cdots \sqcup \R_r (g_{i_0-1})$ by the last part of Remark \ref{rem:property of subset cylinder on tree}.
Hence $T\in O$. Therefore we have
\[ \bigcup_{i=1}^m A_{g_i}=\bigsqcup_{T\in O}\SCyl (T).\]

For each $T\in O$ set
\[ a_T:=\inf_{S\in \SCyl(T) }f(S).\]
We define a step function $\phi$ by
\[ \phi =\sum_{T\in O} a_T\chi_{\SCyl (T)}.\]
Then we have
\begin{align*}
\left| \int f d\mu -\int \phi d\mu \right|
&\leq \int \left| f-\phi \right| d\mu \\
&\leq \frac{\varepsilon}{M} \mu \left( \bigcup_{i=1}^m A_{g_i}\right) \\
&< \varepsilon.
\end{align*}
By the same way, we also have
\[ \left| \int f d\mu_n -\int \phi d\mu_n \right|<\varepsilon.\]
From the assumption, for a sufficiently large $n\in \NN$ we have
\begin{align*}
\left| \int \phi d\mu_n -\int \phi d\mu \right| 
&\leq \sum_{T\in O} |a_T| |\mu_n (\SCyl (T) )-\mu (\SCyl (T))|\\
&< \varepsilon
\end{align*}
Hence
\[ \left| \int f d\mu_n -\int f d\mu \right| < 3\varepsilon.\]
This completes the proof.
\end{proof}

From the proof of the above we have the following corollary:

\begin{corollary}[See {\cite[Proposition 3.7]{KN13}}]\label{cor: Kapovich's open neighborhood of a subset current}
Let $\mu \in \SC (F)$. The family of
\[ \{ \nu\in \SC (F)\mid |\mu (\SCyl (T))-\nu (\SCyl (T))|<\varepsilon \text{ for every }T\in \R_r(\id )\} \]
for $\varepsilon >0$ and $r\in \NN$ forms a fundamental system of open neighborhoods of $\mu$.
\end{corollary}

Let $\mu\in \SC(F), \varepsilon >0$ and $r\in \NN$.
We will construct an SC-graph $\gG$ on $F$ such that there exists $M\in \NN$ such that
\[ \left| \mu (\SCyl (T) )-\frac{1}{M}\eta (\gG )(\SCyl (T))\right| <\varepsilon \]
for any $T\in \R_\id (r)$.
We say that this SC-graph $\gG$ approximates $\mu$.
If we can obtain such an SC-graph $\gG$, then we see that $\mathrm{Span}(\SC_r(F))$ is a dense subset of $\SC(F)$ by Corollary \ref{cor: Kapovich's open neighborhood of a subset current}.
We will write simply $\eta_\gG$ in place of $\eta(\gG)$.

Now, we consider the value $\eta_\gG (\SCyl (T))$ for an SC-graph $\gG$ and $T\in \R_\id (r)$. From the definition of $\eta_\gG$ we have
\begin{align*}
\eta_\gG (\SCyl (T) )
&=\# \{ Y \in \Comp(\gG) \mid Y(\infty )\in \SCyl (T)\}\\
&=\# \{ Y\in \Comp (\gG )\mid Y\cap B(\id ,r)=T\}.
\end{align*}
This equation means that $\eta_\gG (\SCyl (T))$ coincides with the number of components of $\gG$ whose restriction to $B(\id,r)$ equals $T$. This is the most important idea for constructing an SC-graph $\gG$ approximating $\mu$ since we have an information of $\mu(\SCyl(T))$ for every $T\in \R_r(\id )$. Even if $\mu(\SCyl (T))$ is not an integer, we can take $q\in \mathbb{Q}$ approximating $\mu(\SCyl (T))$ and $Mq$ is an integer for some $M\in \NN$.

We also note that for $T\in \R_r (\id )$, $\eta_\gG (\SCyl (T))$ also equals the number of vertices $v$ of the quotient graph $F\backslash \gG$ satisfying the condition that for the connected component $Y$ of $\gG$ containing $\id$ that is projected onto $v$ we have $B(\id, r)\cap Y=T$, which means that the ``$r$-neighborhood'' of $v$ equals $T$. In the case that $\gG=\hat{CH}_H$ for a non-trivial finitely generated subgroup $H$ of $F$, it is easy to calculate $\eta_\gG(\SCyl (T))=\eta_H (\SCyl (T))$ since $F\backslash \hat{CH}_H$ can be identified with $H\backslash CH_H$. See \cite[Subsection 2.3]{Sas15} for a more detailed explanation.

For two vertices $u,v\in V(X)$, we want to combine a round-graph centered at $u$ with a round-graph centered at $v$.
For that purpose, we will use the following definition.

\begin{definition}
Let $r\in \NN$ and $u,v\in V(X)$. We denote by $B(u,v,r)$ the intersection of $B(u,r)$ and $B(v,r)$.
For $T_1\in \R_r(u), T_2\in \R_r(v)$ we say that $T_1$ and $T_2$ are \ti{connectable} if $T_1\cap B(u,v,r)=T_2\cap B(u,v,r)$.
Note that $B(u,v,r)$ can be empty and then $T_1$ and $T_2$ are connectable for any $T_1\in \R_r(u), T_2\in \R_r(v)$.

Assume that $B(u,v,r)$ is not empty.
A subgraph $J$ of $B(u,v,r)$ is called a $(u,v)$-round-graph with radius $r$ if $J\ni u,v$ and there exists $S\in \H (\partial F)$ such that \[J=CH(S)\cap B(u,v,r).\]
We denote by $\R_r(u,v)$ the set of all $(u,v)$-round-graph with radius $r$.
For $J\in \R_r(u,v)$ we define the subset cylinder $\SCyl (J)$ with respect to $J$ by
\[ \SCyl (J):=\{ S\in \H (\partial F)\mid CH(S)\cap B(u,v,r)=J\}.\]
For $T_1\in \R_r(u), T_2\in \R_r(v)$ we say that $T_1$ and $T_2$ are $J$-\ti{connectable} for $J\in \R_r(u,v)$ if $T_1\cap B(u,v,r)=J=T_2\cap B(u,v,r)$.
\end{definition}

\begin{remark}[Property of $(u,v)$-round-graph]\label{rem:property of u,v-round graph}
Let $u,v\in V(X)$ with $B(u,v,r)\not=\emptyset$. For $T\in \R_r(u)$ if $T$ contains $v$, then the intersection $T\cap B(u,v,r)$ belongs to $\R_r(u,v)$.
For any $J\in \R_r(u,v)$ we have
\[ \SCyl (J)=\bigsqcup_{\substack{T\in \R_r(u)\\[1pt] T\cap B(u,v,r)=J}}\SCyl (T)=\bigsqcup_{\substack{T'\in \R_r(v)\\[1pt] T'\cap B(u,v,r)=J}}\SCyl (T). \]
This implies that for any $\mu \in \SC (F)$ we have the equation:
\begin{equation}
\sum_{\substack{T\in \R_r(u)\\[1pt] T\cap B(u,v,r)=J}}\mu (\SCyl (T))=\sum_{\substack{T'\in \R_r(v)\\[1pt] T'\cap B(u,v,r)=J}}\mu (\SCyl (T')). \tag{$\ast_J$}
\end{equation}
This equation will be used for constructing an SC-graph approximating $\mu$.
\end{remark}

\begin{lemma}\label{lem:path connectable on tree}
Let $P$ be a geodesic path from $u\in V(X)$ to $v\in V(X)$, which passes through $v_0=u,v_1,\dots , v_m=v \in V(X)$ in this order.
Take $T_i\in \R_r(v_i)$ for $i=0,1,\dots ,m$. If $T_{i-1}$ and $T_{i}$ are connectable for every $i=1,\dots m$, then $T_0$ and $T_m$ are connectable.
\end{lemma}
\begin{proof}
Since $P$ is a geodesic path in the tree $X$, we have
\[ B(v_0,v_m,r)\subset \bigcap_{i=0}^m B(v_i,r),\]
which implies 
\[ B(v_0, v_m,r)\subset \bigcap_{i=1}^{m}B(v_{i-1},v_i,r).\]
From the assumption,
\[ T_{i-1}\cap B(v_{i-1},v_i,r)=T_i\cap B(v_{i-1},v_i,r)\]
for every $i=1,\dots ,m$. Therefore
\[ T_0\cap B(v_0,v_m,r)=T_1\cap B(v_0,v_m,r)=\cdots =T_m\cap B(v_0,v_m,r).\]
This completes the proof.
\end{proof}

Recall that $B$ is a free basis of $F$.
For $T\in \R_r(\id )$, if $\mu(\SCyl (T))$ is not a rational number, we want to approximate it by a rational number satisfying the equation $(\ast _J)$ in Remark \ref{rem:property of u,v-round graph} for two vertices $\id$ and $u\in B$.
Since $\#\R_r(\id,u)$ is finite and $F$ acts on $\R_r$, the system of the equations $(\ast_ J)$ for all $u\in B$ and $J\in \R_r(\id ,u)$ in Remark \ref{rem:property of u,v-round graph} can be considered as a finite homogeneous system of linear equations with respect to variables $\mu (\SCyl (T) )$ for $T\in \R_r(\id )$.
Hence we can apply the following lemma to the system of the equations $(\ast_J)$ for all $u\in B$ and $J\in \R_r(\id , u)$.

\begin{lemma}\label{lem:approximate by positive rational numbers}
Let $m,n$ be positive integers. Let $u= {}^t\!(u_1,\dots ,u_n)\in \RR^n$ with $u_i\geq 0$ for every $i$. Let $A=[ a_{ij}]$ be an $m\times n$ matrix with $a_{ij}\in \mathbb{Z}$. Assume that $Au=0$. Then for any $\varepsilon>0$ there exists $v\in \RR^n$ such that every coefficient of $v$ is a non-negative rational number, $Av=0$ and $||u-v||<\varepsilon$.
\end{lemma}
\begin{proof}
The proof is by induction on $n$. It is clearly true for $n = 1$.
Assume that $n>1$.
First, we consider the case that every $u_i$ is positive.
Since every entry of $A$ is an integer, we can take $w_1,\dots ,w_k\in \mathbb{Q}^n$ and $c_1,\dots, c_k\in \RR$ such that
\[ u=\sum_{i=1}^k c_iw_i.\]
Then take $d_i \in \mathbb{Q}$ approximating $c_i$ for $i=1,\dots, k$ such that every coefficient of $v:=\sum_i d_iw_i$ is a positive rational number and $||u-v||<\varepsilon$. We also have $Av=0$.

Next, we consider the case that some of $u_i$ equal $0$. We can assume that
\[ u_1,\dots ,u_k>0,\text{ and }u_{k+1}=\cdots =u_n=0.\]
Set $u':= {}^t\!(u_1,\dots ,u_k)$, $A':=[a_{ij}]_{1\leq j\leq k}$. Then $A'u'=0$. By the induction hypothesis, there exists $w= {}^t\!(w_1,\dots ,w_k)\in \RR^k$ such that every $w_i$ is a non-negative rational number, $A'w=0$ and $||u'-w||<\varepsilon$.
Then the vector $v= {}^t\!(w_1,\dots ,w_k, 0,\dots ,0)\in \RR^n$ is a required vector.
\end{proof}

Fix $\mu \in \SC (F)$ and assume that $\mu$ is not the zero measure. Fix $\varepsilon >0$ and $r\in \NN$.
From the above lemma, we can take a map
\[ \theta \: \R_r\rightarrow \mathbb{Z}_{\geq 0}\]
satisfying the following conditions:
\begin{enumerate}
\item $\theta$ is $F$-invariant, that is, for any $T\in \R_r$ and $g\in F$ we have $\theta (T)=\theta (gT)$;
\item there exists $M\in \NN$ for any $T\in \R_r$ we have
\[ \left| \frac{1}{M}\theta(T)-\mu (\SCyl (T) )\right| <\varepsilon ; \]
\item for any $u\in B$ and $J\in \R_r(\id ,u)$ the following equation holds:
\[ \sum_{\substack{T\in \R_r(\id )\\[1pt] T\cap B(\id ,u,r)=J}}\theta (T)=\sum_{\substack{T'\in \R_r(u)\\[1pt] T'\cap B(\id, u,r)=J}}\theta (T').\]
\end{enumerate}
Since $\theta$ is $F$-invariant, for any two adjacent vertices $u,v\in V(X)$ and $J\in \R_r(u,v)$ we have
\[ \sum_{\substack{T\in \R_r(u )\\[1pt] T\cap B(u ,v,r)=J}}\theta (T)=\sum_{\substack{T'\in \R_r(v)\\[1pt] T'\cap B(u, v,r)=J}}\theta (T').\]

The following theorem, which was proved in \cite{Kap17} and named Integral weight realization theorem, is the key for proving the denseness property for $F$. Note that the $\gG$-graph $\gD$ in \cite{Kap17} corresponds to the quotient graph $F\backslash \gG$ for the SC-graph $\gG$ on $F$ in the following theorem.

\begin{theorem}\label{thm:Kapovich's Integral weight realization theorem}
Let $\theta$ be an $F$-invariant map from $\R_r$ to $\mathbb{Z}_{\geq 0}$ satisfying the condition that for any $u\in B$ and $J\in \R_r(\id )$ we have
\[ \sum_{\substack{T\in \R_r(\id )\\[1pt] T\cap B(\id ,u,r)=J}}\theta (T)=\sum_{\substack{T'\in \R_r(u)\\[1pt] T'\cap B(\id, u,r)=J}}\theta (T').\]
Assume that $\theta (T)>0$ for some $T\in \R_r(\id)$.
Then there exists an SC-graph $\gG$ on $F$ such that $\eta_\gG (\SCyl (T))=\theta (T)$ for any $T\in \R_r$.
\end{theorem}
\begin{proof}
We define the vertex set $V(\gG)$ of $\gG$ to be the set
\[ \{ v(g, T, i)\} _{g\in F,\, T\in \R_r(g) ,\, i=1,\dots ,\theta (T)}\]
If $\theta (T)=0$ for $T\in \R_r(g)$, there exists no vertex $v(g, T,i)$.
We regard $v(g,T,i)$ as a copy of $v(g,T,1)$ for $i=2,\dots , \theta (T)$ and we usually write it $v(g,T)$ for short when no confusion can arise.
We define an action of $F$ on $V(\gG )$ by $h v(g,T,i):= v(hg, hT,i)$ for $h\in F $ and $v(g,T,i)\in V(\gG)$.
Note that $\theta(T)=\theta(hT)$ since $\theta$ is $F$-invariant.
A map $\iota$ from $V(\gG)$ to $V(X)=F$ is defined to be the natural projection, that is, $\iota (v(g,T))=g$ for $v(g,T)\in V(\gG)$.

Next, We define the edge set $E(\gG)$ of $\gG$ by connecting two vertices in $V(\gG)$ satisfying certain condition.
Since we require that $F$ acts on $E(\gG)$, we first connect a vertex in $\iota^{-1}(\id )$ to a vertex in $\iota^{-1}( u)$ by an edge for every $u\in B$, and then we copy the edge by using the action of $F$ on $V(\gG)$.
For each $u\in B$ we connect a vertex $v(\id ,T,i)$ to a vertex $v(u, T',i')$ if $T$ and $T'$ is $J$-connectable for some $J\in \R_r(\id, u)$. Since for each $J\in \R_r(\id,u)$ we have
\[ \sum_{\substack{T\in \R_r(\id )\\[1pt] T\cap B(\id ,u,r)=J}}\theta (T)=\sum_{\substack{T'\in \R_r(u)\\[1pt] T'\cap B(\id, u,r)=J}}\theta (T'),\]
the number of vertices $v(\id,T,i)\in \iota^{-1}(\id)$ with $T\cap B(\id, u,r)=J$ equals the number of vertices $v(u, T',i')\in \iota^{-1}(u)$ with $T'\cap B(\id, u ,r )=J$. Hence there exists one-to-one correspondence between $\iota^{-1}(\id )$ and $\iota^{-1}(u)$ satisfying the above condition. Note that
\[ \# \iota^{-1}(\id )=\sum_{T\in \R_r(\id ) }\theta (T) <\infty.\]
From the one-to-one correspondence and the action of $F$ on $V(\gG)$, we obtain the edge set $E(\gG)$.

We see that if $v(\id , T,i)$ is connected to $v(u, T',i')$, then $v(g, T,i )$ is connected to $v(gu,gT',i')$ for every $g\in F$. Moreover, for $v(\id , T ,i )\in V(\gG)$ if $T$ contains $u\in B$, then $J:=T\cap B(\id , u ,r) \in \R_r(\id , u )$ and there exists $v(u, T', i')\in V(\gG )$ such that $T$ and $T'$ are $J$-connectable and $v(\id, T ,i )$ and $v(u, T', i')$ is connected by an edge in $\gG$.
The map $\iota$ sends the edge connecting $v(g, T)$ to $v(gu, T')$ to the edge connecting $g$ to $gu$ for $g\in F$ and $u\in B$.
Then we obtain a graph $(\gG , \iota )$ on $X$ such that $\iota$ is an $F$-equivariant map.
 
Now, we check that $(\gG ,\iota )$ satisfies the condition to be an SC-graph on $F$. It is sufficient to prove that for each connected component $Y$ of $\gG$ the restriction of $\iota $ to $Y$ is injective and $\iota(Y)=CH (\iota(Y)(\infty ))$. Actually, since $\iota$ is locally injective map from the above construction and $X$ is a tree, the restriction of $\iota$ to each connected component is injective.

To see $\iota(Y)=CH (\iota(Y)(\infty ))$, it is sufficient to see that every vertex $v$ of $Y$ has a degree larger than $1$.
Let $v(g,T)$ be a vertex of $Y$. Since there exists $S\in \H (\partial F)$ such that $T=B(g,r)\cap CH(S)$ and $CH(S)\ni g$, the degree of $g$ in $T$ is larger than $1$. 
We prove that $\iota(Y)\cap B(g,1)=T\cap B(g,1)$.
Take a vertex $g'$ of $V(X)$ adjacent to $g$. In the case that $g'\in T$, since $J:=T\cap B(g,g',r)\in \R_r(g,g')$, there exists $v(g',T')\in V(\gG)$ such that $T$ and $T'$ are $J$-connectable and $v(g,T)$ is connected to $v(g',T')$ by an edge. If $g'\not \in T$, then we see that there exists no vertex $v(g',T')\in V(Y)$ connected to $v(g,T)$ by an edge from the construction of $\gG$.
Hence $\iota(Y)\cap B(g,1)=T\cap B(g,1)$.
Therefore the degree of $v(g,T)$ in $Y$ is larger than $1$. Hence $(\gG ,\iota )$ is an SC-graph on $(F,X)$.

Finally, we check that for every $T\in \R_r(\id )$ we have $\eta_{\gG }(\SCyl (T) )=\theta (T)$. From now on, we identify each connected component $Y$ of $\gG$ with $\iota (Y)$. Note that we have $\theta (T)$ copies of $v(\id ,T,1)\in V(\gG)$.
It is sufficient to prove that for $T\in \R_r(\id )$ with $\theta (T)>0$ and for $Y\in \Comp (\gG)$ if $Y$ contains a vertex $v(\id ,T)$, then $Y\cap B(\id ,r)=T$.

Let $T\in \R_r (\id )$ with $\theta (T)>0$ and assume that $v(\id ,T)\in V(Y)$ for $Y\in \Comp (\gG)$.
From the above argument, for any $v(g_1,T_1)\in V(Y)$, there exists $v(g_2,T_2)\in V(Y)$ adjacent to $v(g_1,T_1)$ if and only if $g_1$ and $g_2$ are adjacent vertices of $T_1$. Moreover, we have $Y\cap B(g,1)=T'\cap B(g,1)$ for every $v(g,T')\in V(Y)$.
From Lemma \ref{lem:path connectable on tree}, we can see that for every vertex $v(g, T')\in V(Y)$, $T$ and $T'$ are connectable, that is, $T\cap B(\id, g,r)=T'\cap B(\id ,g,r)$. For each $g\in V(T)\cap B(\id ,r-1)$, by induction on the distance from $\id$ to $g$, we can see that there exists $v(g,T')\in V(Y)$ such that
\begin{align*}
Y\cap B(g,1)=T'\cap B(g,1)=&T'\cap B(\id, g, r)\cap B(g,1)\\
=&T\cap B(\id, g,r )\cap B(g,1)=T\cap B(g,1).
\end{align*}
Therefore $Y\cap B(\id ,r)=T$. This completes the proof.
\end{proof}

By applying Theorem \ref{thm:Kapovich's Integral weight realization theorem} to the map $\theta$ approximating $\mu$, we obtain an SC-graph $(\gG ,\iota )$ on $F$ such that $\eta_\gG (\SCyl (T))=\theta (T)$ for any $T\in \R_r$.
Therefore, for any $T\in \R_r$ we have
\[ \left| \frac{1}{M}\eta_\gG (\SCyl (T) )-\mu (\SCyl (T) )\right| <\varepsilon .\]
This completes the proof of the denseness property of rational subset currents for a free group of finite rank.


\subsection{Approximation by a sequence of subgroups}\label{subsec:approximation by a sequence of sugroups}

In this subsection we assume that the rank of $F$ is $2$ and its free basis $B=\{ x, y\}$ for simplicity.
Every theorem in this subsection can be proved for any free group of finite rank by modifying the definitions and the proofs a little.

\begin{figure}[h]
\begin{center}
\includegraphics[width=9cm]{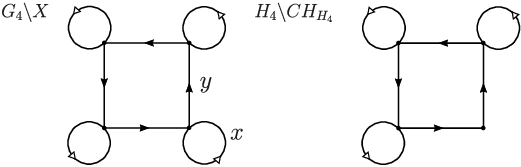}
\vspace{-0.3cm}
\caption{The labels $x,y$ and the orientations are induced from the Cayley graph $X$.}\label{fig:H_4graph}
\end{center}
\end{figure}

For each edge $e$ of the Cayley graph $X=\Cay (F,B)$ we say that $e$ is labeled $\ell \in B$ if $e$ is an edge corresponding to $\ell$ in $X$.
For an integer $n\geq 2$ we consider a normal subgroup $G_n$ of $F$ generated by
\[ \{ x, y^n ,yxy^{-1}, y^2xy^{-2} ,\dots ,y^{n-1}xy^{-n+1} \}. \]
Note that the quotient graph $G_n\backslash X$ is a graph consisting of an $n$-gon each of whose edges are labeled $y$ and each of whose vertices is attached a loop labeled $x$ (see Figure \ref{fig:H_4graph}). The subgroup $G_n$ is an $n$-index subgroup of $F$.
Recall that we have a continuous $\RRR$-linear map $\iota_{G_n}$ from $\SC (G_n)$ to $\SC (F)$ (see Section \ref{sec:relation between subgroups}).
Since $G_n$ is a finite index subgroup of $F$, the map $\iota_{G_n}$ is surjective.

Let $H_n$ be a subgroup of $G_n$ generated by
\[ \{ y^n ,yxy^{-1}, y^2xy^{-2} ,\dots ,y^{n-1}xy^{-n+1} \}. \]
We can obtain the quotient graph $H_n\backslash CH_{H_n}$ by removing certain one loop labeled $x$ in $G_n\backslash X$.
From Proposition \ref{prop:convergence and subset cylinder} and the argument following Corollary \ref{cor: Kapovich's open neighborhood of a subset current}, we can see that the sequence
\[ \left\{ \frac{1}{n}\ \eta_{H_n}\right\}_{n\geq 2}\]
of rational subset currents converges to $\eta_F=\delta_{\partial F}$.
From this property, we can guess that $H_n$ approximates $F$ in some sense.
Actually, we will prove the following theorem in this subsection.

\begin{theorem}\label{thm:approximating theorem for free group}
The union
\[ \bigcup_{n=2}^\infty \iota_{H_n}(\SC_r (H_n) )\]
is a dense subset of $\SC (F)$.
\end{theorem}

\begin{remark}\label{rem:idea of approximating by a sequence of subgroups}
Recall that the map $\iota_{G_n}$ is a surjective continuous $\RRR$-linear map from $\SC (G_n)$ to $\SC (F)$ since $G_n$ is a finite index subgroup of $F$. 
Moreover, we also have
\[ \iota_{H_n} (\SC (H_n) )=\iota_{G_n}\circ \iota_{H_n}^{G_n} (\SC (H_n)) \subset \iota_{G_n}(\SC (G_n ))=\SC (F).\]
Roughly speaking, Theorem \ref{thm:approximating theorem for free group} follows since the ``difference'' between $H_n$ and $G_n$ becomes smaller as $n$ tends to $\infty$. We also remark that for $\mu \in \SC(G)\subset \SC(G_n)$ we have
$\iota_{G_n}(\mu )=n\mu $.

As a corollary to Theorem \ref{rem:idea of approximating by a sequence of subgroups}, we see that $\SC_r(F)$ is a dense subset of $\SC (F)$ since $\iota_{H_n}(\SC_r(H_n))$ is included in $\SC_r(F)$ for every $n\geq 2$.
\end{remark}

Our method of proving Theorem \ref{thm:approximating theorem for free group} is as follows:
Let $\mu \in \SC (F)$. Fix $\varepsilon >0$ and $r\in \NN$, which determine the open neighborhood
\[ \{ \nu\in \SC (F)\mid |\mu (\SCyl (T))-\nu (\SCyl (T)) |<\varepsilon \text{ for every }T\in \R_r(\id )\}\]
of $\mu$. Then for a sufficiently large $n$, we construct a subset current $\nu \in \mathrm{Span}(\SC_r(H_n))$ such that $\iota_{H_n}(\nu )$ belongs to the above neighborhood.
Note that every $H_n$ is a free group of finite rank, and for a free group $H$ of finite rank $\SC_r(H)$ is a dense subset of $\mathrm{Span}(\SC_r(H))$.
During the proof, we do not use the fact that a free group of finite rank has the denseness property of rational subset currents.

In order to obtain $\nu$ we will construct an SC-graph $(Y,f)$ on $(H_n, CH_{H_n})$, which means that $(Y,f)$ satisfies the following conditions:
\begin{enumerate}
\item $f$ is an $H_n$-equivariant graph morphism from $Y$ to $CH_{H_n}$;
\item the restriction of $f$ to each connected component $Y_0$ of $Y$ is injective and the image $f(Y_0)$ coincides with $CH (f(Y_0)(\infty ))$;
\item $\# f^{-1}( \id )<\infty $.
\end{enumerate}
Then we can obtain a subset current $\eta_{Y}\in \mathrm{Span}(\SC_r(H_n))$ by
\[ \eta_Y:=\sum_{Y_0\in \Comp (Y)} \delta_{f(Y_0)(\infty)}.\]
Note that we often identify $Y_0$ with $f(Y_0)$.

\begin{remark}
Theorem \ref{thm:approximating theorem for free group} gives us a new idea to construct an approximating rational subset current. Explicitly, for an infinite hyperbolic group $G$ if we have a sequence $\{ H_n \}$ of quasi-convex subgroups of $G$ such that $a_n\eta_{H_n}$ converges to $\eta_G$ for a sequence $\{ a_n \}$ of $\RRR$ and $H_n$ is a free group of finite rank, then for any $\mu \in \SC (G)$ we may be able to construct $\nu \in \mathrm{Span}({\SC_r (H_n)})$ such that $\iota_{H_n}(\nu)$ approximates $\mu$ for a sufficiently large $n\in \NN$.
In the case that $G$ is a surface group, we will prove the denseness property for $G$ by using this idea (see Theorem \ref{thm:approximation of surface group by free groups}).

Recall that when we proved the denseness property of rational subset currents for $F$, we used the Cayley graph $X$ of $F$, which is a tree, and we constructed an SC-graph on $X$. However, for the Cayley graph of a general hyperbolic group it is much more difficult to construct a subgraph like an SC-graph on the Cayley graph. The difficulty comes from that the Cayley graph of a hyperbolic group is a $\delta$-hyperbolic space for $\delta>0$ and almost everything on a $\delta$-hyperbolic space is vaguely determined in some sense. For example, a geodesic line connecting two points of the boundary is not unique but unique up to some constants.
During the proof of Theorem \ref{thm:approximation of surface group by free groups} we have to prove a lot of lemmas corresponding to such constants. However, the basic idea of the proof of Theorem \ref{thm:approximation of surface group by free groups} is the same as that of Theorem \ref{thm:approximating theorem for free group}.
\end{remark}

Now, we consider the action of $H_n$ on $CH_{H_n}$.
Note that
\[F/G_n=\{ G_n ,yG_n, \dots ,y^{n-1}G_n\} \]
and the finite union of the closed balls
\[ \bigcup_{i=0}^{n-1} B(y^i ,1/2) \]
is a fundamental domain for the action of $G$ on $X$.
We see that
\[ \F :=CH_{H_n}\cap \bigcup_{i=0}^{n-1} B(y^i ,1/2) .\]
is a fundamental domain for the action of $H_n$ on $CH_{H_n}$ and for any non-trivial $h\in H_n$ the intersection of $h\F\cap \F$ is empty or a midpoint of an edge.
Note that
\[ \bigcup_{i=0}^{n-1} B(y^i ,1/2) \setminus \F\]
equals a half-edge labeled $x$ attached to $\id$ since the canonical projection $p_{H_n}$ from $CH_{H_n}$ onto $H_n\backslash CH_{H_n}$ maps $\id$ to the vertex of $H_n\backslash CH_{H_n}$ that is not attached a loop labeled $x$.

Then we see that the set
\[ \{ h\in H_n\setminus \{ \id \} \mid h\F \cap \F\not=\emptyset \}\]
is a generating set of $H_n$. We can take its subset $B_n$ such that
\[ \{ h\in H_n\setminus \{ \id \} \mid h\F \cap \F\not=\emptyset \} =B_n \sqcup B_n^{-1},\]
and then $B_n$ is a free basis of $H_n$.

Now, we consider the Cayley graph $X_n:=\Cay (H_n,B_n)$ of $H_n$ with respect to $B_n$, which is a tree. From the definition of $B_n$ there is one-to-one correspondence between a vertex $h$ of $X_n$ and $h\F \subset X$. Moreover, $h_1,h_2,\in V(X_n)$ are adjacent if and only if $h_1\F\cap h_2\F\not=\emptyset$, which means that $h_1\F$ and $h_2\F$ are also adjacent.

We generalize the notion of a round-graph centered at a vertex with radius $r\in \NN$ and define a round-graph of $r$-neighborhood of a subset of $X$ in order to consider a round-graph of $r$-neighborhood of $h\F$ for $h\in H_n$.

\begin{definition}[Round-graph of $r$-neighborhood]\label{def:round-graph generalized version}
Let $r>0$. For a non-empty bounded subset $Y$ of $X$ we denote by $B(Y,r)$ the $r$-neighborhood of $Y$, explicitly,
\[ B(Y,r):= \{ x\in X\mid d (x, Y)\leq r\} .\]
A subset $T$ of $B(Y,r)$ is called a \ti{round-graph} of $r$-neighborhood of $Y$ if $T\cap Y\not =\emptyset$ and there exists $S\in \H (\partial F)$ such that
\[ T=CH(S)\cap B(Y,r). \]
We denote by $\R_r(Y)$ the set of all round-graphs of $r$-neighborhood of $Y$.
For $T\in \R_r(Y)$ we define the \ti{subset cylinder} $\SCyl (T)$ with respect to $T$ by
\[ \SCyl (T):=\{ S\in \H (\partial F) \mid CH(S)\cap B(Y,r)=T\}.\]

For two non-empty bounded subset $Y$ and $Z$ of $X$ we denote by $B(Y,Z,r)$ the intersection of $B(Y,r)$ and $B(Z,r)$.
For $T_1\in \R_r(Y), T_2\in \R_r(Z)$ we say that $T_1$ and $T_2$ are \ti{connectable} if $T_1\cap B(Y,Z,r)=T_2\cap B(Y,Z,r)$.
Note that $B(Y,Z,r)$ can be empty, and then $T_1$ and $T_2$ are connectable for any $T_1\in \R_r(Y), T_2\in \R_r(Z)$.

Assume that $B(Y,Z,r)$ is not empty.
A subset $J$ of $B(Y,Z,r)$ is called a $(Y,Z)$-round-graph of $r$-neighborhood of $Y,Z$ if $J\cap Y\not=\emptyset ,J\cap Z\not=\emptyset$ and there exists $S\in \H (\partial F)$ such that $J=CH(S)\cap B(Y,Z,r)$.
We denote by $\R_r(Y,Z)$ the set of all $(Y,Z)$-round-graph of $r$-neighborhood of $Y,Z$.
For $J\in \R_r(Y,Z)$ we define the subset cylinder $\SCyl (J)$ with respect to $J$ by
\[ \SCyl (J):=\{ S\in \H (\partial F)\mid CH(S)\cap B(Y,Z,r)=J\}.\]
For $T_1\in \R_r(Y), T_2\in \R_r(Z)$ we say that $T_1$ and $T_2$ are $J$-\ti{connectable} for $J\in \R_r(Y,Z)$ if $T_1\cap B(Y,Z,r)=J=T_2\cap B(Y,Z,r)$.
\end{definition}

As in the case of round-graphs centered at points, for $J\in \R_r(Y,Z)$ the following equation holds:
\[ \SCyl (J)=\bigsqcup_{\substack{T\in \R_r(Y)\\[1pt] T\cap B(Y,Z,r)=J}} \SCyl (T),\]
which implies that for any $\mu \in \SC (F)$ we have
\[ \sum_{\substack{T\in \R_r(Y)\\[1pt] T\cap B(Y,Z,r)=J}}\mu (\SCyl (T))=\sum_{\substack{T'\in \R_r(Z)\\[1pt] T'\cap B(Y,Z,r)=J}}\mu (\SCyl (T')).\]

\begin{lemma}\label{lem:path connectable on tree for fundamental domain}
Let $r>0$. Let $h_0, h_1,\dots ,h_m$ be pairwise disjoint elements of $H_n$ such that $h_{i-1}\F$ is adjacent to $h_{i}\F$ for $i=1,\dots ,m$.
Take $T_i\in \R_r(h_i\F)$ for $i=0,1,\dots ,m$. If $T_{i-1}$ and $T_{i}$ are connectable for every $i=1,\dots m$, then $T_0$ and $T_m$ are connectable.
\end{lemma}
\begin{proof}
The proof is almost the same as that of Lemma \ref{lem:path connectable on tree}.
Since $X$ is a tree, we have
\[ B(h_0\F ,h_m\F ,r)\subset \bigcap_{i=0}^m B(h_i\F ,r),\]
which implies 
\[ B(h_0\F ,h_m\F ,r)\subset \bigcap_{i=1}^{m} B(h_{i-1}\F ,h_i\F ,r).\]
From the assumption,
\[ T_{i-1}\cap B(h_{i-1}\F ,h_i\F ,r)=T_i\cap B(h_{i-1}\F ,h_i\F ,r).\]
for every $i=1,\dots ,m$. Therefore
\[ T_0\cap B(h_0\F ,h_m\F ,r)=T_1\cap B(h_0\F ,h_m\F ,r)=\cdots =T_m\cap B(h_0\F ,h_m\F ,r).\]
This completes the proof.
\end{proof}

Now, we begin to prove Theorem \ref{thm:approximating theorem for free group}.
We divides the proof into 5 steps and we will also follow these steps in the proof of the denseness property for surface groups.

\begin{step}\label{step:1}
Fix $\mu \in \SC (F)$. Fix $\varepsilon >0$ and $r \in \NN$, which determine the open neighborhood $U(r,\varepsilon )$ of $\mu$:
\[ U(r,\varepsilon ):= \{ \nu\in \SC (G)\mid |\mu (\SCyl (T))-\nu (\SCyl (T)) |<\varepsilon \text{ for every }T\in \R_r(\id )\} .\]
Take a sufficiently large $n\in \NN$.
Set $\rho := r+n$.
Recall that
\[ \R_\rho=\bigsqcup_{v\in V(X)} \R_\rho (v).\]
From Lemma \ref{lem:approximate by positive rational numbers}, we can take a map
\[ \theta \: \R_\rho \rightarrow \mathbb{Z}_{\geq 0} \]
satisfying the following conditions:
\begin{enumerate}
\item $\theta$ is $F$-invariant;
\item there exists $M\in \NN$ such that $\frac{1}{M}\theta $ approximates $\mu$, that is, for every $T\in \R_\rho$ 
\[ \left| \frac{1}{M}\theta(T) -\mu (\SCyl (T) )\right| <\varepsilon' ,\]
where $\varepsilon'>0$ depends on $\mu , \varepsilon, r$ and $n$;
\item for any $u\in B=\{ x,y \}$ and any $J \in \R_\rho (\id ,u)$ we have
\[ \sum_{\substack{T\in \R_\rho(\id )\\[1pt] T\cap B(\id ,u,\rho )=J}}\theta (T)=\sum_{\substack{T'\in \R_\rho (u)\\[1pt] T'\cap B(\id, u,\rho )=J}}\theta (T').\]
\end{enumerate}
\end{step}

The above conditions $(1)$ and $(3)$ imply that for any adjacent vertices $v,w \in V(X)$ and $J\in \R_\rho(v,w)$ we have
\[ \sum_{\substack{T\in \R_\rho(v )\\[1pt] T\cap B(v ,w,\rho )=J}}\theta (T)=\sum_{\substack{T'\in \R_\rho (w)\\[1pt] T'\cap B(v, w,\rho )=J}}\theta (T').\]

For each $h\in H_n$ and $T\in \R_r(h\F )$ we define $\theta (T)$ by
\[ \theta (T) :=\sum_{\substack{T'\in \R_\rho (v)\\[1pt] T'\cap B(h\F ,r)=T}} \theta (T' ),\]
where $v$ is a vertex of $T\cap h\F$. Note that the diameter of $\F$ is $n$ and so $B(h\F ,r)$ is included in $B(v, r+n)=B(v,\rho )$ for any vertex $v$ of $T\cap h\F$.

\begin{lemma}\label{lem:restriction to H_n fundamental domain}
The definition of $\theta (T)$ is independent of the choice of $v$ and the map
\[ \theta \: \bigsqcup_{h\in H_n}\R_r(h \F)\rightarrow \mathbb{Z}_{\geq 0}\]
is $H_n$-invariant.
Moreover, for any $u\in B_n$ and any $J\in \R_r(\F ,u\F)$ we have the following equation:
\[ \sum_{\substack{T\in \R_r(\F )\\[1pt] T\cap B(\F, u\F ,r)=J}}\theta (T)=\sum_{\substack{T'\in \R_r(u\F )\\[1pt] T'\cap B(\F ,u\F ,r)=J}}\theta (T').\]
\end{lemma}
\begin{proof}
For $h\in H_n$ and $T\in \R_r(h\F )$ we have
\[ \SCyl (T)=\bigsqcup_{\substack{T'\in \R_\rho (v)\\[1pt] T'\cap B(h\F,r) =T}} \SCyl (T'),\]
Let $v'$ be another vertex of $T\cap h\F$. It is sufficient to consider the case that $v'$ is adjacent to $v$.
Since $B(h\F ,r)\subset B(v,v',\rho )$, we can obtain
\begin{align*}
\sum_{\substack{T'\in \R_\rho (v)\\[1pt] T'\cap B(h\F ,r)=T}} \theta (T' )
&=\sum_{\substack{J\in \R_\rho (v,v')\\[1pt] J\cap B(h\F ,r)=T}} 
	\sum_{\substack{T'\in \R_\rho (v)\\[1pt] T'\cap B(v,v',\rho )=J}} \theta (T' )\\
&=\sum_{\substack{J\in \R_\rho (v,v')\\[1pt] J\cap B(h\F ,r)=T}} 
	\sum_{\substack{T'\in \R_\rho (v')\\[1pt] T'\cap B(v,v',\rho )=J}} \theta (T' )\\
&=\sum_{\substack{T'\in \R_\rho (v')\\[1pt] T'\cap B(h\F ,r)=T}} \theta (T' ).
\end{align*}
Therefore, $\theta (T)$ is independent of the choice of $v$. Moreover, the map
\[ \theta \: \bigsqcup_{h\in H_n}\R_r(h \F)\rightarrow \mathbb{Z}_{\geq 0}\]
is $H_n$-invariant by the definition.

Let $u\in B_n$ and $J\in \R_r(\F ,u\F)$. Since $u\F$ and $\F$ intersect at a midpoint of an edge, there exist two adjacent vertices $v,w$ of $J$ such that $v\in \F $ and $w\in u\F$. Then we have
\begin{align*}
\sum_{\substack{T\in \R_r(\F )\\[1pt] T\cap B(\F, u\F ,r)=J}}\theta (T)
=& \sum_{\substack{T\in \R_r(\F )\\[1pt] T\cap B(\F, u\F ,r)=J}}
		\sum_{\substack{T'\in \R_\rho (v)\\[1pt] T'\cap B(\F ,r)=T}} \theta (T' )\\
=& \sum_{\substack{T\in \R_r(\F )\\[1pt] T\cap B(\F, u\F ,r)=J}}
			\sum_{\substack{J'\in \R_\rho (v,w)\\[1pt] J'\cap B(\F, u\F ,r)=T}} 
				\sum_{\substack{T'\in \R_\rho (v)\\[1pt] T'\cap B(v,w ,\rho )=J'}} \theta (T' )\\
=& \sum_{\substack{T\in \R_r(\F )\\[1pt] T\cap B(\F, u\F ,r)=J}}
			\sum_{\substack{J'\in \R_\rho (v,w)\\[1pt] J'\cap B(\F, u\F ,r)=T}} 
				\sum_{\substack{T'\in \R_\rho (w)\\[1pt] T'\cap B(v,w, \rho )=J'}} \theta (T' )\\			
=& \sum_{\substack{T\in \R_r(\F )\\[1pt] T\cap B(\F, u\F ,r)=J}}
				\sum_{\substack{T'\in \R_\rho (w)\\[1pt] T'\cap B(\F ,r)=T}} \theta (T' )\\
=& \sum_{\substack{T'\in \R_\rho (w)\\[1pt] T'\cap B(\F, u\F ,r)=J}} \theta (T' )\\
=& \sum_{\substack{T'\in \R_r(u\F )\\[1pt] T'\cap B(\F, u\F ,r)=J}}
		\sum_{\substack{T''\in \R_\rho (w)\\[1pt] T''\cap B(u\F ,r)=T'}} \theta (T'' )\\	
=&\sum_{\substack{T'\in \R_r(u\F )\\[1pt] T'\cap B(\F ,u\F ,r)=J}}\theta (T').
\end{align*}
This is the required equation.
\end{proof}

From the above lemma and its proof we can see that $\theta $ can be considered as a measure as long as we consider a value of ``small'' round-graphs by $\theta$. Explicitly, for any subset $Y$ of $X$ and $\ell \in \NN$ satisfying the condition that $B(Y,\ell )\subset B(v ,\rho )$ for a vertex $v\in Y$, we can define
\[ \theta(T) :=\sum_{\substack{T'\in \R_\rho (v)\\[1pt] T'\cap B(Y ,\ell )=T}} \theta (T' )\]
for any $T\in \R_\ell (Y)$. Moreover, if $\frac{1}{M}\theta (T')$ is sufficiently close to $\mu(\SCyl (T'))$ for every $T'\in \R_\rho (\id )$, then $\frac{1}{M} \theta( T)$ is also close to $\mu (\SCyl (T))$ for $T\in \R_\ell (Y)$.

\begin{step}\label{step:2}
By using the map $\theta \: \bigsqcup \R_r (h\F) \rightarrow \mathbb{Z}_{\geq 0}$, we construct a graph $(\gG,\iota )$ on $X_n=\Cay (H_n ,B_n)$ by the same way as we did in the proof of Theorem \ref{thm:Kapovich's Integral weight realization theorem}. Then the graph $(\gG ,\iota )$ has the following properties:
\begin{enumerate}
\item $\iota \: \gG\rightarrow X_n$ is an $H_n$-equivariant map;
\item the restriction of $\iota $ to each connected component of $\gG$ is injective;
\item $\# \iota ^{-1}( \id )<\infty $.
\end{enumerate}
Note that the property (2) in the above is not the same as Condition (SC2) in the definition of an SC-graph on $(X_n ,H_n)$.
\end{step}

We define the vertex set $V(\gG )$ of $\gG$ by
\[ V(\gG ):=\{ v(h ,T ,i)\} _{h\in H_n,\, T\in \R_r (h\F ),\, i=1,\dots ,\theta (T) }.\]
If $\theta(T)=0$ for $T\in \R_r (h\F)$, then we do not have any vertex $v(h,T, i)$. We will write $v(h,T)$ in place of $v(h,T,i )$ when no confusion can arise.
Since for each $u\in B_n$ and $J\in \R_r(\F ,u\F)$ we have 
\[\sum_{\substack{T\in \R_r(\F )\\[1pt] T\cap B(\F, u\F ,r)=J}}\theta (T)=\sum_{\substack{T'\in \R_r(u\F )\\[1pt] T'\cap B(\F ,u\F ,r)=J}}\theta (T'),\]
we can define a one-to-one correspondence between 
\[ \{ v(\id, T )\}_{T\cap B(\F,u\F,r)=J} \text{ and } \{ v(u ,T')\}_{T'\cap B(\F ,u\F ,r)=J}.\]
For this correspondence we connect two vertices by an edge, and we perform this operation for every $u\in B_n$ and every $J\in \R_r(\F, u\F)$. Finally, for every $u\in B_n$ and $h\in H_n$ we connect $v(h,T,i)$ to $v(hu,hT',i')$ by an edge if $v(\id ,T,i)$ and $v(u,T,i')$ are connected by an edge.
In this way we obtain the edge set $E(\gG )$ of $\gG$.

From the construction of $\gG$, we see that $H_n$ acts on $\gG$, and if $v(h_1,T_1) , v(h_2,T_2)\in V(\gG)$ are connected by an edge, then $h_1$ and $h_2$ are adjacent in $X_n$ and $T_1,T_2$ are $J$-connectable for some $J\in \R_r(h_1\F,h_2\F)$. Moreover, for $v(h,T)\in V(\gG)$ if there exists $h'$ adjacent to $h$ in $X_n$ such that $T\cap h'\F\not=\emptyset$, then $T\cap B(h\F,h'\F,r)\in \R_r(h\F, h'\F)$ and there exists $T'\in \R_r(h'\F)$ such that $v(h,T)$ and $v(h',T')$ are connected by an edge.

We also have an $H_n$-equivariant map $\iota$ from $\gG$ to $X_n$ such that $\iota (v(h,T))=h$ for $v(h,T)\in V(\gG)$.
Moreover, the restriction of $\iota$ to each connected component $Y$ of $\gG$ is injective since $X_n$ is a tree. By the definition of $\iota$, we have
\begin{align*}
\#\iota^{-1}(\id )=&\sum_{T\in \R_r (\F )}\theta (T)\\
	\leq &\sum_{v\in V(\F )} \sum_{T\in \R_\rho (v)}\theta (T)=\# V(\F )\sum_{T\in \R_\rho (\id )}\theta (T)<\infty.
\end{align*}

Finally, we note that a connected component $Y$ of $\gG$ may contain a vertex with degree $0$ or $1$ since edges with label $x$ are not attached to the vertex $h\in H_n \subset V(CH_{H_n})$. For example, consider the subgroup $\langle x \rangle $ of $F$ and its limit set $\{ x^{\infty }, x^{-\infty} \}$. Then $T:=CH(\{ x^\infty ,x^{-\infty}\})\cap B(\F ,r) \in \R_r(\F)$ and $T\cap CH_{H_n}=\{ \id \}$.
We see that $v(\id , T)$ will be a vertex with degree $0$ in $\gG$ if $\theta (T)>0$. Therefore even if we define a subset current $\eta_\gG$ on $H_n$ by the same way as we did for an SC-graph on $F$, $\eta_\gG$ loses some information on $\theta (T)$.

\begin{step}\label{step:3}
We construct a graph $(|\gG |, |\iota |)$ on $X$ from $(\gG ,\iota )$
satisfying the following conditions:
\begin{enumerate}
\item $|\iota |\: |\gG|\rightarrow X$ is an $H_n$-equivariant map;
\item the restriction of $|\iota |$ to each connected component of $|\gG|$ is injective;
\item $\# |\iota |^{-1}( \id )<\infty $.
\end{enumerate}
\end{step}

For each connected component $Y$ of $\gG$ we define a subgraph $|Y|$ of $X$ by
\[ |Y|:=\bigcup_{v(h,T)\in V(Y)}T\cap h\F \]
and define $|\gG|$ to be the disjoint union of $|Y|$ over all connected component $Y$ of $\gG$. By the definition, $|Y|$ is included in $\bigsqcup_{h\in H_n }h\F =CH_{H_n}$.
Let $|\iota |$ be the natural projection from $|\gG|$ to $X$, that is, the restriction of $|\iota |$ to $|Y|$ for each $Y\in \Comp (\gG)$ is the inclusion map.

The action of $H_n$ on $\gG$ and on $X$ induce the action of $H_n$ on $|\gG|$ as follows. Let $h_0\in H_n$ and $x\in |Y|$ for $Y\in \Comp (\gG)$. For a moment, we write $(Y,x )$ in place of $x$ to emphasize that $x$ is a point of $|Y|$. Then there exists $v(h,T)\in V(Y)$ such that $x\in T\cap h\F$. Since $H_n$ acts on $\gG$, there exists $v(h_0h,h_0T)\in V(h_0 Y)$ and $h_0x \in h_0T\cap h_0h\F$.
Then we define $h_0(Y,x)$ to be $(h_0Y ,h_0x)$, which is a point of $|h_0Y|$.
We see that the map $|\iota |$ is $H_n$-equivariant by the definition.

\begin{lemma}\label{lem:neighborhood of realization of Gamma}
Let $Y$ be a connected component of $\gG$. Let $v(h,T)\in V(Y)$, $v\in V(T)\cap h\F$. Then we have
\[ |Y|\cap B(v,r)=CH_{H_n}\cap T\cap B(v,r).\]
Moreover, $|Y|$ is connected.
\end{lemma}
\begin{proof}
\underline{Inclusion $\subset$ :} Let $w\in |Y|\cap B(v,r)$. There exists $v(h',T')\in V(Y)$ such that $w\in T'\cap h'\F$.
Since $Y$ is connected there exist a geodesic path of vertices $v(h_0,T_0)=v(h,T),v(h_1,T_1),\dots ,v(h_m,T_m)=v(h',T')\in V(Y)$.
Since $T_{i-1}$ and $T_i$ are connectable for $i=1,\dots ,m$, $T$ and $T'$ are also connectable from Lemma \ref{lem:path connectable on tree for fundamental domain}. Since $B(v,r)\subset B(h\F ,r)$, we have
\begin{align*}
T'\cap h'\F \cap B(v,r)=&T'\cap B(h\F ,h'\F ,r)\cap h'\F \cap B(v,r)\\ 
=&T\cap B(h\F ,h'\F ,r)\cap h'\F \cap B(v,r)
\end{align*}
and so
\[ w\in T'\cap h'\F \cap B(v,r)\subset CH_{H_n}\cap T\cap B(v,r).\]

\underline{Inclusion $\supset$ :} Let $w\in CH_{H_n}\cap T\cap B(v,r)$. Then there exists a geodesic path $P$ from $v$ to $w$ in $CH_{H_n}\cap T\cap B(v,r)$ since all of $CH_{H_n}, T$ and $B(v,r)$ are trees. We can take a sequence $h_0=h, h_1,\dots ,h_m\in H_n$ such that $P$ passes through $h_i\F$ in this order and $w\in h_m\F$. From the construction of the edge set of $\gG$ there exists $T_i\in \R_r(h_i\F)$ for $i=1,\dots ,m$ such that $v(h_1,T_1),\dots , v(h_m ,T_m)\in V(Y)$, $T$ and $T_1$ are connectable, and $T_i$ and $T_{i+1}$ are connectable for $i=1,2,\dots ,m-1$.
From Lemma \ref{lem:path connectable on tree for fundamental domain}, $T$ and $T_m$ are connectable and so
\begin{align*}
w\in &T \cap B(h\F ,h_m\F ,r) \cap h_m\F \cap B(v,r)\\
=&T_m\cap B(h\F ,h_m\F ,r)\cap h_m\F \cap B(v,r).
\end{align*}
This implies that $w\in T_m\cap h_m\F \cap B(v,r)\subset |Y|\cap B(v,r)$.

Finally, we check that $|Y|$ is connected. Take any geodesic path of vertices \[v(h_0,T_0),v(h_1,T_1),\dots ,v(h_m,T_m)\in V(Y).\] 
Since $T_{i-1}$ and $T_i$ are connectable, there exists an edge $e_i$ in $T_{i-1}\cap T_i$ connecting $h_{i-1}\F$ and $h_i\F$ for $i=1,\dots ,m$. Note that $T_i\cap h_i\F $ is connected. Therefore $|Y|$ is connected.
\end{proof}

From the above lemma we have
\begin{align*}
\# |\iota |^{-1}(\id )
=&\# \{ Z\in \Comp (|\gG |)\mid Z\ni \id \} \\
=&\# \{ Y\in \Comp (\gG )\mid v(\id ,T)\in V(Y), T\ni \id \} \\
= &\sum_{\substack{T\in \R_{r}(\F )\\[1pt] T\ni \id}} \theta (T) \\
= &\sum_{\substack{T\in \R_{\rho }(\id )\\[1pt] T\ni \id}} \theta (T) <\infty .
\end{align*}
Therefore we can see that $(|\gG |, |\iota |)$ satisfies all conditions to be an SC-graph on $(H_n, CH_{H_n})$ except the condition that for every connected component $Z$ of $|\gG|$ we have $CH(Z(\infty ))=Z$. The reason is that there exists a vertex $v$ of $Z$ with degree $1$ or $0$ in $Z$. Such a vertex $v$ belongs to $H_n\subset V(CH_{H_n})$ by the construction of $H_n$. This implies that there are finite vertices of $|\iota|^{-1}(\id)$ with degree less than $2$ and any vertex of $|\gG|$ with degree less than $2$ belongs to the $H_n$-orbit of those vertices.

\begin{step}\label{step:4}
We construct an SC-graph $(\hat{\gG},\hat{\iota})$ on $(H_n,CH_{H_n})$ from $(|\gG |,|\iota |)$, that is, $(\hat{\gG},\hat{\iota})$ satisfies the following conditions:
\begin{enumerate}
\item $\hat{\iota} \: \hat{\gG} \rightarrow CH_{H_n}$ is an $H_n$-equivariant map;
\item the restriction of $\hat{\iota}$ to each connected component $Z$ of $\hat{\gG}$ is injective and $\hat{\iota} (Z)=CH (\hat{\iota}(Z)(\infty ))$;
\item $\# \hat{\iota}^{-1}( \id )<\infty $.
\end{enumerate}
\end{step}

Let $v$ be a vertex of $|\iota |^{-1}(\id)$ with degree less than $2$. If the degree of $v$ is 0, then we remove $H_n(v)$ from $|\gG|$.
Now, we consider the case the degree of $v$ is $1$.
Then there exists either an edge connecting a vertex of $|\iota|^{-1}(y)$ to $v$ or an edge connecting a vertex of $|\iota |^{-1}(y^{-1})$ and $v$.
Assume that there exists an edge connecting a vertex of $|\iota|^{-1}(y^{-1})$ to $v$.
Take a subgraph $P$ of $CH_{H_n}$ consisting of two edges connecting $\id$ and $y$, $y$ and $yx$.
Consider the disjoint union
\[ |\gG |\sqcup \bigsqcup_{h\in H_n}h(P)\]
Note that $H_n$ acts on this disjoint union from left.

\begin{figure}[h]
\begin{center}
\includegraphics[width=9cm]{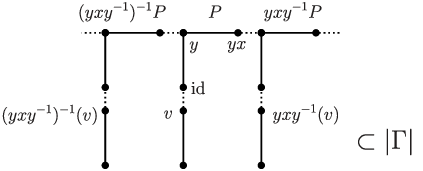}
\vspace{-0.3cm}
\caption{The dotted lines mean the attachment of $H_n(P)$ to $|\gG |$.}\label{fig:attach_P}
\end{center}
\end{figure}

For every $h\in H_n$ we attach the vertex of $hP$ corresponding to $h$ to $hv$ in $|\gG|$, and the vertex of $hP$ corresponding to $h(yx)$ to the vertex of $(hyxy^{-1})P$ corresponding to $(hyxy^{-1})y=hyx$ (see Figure \ref{fig:attach_P}).
Note that if $h\in H_n$, then $hyxy^{-1}\in H_n$. Since this attachment of $H_n(P)$ to $|\gG|$ is $H_n$-invariant, we obtain a graph $|\gG|'$ such that $H_n$ acts on $|\gG|'$ and the degree of $hv$ in $|\gG|'$ equals $2$. For $h\in H_n$ and the vertex $h(y) \in h(P)$ the degree of $h(y)$ in $|\gG |'$ is $3$.
Moreover $|\iota|$ is extended to an $H_n$-equivariant map $|\iota|'$ from $|\gG |'$ to $CH_{H_n}$ such that the restriction of $|\iota |'$ to every connected component of $| \gG |'$ is injective since $CH_{H_n}$ is a tree.

In the case that there exists an edge connecting a vertex of $|\iota|^{-1}(y)$ to $v$, we can perform the same operation by using a subgraph of $CH_{H_n}$ consisting of two edges connecting $\id$ and $y^{-1}$, $y^{-1}$ and $y^{-1}x$.

We perform this operation until every vertex of $|\gG|$ has a degree larger than or equal to $2$.
The resulting graph, which we denoted by $(\hat{\gG}, \hat{\iota})$, is an SC-graph on $(H_n ,CH_{H_n})$.
Let $C_1$ be the number of vertices of $|\iota|^{-1}(\id )$ with degree $1$ in $|\gG|$.
Then we need to perform the above operation exactly $C_1$ times in order to obtain $\hat{\gG}$.

\begin{step}\label{step:5}
Set 
\[ \eta_{\hat{\gG}}:=\sum_{Z\in \Comp (\hat{\gG})} \delta_{Z(\infty )}\in \mathrm{Span}(\SC_r (H_n)).\]
We prove that $\frac{1}{nM}\iota_{H_n}(\eta_{\hat{\gG}} )$ belongs to the open neighborhood $U(r,\varepsilon )$ of $\mu$ for a sufficiently large $n$.
\end{step}

Take $T\in \R_r(\id )$.
Then we have
\begin{align*}
\iota_{H_n}(\eta_{\hat{\gG}} )(\SCyl (T))=&\sum _{gH_n\in F/H_n}g_{\ast}(\eta_{\hat{\gG}}) (\SCyl (T) )\\
=&\sum _{gH_n\in F/H_n}\eta_{\hat{\gG}} (\SCyl (g^{-1}T) ).
\end{align*}
If $g^{-1}T$ is not included in $CH_{H_n}$ for $g\in F$, then $\SCyl (T)\cap \H (\partial \gL (H_n))=\emptyset$ and so $\eta_{\hat{\gG}} (\SCyl (g^{-1}T) )=0$.
Since the fundamental domain $\F$ for the action of $H_n$ on $CH_{H_n}$ includes the vertices $\id ,y ,\dots ,y^{n-1}$, we have
\[ V(CH_{H_n})=H_n\sqcup H_n y\sqcup \cdots \sqcup H_ny^{n-1}.\]
This implies that if $gH_n\in G/H_n$ is different from every $y^{-i}H_n$ for $i=0,1,\dots ,n-1$, then $g^{-1}\not \in V(CH_{H_n})$, which implies that $\eta_{\hat{\gG}} (\SCyl (g^{-1}T) )=0$. Note that $T\ni \id$ and so $g^{-1}T\ni g^{-1}$.
Therefore we have
\[ \iota_{H_n}(\eta_{\hat{\gG}} )(\SCyl (T))=\sum_{i=0}^{n-1}\eta_{\hat{\gG}} (\SCyl (y^i T) ).\]

Now, we can assume that $n$ is much larger than $r$. For each $i=0,1,\dots ,n-1$ we calculate and evaluate $\eta_{\hat{\gG}} (\SCyl (y^iT) )$.
The point is that any connected component $Z$ of $\hat{\gG}$ satisfies the condition that $CH(Z(\infty ))=Z$, which implies that $Z(\infty)$ belongs to $\SCyl (T)$ for $v\in V(X)$ and $T\in \R_r (v) $ if and only if $Z\cap B(v,r)=T$.
Hence we have

\begin{align*}
\eta_{\hat{\gG}} (\SCyl (y^iT) )
=&\#\{ Z\in \Comp (\hat{\gG})\mid Z(\infty )\in \SCyl (y^iT) \} \\
=&\#\{ Z\in \Comp (\hat{\gG})\mid Z\cap B(y^i,r)=y^iT\} .
\end{align*}

\underline{Case 1:} The number $i$ belongs to $\{ r, \dots , n-r\} $.

In this case, the point is that the closed ball $B(y^i, r)$ in $X$ is included in $CH_{H_n}$. Consider a connected component $Y$ of $\gG$ with $|Y|\ni y^i$.
Since $y^i\in \F$, there exists $v(\id , T')\in V(Y)$ and we have
\[ |Y|\cap B(y^i ,r)=CH_{H_n}\cap T'\cap B(y^i ,r )=T'\cap B(y^i ,r ) \]
by Lemma \ref{lem:neighborhood of realization of Gamma}.
Hence for the connected component $Z$ of $\hat{\gG}$ containing $|Y|$, we also see that
$Z \cap B(y^i,r)=|Y|\cap B(y^i,r)$.
Note that for a connected component $Z$ of $\hat{\gG}$ containing $y^i$, $Z$ must include a subgraph $|Y|$ for a connected component $Y$ of $\gG$ and $|Y|\ni y^i$. 
Hence we have
\begin{align*}
&\eta_{\hat{\gG}} (\SCyl (y^iT) )\\
=&\#\{ Z\in \Comp (\hat{\gG})\mid Z\cap B(y^i,r)=y^iT\} \\
=&\#\{ Y\in \Comp (\gG)\mid |Y|\cap B(y^i,r)=y^iT\} \\
=&\#\{ Y\in \Comp (\gG)\mid v(\id ,T')\in V(Y), T'\cap B(y^i,r)=y^iT\} \\
=&\sum_{\substack{T'\in \R_{r}(\F)\\[1pt] T'\cap B(y^i,r)=y^iT}}\theta (T')\\
=&\sum_{\substack{T'\in \R_{\rho }(y^i )\\[1pt] T'\cap B(y^i,r)=y^iT}}\theta (T')\\
=&\sum_{\substack{T'\in \R_{\rho }(\id )\\[1pt] T'\cap B(\id ,r)=T}}\theta (T')=\theta (T).
\end{align*}
Note that
\[ \SCyl (T)=\bigsqcup_{\substack{T'\in \R_{\rho }(\id )\\[1pt] T'\cap B(\id ,r)=T}}\SCyl (T').\]
Recall that we took $\theta $ after fixing $n$ and $\rho=r+n$. Since $\frac{1}{M}\theta (T')$ is close to $\mu (\SCyl (T'))$ for $T'\in \R_{\rho}(\id )$ and the cardinality of $\R_{\rho }(\id )$ is finite and depends on $\rho$, $\frac{1}{M}\theta (T)$ is also close to $\mu (\SCyl (T))$.

\underline{Case 2:} The number $i$ belongs to $\{ 0,\dots ,r-1, n-r+1,\dots ,n-1\} $.

For a connected component $Z$ of $\hat{\gG}$ containing $y^i$ the intersection of $Z$ and $B(y^i,r)$ is influenced by our construction of $\hat{\gG}$ from $|\gG|$. The point is that we can make $r/n$ as small as we like since we choose $n$ after $r$.
Recall that $C_1$ is the number of vertices of $|\iota|^{-1}(\id )$ with degree $1$ in $|\gG|$.
Then we have
\begin{align*}
&\eta_{\hat{\gG}} (\SCyl (y^iT) )\\
=&\#\{ Z\in \Comp (\hat{\gG})\mid Z\cap B(y^i,r)=y^iT\} \\
\leq &\#\{ Z\in \Comp (\hat{\gG})\mid Z\ni y^i\} \\
\leq &\#\{ Y\in \Comp (\gG )\mid v(\id ,T')\in V(Y), T'\ni y^i\} +C_1\\
\leq &\sum_{T'\in \R_{\rho }(y^i )}\theta (T')+C_1\\
\leq &\sum_{T'\in \R_{\rho }(\id )}\theta (T')+C_1
\end{align*}
Note that
\[ C_1\leq \#|\iota|^{-1}(\id )\leq \sum_{T'\in \R_{\rho }(\id )} \theta (T') \]
and
\[ \bigsqcup_{T'\in \R_{\rho}(\id )}\SCyl (T)=A_\id .\]
Hence for
\[ \theta (\id ):=\sum_{T'\in \R_{\rho }(\id )}\theta (T'),\]
$\frac{1}{M}\theta(\id)$ is also close to $\mu(A_\id)$ and there exists a constant $C$ depending on $\mu (A_\id )$ such that
\[ \frac{1}{M}\theta(\id )\leq C.\]
Then we see that 
\[ \eta_{\hat{\gG}} (\SCyl (y^iT) )\leq 2CM.\]
Note that
\[ \theta(T)\leq \theta(\id )\leq CM.\]

From Case 1 and Case 2 we have
\begin{align*}
&\left| \frac{1}{nM}\iota_{H_n}(\eta_{\hat{\gG}} )(\SCyl (T))-\mu (\SCyl (T))\right| \\
\leq &\left| \frac{n-2r+1}{nM}\theta (T)-\mu (\SCyl (T))\right| + \frac{2r-1}{nM}\cdot 2CM\\
\leq &\left| \frac{1}{M}\theta (T)-\mu (\SCyl (T))\right| + \frac{2r-1}{nM}\theta (T)+ \frac{2(2r-1)C}{n}\\
\leq &\left| \frac{1}{M}\theta (T)-\mu (\SCyl (T))\right| + \frac{3(2r-1)C}{n}.
\end{align*}
Therefore, if we take $n$ sufficiently large and take $\theta $ such that $\frac{1}{M} \theta $ is sufficiently close to $\mu$, then we have
\[ \left| \frac{1}{nM}\iota_{H_n}(\eta_{\hat{\gG}} )(\SCyl (T))-\mu (\SCyl (T))\right| <\varepsilon\]
for every $T\in \R_r(\id)$.
This completes the proof of Theorem \ref{thm:approximating theorem for free group}.

\subsection{Denseness property of surface groups}\label{subsec:denseness property of surface groups}

We prove the following theorem in this subsection, which is our main result:

\begin{theorem}\label{thm:surface group has denseness property}
Surface groups have the denseness property of rational subset currents.
\end{theorem}

Note that the fundamental group of a compact hyperbolic surface is a free group of finite rank or a surface group.
Hence we also have the following theorem:

\begin{theorem}\label{thm:subset currents on hyperbolic surface has denseness property}
The fundamental group $\pi_1(\gS)$ of a compact hyperbolic surface $\gS$ has the denseness property of rational subset currents.
\end{theorem}

From now on, let $\gS$ be a closed hyperbolic surface and $G$ the fundamental group of $\gS$.
In this subsection we write $\SC (G)$ to denote the space of subset currents on $G$ since we consider both the universal cover $\HH$ of $\gS$ and the Cayley graph of $G$ with respect to a finite generating set.

The strategy to prove Theorem \ref{thm:surface group has denseness property} is based on the proof of Theorem \ref{thm:approximating theorem for free group} in the previous subsection.
However, in this case our proof will be more complicated.
We first take a certain sequence of finitely generated subgroups $\{ H_n \}$ of $G$, which are free groups of finite rank, but we need to modify $H_n$ during the proof.
Recall that in Step \ref{step:4} of the proof of Theorem \ref{thm:approximating theorem for free group} we constructed the graph $(\hat{\gG},\hat{\iota})$ from $(|\gG |, |\iota|)$.
We need to modify $H_n$ in this context.
Explicitly, we take $u_0\in G$ independent of $n$ such that $\hat{H_n}:=\langle H_n \cup \{ u_0 \} \rangle$ is isomorphic to the free product of $H_n$ and $\langle u_0 \rangle$, and satisfies several conditions in addition. Then we construct $\nu \in \mathrm{Span}(\SC_r(\hat{H_n}))$ such that $\iota_{\hat{H_n}}(\nu)$ is sufficiently close to a given subset current $\mu \in \SC (G)$.
Note that $\iota_{\hat{H_n}}(\SC (\hat{H_n}))$ includes $\iota_{H_n}(\SC(H_n))$ since $\iota_{H_n}=\iota_{\hat{H_n}}\circ \iota_{H_n}^{\hat{H_n}}$.

We can obtain Theorem \ref{thm:surface group has denseness property} as a corollary to the following theorem:

\begin{theorem}\label{thm:approximation of surface group by free groups}
There exists a sequence of finitely generated subgroups $\{ J_n \}_{n\in \NN}$ of $G$ such that each $J_n$ is a free group of finite rank and the union
\[ \bigcup_{n\in \NN} \iota_{J_n}(\mathrm{Span}(\SC_r (J_n)) )\]
is a dense subset of $\SC (G )$.
\end{theorem}

From the above theorem, for any $\mu \in \SC (G)$ and its open neighborhood $U$ there exist $n\in \NN$ and $\nu \in \mathrm{Span}(\SC_r(J_n))$ such that $\iota_{J_n}(\nu )$ belongs to $U$. Since $\iota_{J_n}$ is continuous and $J_n$ has the denseness property of rational subset currents, there exists $\nu' \in \SC_r(J_n)$ such that $\iota_{J_n}(\nu' )$ belongs to $U$. Note that $\iota_{J_n}(\SC_r(J_n))\subset \SC_r (G)$. Hence $G$ has the denseness property of rational subset currents.

\begin{figure}[h]
\begin{center}
\includegraphics[width=9cm]{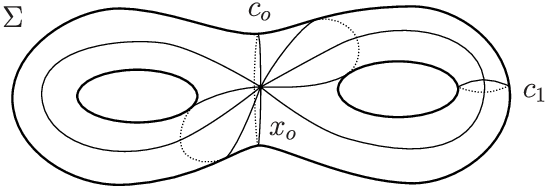}
\vspace{-0.3cm}
\caption{The four closed curves except $c_0$ and $c_1$ represent the $1$-skeleton of the CW-complex structure of $\gS$ and also represent the generating set $B_G$ of $G$.}\label{fig:genus 2 surface}
\end{center}
\end{figure}

For the simplicity of describing subgroups of $G$, we assume that the genus of $\gS$ is $2$.
We construct $\gS$ by gluing edges of an octagon by the fundamental way.
This construction gives $\gS$ a CW-complex structure, a base point $x_0$ of $G=\pi_1(\gS ,x_0)$ and a generating set $B_G$ of $G$.
Set $X:=\Cay (G,B_G)$.
We also fix a hyperbolic metric on $\gS$ and assume that there exists a closed geodesic $c_0$ passing through the base point $x_0$ and dividing $\gS$ into two compact surfaces, each of which is a torus with one boundary component and contains two generators of $G$ (see Figure \ref{fig:genus 2 surface}).

The CW-complex structure on $\gS$ induces the CW-complex structure on the universal cover $\HH$ of $\gS$. 
Fix a lift $\tilde{x_0}$ of $x_0$ in $\HH$. Then we can see that there exists a $G$-equivariant homeomorphism $\Phi$ from the Cayley graph $X$ to the $1$-skeleton $\HH^{(1)}$ of $\HH$ such that $\Phi (g)=g\tilde{x_0}$ for every $g\in G$. 
Moreover, the map $\Phi$ is a quasi-isometry from the \v{S}varc-Milnor Lemma.

\begin{figure}[h]
\begin{center}
\includegraphics[width=9cm]{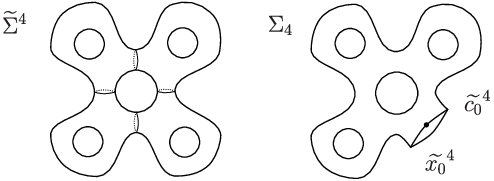}
\vspace{-0.3cm}
\caption{The four closed curves on $\tilde{\gS}^4$ are the copies of $c_1$. This figure corresponds to Figure \ref{fig:H_4graph} in the previous subsection.}\label{fig:4-fold covering}
\end{center}
\end{figure}

Take a closed geodesic $c_1$ cutting one of the handles of $\gS$ (see Figure \ref{fig:genus 2 surface}).
For $n\geq 2$ we can obtain an $n$-fold covering space $\tilde{\gS}^n$ of $\gS$ by cutting $\gS$ along $c_1$ and gluing $n$-copies of $\gS\setminus c_1$ along $c_1$ (see the left of Figure \ref{fig:4-fold covering} for $\tilde{\gS}^4$).
Let $p_{\tilde{\gS}^n}$ be the covering map from $\tilde{\gS}^n$ to $\gS$ and $\tilde{x_0}^n$ a lift of $x_0$ in $\tilde{\gS}^n$.
Let $G_n$ be the image of the homomorphism $(p_{\tilde{\gS}^n})_\#$ from $\pi_1(\tilde{\gS}^n, \tilde{x_0}^n)$ to $G=\pi_1(\gS ,x_0)$.
Consider a lift $\tilde{c_0}^n$ of $c_0$ passing through $\tilde{x_0}^n$ in $\tilde{\gS}^n$. Then $\tilde{c_0}^n$ divides $\tilde{\gS}^n$ into two connected components, one of which is a torus with one boundary component and the other of which is an $n$-genus surface with one boundary component, denoted by $\gS_n$ (see the right of Figure \ref{fig:4-fold covering} for $\gS_4$). The point is that $\gS_n$ ``approximates'' $\tilde{\gS}^n$ if $n$ is large.

Set $H_n:=(p_{\tilde{\gS}^n}|_{\gS_n})_\# (\pi_1(\gS_n, \tilde{x_0}^n))<G$.
Since $\tilde{c_0}^n$ is a closed geodesic of $\tilde{\gS}^n$, the convex core $C_{H_n}$ is identified with $\gS_n$.
Then $CH_{H_n}$ contains $\tilde{x_0}$.
We can see that $\frac{1}{n}\eta_{H_n}$ converges to $\eta_G$ in $\SC (G)$ from the proof of Theorem \ref{thm:approximation of surface group by free groups}.

Let $d=d_X$ be the path metric on $X$ such that each edge of $X$ has length one.
We fix $\delta > 0$ such that $X=\Cay (G, B_G)$ is a $\delta$-hyperbolic space, which means that for every geodesic triangle $\gD$ in $X$, every side of $\gD$ is included in the union of the $\delta$-neighborhood of the other two sides.
Recall that a map $f$ from an interval $I$ of $\RR$ to $X$ is called an $(a,b)$-\ti{quasi-geodesic} for $a\geq 1 ,b\geq 0$ if for any $t_1,t_2\in I$ we have
\[ \frac{1}{a}|t_1-t_2|-b\leq d(f(t_1),f(t_2)) \leq a |t_1-t_2|+b.\]
In the case that $I=\RR$, we call it an $(a,b)$-quasi-geodesic line.
A quasi-geodesic (line) is an $(a,b)$-quasi-geodesic (line) for some $a\geq 1$ and $b\geq 0$. 
We consider only continuous quasi-geodesics.

\begin{remark}[Constants related to $\delta$]
This remark is the most important remark in this subsection.
In the case of a free group $F$ of finite rank, the Cayley graph of $F$ with respect to a free basis is a tree, which is a $0$-hyperbolic space. Then we could construct subtrees or geodesics of the Cayley graph clearly. However, when we construct something on $X$, we will be always annoyed with some positive constants coming from $\delta$.
Here, we introduce some conventions in order to reduce the complicatedness of such constants.
Even if a constant depends on not only $\delta$ but also objects not depending on $n$ appeared in the above, such as $\HH$, $\Phi$, and the degree of a vertex of $X$, we will say that the constant depends only on $\delta$.
We will use a symbol $\delta '$ to represent a constant depending only on $\delta$, which can be different in each situation.
We say that a quasi-geodesic $\gamma$ is a $\delta$-\ti{quasi-geodesic} if there exist $a\geq 1, b\geq 0$ depending only on $\delta$ such that $\gamma$ is an $(a,b)$-quasi-geodesic. Note that $a$ and $b$ can be different in each situation.
\end{remark}

\begin{definition}[Convex hull in $X$]
We identify $\partial G$ with the boundary $\partial X$ of $X$.
Recall that for $S\in \H (\partial G)$ the weak convex hull $WC(S)$ of $S$ in $X$ is the union of all geodesic lines connecting two points of $S$. It is known that $WC(S)$ is not necessarily a convex subset of $X$ but $\delta'$-quasi-convex subset of $X$, that is, any geodesic path connecting two points of $WC(S)$ is included in the (closed) $\delta'$-neighborhood $B(WC(S),\delta')$ of $WC(S)$.
We introduce a notion of the \ti{convex hull} $CH(S)$ of $S$ \ti{in} $X$. Note that $X$ is a planar graph since $X$ is homeomorphic to the $1$-skeleton $\HH^{(1)}$.

Let $\xi, \eta \in \partial G$ with $\xi\not=\eta$. We give an orientation to a geodesic line $\ell$ joining $\xi$ to $\eta$.
Since $X$ is planar, we can define the left side $\mathrm{Left}(\ell)$ of $\ell$ and the right side $\mathrm{Right}(\ell)$ of $\ell$, each of which includes $\ell$. We say that an edge $e$ of $WC(\{ \xi ,\eta \})$ is leftmost if $e$ is included in the left side of $\ell$ for every geodesic line $\ell$ from $\xi $ to $\eta$. We can define a rightmost edge of $WC (\{ \xi ,\eta\})$ by the same way.
Then we can see that the union of all leftmost edges forms a quasi-geodesic line joining $\xi$ to $\eta$, which is denoted by $\mathrm{Left}(\xi,\eta)$. 
The union of all rightmost edges also forms a quasi-geodesic line $\mathrm{Right}(\xi ,\eta )$ joining $\xi$ to $\eta$.
Note that $\mathrm{Left}(\xi,\eta)=\mathrm{Right}(\eta, \xi)$.
We define the convex hull $CH(\{ \xi , \eta \})$ of $\{ \xi ,\eta \}$ to be the intersection of the right side of $\mathrm{Left}(\xi,\eta)$ and the left side of $\mathrm{Right}(\xi ,\eta)$.

Let $S\in \H (\partial G)$. If $S=\partial G$, then we define $CH(S)$ to be $X$.
Assume that $S\not=\partial G$.
Recall that $\partial G\setminus S$ is the union of at most countably many open intervals $\{ I_\lambda \}_{\lambda \in \gL}$.
We give an orientation to $\partial G$, which induces an orientation on each $I_\lambda$. For the orientation of $I_\lambda$ we give an orientation to $\partial I_\Lambda=\{ \xi_\lambda ,\eta_\lambda\}$ such that the limit set of the right side of a geodesic line joining $\xi_\lambda $ to $\eta_\lambda$ equals $\ol{I_\lambda}$. Now, we define $CH(S)$ to be the intersection of the left side of $\mathrm{Right}(\xi _\lambda,\eta_\lambda)$ taken over all $\lambda \in \gL$.
We call each $\mathrm{Right}(\xi_\lambda ,\eta_\lambda )$ a boundary component of $CH(S)$ and the union of $\mathrm{Right}(\xi_\lambda, \eta_\lambda)$ taken over all $\lambda \in \gL$ the boundary of $CH(S)$.
Note that every boundary component of $CH(S)$ is a $\delta$-quasi-geodesic line.
\end{definition}

From the definition of $CH(S)$ for $S\in \H (\partial G)$ we can see that $CH(S)$ has the following properties:
\begin{enumerate}
\item $CH(S)$ is a $\delta'$-quasi-convex connected subgraph of $X$;
\item $CH(S)\supset WC(S)$;
\item $CH(S)$ is included in the $\delta'$-neighborhood of $WC(S)$;
\item for every $x,y\in CH(S)$ there exists a $\delta$-quasi-geodesic joining $x$ to $y$ in $CH(S)$.
\end{enumerate}
Recall that when we consider a quasi-geodesic, we always assume that it is continuous, and that is the point in the property (4).
If a geodesic joining $x$ to $y$ goes out from $CH(S)$, then we can consider a $\delta$-quasi-geodesic traveling along the boundary of $CH(S)$. As a result, the property (4) holds, which will play an important role later. 

Now, we can define the notion of round-graphs and subset cylinders with respect to round-graphs by using the convex hull defined in the above by the same way as Definition \ref{def:round-graph generalized version}.
Note that for a round-graph $T\in \R_r(\id)$ we can see that $\SCyl (T)$ is a Borel subset of $\H (\partial G)$ but neither open nor closed.
Therefore, we need to develop a new neighborhood of $\mu \in \SC (G)$ instead of Corollary \ref{cor: Kapovich's open neighborhood of a subset current}.

\begin{notation}
Let $Y$ be a non-empty bounded subset of $X$.
Let $a,r>0$. For $T_1,T_2\in \R_r(Y)$ we denote by $T_1\underset{a}{\sim}T_2$ if $T_1\subset B(T_2,a)$ and $T_2\subset B(T_1,a)$.
Let $d_{\ol{X}}$ be a visual metric on $\ol{X}:=X\cup \partial X$. 
Let $d_{\mathrm{Haus}}$ be the Hausdorff distance on $\H (\partial G)$ induced by the restriction of $d_{\ol{X}}$ to $\partial X=\partial G$.
\end{notation}

\begin{lemma}\label{lem:diameter of subset cylinder tends to zero}
Let $a>0$. Let $Y$ be a non-empty bounded subset of $X$. The supremum of $d_{\mathrm{Haus}}(S_1,S_2)$ taken over all $S_1\in \SCyl (T_1)$, $S_2\in \SCyl (T_2)$ for all $T_1,T_2\in \R_r(Y)$ with $T_1\underset{a}{\sim}T_2$ converges to $0$ when $r\rightarrow \infty $. 
\end{lemma}
\begin{proof}
To obtain a contradiction, suppose that there exists $\varepsilon>0$ such that for any $r_0>0$ there exist $r\geq r_0$ and $T_1,T_2\in \R_r(Y)$ with $T_1\underset{a}{\sim}T_2$ and $S_1\in \SCyl (T_1), S_2\in \SCyl (T_2)$ such that $d_{\mathrm{Haus}}(S_1,S_2)>\varepsilon$.
For such $S_1,S_2$ we can assume that there exists $\xi \in S_1$ such that $d_{\ol{X}}(\xi ,S_2)>\varepsilon$ without loss of generality.
Then there exists $\varepsilon '>0$ depending only on $\varepsilon$ such that $d_{\ol{X}}(\xi , CH(S_2))>\varepsilon'$.
Let $B_{\ol{X}}(\xi ,\varepsilon')$ be the closed ball centered at $\xi$ with radius $\varepsilon '$ with respect to $d_{\ol{X}}$.
Then $B_{\ol{X}}(\xi ,\varepsilon')\cap CH(S_2)=\emptyset$.

Now, we assume that $r_0$ is sufficiently large. Then $B_{\ol{X}}(\xi ,\varepsilon')\cap B(Y,r)\not=\emptyset$, and $B_{\ol{X}}(\xi ,\varepsilon')$ also intersects $T_1$ since $CH(S_1)\cap B(Y,r)=T_1$.
Moreover, $B_{\ol{X}}(\xi ,\varepsilon')$ also intersects $T_2$ since $T_1\subset B(T_2,a)$ for the fixed constant $a>0$.
Therefore $B_{\ol{X}}(\xi ,\varepsilon')$ intersects $CH(S_2)$, a contradiction.
\end{proof}

From the above lemma we can see that the supremum of the diameter of $\SCyl (T)$ in $\H (\partial G)$ taken over $T\in \R_r(Y)$ tends to $0$ when $r\rightarrow \infty$.

The following lemma is a technical lemma that will be used in the proofs of Lemmas \ref{lem:difference between Y and WC(Y infty)} and \ref{lem:quasi-geodesic line and quasi-convex implies good condition}.

\begin{lemma}\label{lem:neighborhood and distance ineq}
Let $a,b>0$.
Let $Y$ be a $b$-quasi-convex subset of $X$.
Let $v_0\in X$ and $y_0\in Y$.
Assume that $r>0$ is much larger than $a, b, d(v_0,y_0)$.
If $x$ belongs to $B(Y,a)\cap B(v_0,r)$, then $x$ also belongs to $B(Y\cap B(v_0,r) ,2(a+b+d(v_0,y_0)) )$.
\end{lemma}
\begin{proof}
Suppose that $x$ belongs to $B(Y,a)\cap B(v_0,r)$.
Take $y\in Y$ such that $d(x,y)\leq a$.
If $y$ belongs to $B(v_0,r)$, then our claim follows.
Hence we assume that $d(v_0,y)>r$.
Take a geodesic $\ell$ joining $y_0$ to $y$.
Note that $d(y_0,y)\geq d(v_0,y)-d(y_0,v_0)\geq r-d(y_0,v_0)$.
Hence we can take $p\in \ell$ such that $d(y_0,p)=r-b-d(v_0,y_0)$.
Then we have
\begin{align*}
d(p,y)=d(y_0,y)-d(y_0,p)
&\leq d(y_0,v_0)+d(v_0,y)-(r-b-d(v_0,y_0))\\
&\leq d(v_0,y)-r+b+2d(v_0,y_0)  \\
&\leq d(v_0,x)+d(x,y)-r+b+2d(v_0,y_0)\\
&\leq a+b+2d(v_0,y_0).
\end{align*}
Since $Y$ is $b$-quasi-convex, there exists $p'\in Y$ such that $d(p,p')\leq b$.
Then
\[ d(v_0,p')\leq d(v_0,y_0)+d(y_0,p)+d(p,p')\leq r,\]
which implies that $p'\in B(v_0,r)$.
Moreover, we have
\[ d(x,p')\leq d(x,y)+d(y,p)+d(p,p')\leq a + (a+b+ 2d(v_0,y_0) ) + b.\]
Therefore $x$ belongs to $B(Y\cap B(v_0,r) ,2(a+b+d(v_0,y_0)) )$.
\end{proof}

For $U\subset \H (\partial G )$ and $a>0$ set
\[ B_\H (U ,a):=\{ S\in \H (\partial G )\mid d_{\mathrm{Haus}}(U, S )\leq a\} ,\]
the $a$-neighborhood of $U$ in $\H (\partial G)$.
Then we have the following lemma:

\begin{lemma}\label{lem:difference between Y and WC(Y infty)}
Let $\varepsilon, a >0$. Let $Y$ be a non-bounded subset of $X$ with $Y(\infty )\in \H (\partial G)$. Let $y\in Y$.
For a sufficiently large $r>0$ depending on $\varepsilon, a, \delta$, if $Y\cap B(y, r)\in \R_r(y)$ and
\[ Y \cap B(y,r )\underset{a}{\sim }WC (Y(\infty ))\cap B(y,r ),\]
then $Y(\infty )$ belongs to $B_\H (\SCyl (Y\cap B(y ,r )), \varepsilon )$.
\end{lemma}
\begin{proof}
Take $S\in \SCyl (Y\cap B(y,r))$, which implies that $CH (S)\cap B(y,r)= Y\cap B(y,r)$.
Note that $Y(\infty) \in \SCyl (CH(Y(\infty))\cap B(y,r))$ by definition.
Take $\delta'>0$ such that $CH (Y(\infty ))$ is included in $B(WC (Y(\infty)), \delta')$ and $WC(Y(\infty))$ is $\delta'$-quasi-convex.
We are going to prove that
\[ CH(Y(\infty ))\cap B(y,r )\underset{4\delta' +3a}{\sim} CH(S)\cap B(y,r ).\]
Then from Lemma \ref{lem:diameter of subset cylinder tends to zero}, we see that the Hausdorff distance between $Y(\infty )$ and $S$ is smaller than $\varepsilon$ for a sufficiently large $r>0$.

Since $Y \cap B(y,r )\underset{a}{\sim }WC (Y(\infty ))\cap B(y,r )$, we can take $y_0\in WC(Y)$ such that $d(y,y_0)\leq a$.
Take $x\in CH(Y(\infty ))\cap B(y,r )$.
Then $x$ belongs to $B(WC(Y(\infty)),\delta' )\cap B(y,r)$.
From Lemma \ref{lem:neighborhood and distance ineq}, $x$ belongs to $B(WC(Y(\infty))\cap B(y,r) ,2(\delta' +\delta'+a))$.
Moreover, from the assumption, we have
\begin{align*}
B(WC(Y(\infty ))\cap B(y,r),4\delta' +2a)
&\subset B(Y\cap B(y,r),4\delta' +3a)\\
&=	B(CH(S)\cap B(y,r),4\delta' +3a ).
\end{align*}
Hence
\[ CH(Y(\infty ))\cap B(y,r )\subset B(CH(S)\cap B(y,r),4\delta' +3a ).\]

Since $WC (Y(\infty ))\subset CH (Y(\infty ))$, we have
\begin{align*}
CH(S)\cap B(y,r )=Y\cap B(y,r) 
&\subset B(WC(Y(\infty ))\cap B(y,r) ,a) \\
&\subset B(CH(Y(\infty ))\cap B(y,r) ,a).
\end{align*}
Therefore 
\[ CH(Y(\infty ))\cap B(y,r )\underset{4\delta'+ 3a}{\sim} CH(S)\cap B(y,r).\]
This completes the proof.
\end{proof}

Let $\mu \in \SC (G)$. For compactly supported continuous functions $f_1,\dots ,f_k$ on $\H (\partial G)$ and $\varepsilon>0$ we have an open neighborhood $U(f_1,\dots ,f_k ; \varepsilon )$ of $\mu$ defined by
\[ \{ \nu\in \SC(G)\mid \left|\int f_i d \mu -\int f_i d\nu \right| <\varepsilon \text{ for every } i=1,\dots ,k\}, \]
and the family of all such open neighborhoods of $\mu$ forms a fundamental system of open neighborhoods of $\mu$.

Since the proof of Theorem \ref{thm:approximation of surface group by free groups} is long and includes many constants, we will write \underline{\textbf{Setting}} when we fix something.

\begin{setting8}\label{set:1}
Fix $\mu \in \SC (G)$ and compactly supported continuous functions $f_1,\dots , f_k$ on $\H (\partial G)$ and $\varepsilon_\mu >0$. We assume that $\mu$ is not the zero measure.
Take $r_\mu\in \NN$ such that 
\[ A(\id ,r_\mu ):=\{ S\in \H (\partial G)\mid CH(S)\cap B(\id ,r_\mu )\not= \emptyset \}\]
includes the support of $f_i$ for every $i=1,\dots ,k$.
\end{setting8}

The set $A(\id ,r_\mu )$ is a compact subset of $\H (\partial G)$.
Since each $f_i$ is compactly supported, $f_i$ is a uniformly continuous function.

Let $m$ be a Borel measure on a topological space $\Omega$. Set $|m|:=m(\Omega)$. For a non-empty Borel subset $A$ of $\Omega$ we denote by $m|_A$ the restriction of $m$ to $A$. The support of $m$, denoted by $\mathrm{supp}\, m$, is the smallest closed subset $A$ of $\Omega$ such that $m(A^c)=0$. Then $|m|=m(\Omega )=m (\mathrm{supp}\,m)$.

The following lemma describes a sufficient condition for a subset current to belong to the open neighborhood $U(f_1,\dots ,f_k;\varepsilon_\mu)$ of $\mu$.

\begin{lemma}\label{lem:describe a neighborhood of a given subset current}
Let $r_\mu' \geq r_\mu $.
There exist $\rho>0, \varepsilon_1>0, \varepsilon_2 >0$ such that if $\nu\in \SC(F)$ satisfies the following conditions, then $\nu \in U(f_1,\dots ,f_k;\varepsilon_\mu)$:
\begin{enumerate}
\item there exist Borel measures $\nu ', \nu_T$ on $\H (\partial G)$ for $T\in \R_\rho(B(\id ,r_\mu' ))$ such that
\[ \nu|_{A(\id ,r_\mu )}=\sum_{T\in \R_\rho ( B(\id, r_\mu '))}\nu _T|_{A(\id ,r_\mu )}+\nu' ;\] 
\item $\mathrm{supp}\,\nu_T \subset \ol{B_\H (\SCyl (T), \varepsilon_1 )}$ for every $T\in \R_\rho ( B(\id, r_\mu' ))$;
\item $|\nu '|<\varepsilon_2$;
\item $\left| |\nu_T |-\mu (\SCyl (T) ) \right| <\varepsilon_2$ for every $T\in \R_\rho ( B(\id, r_\mu' ))$.
\end{enumerate}
\end{lemma}
\begin{proof}
Note that
\[ A(\id ,r_\mu )\subset A(\id ,r_\mu' )=\bigsqcup_{T\in \R_\rho (B(\id, r_\mu' ))}\SCyl (T).\]
Let $f$ be an element of $\{ f_1,\dots ,f_k\}$.
Since $\mathrm{supp}\,f$ is included in $A(\id ,r_\mu)$, we have
\begin{align*}
&\left| \int f d\nu -\int f d\mu \right|\\
= &\left| \sum_{T\in \R_\rho (B(\id, r_\mu' ))} \int f d \nu_T + \int f d \nu' -\sum_{T\in \R_\rho (B(\id, r_\mu' ))} \int_{\SCyl (T)}f d\mu \right| \\
\leq &\sum_{T\in \R_\rho (B(\id, r_\mu' ))} \left|\int f d\nu_T -\int_{\SCyl (T)}f d\mu \right| +|\nu '|\max |f|.
\end{align*}

Let $\varepsilon_3>0$. From Lemma \ref{lem:diameter of subset cylinder tends to zero}, for a sufficiently large $\rho$ and small $\varepsilon_1>0$ the diameter of $K_T:=\ol{B_\H (\SCyl (T), \varepsilon_1 )}$ is sufficiently small, and then we have
\[ \sup_{S\in K_T}f(S)-\inf_{S\in K_T}f(S)<\varepsilon_3 \]
for $T\in \R_\rho (B(\id, r_\mu' ))$. Set
\[ M_T:=\sup_{S\in K_T}f(S).\]
Then for each $T\in \R_\rho (B(\id, r_\mu' ))$
\begin{align*}
&\left|\int f d\nu_T -\int_{\SCyl (T)}f d\mu \right| \\
= &\Bigg| \int f d\nu_T - M_T|\nu_T | +M_T|\nu_T | \\
	&\quad -M_T \mu (\SCyl (T))+M_T\mu (\SCyl (T))-\int_{\SCyl (T)}f d\mu \Bigg| \\
\leq &\varepsilon_3 |\nu_T |+|M_T|\varepsilon_2 +\varepsilon_3 \mu (\SCyl (T))\\
< &\varepsilon_3 (\mu (\SCyl (T))+\varepsilon_2 )+|M_T|\varepsilon_2 +\varepsilon_3 \mu (\SCyl (T)).
\end{align*}
Hence
\begin{align*}
&\left| \int f d\nu -\int f d\mu \right| \\
< & \varepsilon_3 \mu (A(\id, r_\mu '))+ \varepsilon_2 \varepsilon_3\# \R_\rho (B( \id,r_\mu '))\\
		&+\# \R_\rho (B( \id,r_\mu ')) \cdot \varepsilon_2 \cdot \max |f|+\varepsilon_3 \mu (A(\id ,r_\mu '))+\varepsilon_2 \max |f|.
\end{align*}
Now, we assume that $\varepsilon_3$ is sufficiently small. Then we need to take small $\varepsilon_1$ and large $\rho$.
Hence $\# \R_\rho (B( \id,r_\mu' ))$ will be large.
Finally, we take $\varepsilon_2$ sufficiently small. Then we can obtain 
\[ \left| \int f d\nu -\int f d\mu \right| <\varepsilon_\mu. \]
This completes the proof.
\end{proof}

\begin{setting8}\label{set:2}
The gap between $r_\mu'$ and $r_\mu$ depends on $\delta$, and $r_\mu'$ will be determined later. 
We fix $\rho ,\varepsilon_1,\varepsilon_2>0$ satisfying the conditions in the above lemma.
We assume that $\rho$ is much larger than constants depending on $\delta$.
\end{setting8}

We are going to construct $\nu \in \mathrm{Span} (\SC _r(G))$ satisfying the conditions in the above lemma.
When we check the condition (2), we will use Lemma \ref{lem:difference between Y and WC(Y infty)}.
Recall that we constructed the SC-graph $(\hat{\gG },\hat{\iota})$ on $F$ in Step \ref{step:4} in the previous subsection such that each connected component $Z$ of $\hat{\gG}$ satisfying the condition that $Z=CH(Z(\infty))$. Since the Cayley graph of $F$ with respect to a free basis is a tree, the condition that every vertex of $Z$ has degree larger than $1$ implies that $Z=CH(Z (\infty ))$.
In the case of the Cayley graph $X$ of $G$ we need to give a new criterion in order to use Lemma \ref{lem:difference between Y and WC(Y infty)}.

\begin{lemma}\label{lem:quasi-geodesic line and quasi-convex implies good condition}
Let $Y$ be a non-bounded subset of $X$ and $y\in Y$. Assume that $Y$ is $c$-quasi-convex in $X$ for a constant $c\geq 0$.
Take $r>0$ much larger than $c$ and $\delta$. 
If for every $z\in Y\cap B(y,r)$ there exists a $\delta$-quasi-geodesic line $\ell$ in $Y$ such that $d(z,\ell )\leq c$,
then there exists $a>0$ depending only on $c$ and $\delta$ such that
\[ Y\cap B(y,r )\underset{a}{\sim }WC(Y(\infty ))\cap B(y,r ) .\]
\end{lemma}
\begin{proof}
Take $z\in Y\cap B(y,r)$. From the assumption there exists a $\delta$-quasi-geodesic line $\ell$ in $Y$ such that $z\in B(\ell,c)$.
Then the $\delta'$-neighborhood of a geodesic line $\ell'$ connecting two endpoints of $\ell$ includes $\ell$, which implies that
\[ z\in B(\ell ,c)\subset B(\ell ', \delta'+c )\subset B(WC (Y(\infty )),\delta '+c).\]
Note that $WC(Y(\infty))$ is $\delta''$-quasi-convex for $\delta''>0$ depending only on $\delta$ and there exists $y_0\in WC(Y(\infty))$ such that $d(y,y_0)\leq \delta' +c$.
Then from Lemma \ref{lem:neighborhood and distance ineq} we see that
\[ z\in B(WC(Y(\infty ))\cap B(y,r), 2(\delta' +c +\delta'' +\delta' +c)).\]

Take $z\in WC(Y(\infty))\cap B(y,r)$.
Let $\ell $ be a geodesic line connecting two points of $Y(\infty)$ passing through $z$. Since $Y$ is $c$-quasi-convex, $\ell$ is included in $B(Y, \delta'+c )$, which implies that $z\in B(Y,\delta'+c)\cap B(y,r)$.
Hence $z$ belongs to $B(Y\cap B(y,r), 2(\delta' +c + c))$ by Lemma \ref{lem:neighborhood and distance ineq}.
From the above, $a:=2(2\delta'  +\delta '' +2c) $ satisfies the condition in our claim.
\end{proof}

In order to use Lemma \ref{lem:quasi-geodesic line and quasi-convex implies good condition} we need to see the existence of $\delta$-quasi-geodesic lines in $Y$.
Hence when we construct a graph by ``connecting'' round-graphs, we need to construct a quasi-geodesic line in each connected component of the graph.
For the purpose, we modify the definition of a round-graph in Definition \ref{def:round-graph generalized version}.

\begin{definition}[Round-graph with information of geodesics]\label{def:round-graph with information of geodesics}
Let $r>0$. Let $Y$ be a non-empty bounded subset of $X$ and $T\in \R_r(Y)$.
Let $\gamma_1,\dots ,\gamma_m$ be subsets of $B(Y,r)$ such that for every $\gamma_i$ there exists a geodesic line $\ell$ such that $\ell \cap B(Y,r)=\gamma_i$. Note that $\gamma_i$ can be non-connected, but we call $\gamma_i$ a geodesic in $B(Y,r)$.
We call a pair $(T, \{ \gamma_1,\dots ,\gamma_m\} )$ \ti{a round-graph of $r$-neighborhood of $Y$ with information of geodesics} if 
there exists $S\in \H (\partial G)$ satisfying the following conditions:
\begin{enumerate}
\item $T\cap Y\not=\emptyset$;
\item $T=CH(S)\cap B(Y,r)$;
\item for every $\gamma_i$ there exists a geodesic line $\ell$ connecting two points of $S$ such that $\ell\cap B(Y,r)=\gamma_i$;
\item for every geodesic line $\ell$ connecting two points of $S$ there exists $\gamma_i$ such that $\ell\cap B(Y,r)=\gamma_i$.
\end{enumerate}
From the conditions $(3)$ and $(4)$, we see that $WC(S)\cap B(Y,r)=\bigcup_i \gamma_i$.
We denote by $\R_r^\ast(Y)$ the set of all round-graphs of $r$-neighborhood of $Y$ with information of geodesics.
For $T_\ast=(T,\gamma_T)\in \R_r^\ast(Y)$, we define $|T_\ast|$ to be $T$
and we will write the pair $(T,\gamma_T)$ simply as $T$.
In this notation $T\in \R_r^\ast (Y)$ means that $T=(|T|,\gamma_T)$.
We call an element of $\gamma_T$ a geodesic of $T$.

For $T\in \R_r^\ast(Y)$ we define the \ti{subset cylinder} $\SCyl (T)$ with respect to $T$ to be a subset of $\H (\partial G)$ consisting of $S$ satisfying the conditions $(2),(3),(4)$ in the above.
For a subset $Z$ of $B(Y,r)$ the restriction of $T$ to $Z$, denoted by $T|_Z$, is defined to be the pair of $|T|\cap Z$ and the set consisting of $Z\cap \gamma$ for every $\gamma \in \gamma_T$.

Let $Y,Z$ be non-empty bounded subsets of $X$.
For $T_1\in \R_r^\ast(Y), T_2\in \R_r^\ast(Z)$ we say that $T_1$ and $T_2$ are \ti{connectable} if $T_1|_{B(Y,Z,r)}=T_2|_{B(Y,Z,r)}$.
Note that $B(Y,Z,r)$ can be empty and then $T_1$ and $T_2$ are connectable for any $T_1\in \R_r^\ast(Y), T_2\in \R_r^\ast(Z)$.

Assume that $B(Y,Z,r)$ is not empty.
A pair of a subset $J$ of $B(Y,Z,r)$ and a set of geodesics $\gamma_1\dots ,\gamma_m$ in $B(Y,Z,r)$ is called a $(Y,Z)$-round-graph of $r$-neighborhood of $Y,Z$ with information of geodesics if there exists $S\in \H (\partial G)$ satisfying the following conditions:
\begin{enumerate}
\item $J\cap Y\not=\emptyset ,J\cap Z\not=\emptyset$;
\item $J=CH(S)\cap B(Y,Z,r)$;
\item for every $\gamma_i$ there exists a geodesic line $\ell$ connecting two points of $S$ such that $\ell\cap B(Y,Z,r)=\gamma_i$.
\item for every geodesic line $\ell$ connecting two points of $S$ there exists $\gamma_i$ such that $\ell\cap B(Y,Z,r)=\gamma_i$.
\end{enumerate}
We denote by $\R_r^\ast(Y,Z)$ the set of all $(Y,Z)$-round-graph of $r$-neighborhood of $Y,Z$ with information of geodesics.
For $J\in \R_r^\ast(Y,Z)$ we define the subset cylinder $\SCyl (J)$ with respect to $J$ to be a subset of $\H (\partial G)$ consisting of $S$ satisfying the conditions $(2),(3),(4)$ in the above.
For $T_1\in \R_r^\ast(Y), T_2\in \R_r^\ast(Z)$ we say that $T_1$ and $T_2$ are $J$-\ti{connectable} for $J\in \R_r^\ast(Y,Z)$ if $T_1|_{B(Y,Z,r)}=J=T_2|_{B(Y,Z,r)}$.
\end{definition}

\begin{remark}
For $T\in \R_r^\ast(Y)$ we can see that the subset cylinder with respect to $T$ is included in the subset cylinder with respect to $|T|$ since $T$ has more information than $|T|$. Actually, for every $T_0\in \R_r (Y)$ we have
\[ \SCyl (T_0)=\bigsqcup_{\substack{T\in \R_r^\ast(Y) \\[1pt] |T|=T_0}}\SCyl (T).\]

For $J\in \R_r^\ast(Y,Z)$ the following equation holds:
\[ \SCyl (J)=\bigsqcup_{\substack{T\in \R_r^\ast(Y)\\[1pt] T|_{ B(Y,Z,r)}=J}} \SCyl (T),\]
which implies that for any $\nu \in \SC (F)$ we have
\[ \sum_{\substack{T\in \R_r^\ast(Y)\\[1pt] T|_{B(Y,Z,r)}=J}}\nu (\SCyl (T))=\sum_{\substack{T'\in \R_r^\ast(Z)\\[1pt] T'|_{B(Y,Z,r)}=J}}\nu (\SCyl (T')).\]
\end{remark}

\begin{setting8}\label{set:3}
Fix $n\in \NN $ with $n\geq 2$. We will assume that $n$ is sufficiently large.
\end{setting8}

Recall that we have the homeomorphism $\Phi$ from $X$ to the 1-skeleton $\HH^{(1)}$ of $\HH$.
Set $X_{H_n}:=\Phi^{-1}(CH_{H_n}\cap \HH^{(1)})$. Then $X_{H_n}$ is an $H_n$-invariant subgraph of $X$. Moreover, we can see that for any two points $x,y\in X_{H_n}$ there exists a geodesic joining $x$ to $y$ in $X_{H_n}$ since $X\cong \HH^{(1)}$ is a planar graph and $X_{H_n}$ is surrounded by geodesic lines in $X$, which are called boundary components of $X_{H_n}$. We denote by $\partial X_{H_n}$ the union of boundary components of $X_{H_n}$ and call it the boundary of $X_{H_n}$.
Note that the CW-complex structure on $\gS$ induces a CW-complex structure on $\tilde{\gS}^n$ and $\gS_n=C_{H_n}$ includes all vertices of $\tilde{\gS}^n$.
The intersection of $\gS_n$ and the $1$-skeleton of $\tilde{\gS}^n$ is the $1$-skeleton of $\gS_n$, which can be identified with the quotient graph $H_n\backslash X_{H_n}$.

Consider the action of $H_n$ on $CH_{H_n}\subset \HH$. By cutting $\gS_n$ along appropriate piecewise geodesics, we can obtain a bounded connected fundamental domain $\F_0$ for the action of $H_n$ on $CH_{H_n}$ satisfying the following conditions:
\begin{enumerate}
\item $H_n(\F)=CH_{H_n}$;
\item $h\F\cap \F=\emptyset$ for non-trivial $h\in H_n$;
\item $\ol{\F_0}$ is a polygon, which is not necessarily convex;
\item every edge of $\partial \F_0$ intersects an edge of $\HH^{(1)}$ transversally;
\item we can obtain a free basis $B_n$ of $H_n$ as side-pairing transformations of $\F_0$, that is,
$B_n\sqcup B_n^{-1}=\{ h\in H_n\setminus \{ \id \} \mid h \ol{\F_0} \cap \ol{\F_0} \not=\emptyset \}$.
\end{enumerate}
Set $\F :=\Phi^{-1}(\F_0\cap \HH^{(1)})$. Then $\F$ is a fundamental domain for the action of $H_n$ on $X_{H_n}$ and we also have
$B_n\sqcup B_n^{-1}=\{ h\in H_n\setminus \{ \id \} \mid h \ol{\F} \cap \ol{\F} \not=\emptyset \}$.
The fundamental domain $\F$ is a non-connected subset of $X_{H_n}$ in general. We can assume that $\F \ni \id$ and $\ol{\F}$ contains exactly $n$ vertices since the $0$-skeleton of $\tilde{\gS}^n$ consists of $n$ vertices. 
 
Set $X_n:=\Cay (H_n ,B_n)$.
Then $X_n$ is a tree, and each vertex $h\in V(X_n)$ corresponds to $h\F \subset X_{H_n}$. 
From the property of $B_n$, we can see that two vertices $h_1,h_2\in V(X_n)$ are adjacent if and only if $h_1\not=h_2$ and $h_1\ol{\F}\cap h_2\ol{\F}\not=\emptyset$.

\begin{setting8}\label{set:4}
Fix a sufficiently large $\rho_0 $. We will take $\rho_1,\rho_2, \rho_3$ later such that $\rho_3\leq \rho_2\leq \rho_1\leq \rho_0$, where the gaps depend on some constants depending on $n$ and $\delta$. We assume that all of $\rho_0,\rho_1,\rho_2 ,\rho_3$ are much larger than any constants depending on $\delta$.
\end{setting8}

By the same way as Step \ref{step:1} in the previous subsection, we can take a map
\[ \theta \: \bigsqcup_{v\in V(X)}\R_{\rho_0}^\ast (v) \rightarrow \mathbb{Z}_{\geq 0} \]
satisfying the following conditions:
\begin{enumerate}
\item $\theta$ is $G$-invariant;
\item there exist $M\in \NN$ such that $\frac{1}{M}\theta $ approximates $\mu$, that is, $\frac{1}{M}\theta(T)$ is sufficiently close to $\mu (\SCyl (T))$ for every $T\in \R_{\rho}^\ast (v)$;
\item for any $u\in B_G$ and any $J \in \R_{\rho_0}^\ast(\id ,u)$ we have
\[ \sum_{\substack{T\in \R_{\rho_0}^\ast (\id )\\[1pt] T|_{B(\id ,u,\rho_0 )}=J}}\theta (T)=\sum_{\substack{T'\in \R_{\rho_0}^\ast(u)\\[1pt] T'|_{B(\id, u,\rho_0 )}=J}}\theta (T').\]
\end{enumerate}
We note that the same equation as the above follows for any adjacent $u,v\in V(X)$ and $J\in \R_{\rho_0}^\ast (u,v)$.

In addition, we can define $\theta (T)$ for every round-graph $T$ (with information of geodesics) included in $B(v, \rho_0)$ for some $v\in V(X)$ and we can assume that $\frac{1}{M}\theta (T)$ is also close to $\mu (\SCyl (T))$.

For an appropriate $r>0$ we want to define $\theta (T)$ for $h\in H_n$ and $T\in \R_{r}^\ast (h\F)$.
We note that $|T|\cap h\F\not=\emptyset$ by the definition but $|T| \cap h\F$ may contain no vertex.
Nevertheless we can take a vertex $v\in |T|\cap B(h\F ,1)$. Hence we need to see that $B(h\F ,r)$ is included in $B(v,\rho_0)$.
Moreover, in order to see that the $\theta(T)$ is independent of the choice of $v$ we need to consider a geodesic connecting two vertices of $B(h\F ,1)$, and for every vertex $w$ on the geodesic $B(h\F ,r)$ should be included in $B(w,\rho_0)$.

\begin{setting8}\label{set:5}
Assume that $\rho_0$ is sufficiently larger than the diameter of $\F$, which depends on $n$.
Since $\F$ is bounded, there exists a constant $c_\F>0$ depending on $\F$ such that $B(\F ,1)$ is $c_\F$-quasi-convex.
Set \[\rho_1:=\rho_0-\mathrm{diam}\F -c_\F -1.\]
\end{setting8}

For two vertices $v,v'\in B(\F,1)$ and any vertex $w$ on a geodesic $\ell$ joining $v$ to $v'$, we see that
$ B(w,\rho_0 )\supset B(\F , \rho_1 )$.
Therefore we can prove the following lemma by the same way as the proof of Lemma \ref{lem:restriction to H_n fundamental domain}.

\begin{lemma}\label{lem:restrict theta to rho_1}
For each $h\in H_n$ and $T\in \R_{\rho_1}^\ast(h\F )$ we define $\theta (T)$ by
\[ \theta (T) :=\sum_{\substack{T'\in \R_{\rho_0}^\ast (v)\\[1pt] T'|_{B(h\F ,\rho_1 )}=T}} \theta (T' ),\]
where $v$ is a vertex of $|T|\cap B(h\F, 1)$. Then the definition of $\theta (T)$ is independent of the choice of $v$ and we obtain an $H_n$-invariant map
\[ \theta \: \bigsqcup_{h\in H_n}\R_{\rho_1}^\ast (h \F)\rightarrow \mathbb{Z}_{\geq 0}.\]
Moreover, for any $u\in B_n$ and any $J\in \R_{\rho_1}(\F ,u\F)$ we have the following equation:
\[ \sum_{\substack{T\in \R_{\rho_1}^\ast (\F )\\[1pt] T|_{ B(\F, u\F ,\rho_1)=J} }}\theta (T)=\sum_{\substack{T'\in \R_{\rho_1}^\ast(u\F )\\[1pt] T'|_{B(\F ,u\F ,\rho_1)=J} }}\theta (T').\]
\end{lemma}

Following Step \ref{step:2} in the previous subsection, we construct a graph $(\gG,\iota )$ on $(H_n,X_n)$. Then the graph $(\gG ,\iota )$ satisfies the following conditions:
\begin{enumerate}
\item $\iota \: \gG\rightarrow X_n$ is an $H_n$-equivariant map;
\item the restriction of $\iota $ to each connected component of $\gG$ is injective;
\item $\# \iota ^{-1}( \id )<\infty $.
\end{enumerate}
Explicitly, 
\[ V(\gG ):= \{ v(h ,T, i)\} _{h\in H_n,\, T\in \R_{\rho_1}^\ast (h\F ),\, i=1,\dots ,\theta (T)}.\]
If two vertices $v(h_1,T_1) , v(h_2,T_2)$ of $V(\gG)$ are connected by an edge, then $h_1$ and $h_2$ are adjacent in $X_n$ and $T_1,T_2$ are $J$-connectable for some $J\in \R_{\rho_1}^\ast(h_1\F,h_2\F)$.
For $v(h,T)\in V(\gG)$ if there exists $h'$ adjacent to $h$ in $X_n$ such that $T\cap h'\F\not=\emptyset$, then $T|_{B(h\F,h'\F,\rho_1)}\in \R_r^\ast(h\F, h'\F)$ and there exists $T'\in \R_{\rho_1}^\ast(h'\F)$ such that $v(h,T)$ and $v(h',T')$ are connected by an edge in $\gG$.
The map $\iota$ maps $v(h,T)\in V(\gG)$ to $h\in X_n$.
Finally, we check that
\begin{align*}
\#\iota^{-1}(\id )=&\sum_{T\in \R_{\rho_1}^\ast (\F )}\theta (T)\\
	\leq &\sum_{v\in V(B(\F ,1))} \sum_{T\in \R_{\rho_0}^\ast (v)}\theta (T)=\# V(B(\F ,1) )\sum_{T\in \R_{\rho_0}^\ast (\id )}\theta (T)<\infty.
\end{align*}

We construct a graph $(|\gG |, |\iota |)$ on $X_{H_n}=\Phi^{-1}(CH_{H_n}\cap \HH^{(1)})(\subset X)$ from $(\gG ,\iota )$ by the same way as we did in Step \ref{step:3} in the previous subsection.
Explicitly, for each connected component $Y$ of $\gG$ we define a subgraph $|Y|$ of $X$ by
\[ |Y|:=\bigcup_{v(h,T)\in V(Y)}|T|\cap h\F \] 
and define $|\gG|$ to be the disjoint union of $|Y|$ over all connected component $Y$ of $\gG$. Note that $|Y|$ could be non-connected but $|Y|$ is a subgraph of $X$ although $h\F$ is just a subset of $X$ for $h\in H_n$.
Consider the case that an edge $e$ of $X_n$ is covered by $h_1\F, \dots , h_k\F$ for $h_1,\dots h_k \in H_n$. Then we can assume that $h_{i}$ and $h_{i+1}$ are adjacent for $i=1,\dots ,k-1$. If $Y$ contains a vertex $v(h_j, T_j) \in V(Y)$ with $|T_j| \supset e$, then there exists $v(h_i, T_i)\in V(Y)$ for $i=1,\dots ,j-1, j+1 ,\dots ,k$ such that $v(h_i ,T_i)$ and $v(h_{i+1}, T_{i+1})$ are adjacent in $Y$ for every $i=1,\dots ,k-1$.
Since $T_i$ and $T_{i+1}$ are connectable for every $i=1,\dots ,k-1$, $|T_i|$ includes $e$ for every $i$. Therefore $|Y|$ includes $e$.

The map $|\iota |$ is an $H_n$-equivariant map from $|\gG|$ to $X_{H_n}$ and the restriction of $|\iota |$ to $|Y|$ for each connected component $Y$ of $\gG$ is the inclusion map. Hence we will identify $|Y|$ with $|\iota|(|Y|)$.

Now, we want to prove a certain lemma corresponding to Lemma \ref{lem:path connectable on tree for fundamental domain}.
Note that Lemma \ref{lem:path connectable on tree for fundamental domain} deeply depends on the property that the space $X$ is a tree in the previous subsection.

Let $\phi$ be the inclusion map from $H_n$ to $X$ sending $h\in H_n$ to $h\in V(X)=G$. Since $H_n$ is a quasi-convex subgroup of $G$, we can extend $\phi$ to a quasi-isometric embedding from $X_n$ to $X$.

\begin{lemma}\label{lem:construct rho 2 from rho 1}
Assume that $\phi$ is an $(a ,c)$-quasi-isometric embedding for constants $a \geq 1, c\geq 0$, which depend on $n$. Let $Y$ be a connected component of $\gG$. Let $v=v(h,T),v'=v(h',T')\in V(Y)$.
Set
\[ \rho_2 :=\frac{\rho_1 -a(2\mathrm{diam} \F +c)(\mathrm{diam}\F +1)}{1+2a(\mathrm{diam}\F +1)}\]
and assume that $\rho_2>0$.
Then $T|_{B(h\F, \rho_2)}\in \R_{\rho_2}^\ast (h\F)$ and $T'|_{B(h'\F,\rho_2)}\in \R_{\rho_2}^\ast (h'\F)$ are connectable.
\end{lemma}
\begin{proof}
We denote by $d_{B_n}$ the path metric on $X_n=\Cay (H_n, B_n)$. We identify $Y$ with $\iota (Y)$, which is a subtree of $X_n$.
Take the geodesic $\ell$ from $v$ to $v'$ in $Y$, which passes through vertices $v_0=v, v_1,\dots ,v_m=v'$ in this order.
Note that $m=d_{B_n}(h,h')$.
Since $v_{i-1}=v(h_{i-1},T_{i-1}) ,v_i=v(h_i,T_i)$ are connected by an edge, $T_{i-1}$ and $T_i$ are $J_i$-connectable for some $J_i\in \R_{\rho_1}^\ast (h_{i-1}\F,h_i\F)$ for $i=1,\dots ,m$. This implies that the restriction of $T$ to
\[ U:=B(h_0\F, \rho_1)\cap B(h_1\F,\rho_1)\cap \cdots \cap B(h_m\F ,\rho_1 )\]
coincides with that of $T'$ to $U$.
Therefore it is sufficient to see that $B(h\F ,h'\F ,\rho_2)$ is included in $U$.

From the assumption we have
\[ \frac{1}{a}m-c\leq d(h,h') \leq a m +c.\]
Since $\F\ni \id$, $h\F$ and $h'\F$ contain $h$ and $h'$ respectively and so
\[ d(h\F ,h'\F )\geq \frac{1}{a}m-c-2\mathrm{diam}\F.\]
If $d(h\F ,h'\F) >2\rho_2$, then $B(h\F,h'\F, \rho_2)=\emptyset$ and $T|_{B(h\F, \rho_2)}$ and $T'|_{B(h'\F,\rho_2)}$ are connectable.
Therefore it is sufficient to consider the case that
\[ \frac{1}{a}m-c-2\mathrm{diam}\F \leq 2\rho_2 ,\]
that is, $m\leq a(2\rho_2 +2\mathrm{diam}\F +c)$.

Since $h_{i-1},h_i$ are adjacent, for any $\alpha >0$ we have
\[ B(h_{i-1}\F, \alpha-\mathrm{diam}\F -1)\subset B(h_i\F ,\alpha)\]
for every $i=1,\dots ,m$. 
Hence
\begin{align*}
B(h_0\F, \rho_1 -m(\mathrm{diam}\F +1) )
&\subset B(h_1\F , \rho_1 -(m-1)(\mathrm{diam}\F +1) )\\
&\ \, \vdots \\
&\subset B(h_m\F , \rho_1 ),
\end{align*}
which implies that
\[ B(h_0\F ,\rho_1 -m\mathrm{diam}\F -m )\subset U.\]
Since $m\leq a(2\rho_2 +2\mathrm{diam}\F +c)$, we have
\[ \rho_1-m(\mathrm{diam}\F+1 )\geq \rho_1-a (2 \rho_2 +2\mathrm{diam}\F +c)(\mathrm{diam}\F +1).\]
We can see that 
\[ \rho_1-a (2 \rho_2 +2\mathrm{diam}\F +c)(\mathrm{diam}\F +1)=\rho_2.\]
In fact,
\begin{align*}
&\rho_1-a (2 \rho_2 +2\mathrm{diam}\F +c)(\mathrm{diam}\F +1)-\rho_2\\
=&\rho_1-a (2\mathrm{diam}\F +c)(\mathrm{diam}\F +1)-\rho_2 (1+2a(\mathrm{diam}\F +1)) \\
=&0.
\end{align*}
Hence 
\[ B(h\F,h'\F, \rho_2)\subset B(h_0\F ,\rho_2 )\subset U.\]
This completes the proof.
\end{proof}

\begin{setting8}\label{set:6}
We take $\rho_2$ in the above lemma.
Recall that the length of a $\delta$-quasi-geodesic connecting two points with distance $d$ is smaller than or equal to $\delta' d+ \delta'$.
We also assume that $\rho_2':=\delta' \rho_2+\delta' \leq \rho_1-1$.
\end{setting8}

Now, we prove the following lemma corresponding to Lemma \ref{lem:neighborhood of realization of Gamma}.

\begin{lemma}\label{lem:rho_2 neighborhood of an interior point is good}
Let $Y$ be a connected component of $\gG$. Let $v(h,T)\in V(Y)$, $v\in |T|\cap B(h\F,1)$. 
Assume that $B(v,\rho_2')\subset X_{H_n}=\Phi^{-1}(CH_{H_n}\cap \HH^{(1)})$.
Then we have
\[ |Y|\cap B(v,\rho_2 )=|T|\cap B(v,\rho_2 ).\]
Moreover, for the connected component $Z$ of $|Y|$ containing $v$,
\[ Z\cap B(v,\rho_2 )=|Y|\cap B(v,\rho_2 )=|T|\cap B(v,\rho_2 ).\]
\end{lemma}
\begin{proof}
Take $x\in |Y|\cap B(v,\rho_2)$. There exists $v(h_0,T_0)\in V(Y)$ such that $x \in |T_0|\cap h_0\F$.
Since $T|_{B(h\F ,\rho_2)}$ and $T_0|_{B(h_0\F,\rho_2)}$ are connectable, we have
\[ x\in |T_0|\cap B(h\F , h_0\F, \rho_2 )=|T|\cap B(h\F , h_0\F , \rho_2 ).\]
Hence $x\in |T|\cap B(v,\rho_2)$.

Take $x\in |T|\cap B(v,\rho_2 )$ and $S\in \SCyl (T)$. Then $x\in CH(S)\cap B(v,\rho_2)$.
The point is that we can take a $\delta$-quasi-geodesic $\ell$ joining $v$ to $x$ in $CH(S)$.
Hence $\ell$ is included in $B(v,\rho_2') (\subset B(h\F ,\rho_2'+1))$, which implies that 
\[ \ell\subset |T|\cap B(h\F ,\rho_2'+1)=CH(S)\cap B(h\F ,\rho_2'+1) . \] 
From the construction of $\gG$ there exists a path of vertices $v(h_0, T_0)=v(h,T),\dots ,v(h_m,T_m)$ in $Y$ such that $\ell$ passes through $h_i\F$ in this order and $x\in h_m\F$.
Since $T|_{B(h\F ,\rho_2)}$ and $T_m|_{B(h_m\F,\rho_2)}$ are connectable, we have
\[ x \in |T|\cap B(h\F, h_m\F ,\rho_2) =|T_m|\cap B(h\F, h_m\F ,\rho_2). \]
This implies that $x\in |T_m|\cap h_m\F \subset |Y|$.

Note that $v\in |T|\cap B(v,\rho_2 )=|Y|\cap B(v,\rho_2 )$, particularly, $v\in |Y|$.
From the above, for any $x\in |T|\cap B(v,\rho_2)$ there exists a path $\ell$ joining $v$ to $x$ in $|Y|$, which implies that $x\in Z\cap B(v,\rho_2)$ for the connected component $Z$ of $|Y|$ containing $v$. Hence $Z\cap B(v,\rho_2 )=|Y|\cap B(v,\rho_2 )$.
\end{proof}

In the above proof, we took a $\delta$-quasi-geodesic $\ell$ in $CH(S)$ connecting two points of $CH(S)$. This is the reason why we introduce the notion of the convex hull and define the round-graph by using the convex hull instead of the weak convex hull.

Let $Y$ be a connected component of $\gG$.
Take $v(h,T)\in V(Y)$ and assume that there exists $\gamma \in \gamma_T$ such that $\gamma \cap h\F \not=\emptyset$.
Take $v(h',T') \in V(Y)$ adjacent $v(h,T)$ in $Y$ such that $\gamma \cap h' \F \not=\emptyset$ .
Then $T|_{B(h\F ,\rho_2)}$ and $T'|_{B(h'\F ,\rho_2)}$ are $J$-connectable for $J=T|_{B(h\F ,h'\F ,\rho_2)}$.
This implies that there exists $\gamma'\in \gamma_{T'}$ such that
\[ \gamma \cap B(h\F, h'\F ,\rho_2 )=\gamma' \cap B(h\F ,h'\F ,\rho_2 ) \ (\not=\emptyset ).\]
Therefore we can extend $\gamma \cap B(h\F ,\rho_2)$ by connecting $\gamma \cap B(h\F,\rho_2 )$ to $\gamma \cap B(h'\F ,\rho_2 )$ and we can perform this operation over and over until the extension of $\gamma \cap B(h\F ,\rho_2)$ meets the boundary of $CH_{H_n}$.

By the definition, there exists a geodesic line $\ell$ such that $\ell \cap B(h\F ,\rho_1)=\gamma$, which implies that the extension of $\gamma \cap B(h\F ,\rho_2)$ is a $2\rho_2$-local geodesic, that is, every sub-arc with length less than or equal to $2\rho_2$ is a geodesic segment. It is known that an $L$-local geodesic for $L>0$ will be a $\delta$-quasi-geodesic if $L$ is larger than a constant depending on $\delta$.
We can assume that $\rho_2$ is sufficiently large such that the extension of $\gamma \cap B(h\F ,\rho_2)$ is a $\delta$-quasi-geodesic.
Note that the extension of $\gamma \cap B(h\F ,\rho_2)$ will be a $\delta$-quasi-geodesic line if it does not meet the boundary of $X_{H_n}$. We call each extension of $\gamma \cap B(h\F ,\rho_2)$ for $\gamma \in \gamma_T$ a \ti{$Y$-quasi-geodesic}.
If the extension of $\gamma \cap B(h\F ,\rho_2)$ is a $\delta$-quasi-geodesic line, then we call it a $Y$-quasi-geodesic line.

In order to apply Lemma \ref{lem:quasi-geodesic line and quasi-convex implies good condition} to each connected component of $|\gG|$, we prove that every connected component of $|\gG|$ is $\delta'$-quasi-convex.

\begin{lemma}\label{lem:a connected component of Y is quasi-convex}
Let $Y$ be a connected component of $\gG$ and $Z$ a connected component of $|Y|$.
Then $Z$ is a $\delta'$-quasi-convex subgraph of $X$.
\end{lemma}
\begin{proof}
Let $x, y\in Z$. Consider a shortest path $\ell$ joining $x$ to $y$ in $Z$. We prove that $\ell$ is a $\delta$-quasi-geodesic in $X$ and then 
$Z$ is $\delta'$-quasi-convex from the stability of quasi-geodesics.
In order to see that $\ell$ is $\delta$-quasi-geodesic, it is sufficient to see that for a large constant $L>0$ depending on $\delta$, $\ell$ is $L$-local $\delta$-quasi-geodesic, that is, every sub-arc of $\ell$ with length less than or equal to $L$ is $\delta$-quasi-geodesic. We can assume that $\rho_2$ is much larger than $L$.
Then it is sufficient to consider the case that $d(x,y)\leq L (< \rho_2)$ and prove that there exists a $\delta$-quasi-geodesic joining $x$ to $y$.

Take $v(h,T)\in Y$ such that $x\in h\F$, which implies that $y\in B(h\F ,L) $.
Take $S\in \SCyl (T)$, which implies that $CH(S)\cap B(h\F ,\rho_1)=|T|$. Then we can take a $\delta$-quasi-geodesic $\gamma$ joining $x$ to $y$ in $CH(S)$ and $\gamma$ is included in $CH(S)\cap B(x,\rho_2')$. 

Now it is sufficient to see that $\gamma $ is included in $X_{H_n}$.
Actually, if $\gamma$ is included in $X_{H_n}$, then we can see that $\gamma$ is included in $Z$ by the same argument in Lemma \ref{lem:rho_2 neighborhood of an interior point is good}, which is the desired conclusion.

We prove that $x$ and $y$ are included in the same connected component of the intersection of $CH(S)$ and $X_{H_n}$.
Then we can take $\gamma$ such that $\gamma$ is included in $CH(S)\cap X_{H_n}$ since every boundary component of $CH(S)$ is a $\delta$-quasi-geodesic line and we can consider a $\delta$-quasi-geodesic traveling along the boundary of $CH(S)$. Hence $\gamma$ is included in $Z$, which is the desired conclusion.

To obtain a contradiction, suppose that the connected component of $CH(S)\cap X_{H_n}$ containing $x$ is different from that containing $y$. Then a path $\ell$ joining $x$ to $y$ need to ``take a detour'', that is, a geodesic $[x,y]$ joining $x$ to $y$ in $X_{H_n}$ must meets a boundary component $b$ of $CH(S)$ at $z$.
Take $v(h',T')\in V(Y)$ such that $z\in h'\F$, which implies that $b\cap h'\F\not=\emptyset$.
We consider the extension of $b\cap B(h'\F,\rho_2)$ by the same way as we did in the above in order to obtain a $Y$-quasi-geodesic.
The extension of $b\cap B(h'\F,\rho_2)$ is a $\delta$-quasi-geodesic and can be considered as a boundary component of $Z$.
Then we see that a path joining $x$ to $y$ in $X_{H_n}$ must cross the extension of $b\cap B(h'\F,\rho_2)$, which implies that there exists no path joining $x$ to $y$ in $Z$, a contradiction.
\end{proof}

Now, we assume the following condition for a while:

\textbf{Assumption} $(\ast )$:
For every $v(h,T)\in V(\gG)$ and every $\gamma \in \gamma_T$ with $\gamma \cap h\F \not=\emptyset$, the extension of $\gamma \cap B(h\F ,\rho_2)$ is a $\delta$-quasi-geodesic line if $\gamma$ (or its extension) contains a point $x$ with $B(x,C_0) \subset X_{H_n}$ for a constant $C_0>0$ independent of $n$.

Set
\[ \eta_\gG:= \sum_{Z\in \Comp (|\gG | )}\delta_{Z(\infty)} \in \SC (H_n).\]
Under Assumption $(\ast )$ we prove that
\[ \nu:=\frac{1}{nM} \iota_{H_n}(\eta_\gG ) \]
belongs to the open neighborhood $U(f_1,\dots ,f_k; \varepsilon_\mu )$ of $\mu$ for a sufficiently large $n$ by using Lemma \ref{lem:describe a neighborhood of a given subset current}.
Note that $\eta_\gG$ is a subset current on $H_n$.
In the case that $|\gG|$ does not satisfy the condition in Assumption $(\ast )$, we construct $\hat{\gG}$ from $|\gG|$ by a similar way as we did in Step \ref{step:4} in the previous subsection such that $\hat{\gG}$ satisfies the condition in Assumption $(\ast )$.
During the construction of $\hat{\gG}$ the constant $C_0$ will plays an essential role.

Let $Y\in \Comp (\gG)$, $Z\in \Comp (|Y|)$, $g \in V(Z)$. Take $v(h,T)\in V(Y)$ such that $g \in |T|\cap h\F$.
Assume that $B(g,\rho_2')\subset X_{H_n}$.
Then $Z\cap B(g,\rho_2)=|T|\cap B(g,\rho_2)$ from Lemma \ref{lem:rho_2 neighborhood of an interior point is good}.
By the definition of $T$, there exists a constant $\delta_1>0$ depending on $\delta$ such that for every $x\in B(g,\rho_2- \delta_1)$ there exists $\gamma \in \gamma_T$ such that $d(x,\gamma ) \leq \delta_1 $.

\begin{setting8}\label{set:7}
Set $\rho_3:=\rho_2-\delta_1-C_0$ and assume that $\rho_3>0$.
\end{setting8}

\begin{lemma}[Under Assumption $(\ast )$]\label{lem:with assumption and Z infty in a neighborhood of SCyl}
Assume that $\rho_3$ is sufficiently large. Let $Z\in \Comp (|\gG |)$, $g\in V(Z)$.
If $B(g,\rho_2')$ is included in $X_{H_n}$, then 
\[ g^{-1}Z(\infty )\in B_{\H}(\SCyl ( g^{-1}Z\cap B(\id ,\rho_3 )) ,\varepsilon_1 ).\]
\end{lemma}
\begin{proof}
Take $Y\in \Comp (\gG)$ such that $Z\in \Comp (|Y|)$.
Take $v(h,T)\in V(Y)$ such that $g \in |T|\cap h\F$.
Then, for $x\in Z\cap B(g,\rho_3)$ we can take $\gamma \in \gamma_T$ such that $d(x,\gamma )\leq \delta_1$, and then $\gamma$ contains a point $y$ such that $d(x,y)\leq \delta_1$ and $B(y, C_0) \subset X_{H_n}$.
By considering a path from $g$ to $y$ included in $|T|\cap B(h\F, \rho_2)$, we can take $v(h',T')\in V(Y)$ and $\gamma'\in \gamma_{T'}$ such that $y\in \gamma \cap h'\F \cap |T'|$ and 
\[ y\in \gamma \cap B(h\F , h'\F ,\rho_2 )=\gamma' \cap B(h\F ,h' \F ,\rho_2 ).\]
Hence from Assumption $(\ast )$, the extension of $\gamma'\cap B(h'\F,\rho_2)$ will be a $Y$-quasi-geodesic line $\ell$ in $Z$, and $d(x,\ell)\leq \delta_1$.

As a result, we see that for every $x\in Z\cap B(g,\rho_3)$ there exists a $Y$-quasi-geodesic line $\ell$ in $Z$ such that $d(x,\ell )\leq \delta_1$.
Now, we can apply Lemma \ref{lem:quasi-geodesic line and quasi-convex implies good condition} to $Z$ and we can see that for a constant $a>0$ depending only on $\delta$ (and $\delta_1$), we have
\[ Z\cap B(g,\rho_3 )\underset{a}{\sim }WC(Z(\infty ))\cap B(g,\rho_3 ) .\]
Note that $Z\cap B(g,\rho_3)=|T|\cap B(g,\rho_3)\in \R_{\rho_3}(g)$ by Lemma \ref{lem:rho_2 neighborhood of an interior point is good}.
Now, we assume that $\rho_3$ is sufficiently large to apply Lemma \ref{lem:difference between Y and WC(Y infty)} to the constant $\varepsilon_1>0$ and $g^{-1}Z\cap B(\id, \rho_3) $.
The constant $r$ in Lemma \ref{lem:difference between Y and WC(Y infty)} depends on the base point $y$ but as long as we use $\id $ as the base point we do not need to consider the problem.
Therefore we have
\[ g^{-1}Z(\infty )\in B_{\H}(\SCyl ( g^{-1}Z\cap B(\id, \rho_3 )) ,\varepsilon_1 ).\]
This completes the proof.
\end{proof}

Take a complete system of representatives $\gL_0$ of $G/H_n$.
To apply Lemma \ref{lem:describe a neighborhood of a given subset current} to $\nu$ we consider the restriction of $\nu$ to $A(\id ,r_\mu)$.
Set
\[ \gL_1:= \{ g\in \gL_0 \mid g CH(\gL(H_n) )\cap B(\id ,r_\mu )\not =\emptyset \}, \]
which is a finite set.
Note that $CH(\gL (H_n))$ is the convex hull of $\gL (H_n)$ in $X$ and we write $CH_{H_n}$ to represent the convex hull of $\gL (H_n)$ in $\HH$.
Then
\[ \iota_{H_n}(\eta_\gG )|_{A(\id, r_\mu )}=\sum_{g\in \gL_1 }g_{\ast}(\eta_\gG )|_{A(\id ,r_\mu )}.\]

Recall that every boundary component of $X_{H_n}$ is a geodesic line, and so $CH(\gL (H_n))$ includes $X_{H_n}$. Hence for $g\in \gL_0$ if $gX_{H_n}\cap B(\id ,r_\mu)\not=\emptyset$, then $g\in \gL_1$.
Recall that $\F$ includes exactly $n$ vertices of $X_{H_n}$, each of which corresponds to a vertex of $\gS_n$. Let $g_1=\id ,g_2^{-1},\dots , g_n^{-1}$ be the vertices of $\F$. We can assume that $g_1,\dots ,g_n\in \gL_1$.

\begin{lemma}\label{lem:gl 1 and g1 dots gn estimate}
The sequence
\[ \frac{1}{n}\# \left( \gL_1 \setminus \{ g_1,\dots ,g_n\}\right) \]
tends to $0$ when $n\rightarrow \infty$.
\end{lemma}
\begin{proof}
First, we have
\begin{align*}
\#\gL_1 =&\# \{ gH_n\in G/H_n |\ gCH(\gL(H_n))\cap B(\id ,r_\mu )\not =\emptyset \} \\
=&\# \{ gH_n\in G/H_n |\ gB(CH(\gL(H_n)),r_\mu )\ni \id \} \\
=&\# \{ gH_n\in G/H_n |\ B(CH(\gL(H_n)),r_\mu )\ni g^{-1} \} \\
=&\# V(H_n\backslash B(CH(\gL(H_n)), r_\mu )).
\end{align*}

Note that $CH(\gL(H_n))\underset{\delta'}{\sim} X_{H_n}$. 
From the definition of $X_{H_n}$ the quotient $H_n\backslash X_{H_n}$ is isomorphic to the $1$-skeleton of $\gS_n$ and includes $n$ vertices.
Moreover, the degree of every vertex of $\gS_n$ except $\tilde{x_0}^n$ coincides with the degree of $\id$ in $X$, denoted by $\mathrm{deg}_X(\id )$.
Since 
\[ B(CH(\gL(H_n)), r_\mu ))\subset B(X_{H_n},r_\mu +\delta') = X_{H_n}\cup B(\partial X_{H_n}, r_\mu +\delta') \]
and, we have
\begin{align*}
&\# V(H_n\backslash B(CH(\gL(H_n)), r_\mu ))-\# V(H_n\backslash X_{H_n}) \\
\leq &\# V(H_n\backslash B(X_{H_n}, r_\mu +\delta ' ))-\# V(H_n\backslash X_{H_n}) \\
\leq &\big( \mathrm {deg}_X( \id ) \big)^{r_\mu +\delta '},
\end{align*}
which implies
\[ \frac{1}{n}\# \left( \gL_1 \setminus \{ g_1,\dots ,g_n\}\right)\leq \frac{1}{n}\big( \mathrm {deg}_X( \id ) \big)^{r_\mu +\delta '}.\]
This proves our claim.
\end{proof}

\begin{setting8}\label{set:8}
Set $\gL :=\{ g_i \mid B(g_i^{-1} ,\rho_2' )\subset X_{H_n} (i=1,\dots, n)\} $.
\end{setting8}

\begin{remark}[About constants $\rho_0,\rho_1,\rho_2,\rho_2',\rho_3$]
Since we need to take sufficiently large $\rho_3$ in Setting \ref{set:7}, which depends on the neighborhood $U(f_1,\dots ,f_k;\varepsilon_\mu)$ of $\mu$ and on constants related to $\delta$, we determine $\rho_3,\rho_2, \rho_1 $ and $\rho_0$ in this order.
The point is that $\rho_3, \rho_2, \rho_2'$ are independent of $n$.
\end{remark}

\begin{lemma}\label{lem:gL 1 and gL estimate}
The sequence
\[ \frac{1}{n}\# \left( \gL_1 \setminus \gL \right) \]
tends to $0$ when $n\rightarrow \infty$.
\end{lemma}
\begin{proof}
Recall that $\Phi$ is a quasi-isometric map from $X$ to $\HH$. Then the restriction of $\Phi$ to $X_{H_n}$ is also a quasi-isometric map to $CH_{H_n}$. There exists a constant $c$ depending on $\rho_2'$ and $\Phi$ such that if $B(g_i^{-1} ,\rho_2 ) \not\subset X_{H_n}$, then $\Phi (g_i^{-1})$ is contained in the $c$-neighborhood of the boundary of $CH_{H_n}$. By considering the quotient space $\gS_n=H_n\backslash CH_{H_n}$ and the $c$-neighborhood of the boundary component $\tilde{c_0}^n$ of $\gS_n$, the number of $g_i$ such that $B(g_i^{-1},\rho_2')\not\subset X_{H_n}$ is bounded by a constant depending on $c$ and independent of $n$. This proves our claim.
\end{proof}

From the above setting, we have
\begin{align*}
\iota_{H_n}(\eta_\gG )|_{A(\id, r_\mu )} 
=&\sum_{g\in \gL_1 }g_{\ast}(\eta_\gG )|_{A(\id ,r_\mu )}\\
=&\left( \sum_{g\in \gL_1\setminus \gL } +\sum_{g\in \gL} \right) g_{\ast}(\eta_\gG )|_{A(\id ,r_\mu )}.
\end{align*}
We mainly consider the sum taken over $g\in \gL$ and 
\begin{align*}
\sum_{g\in \gL} g_{\ast}(\eta_\gG )|_{A(\id ,r_\mu )}
=&\sum_{g\in \gL }\left( \sum_{Z\in \Comp (|\gG |)} \delta_{gZ(\infty )}\right) \Bigg|_{A(\id , r_\mu )}\\
=&\sum_{g\in \gL } \sum_{\substack{Z\in \Comp (|\gG |) \\[1pt] gZ(\infty )\in A(\id ,r_\mu )}} \delta_{gZ(\infty )}.
\end{align*}

Now, we consider $Z\in \Comp (|\gG |)$ with $gZ(\infty )\in A(\id, r_\mu )$ for $g\in \gL$. For $Z$ we denote by $Y_Z$ the connected component $Y$ of $\gG$ such that $Z$ is a connected component of $|Y|$.
Recall that from Lemma \ref{lem:a connected component of Y is quasi-convex}, $Z$ is $\delta '$-quasi-convex.
Hence 
\[ WC (Z(\infty ))\subset B(Z, \delta' ),\]
and we can take a constant $\alpha>0$ depending on $\delta$ such that
\[ CH(Z(\infty ))\subset B(Z, \alpha ).\]
Since $gCH(Z(\infty ))\cap B(\id ,r_\mu )\not=\emptyset $, we see that $Z\cap B(g^{-1}, r_\mu +\alpha )\not=\emptyset$.

\begin{setting8}[Under Assumption $(\ast )$]\label{set:9}
For the constant $\alpha$ in the above, we set $r_\mu':=r_\mu +\alpha$, which corresponds to the constant $r_\mu'$ in Lemma \ref{lem:describe a neighborhood of a given subset current}.
We assume that $\rho_2 \geq 2r_\mu '+\rho$. Recall that we fixed $\rho$ in Setting \ref{set:2}.
\end{setting8}

The following lemma does not depend on Assumption $(\ast )$.

\begin{lemma}\label{lem: intersection will be round graph}
Let $g\in \gL$ and $Z\in \Comp (|\gG |)$ with $gZ(\infty )\in A(\id, r_\mu )$.
Then $gZ \cap B(\id , r_\mu' +\rho )$ is an element of $\R_\rho (B(\id ,r_\mu' ))$.
\end{lemma}
\begin{proof}
Note that $Z\cap B(g^{-1}, r_\mu')$ contains a vertex $g_0$ since $Z$ is a subgraph of $X$. Then there exists $v(h_0,T_0)\in V(Y_Z)$ such that $g_0\in h_0\F \cap |T_0|$. Hence $gg_0\in g|T_0|\cap B(\id, r_\mu')$.

Since $\rho_2 \geq 2r_\mu '+\rho$, we have
\[ B(g^{-1},r_\mu' +\rho )\subset B(g_0 ,2r_\mu '+\rho )\subset B(h_0\F, \rho_2 ). \]
Since $g\in \gL$, we have $B(g^{-1}, \rho_2')\subset X_{H_n}$.
By Lemma \ref{lem:rho_2 neighborhood of an interior point is good} we have
\[ Z\cap B(g^{-1} ,r_\mu' +\rho ) =|T_0|\cap B(g^{-1}, r_\mu'+\rho ).\]
Hence
\[ gZ \cap B(\id , r_\mu' +\rho )=g|T_0|\cap B(\id ,r_\mu' +\rho ),\]
which is an element of $\R_{\rho}(B(\id ,r_\mu '))$.
\end{proof}

From the above lemma, we have
\begin{align*}
&\sum_{g\in \gL} g_{\ast}(\eta_\gG )|_{A(\id ,r_\mu )} \\
=&\sum_{g\in \gL } \sum_{\substack{Z\in \Comp (|\gG |) \\[1pt] Z\cap B(\id , r_\mu')\not=\emptyset}} \delta_{gZ(\infty )}|_{A(\id ,r_\mu )} \\
=&\sum_{g\in \gL } \sum_{T\in \R_\rho (B(\id ,r_\mu' ))}\sum_{\substack{Z\in \Comp (|\gG |) \\[1pt] gZ \cap B(\id , r_\mu' +\rho )=T}} \delta_{gZ(\infty )}|_{A(\id ,r_\mu )}.
\end{align*}
For each $T\in \R_{\rho}(B(\id, r_\mu'))$ set
\[ \iota_{H_n}(\eta_\gG )_T:=\sum_{g\in \gL } \sum_{\substack{Z\in \Comp (|\gG |) \\[1pt] gZ \cap B(\id , r_\mu' +\rho )=T}} \delta_{gZ(\infty )}. \]
Then
\[ \sum_{g\in \gL} g_{\ast}(\eta_\gG )|_{A(\id ,r_\mu )} =\sum_{T\in \R_\rho (B(\id , r_\mu '))} \iota_{H_n}(\eta_\gG )_T|_{A(\id ,r_\mu )}.\]

For every $T\in \R_{\rho}(B(\id , r_\mu'))$ we can define $\theta(T)$ by the same way as we did in Lemma \ref{lem:restriction to H_n fundamental domain} (and \ref{lem:restrict theta to rho_1}). Explicitly, for some vertex $u\in T\cap B(\id, r_\mu') $
\[ \theta (T)=\sum_{\substack{T'\in \R_{\rho_0}^\ast (u)\\[1pt] |T'|\cap B(\id , r_\mu' +\rho ) =T}}\theta (T'), \]
which is independent of the choice of $u$. Moreover, for every $g\in G$ we can define $\theta (gT)$ by the same way, and we have $\theta(gT)=\theta(T)$.
Note that $T\in \R_{\rho}(B(\id , r_\mu'))$ does not include the information of geodesics and $T\cap B(\id, r_\mu')\not=\emptyset$ implies that $T\cap B(\id, r_\mu')$ contains a vertex.
We can see that $\frac{1}{M}\theta (T)$ is also close to $\mu (\SCyl (T))$ for $T\in \R_{\rho}(B(\id , r_\mu'))$, since we take $\theta $ after $r_\mu', \rho$.

\begin{lemma}[Under Assumption $(\ast )$]\label{lem:under assumption iota Hn eta gG is good}
For each $T\in \R_{\rho}(B(\id ,r_\mu'))$ we have
\[ \mathrm{supp} (\iota_{H_n}(\eta_\gG )_T)\subset \ol{B_{\H}(\SCyl (T) ,\varepsilon_1 )}.\]
Moreover,
\[ | \iota_{H_n}(\eta_\gG )_T |=\# \gL \cdot \theta (T).\]
\end{lemma}
\begin{proof}
Fix a vertex $u \in T\cap B(\id, r_\mu')$. For $g\in \gL$ we consider $Z\in \Comp (|\gG |)$ satisfying the condition that $gZ \cap B(\id , r_\mu' +\rho )=T$.
Note that $u\in T \cap B(\id,r_\mu')=gZ \cap B(\id, r_\mu')$.
Then $g^{-1}u\in Z\cap B(g^{-1}, r_\mu')$ and take $v(h_0,T_0)\in V(Y_Z)$ such that $h_0\F \cap |T_0|\ni g^{-1}u$.
We see that 
\[ B(g^{-1}u, r_\mu'+\rho )\subset B(g^{-1}, 2r_\mu' +\rho)\subset B(g^{-1}, \rho_2)\subset X_{H_n}\]
since $g\in \gL$. Hence by Lemma \ref{lem:rho_2 neighborhood of an interior point is good}
\[ |T_0| \cap B(g^{-1} , r_\mu '+ \rho )=Z\cap B(g^{-1}, r_\mu ' +\rho )=g^{-1}T,\]
which implies
\[ g|T_0|\cap B(\id ,r_\mu' +\rho )= gZ \cap B(\id , r_\mu' +\rho )=T.\]
From Lemma \ref{lem:with assumption and Z infty in a neighborhood of SCyl}
\[ gZ(\infty )\in B_{\H}(\SCyl ( gZ\cap B(\id,\rho_3 )) ,\varepsilon_1 ).\]
We can assume that $\rho_3 \geq r_\mu '+\rho$.
Then $T=(gZ \cap B(\id ,\rho_3))\cap B(\id , r_\mu' +\rho )$, and so we have
\[ \SCyl (gZ \cap B(\id ,\rho_3))\subset \SCyl (T) ,\]
which implies that
\[ gZ(\infty )\in B_{\H}(\SCyl (T) ,\varepsilon_1 ).\]
Therefore we obtain
\[ \mathrm{supp} (\iota_{H_n}(\eta_\gG )_T)\subset \ol{B_{\H}(\SCyl (T) ,\varepsilon_1 )}.\]

Now, we calculate $| \iota_{H_n}(\eta_\gG )_T |$.
From the above argument, for $g\in \gL$ and $Z\in \Comp (|\gG |)$, we have $gZ \cap B(\id , r_\mu' +\rho )=T$ if and only if
for $h_0\in H_n$ with $h_0\F \ni g^{-1}u$ there exists $v(h_0, T_0)\in V(Y_Z)$ such that 
\[ |T_0|\cap B(g^{-1},r_\mu' +\rho)=g^{-1}T.\]
Actually, by Lemma \ref{lem:rho_2 neighborhood of an interior point is good} $|T_0|\cap B(g^{-1},r_\mu' +\rho)=g^{-1}T$ implies $gZ \cap B(\id , r_\mu' +\rho )=T$.
Note that $h_0$ depends on $g$.
Therefore we have
\begin{align*}
&| \iota_{H_n}(\eta_\gG )_T |\\
=&\sum_{g\in \gL }\# \{ Z\in \Comp (|\gG |) \mid gZ \cap B(\id , r_\mu' +\rho )=T\} \\
=&\sum_{g\in \gL }\# \{ Z\in \Comp (|\gG |) \ |\\
	&\quad \exists v(h_0, T_0)\in V(Y_Z) \text{ s.t. } h_0\F \ni g^{-1}u \text{ and }|T_0|\cap B(g^{-1}, r_\mu '+\rho )=g^{-1}T \} \\
=&\sum_{g\in \gL }\ \sum_{\substack{g^{-1}u\in h_0\F, T_0\in \R_{\rho_1}^\ast (h_0\F ) \\[1pt] |T_0|\cap B(\id , r_\mu' +\rho )=g^{-1}T}}\theta (T_0)\\
=&\sum_{g\in \gL }\ \sum_{\substack{g^{-1}u\in h_0\F, T_0\in \R_{\rho_1}^\ast (h_0\F ) \\[1pt] |T_0|\cap B(\id , r_\mu' +\rho )=g^{-1}T}}\ 
\sum_{\substack{T'\in \R_{\rho_0}^\ast (g^{-1}u ) \\[1pt] T'|_{B(h_0\F ,\rho_1 )}=T_0} } \theta (T')\\
=&\sum_{g\in \gL }\ \sum_{\substack{T'\in \R_{\rho_0}^\ast (g^{-1}u ) \\[1pt] |T'|\cap B(\id , r_\mu' +\rho )=g^{-1}T}}\theta (T')\\
=&\sum_{g\in \gL }\theta (g^{-1}T)=\# \gL \cdot \theta (T).
\end{align*}
This completes the proof.
\end{proof}

For $g\in V(X)=G$ we set
\[ \theta (g)=\sum_{T\in \R_{\rho_0}^\ast (g)}\theta (T).\]
Then we can see that $\theta (g ) =\theta (\id )$ for every $g\in V(X)$.
Note that
\[ \bigsqcup_{T\in \R_{\rho_0}^\ast (g)}\SCyl (T) =A(g, 0)=\{ S\in \H (\partial G)\mid CH(S)\ni g\}. \]

Now, we consider the other part of $\iota_{H_n}(\eta_\gG )$.
Let $g\in \gL_1\setminus \gL$. Then
\begin{align*}
&\big| (g_\ast (\eta_\gG ) )|_{A(\id ,r_\mu )} \big|\\
=&\eta_\gG (A(g^{-1} ,r_\mu ))\\
=&\# \{ Z\in \Comp (|\gG |) \mid CH(Z(\infty ))\cap B(g^{-1}, r_\mu )\not=\emptyset \} \\
\leq &\# \{ Z\in \Comp (|\gG |) \mid Z\cap B(g^{-1}, r_\mu' )\not=\emptyset \} \\
\leq &\sum_{v\in V(B(g^{-1},r_\mu '))}\# \{ Z\in \Comp (|\gG |) \mid Z\ni v \} \\
= &\sum_{v\in V(B(g^{-1},r_\mu '))}\# \{ v(h_v,T )\in V(\gG ) \mid |T|\cap h_v\F \ni v \} \\
\leq &\sum_{v\in V(B(g^{-1},r_\mu '))}\sum_ {\substack{T\in \R_{\rho_1}^\ast (h_v\F ) \\[1pt] |T|\cap h_v\F \ni v}}\theta (T)\\
=&\sum_{v\in V(B(g^{-1},r_\mu '))}\sum_ {T\in \R_{\rho_0}^\ast (v)}\theta (T)\\
=&\# V(B(\id ,r_\mu ')) \theta (\id ).
\end{align*}
Since $\theta (\id )$ is close to $\mu (A(g,0))$, we can see that $\big| (g_\ast (\eta_\gG ) )|_{A(\id ,r_\mu )} \big|$ is bounded by a constant independent of $n$.

For $T\in \R_\rho (B(\id ,r_\mu '))$ set
\[ \nu_T :=\frac{1}{nM}\iota_{H_n}(\eta_\gG )_T\]
and
\[ \nu ':=\frac{1}{nM} \sum_{g\in \gL_1\setminus \gL} g_\ast (\eta_\gG )|_{A(\id ,r_\mu )}. \]
Then we have
\begin{align*}
\nu|_{A(\id ,r_\mu )}=&\frac{1}{nM} \iota_{H_n}(\eta_\gG )|_{A(\id ,r_\mu )} \\
=&\sum_{T\in \R_\rho ( B(\id, r_\mu' ))}\nu _T|_{A(\id ,r_\mu )}+\nu' .
\end{align*}
Now, we prove that for a sufficiently large $n\in \NN$, $\nu \in U(f_1,\dots ,f_k;\varepsilon_\mu)$ by using Lemma \ref{lem:describe a neighborhood of a given subset current}.
From Lemma \ref{lem:under assumption iota Hn eta gG is good} for every $T\in \R_\rho ( B(\id, r_\mu' ))$ we have
\[ \mathrm{supp }\nu_T \subset \ol{B_{\H}(\SCyl (T) ,\varepsilon_1 )}\]
and
\begin{align*}
&\big| |\nu_T |-\mu (\SCyl (T)) \big| \\
=&\left| \frac{1}{nM} \# \gL \theta (T) -\mu (\SCyl (T)) \right| \\
=&\left| \frac{1}{M} \frac{\# \gL}{n} \theta (T) -\frac{1}{M}\theta(T) \right| +\left| \frac{1}{M}\theta (T)-\mu (\SCyl (T)) \right| \\
=& \frac{n-\# \gL}{n} \frac{1}{M}\theta (T) +\left| \frac{1}{M}\theta (T)-\mu (\SCyl (T)) \right| .
\end{align*}
Since $\frac{1}{M}\theta (T)$ is close to $\mu (\SCyl (T))$, from Lemma \ref{lem:gl 1 and g1 dots gn estimate} and \ref{lem:gL 1 and gL estimate}, if $n$ is sufficiently large, then we have
\[ \big| |\nu_T |-\mu (\SCyl (T)) \big| <\varepsilon_2 .\]
Finally, 
\[ |\nu' | \leq \frac{\# (\gL_1\setminus \gL )}{nM} \# V(B(\id,r_\mu ')) \theta (\id )=\frac{\# (\gL_1\setminus \gL )}{n} \# V(B(\id ,r_\mu '))\frac{\theta (\id )}{M} .\]
Hence if $n$ is sufficiently large, then we have
\[ |\nu' | <\varepsilon_2.\]
Therefore we see that $\nu$ belongs to $U(f_1,\dots ,f_k; \varepsilon_\mu )$ under Assumption $(\ast )$.
\medskip

Now, we consider the case that the condition in Assumption $(\ast )$ does not hold.
Let $Y\in \Comp (\gG )$. Consider a $Y$-quasi-geodesic $\ell $ in $|Y|$. From the construction of $X_{H_n}$ the degree of a vertex $v$ in $X_{H_n}$ is less than the degree of $v$ in $X$ if and only if $v$ belongs to $H_n\subset V(X)$. This implies that we can not extend the $Y$-quasi-geodesic $\ell$ to a $Y$-quasi-geodesic line if and only if $\ell$ meets a vertex of $H_n$. This situation corresponds to the situation that $|\gG|$ has a vertex with degree less than $2$ in the previous subsection. Recall that in that case we constructed the SC-graph $(\hat{\gG}, \hat \iota )$ on $(H_n ,CH_{H_n})$ from $(|\gG |, |\iota |)$ in Step \ref{step:4}.

We also construct such a graph $\hat \gG$ from $|\gG|$ so that Assumption $(\ast )$ holds in $\hat \gG$. Explicitly,
$\hat{\gG }$ will include $|\gG |$ and if we have a $Y$-quasi-geodesic $\ell$ in $|Y|$ containing a point $x$ such that $B(x, C_0)\subset X_{H_n}$ for a constant $C_0>0$, then we can extend $\ell$ to a $\delta$-quasi-geodesic line in the connected component $W$ of $\hat \gG$ including $\ell$.
The point is that we need to modify the subgroup $H_n$ in contrary to the previous subsection.

In order to extend a $\delta$-quasi-geodesic segment $\gamma$ to a $\delta$-quasi-geodesic line, we consider a \ti{piecewise quasi-geodesic curve} in $\HH$, which is a curve consisting of at most countably many quasi-geodesic pieces. From the fundamental hyperbolic geometry in $\HH$, we can see that if a piecewise geodesic curve $\ell$ satisfies the following conditions, then $\ell$ is an $(a,c)$-quasi-geodesic for constants $a\geq 1,c\geq 0$ depending on the following constants $\theta_0 >0, L>0$:
\begin{enumerate}
\item every interior angle of $\ell$ is bounded below by some $\theta_0>0$;
\item the length of every geodesic piece of $\ell$ is larger than $L>0$, which depends on $\theta_0$.
\end{enumerate}

For a piecewise ``quasi-geodesic'' curve $\ell$, we can obtain a piecewise geodesic curve $\ell'$ by connecting endpoints of each quasi-geodesic piece of $\ell$ by a geodesic segment. Then we can see that if $\ell$ satisfies the following conditions, then $\ell$ is an $(a,c)$-quasi-geodesic for constants $a\geq 1,c\geq 0$ depending on the following constants $s\geq 1, t\geq 0, \theta_0 >0$:
\begin{enumerate}
\item there exist $s\geq 1, t\geq 0$ such that every quasi-geodesic piece of $\ell$ is a $(s,t)$-quasi-geodesic.
\item every interior angle of $\ell'$ is bounded below by some $\theta_0>0$;
\item the length of every geodesic piece of $\ell'$ is larger than $L_0>0$ depending on $s,t$ and $\theta_0$.
\end{enumerate}

Since we need to consider a quasi-geodesic line in $X$, we want to check whether a piecewise quasi-geodesic in $X$ is a quasi-geodesic or not. By using the quasi-isometry $\Phi$ from $X$ to $\HH$ we can see that a curve $\ell$ in $X$ is an $(a',c')$-quasi-geodesic if $\Phi (\ell )$ is an $(a,c)$-quasi-geodesic in $\HH$. The constants $a',c'$ depend on $a,c$ and $\Phi$.
From the above, we obtain the following lemma, which gives a sufficient condition for a piecewise quasi-geodesic curve $\ell$ in $X$ to be a quasi-geodesic in $X$.

\begin{lemma}\label{lem:piecewise quasi-geodesic will be quasi-geodesic algorithm}
Let $\ell$ be a piecewise quasi-geodesic curve in $X$. Assume that every quasi-geodesic piece of $\ell$ is an $(a, c)$-quasi-geodesic for $a\geq 1, c\geq 0$.
Let $\ell'$ be the piecewise geodesic of $\HH$ consisting of geodesic segments connecting endpoints of $\Phi(\gamma )$ for each quasi-geodesic piece $\gamma$ of $\ell$. Fix $\theta_0>0$. If $\ell$ satisfies the following conditions, then $\ell$ is an $(a',c')$-quasi-geodesic in $X$ for constants $a'\geq 1, c'\geq 0$:
\begin{enumerate}
\item every interior angle of $\ell'$ is bounded below by $\theta_0$;
\item the length of every quasi-geodesic piece of $\ell$ is larger than $L_0>0$, which depends on $a, c, \theta_0, \Phi$.
\end{enumerate}
The constants $a',c'$ depend on $a,c, \theta_0, \Phi$. 
\end{lemma}

We will use the above lemma for the case that $\theta_0$ is close to $\pi/2$.
Note that if $a,c$ depend only on $\delta$, then $a',c'$ depend only on $\delta, \theta_0, \Phi$, which implies that $a',c'$ are independent of $n$.

Now, we prepare for modifying $H_n$ and construct a graph $\hat \gG$ from $|\gG|$. 
Recall the construction of $\gS_n$. Let $\tilde{B}$ be the boundary component of $CH_{H_n}$ passing through $\tilde{x_0}$.
Then $\tilde{B}$ is a lift of the closed curve $c_0$, and $h_0:=[c_0]\in G=\pi_1(\gS,x_0)$ acts on $\tilde{B}$.
Note that $h_0\in H_n$ for every $n\geq 2$ and $\tilde{B}$ coincides with the axis $\mathrm{Ax}_\HH(h_0)$ of $h_0$ in $\HH$.
The point is that $\tilde{B}$ and $h_0$ do not depend on $n$.

We give an orientation to $\tilde{B}$ such that the left side of $\tilde{B}$ is the interior of $CH_{H_n}$.
Then we take a non-trivial element $u_0 \in G$ satisfying the following conditions:
\begin{enumerate}
\item the axis $\mathrm{Ax}_{\HH}(u_0)$ of $u_0$ in $\HH$ is included in the right side of $\tilde{B}$;
\item the hyperbolic distance $d_{\HH}(\tilde{B},\mathrm{Ax}_{\HH} (u_0))$ between $\tilde{B}$ and $\mathrm{Ax}_{\HH} (u_0)$ is sufficiently large;
\item the translation length $\tau_{\HH} (u_0)$ of $u_0$ in $\HH$ is also sufficiently large.
\end{enumerate}
Note that $d_{\HH}(\tilde{B},\mathrm{Ax}_{\HH} (u_0))$ and $\tau_\HH (u_0)$ depend on constants related to $\delta$ but do not depend on $n$.

For $u_0$ in the above, we can take a $\delta$-quasi-geodesic line $\mathrm{Ax}(u_0)$ in $X$ connecting the two points of $\gL (\langle u \rangle )$ such that $\mathrm{Ax}(u_0)$ is $\langle u \rangle $-invariant, which can be considered as an axis of $u_0$ in $X$.
For $h_0$ there exists a unique geodesic line $\mathrm{Ax}(h_0)$ in $X$ connecting the two points of $\tilde{B}(\infty )=\gL (\langle h_0 \rangle )$ such that $\mathrm{Ax}(h_0)$ includes $\langle h_0 \rangle (\subset V(X) )$. Note that $\mathrm{Ax}(h_0)$ coincides with the boundary component of $X_{H_n}$ passing through $\id$.
Then we can see that $d(\mathrm{Ax}(u_0), \mathrm{Ax}(h_0) )$ is sufficiently large and the translation length $\tau_X (u_0)$ in $X$ is also sufficiently large.

Take a geodesic $\ell_X (u_0)$ joining a point $p_{u_0}$ of $\mathrm{Ax} (u_0)$ to a point $h$ of $\langle h_ 0 \rangle $ such that the length of $\ell_X (u_0)$ equals $d(\tilde{B}, \mathrm{Ax}(u_0))$. Here, we can assume that $h =\id $ by using $h^{-1} uh$ instead of $u$.
See Figure \ref{fig:geodesic connecting axis} for the setting.
Then we can obtain the following lemma:

\begin{figure}[h]
\begin{center}
\includegraphics[width=8cm]{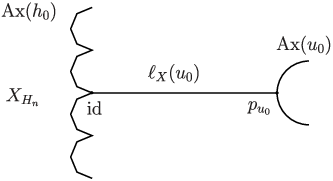}
\vspace{-0.3cm}
\caption{This figure shows the positional relationship between $X_{H_n}$ and $\mathrm{Ax}(u_0)$. By Lemma \ref{lem:extending a quasi-geodesic to piecewise quasi-geodesic}, a $\delta$-quasi-geodesic segment in $X_{H_n}$ one of whose endpoints is $\id$ can be extended to a $\delta$-quasi-geodesic half-line in $X_{H_n}\cup \ell_X(u_0)\cup \mathrm{Ax}(u_0)$.}\label{fig:geodesic connecting axis}
\end{center}
\end{figure}

\begin{lemma}\label{lem:extending a quasi-geodesic to piecewise quasi-geodesic}
Let $\gamma $ be a $\delta$-quasi-geodesic in $X_{H_n}$ from $v\in X_{H_n}$ to $\id$.
Consider a piecewise quasi-geodesic $\gamma'$ by connecting $\gamma$ to $\ell_X (u_0)$ at $\id$, and connecting $\ell_X(u_0)$ to a quasi-geodesic half-line of $\mathrm{Ax}(u_0)$ at $p_{u_0}$.
If the length of $\gamma$ is sufficiently large, then $\gamma '$ is a $\delta$-quasi-geodesic half-line.
\end{lemma}
\begin{proof}
The point is that $\Phi(\ell_X(u_0))$ is close to the common perpendicular of $\tilde{B}$ and $\mathrm{Ax}_{\HH }(u_0)$.
Then we can apply Lemma \ref{lem:piecewise quasi-geodesic will be quasi-geodesic algorithm} to $\gamma'$ and our claim follows.
\end{proof}

Set
\[ \hat{H_n}:= \langle H_n \cup \{ u_0 \} \rangle .\]
We assume that the translation length $\tau_X(u_0)$ and $d(\mathrm{Ax}(h_0), \mathrm{Ax}(u_0))$ are sufficiently large such that $h_0$ and $u_0$ generate a Schottky subgroup of $G$.
Then by the Ping-Pong argument, $\hat{H_n}$ satisfies the following properties:
\begin{enumerate}
\item $\hat{H_n}$ is isomorphic to the free product of $H_n$ and $\langle u_0 \rangle $;
\item for any $g\in \hat{H_n}\setminus H_n$ we have
\[ g(X_{H_n})\cap X_{H_n}=\emptyset ;\]
\item for every non-trivial $h\in H_n$ we have
\[ h (\mathrm{Ax}(u_0))\cap \mathrm{Ax}(u_0)=\emptyset . \]
\end{enumerate}

We consider each connected component of $|\gG|$ as a subgraph of $X$, and then for $g\in \hat{H_n}$ we define $g| \gG |$ to be the disjoint union of the image of connected components of $|\gG|$ by $g$.
Set
\[ | \gG |^\ast := \{ (gH_n ,x )\mid gH_n \in \hat{H_n}/H_n,\ x \in g |\gG | \}.\]
Then $|\gG|^\ast$ is homeomorphic to the disjoint union $\bigsqcup_{gH\in \hat{H_n}/H_n} g|\gG |$.
Note that this way of constructing $|\gG |^\ast$ corresponds to the map $\iota_{H_n}^{\hat{H_n}}$ from $\SC(H_n)$ to $\SC (\hat{H_n})$.
Then $\hat{H_n}$ acts on $|\gG|^\ast$ by
\[ g (g'H_n,x ):=(gg'H_n, gx)\]
for $g\in \hat{H_n}$ and $(g'H_n,x )\in |\gG |^\ast$.

Take the sub-arc $[p_{u_0}, u_0(p_{u_0})]$ of $\mathrm{Ax}(u_0)$ joining $p_{u_0}$ to $u_0 (p_{u_0})$.
Set
\[ P:= \ell_X (u_0 )\cup [p_{u_0}, u_0(p_{u_0})] .\]
Note that this subgraph $P$ of $X$ corresponds to the subgraph $P$ for constructing $\hat {\gG}$ in the previous subsection.
By the Ping-Pong argument, we can see that for every non-trivial $h\in \hat{H_n}$, $hP\cap P\not=\emptyset$ if and only if $h=u_0$ or $u_0^{-1}$ and $hP\cap P=\{ u_0 (p_{u_0})\}$ or $\{ p_{u_0}\}$, respectively.

Let $v(h,T)\in V(\gG )$ and take $\gamma \in \gamma_{T}$ with $\gamma \cap h\F \not=\emptyset$.
Fix a constant $C_0>0$ independent of $n$.
Consider the case that $\gamma$ contains a point $x$ with $B(x,C_0)\subset X_{H_n}$ and the extension $\gamma '$ of $\gamma \cap B(h\F ,\rho_2)$ is not a $\delta$-quasi-geodesic line.
In this case $\gamma '$ must meet a vertex $g$ of $H_n$. By considering $g^{-1}(\gamma ')$ instead of $\gamma '$, we can assume that $\gamma '$ meets $\id$. In this setting, the length of $\gamma'$ is larger than or equal to $C_0$, and so we can assume that the length of $\gamma'$ is sufficiently large to apply Lemma \ref{lem:extending a quasi-geodesic to piecewise quasi-geodesic} to $\gamma'$.

Now, we consider the disjoint union
\[ |\gG |^\ast \sqcup \bigsqcup_{h\in \hat{H_n}}h(P ).\]
Note that $\hat{H_n}$ acts on this union from left.
First, for every $h\in \hat{H_n}$ we attach the vertex $h$ of $hP$ to the vertex $h$ of $h\gamma '\subset |\gG|^\ast$.
Then for every $h\in \hat{H_n}$ we attach the vertex $h(u_0(p_{u_0}))$ of $hP$ to the vertex $hu_0(p_{u_0})$ of $hu_0P$.
By this operation of the attachment we obtain $|\gG |'$ such that $\hat{H_n}$ acts on $|\gG |'$ and the connected component of $|\gG|'$ including $\gamma'$ includes $\ell_X (u_0)$ and $\mathrm{Ax}(u_0)$.
Hence for every $h\in \hat{H_n}$ we can extend $h\gamma'$ to a $\delta$-quasi-geodesic line by using Lemma \ref{lem:extending a quasi-geodesic to piecewise quasi-geodesic}.

We can perform this operation for the disjoint union $|\gG |'\sqcup \bigsqcup_{h\in \hat{H_n}}h(P)$ and repeat the same operation until $|\gG |'$ satisfies the condition that for every $v(h,T)\in V(\gG) $ and $\gamma \in \gamma_T$ with $\gamma h\F \not=\emptyset$ if $\gamma$ contains a point $x$ with $B(x,C_0)\subset X_{H_n}$, then there exists a $\delta$-quasi-geodesic line $\ell$ in $|\gG |'$ such that $\ell$ includes $\gamma \cap B(h\F ,\rho_2)$.
Then we denote by $\hat{\gG}$ the graph that we obtain as the result of the above operation.
Note that in order to obtain $\hat{\gG}$ we perform the above operation at most $\# |\iota|^{-1}(\id )$ times since in the case that two quasi-geodesics $\gamma_1$ and $\gamma_2$ meet $\id$ in the same connected component of $|\gG|$, it is sufficient to perform the above operation only once.
We have
\begin{align*}
\# |\iota |^{-1}(\id )
=&\# \{ Z\in \Comp (|\gG |)\mid Z\ni \id \} \\
=&\# \{ v(\id ,T )\in V(\gG )\mid |T|\ni \id \} \\
=&\sum_{\substack{T\in \R_{\rho_1}^\ast (\F ) \\[1pt] |T|\ni \id}}\theta (T) \\
=&\theta (\id ).
\end{align*}
Let $\hat{m}$ be the number of times we perform the above operation. Then $\hat{m}\leq\theta (\id )$.
Denote by $P_j$ the copy of $P$ that we used in the $j$-th operation for $j=1,\dots ,\hat{m}$.

The projection from the disjoint union $|\gG |^\ast \sqcup \bigsqcup_{h\in \hat{H_n}}h(P)$ to $X$ induces a map $\hat{\iota}$ from $\hat{\gG}$ to $X$.
We can see that the restriction of $\hat{\iota}$ to each connected component $W$ of $\hat{\gG}$ is injective from the Ping-Pong argument.
We identify each connected component $W$ of $\hat{\gG}$ with $\hat{\iota}(W)$.
From the above, $(\hat{\gG} , \hat{\iota})$ satisfies the condition in Assumption $(\ast)$ essentially.

Now, we define $\eta_{\hat \gG }$ by
\[ \eta_{\hat \gG }:= \sum_{W\in \Comp (\hat{\gG })} \delta_{W(\infty )}.\]
Then we can see that $\eta_{\hat \gG} \in \SC (\hat {H_n})$. The local finiteness of $\eta_{\hat \gG}$ follows by the argument below.
Set
\[ \nu: =\frac{1}{nM}\iota_{\hat{H_n}}(\eta_{\hat{\gG }}) \in \SC (G ).\]
We prove that $\nu$ belongs to the open neighborhood $U(f_1,\dots ,f_k; \varepsilon_\mu )$ of $\mu$ for a large $n$ by using Lemma \ref{lem:describe a neighborhood of a given subset current}.

\begin{lemma}\label{lem: conn comp of hat gG is quasi convex}
Every connected component of $\hat{\gG}$ is a $\delta'$-quasi-convex subgraph of $X$.
\end{lemma}
\begin{proof}
Let $W$ be a connected component of $\hat{\gG}$. Take $x , y \in W$. We prove that there exists a $\delta$-quasi-geodesic joining $x$ to $y$ included in the $\delta'$-neighborhood of $W$.
Then by the stability of quasi-geodesics, $W$ is $\delta'$-quasi-convex.
If $x,y$ belong to $Z$ for a connected component $Z$ of $|\gG|^\ast$, then $W$ includes $Z$ and there exists a $\delta$-quasi-geodesic joining $x$ to $y$ in $Z$ by Lemma \ref{lem:a connected component of Y is quasi-convex}.

Hence we consider the case that for different connected components $Z,Z'$ of $|\gG|^\ast$, $x\in Z$ and $y\in Z'$.
Take a shortest path $\ell$ from $x$ to $y$ in $W$. From the construction of $\hat{\gG}$ there exists a sequence of connected components $Z_0=Z, Z_1,\dots ,Z_k=Z'$ of $|\gG |^\ast$ such that $\ell$ passes through these components in this order.
From $Z_{i-1}$ to $Z_{i}$, the path $\ell$ passes through $h_iP$ when $\ell$ goes out from $Z_{i-1}$, and passes through $h_i'P$ when $\ell$ goes into $Z_{i}$ for some $h_i,h_i'\in \hat{H_n}$. Since the translation length $\tau(u_0)$ and the length of $\ell_X(u_0)$ are sufficiently large, the restriction of $\ell$ to this part is a $\delta$-quasi-geodesic in $X$.

Now, for each $i=1,\dots ,k$ we take the mid-point $m_i$ of $h_i (\ell_X(u_0))$ and $m_i'$ of $h_i' (\ell_X (u_0))$ and consider a geodesic $[m_{i}', m_{i+1}]$ joining $m_{i}'$ to $m_{i+1}$ in $X$, which is included in the $\delta'$-neighborhood of the union of $Z_i$, $h_i'(\ell_X(u_0))$ and $h_{i+1}(\ell_X (u_0))$.
Then we consider the following path $\ell'$ from $x$ to $y$:
\begin{enumerate}
\item starts from $x$ and bounds for $m_1'$ along $\ell$;
\item from $m_i'$ to $m_{i+1}$ travels along the geodesic $[m_i',m_{i+1}]$, and from $m_{i+1}$ to $m_{i+1}'$ travels along $\ell$ for $i=1,\dots ,k$;
\item from $m_{k}'$ to $y$ travel along $\ell$.
\end{enumerate}
The path $\ell'$ is a piecewise quasi-geodesic in $X$ and if the translation length $\tau(u_0)$ and the length of $\ell_X(u_0)$ are sufficiently large, then $\ell'$ is a $\delta$-quasi-geodesic in $X$.

In other cases we can construct the almost same piecewise quasi-geodesic joining $x$ to $y$.
\end{proof}

Then we obtain the following lemma for the constant $\rho_3=\rho_2-\delta_1-C_0$, which corresponds to Lemma \ref{lem:with assumption and Z infty in a neighborhood of SCyl} under Assumption $(\ast )$.

\begin{lemma}\label{lem:no assumption and W infty in a neighborhood of SCyl}
Assume that $\rho_3$ is sufficiently large. Let $W\in \Comp (\hat \gG )$, $g \in V(W)$.
If $B(g,\rho_2')$ is included in $X_{H_n}$, then 
\[ g^{-1}W(\infty )\in B_{\H}(\SCyl ( g^{-1}W\cap B(\id ,\rho_3 )) ,\varepsilon_1 ).\]
\end{lemma}
\begin{proof}
Since $B(g,\rho_2')\subset X_{H_n}$, there exists a connected component $Z$ of $|\gG|$ such that
\[ W\cap B(g, \rho_2' )=Z\cap B(g,\rho_2' )\]
by the construction of $\hat{\gG}$.
Then by the same argument as that in the proof of Lemma \ref{lem:with assumption and Z infty in a neighborhood of SCyl}, we see that $W\cap B(g, \rho_3)=Z\cap B(g ,\rho_3)$ belongs to $\R_{\rho_3}(g)$, and
\[ g^{-1}W(\infty )\in B_{\H}(\SCyl ( g^{-1}W\cap B(\id ,\rho_3 )) ,\varepsilon_1 )\]
if $\rho_3$ is sufficiently large.
\end{proof}

Now, we construct a subgraph $X_{\hat{H_n}}$ of $X$ such that every connected component of $\hat \gG$ is included in $X_{\hat{H_n}}$. 
By the same way as we did for $|\gG|$, we set
\[ X_{H_n}^\ast :=\{ (gH_n, x)\in \hat{H_n}/H_n \times X\mid x \in gX_{H_n} \} \]
and consider the disjoint union
\[ X_{H_n}^\ast \sqcup \bigsqcup_{h\in \hat{H_n}}h(P).\]
For every $h\in \hat{H_n}$ we attach the vertex $h$ of $hP$ to the vertex of $h$ of $X_{H_n}^\ast$ and attach the vertex $h(u_0(p_{u_0}))$ of $hP$ to the vertex $hu_0(p_{u_0})$ of $hu_0P$.
By this attachment we obtain a connected graph $X_{\hat{H_n}}$ and the inclusion map from $X_{H_n}^\ast \sqcup \bigsqcup_{h\in \hat{H_n}}h(P)$ to $X$ induces an injective map from $X_{\hat{H_n}}$ to $X$ from the property of $\hat{H_n}$.
Hence we can consider $X_{\hat{H_n}}$ as a subgraph of $X$, which is $\hat{H_n}$-invariant.
Moreover, by the same argument as that in Lemma \ref{lem: conn comp of hat gG is quasi convex}, we see that $X_{\hat{H_n}}$ is a $\delta'$-quasi-convex subgraph of $X$ and for every $x\in X_{\hat{H_n}}$ there exists a $\delta$-quasi-geodesic line passing through $x$. Hence we see
\[ X_{\hat{H_n}}\underset{\delta '}{\sim }CH (\gL (\hat{H_n})).\]

Note that the quotient graph $\hat{H_n}\backslash X_{\hat{H_n}}$ can be described as follows.
Recall that $H_n\backslash X_{H_n}$ can be identified with the $1$-skeleton of $\gS_n$.
By attaching the vertex $p_{u_0}$ of $P$ to the vertex $u_0(p_{u_0})$, we obtain a graph $P'$, which is homotopic to a circle.
Then we attach the vertex $\id$ of $P'$ to the vertex $\tilde{x_0}^n$ of $H_n\backslash X_{H_n}$.
The resulting graph is isomorphic to $\hat{H_n}\backslash X_{\hat{H_n}}$.

Take a complete system of representatives $\hat{\gL}_0$ of $G/\hat{H_n}$.
Set
\[ \hat{\gL}_1:= \{ g\in \hat{\gL}_0 \mid g CH(\gL(\hat{H_n}) )\cap B(\id ,r_\mu )\not =\emptyset \}, \]
which is a finite set.
Then
\[ \iota_{H_n}(\eta_{\hat \gG })|_{A(\id, r_\mu )}=\sum_{g\in \hat{\gL}_1 }g_{\ast}(\eta_\gG )|_{A(\id ,r_\mu )}.\]

Recall that $\F$ includes exactly $n$ vertices $g_1=\id ,g_2^{-1},\dots , g_n^{-1}$ of $X_{H_n}$.
By considering the action of $\hat{H_n}$ on $X_{\hat{H_n}}$, we see that $g_1\hat{H_n}, \dots ,g_n\hat{H_n}$ are mutually disjoint.
Hence we can assume that $g_1,\dots ,g_n\in \hat{\gL}_1$.

The following lemma corresponds to Lemma \ref{lem:gl 1 and g1 dots gn estimate}.

\begin{lemma}\label{lem:in hat gG gl 1 and g1 dots gn estimate}
The sequence
\[ \frac{1}{n}\# \left( \hat{\gL}_1 \setminus \{ g_1,\dots ,g_n\}\right) \]
tends to $0$ when $n\rightarrow \infty$.
\end{lemma}
\begin{proof}
Note that the translation length $\tau(u_0)$ and the length of $\ell_X(u_0)$ are independent of $n$.
Hence $\# V(P)$ is independent of $n$.

First, we have
\begin{align*}
\#\hat{\gL}_1 =&\# \{ g\hat{H_n}\in G/\hat{H_n} |\ gCH(\gL(\hat{H_n}))\cap B(\id ,r_\mu )\not =\emptyset \} \\
=&\# V(\hat{H_n}\backslash B(CH(\gL (\hat{H_n}) ), r_\mu )).
\end{align*}

Since $X_{\hat{H_n}}\underset{\delta '}{\sim }CH (\gL (\hat{H_n}))$, we have
\[ \# \hat{\gL}_1 \leq \# V(\hat{H_n}\backslash B(X_{\hat{H_n}}, r_\mu +\delta ')) .\]
From the definition of $X_{\hat{H_n}}$ we have
\[ \# V(\hat{H_n}\backslash X_{\hat{H_n}})=\#V (H_n\backslash X_{H_n})+ \# V(P)-2.\]
Note that $V (H_n\backslash X_{H_n})$ corresponds to $\{ H_ng_1^{-1}, \dots , H_n g_n^{-1} \}$.
By considering the degree of each vertex of $\hat{H_n}\backslash X_{\hat{H_n}}$ we have
\begin{align*}
&\# V(\hat{H_n}\backslash B(CH(\gL (\hat{H_n}) ), r_\mu ))-\# V(H_n\backslash X_{H_n}) \\
\leq &\# V(\hat{H_n}\backslash B(X_{\hat{H_n}}, r_\mu +\delta ' ))-\# V(H_n\backslash X_{H_n}) \\
\leq &\# V(P)\big( \mathrm {deg}_X( \id ) \big)^{ r_\mu +\delta '},
\end{align*}
which implies
\[ \frac{1}{n}\# \left( \hat{\gL}_1 \setminus \{ g_1,\dots ,g_n\}\right)\leq \frac{1}{n} \# V(P) \big( \mathrm {deg}_X( \id ) \big)^{r_\mu +\delta '}.\]
This proves our claim.
\end{proof}

From the above proof it is easy to see that the argument for $\hat{\gG}$ is almost the same as that for $|\gG |$ under Assumption $(\ast )$. Moreover, since $\# V(P)$ is a constant not depending on $n$, $\# V(P)$ does not influence our argument.
For the completeness of the proof, we continue the almost same argument as that under Assumption $(\ast )$.

Recall that $\gL =\{ g_i \mid B(g_i^{-1} ,\rho_2' )\subset X_{H_n} \}$.
We also see that $\frac{1}{n} \# (\hat{\gL}_1\setminus \gL )$ tends to $0$ when $n\rightarrow \infty $ by the same argument as that in Lemma \ref{lem:gL 1 and gL estimate}.
Then
\begin{align*}
\iota_{\hat{H_n}}(\eta_{\hat \gG} )|_{A(\id, r_\mu )} 
=&\sum_{g\in \hat{\gL}_1 }g_{\ast}(\eta_{\hat \gG} )|_{A(\id ,r_\mu )}\\
=&\left( \sum_{g\in \hat{\gL}_1\setminus \gL } +\sum_{g\in \gL} \right) g_{\ast}(\eta_{\hat \gG} )|_{A(\id ,r_\mu )},
\end{align*}
and we mainly consider the sum taken over $g\in \gL$.

First we have
\[ \sum_{g\in \gL} g_{\ast}(\eta_{\hat{\gG}} )|_{A(\id ,r_\mu )}
=\sum_{g\in \gL } \sum_{\substack{W\in \Comp (\hat \gG ) \\[1pt] gW(\infty )\in A(\id ,r_\mu )}} \delta_{gW(\infty )}.\]
Note that every connected component $W$ of $\hat \gG$ is $\delta '$-quasi-convex. By the same argument as before, for a constant $\beta >0$ depending on $\delta$, we see that if $gW(\infty ) \in A(\id ,r_\mu )$, then $gW\cap B(\id, r_\mu +\beta )\not=\emptyset$.

The following setting corresponds to Setting \ref{set:9} under Assumption $(\ast )$:

\begin{setting8}\label{set:10}
For the constant $\beta$ in the above, we set $r_\mu':=r_\mu +\beta$, which corresponds to the constant $r_\mu'$ in Lemma \ref{lem:describe a neighborhood of a given subset current}.
We assume that $\rho_2 \geq 2r_\mu' +\rho$.
\end{setting8}
\begin{lemma}
Let $g\in \gL$ and $W\in \Comp (\hat \gG )$ with $gW(\infty )\in A(\id ,r_\mu )$.
Then $gW \cap B(\id , r_\mu' +\rho )$ is an element of $\R_\rho (B(\id ,r_\mu' ))$.
\end{lemma}
\begin{proof}
The point is that $B(g^{-1}, r_\mu'+\rho )\subset B(g^{-1}, \rho_2)\subset X_{H_n}$ implies that there exists a connected component $Z$ of $|\gG|$ such that $Z\subset W$ and
\[ W\cap B(g^{-1}, \rho_2 )=Z\cap B(g^{-1}, \rho_2 ).\]
Hence we have
\[ gW \cap B(\id ,r_\mu' +\rho )=gZ \cap B(\id ,r_\mu' +\rho ) ,\]
which is an element of $\R_\rho (B(\id ,r_\mu' ))$ from Lemma \ref{lem: intersection will be round graph}.
\end{proof}

For $T\in \R_\rho (B(\id ,r_\mu ')) $ we set
\[ \iota_{\hat{H_n}}(\eta_{\hat{\gG}})_T 
:=\sum_{g\in \gL } \sum_{\substack{W\in \Comp (\hat \gG ) \\[1pt] gW(\infty )\cap B(\id ,r_\mu' +\rho )= T}}\delta_{ g W(\infty )}. \] 
Then
\[ \sum_{g\in \gL} g_{\ast}(\eta_{\hat{\gG}} )|_{A(\id ,r_\mu )}= \sum_{T\in \R_\rho (B(\id , r_\mu '))} \iota_{\hat{H_n}}(\eta_{\hat{\gG}})_T .\]

Now, we prove the following lemma, which corresponds to Lemma \ref{lem:under assumption iota Hn eta gG is good}:

\begin{lemma}\label{lem: iota Hn eta hat gG is good}
For each $T\in \R_{\rho}(B(\id ,r_\mu'))$ we have
\[ \mathrm{supp} (\iota_{\hat{H_n}}(\eta_{\hat{\gG}})_T )\subset \ol{B_{\H}(\SCyl (T) ,\varepsilon_1 )}.\]
Moreover,
\[ | \iota_{\hat{H_n}}(\eta_{\hat{\gG}})_T |=\# \gL \cdot \theta (T).\]
\end{lemma}
\begin{proof}
Fix $T\in \R_{\rho}(B(\id ,r_\mu '))$.
For $g\in \gL$ consider $W\in \Comp (\hat{\gG})$ with $gW\cap B(\id, r_\mu'+\rho )=T$.
From Lemma \ref{lem:no assumption and W infty in a neighborhood of SCyl}, we have
\[ gW(\infty )\in B_{\H}(\SCyl ( gW\cap B(\id ,\rho_3 )) ,\varepsilon_1 ).\]
We assume that $\rho_3 >r_\mu' +\rho$.
Since $T=(gW\cap B(\id ,\rho_3 ))\cap B(\id, r_\mu'+\rho)$,
\[ \SCyl ( gW\cap B(\id ,\rho_3 ))\subset \SCyl (T),\]
and so
\[ gW(\infty )\in B_{\H}(\SCyl (T) ,\varepsilon_1 ).\]
Therefore
\[ \mathrm{supp} (\iota_{\hat{H_n}}(\eta_{\hat{\gG}})_T )\subset \ol{B_{\H}(\SCyl (T) ,\varepsilon_1 )}.\]

Now, we calculate $| \iota_{\hat{H_n}}(\eta_{\hat{\gG}})_T |$.
Fix $g\in \gL$. Take a vertex $u\in T\cap B(\id, r_\mu')$.
Take $W\in \Comp (\hat{\gG})$.

Suppose that $gW\cap B(\id, r_\mu'+\rho )=T$.
Then there exists a connected component $Z$ of $|\gG|$ such that $Z\subset W$ and
\[ W\cap B(g^{-1}, \rho_2 )=Z\cap B(g^{-1}, \rho_2 ).\]
Moreover, for $v(h', T')\in V(Y_Z)$ with $g^{-1}u\in h'\F \cap |T'|$, we have
\[ T=gW\cap B( \id, r_\mu' +\rho) =gZ \cap B(\id ,r_\mu '+ \rho )=g|T'| \cap B(\id ,r_\mu' +\rho )\]
by the same argument as that in the proof of Lemma \ref{lem:under assumption iota Hn eta gG is good}.
Hence
\[ |T'|\cap B(g^{-1},r_\mu' +\rho )=g^{-1}T.\]

Conversely, suppose that there exists a connected component $Z$ of $|\gG|$ and $v(h',T')\in V(Y_Z)$ with $g^{-1}u\in h'\F \cap |T'|$ such that $Z\subset W$ and
\[ |T'|\cap B(g^{-1},r_\mu' +\rho )=g^{-1}T.\]
Then 
\[ W\cap B(g^{-1}, \rho_2 )=Z\cap B(g^{-1}, \rho_2 ),\]
and so
\[ gW\cap B( \id, r_\mu' +\rho) =gZ \cap B(\id ,r_\mu '+ \rho )=g|T'| \cap B(\id ,r_\mu' +\rho )=T.\]

Hence the number of $W\in \Comp (\hat{\gG})$ satisfying the condition that $gW\cap B(\id, r_\mu'+\rho )=T$ equals the number of $Z\in \Comp(|\gG|)$ satisfying the condition that
there exists $v(h', T')\in V(Y_Z)$ with $g^{-1}u\in h'\F \cap |T'|$ such that
\[ |T'|\cap B(g^{-1},r_\mu' +\rho )=g^{-1}T.\]

Therefore, from the proof of Lemma \ref{lem:under assumption iota Hn eta gG is good} we have
\begin{align*}
&| \iota_{\hat{H_n}}(\eta_{\hat{\gG}})_T |\\
=&\sum_{g\in \gL }\# \{ W\in \Comp (\hat \gG) \mid gW \cap B(\id , r_\mu' +\rho )=T\} \\
=&\sum_{g\in \gL }\# \{ Z\in \Comp (|\gG |) \ |\\
	&\quad \exists v(h', T')\in V(Y_Z) \text{ s.t. } h'\F \ni g^{-1}u \text{ and }T'\cap B(g^{-1}, r_\mu '+\rho )=g^{-1}T \} \\
=&\sum_{g\in \gL }\theta (g^{-1}T)=\# \gL \cdot \theta (T).
\end{align*}
This completes the proof.
\end{proof}

Now, we consider the other part of $\iota_{\hat{H_n}}(\eta_{\hat \gG})$. Let $g\in \hat{\gL}_1\setminus \gL$.
Then we have
\begin{align*}
&\big| (g_\ast (\eta_{\hat \gG }) )|_{A(\id ,r_\mu )} \big|\\
=&\eta_{\hat{\gG}} (A(g^{-1} ,r_\mu ))\\
=&\# \{ W\in \Comp (\hat \gG ) \mid CH(W(\infty ))\cap B(g^{-1}, r_\mu )\not=\emptyset \} \\
\leq &\# \{ W\in \Comp (\hat \gG ) \mid W\cap B(g^{-1}, r_\mu' )\not=\emptyset \}.
\end{align*}
If $W\cap B(g^{-1}, r_\mu ')\not=\emptyset $ for $W\in \Comp (\hat \gG)$, then there exists $Z\in \Comp (|\gG| )$ such that $Z\subset W$ and $Z\cap B(g^{-1}, r_\mu ')\not=\emptyset$, or there exist $j\in \{ 1,\dots ,\hat{m} \}$ and $g_0 \in \hat{H_n}$ such that $g_0P_j \subset W$ and $g_0P_j\cap B(g^{-1}, r_\mu ')\not=\emptyset $.
Note that $g_0P_j\cap B(g^{-1}, r_\mu ')\not=\emptyset $ implies that $B(gg_0P_j, r_\mu') \ni \id$.
Hence for each $j\in \{ 1, \dots ,\hat{m} \}$ the number of $g_0\in \hat{H_n}$ satisfying the condition that $g_0P_j\cap B(g^{-1}, r_\mu ')\not=\emptyset $ is less than or equal to the number of vertices of $B(gg_0P_j, r_\mu')$, which is less than
\[ D:=\# V(P)\big( \mathrm {deg}_X( \id ) \big)^{ r_\mu'} .\]
Therefore,
\begin{align*}
&\big| (g_\ast (\eta_{\hat \gG }) )|_{A(\id ,r_\mu )} \big|\\
< &\sum_{v\in V(B(g^{-1},r_\mu '))}\# \{ Z\in \Comp (|\gG |) \mid Z\ni v \} + \hat{m}D\\
\leq & \# V(B(g^{-1},r_\mu ')) \theta (\id )+ \hat{m}D\\
= & \# V(B(\id ,r_\mu ')) \theta (\id )+ \hat{m}D.
\end{align*}

For $T\in \R_\rho (B(\id ,r_\mu '))$ set
\[ \nu_T:= \frac{1}{nM}\iota_{\hat{H_n}}(\eta_{\hat \gG })_T \]
and set
\[ \nu' :=\frac{1}{nM}\sum_{g\in \hat{\gL}_1\setminus \gL }g_\ast (\eta_{\hat{\gG}})|_{A(\id, r_\mu )}.\]
Then we have
\begin{align*}
\nu |_{A(\id ,r_\mu )}
=&\frac{1}{nM}\iota_{\hat{H_n}}(\eta_{\hat \gG })|_{A(\id ,r_\mu )}\\
=&\sum_{T\in \R_\rho (B(\id ,r_\mu '))} \nu_T |_{A(\id ,r_\mu )}+ \nu '.
\end{align*}
We prove that for a sufficiently large $n\in \NN$ we have $\nu \in U(f_1, \dots , f_k; \varepsilon _\mu)$ by using Lemma \ref{lem:describe a neighborhood of a given subset current}.
From Lemma \ref{lem: iota Hn eta hat gG is good}, for every $T\in \R_\rho (B(\id ,r_\mu '))$ we have

\begin{align*}
&\big| |\nu_T |-\mu (\SCyl (T)) \big| \\
=&\left| \frac{1}{nM} \# \gL \theta (T) -\mu (\SCyl (T)) \right| \\
=&\left| \frac{1}{M} \frac{\# \gL}{n} \theta (T) -\frac{1}{M}\theta(T) \right| +\left| \frac{1}{M}\theta (T)-\mu (\SCyl (T)) \right| \\
=& \frac{n-\# \gL}{n} \frac{1}{M}\theta (T) +\left| \frac{1}{M}\theta (T)-\mu (\SCyl (T)) \right| .
\end{align*}
Therefore if $n$ is sufficiently large and $\frac{1}{M} \theta(T)$ is close to $\mu (\SCyl (T))$, then
\[ \big| |\nu_T | - \mu ( \SCyl (T) ) \big| <\varepsilon_2.\]

Finally,
\begin{align*}
|\nu'| 
<& \frac{\# (\hat{\gL}_1\setminus \gL )}{nM} (\# V(B(\id ,r_\mu ')) \theta (\id )+ \hat{m}D)\\
\leq & \frac{\# (\hat{\gL}_1\setminus \gL )}{nM} (\# V(B(\id ,r_\mu ')) \theta (\id )+ \theta (\id )D)\\
=& \frac{\# (\hat{\gL}_1\setminus \gL )}{n} \left( \# V(B(\id ,r_\mu ')) +D \right) \frac{\theta (\id )}{M}.
\end{align*}
Since $\# (\hat{\gL}\setminus \gL )/ n$ tends to $0$ when $n\rightarrow \infty$, for a sufficiently large $n\in \NN$ we have
\[ |\nu ' |<\varepsilon_2 .\]
Therefore $\nu$ belongs to the open neighborhood $U(f_1,\dots ,f_k; \varepsilon_\mu)$ of $\mu$.
This completes the proof of Theorem \ref{thm:surface group has denseness property} and \ref{thm:approximation of surface group by free groups}.\hspace{\fill} Q.E.D.

\medskip
From Theorem \ref{thm:approximation of surface group by free groups}, it is natural to propose the following problem:
\begin{problem}
Let $G$ be an infinite hyperbolic group. Is there a sequence of quasi-convex subgroups $\{ H_n \}_{n\in \NN}$ of $G$ such that each $H_n$ is a free group of finite rank and the union 
\[ \bigcup_{n\in \NN} \iota_{H_n}(\SC (H_n ) )\]
is a dense subset of $\SC (G)$?
\end{problem}

This problem gives us an approach to Problem \ref{problem: denseness property for hyperbolic group}.
Moreover, if this problem is solved affirmatively, then we can say that an infinite hyperbolic group can be approximated by free quasi-convex subgroups in the meaning of subset currents.
We note that for an infinite hyperbolic group $G$ and for every $S\in \H (\partial G)$ there exists a quasi-convex subgroup $H$ of $G$ such that $H$ is a free group of finite rank and the limit set $\gL (H)$ is close to $S$.


\end{document}